\documentclass[10pt]{article}
\RequirePackage[%
    papersize={15cm,24cm},%A4
 hmargin=1cm,%
  vmargin=1.8cm,%
  head=0.2cm,% will be changed later
  headsep=0pt,%
  foot=1.1cm% will be changed later
  ]{geometry}

\usepackage{amsmath,  amssymb,  amsthm}
\usepackage{color}
\usepackage{float}
\usepackage{cases}
\usepackage[hang]{subfigure}
\usepackage{mathrsfs}
\usepackage{graphicx}

\renewcommand{\vec}[1]{\mbox{\boldmath \small $#1$}}

\newcommand\bbR{\mathbb{R}}
\newcommand\bbN{\mathbb{N}}

\newcommand\mcP{\mathcal{\bP}}

\newcommand\dd{\mathrm{d}}

\newcommand\bx{\boldsymbol{x}}

\newcommand\by{\boldsymbol{y}}
\newcommand\bF{\boldsymbol{F}}
\newcommand\bE{\boldsymbol{E}}
\newcommand\bH{\boldsymbol{H}}
\newcommand\bP{\boldsymbol{P}}

\newcommand\pd[2]{\dfrac{\partial {#1}}{\partial {#2}}}

\newcommand\abs[1]{\lvert #1 \rvert}

\newcommand\bt{\bar\theta}

\newcommand\pro[2]{\langle{#1}, {#2}\rangle}

\newtheorem{lemma}{Lemma}

\newtheorem{remark}{Remark}
\newtheorem{example}{Example}
\numberwithin{equation}{section}
\numberwithin{figure}{section}
\numberwithin{table}{section}
\numberwithin{theorem}{section}
\numberwithin{remark}{section}
\numberwithin{example}{section}
\begin{document}%\begin{CJK*}{GBK}{fs}

\title{Model reduction of  a kinetic swarming model by operator projection}
\author{Junming Duan$^*$,
%\thanks{LMAM, School of Mathematical Sciences,
%Peking University, Beijing 100871, P.R. China.
 %Email: {\tt duanjm@pku.edu.cn}},
 %
 Yangyu Kuang\thanks{LMAM, School of Mathematical Sciences,
Peking University, Beijing 100871, P.R. China.
 Email: {\tt duanjm@pku.edu.cn}; {\tt kyy@pku.edu.cn}},
 % \and
Huazhong Tang\thanks{
            HEDPS, CAPT \& LMAM, School of Mathematical Sciences, Peking University,
            Beijing 100871, P.R. China;
            School of Mathematics and Computational Science,
           Xiangtan University,  Xiangtan 411105,  Hunan Province, P.R. China.
               Email: {\tt hztang@math.pku.edu.cn}}
}
\maketitle
%\tableofcontents %If the document is very long
\begin{abstract}
  This paper derives the arbitrary order globally hyperbolic moment system
  for a non-linear kinetic description of the Vicsek swarming model
  by using the operator projection. It is built on our
  careful study of  a family of the complicate Grad type
  orthogonal functions  depending on a parameter (angle of macroscopic velocity).
  We  calculate their derivatives with respect to the independent variable,
  and projection of those derivatives, the product of velocity and basis, and collision term.
 The moment system is also proved to be
  hyperbolic,      rotational invariant, and mass-conservative.
The relationship between Grad type expansions in different parameter is also established.
  A semi-implicit numerical scheme is presented to solve a Cauchy problem of our
  hyperbolic moment system in order to verify the convergence behavior of the moment method.
  It is also compared to the spectral method for the kinetic equation.
  The results show that the solutions of our hyperbolic moment system {converge}
  to the solutions of the kinetic equation for the Vicsek   model as the order of the  moment system increases, and
  the moment method can capture key features such as vortex formation and traveling waves.

\end{abstract}
\section{Introduction}
Swarm behaviour, or swarming, is a collective behaviour exhibited by entities, particularly animals, of similar size which aggregate together, perhaps milling about the same spot or perhaps moving en masse or migrating in some direction \cite{Bouffanais}.
It is a highly interdisciplinary topic.
Some works studied kinetic models for swarming \cite{Canizo2011,Degond-Liu2013,Ha-Tadmor2008},
but few  did a numerical investigation.
The first numerical method was presented in \cite{Gamba-Haack2015} for a kinetic description of the Vicsek swarming model.
The main contribution was to use a spectral representation linked with a discrete constrained optimization to compute those interactions. Unfortunately, only first-order accurate upwind method was used to
approximated the transport term.
  %To test the numerical scheme we investigate the kinetic model at different scales and compare the solution with the microscopic and macroscopic descriptions of the Vicsek model. We observe that the kinetic model captures key features such as vortex formation and traveling waves.

The kinetic theory has been widely studied   and
played an important role in many fields  during several decades, see e.g. \cite{RB:1988,Cha:1970}.
The kinetic equation  can determine the  distribution function
hence the transport coefficients, however such task is not so easy.
% Hilbert showed that an  approximate solution of the integro-differential
%equation  could be obtained from a power series expansion of a parameter
%(being proportional to the mean free path).
%Chapman and Enskog calculated   independently
% the transport coefficients for gases whose molecules interacted according to any kind of spherically symmetric potential function.
The moment method \cite{GRAD:1949,GRAD2:1949} is a model reduction for the kinetic equation
by expanding the distribution
function in terms of tensorial Hermite polynomials and introducing
the balance equations corresponding to higher order moments of the distribution function.
One major disadvantage of the Grad moment method is the loss of hyperbolicity, which
 will   cause the solution blow-up when the distribution is far away
from the equilibrium state. Increasing the number of moments could not  avoid  such blow-up \cite{FAR:2012}.
%
%The crucial ingredient of the  Chapman-Enskog  method is
%the assumption that in the hydrodynamic regime
%the distribution function can be expressed as a function of the hydrodynamic variables and their gradients.
%
%
Up to now, there has been some latest progress on the Grad moment method for the kinetic equation.
Numerical regularized moment method of arbitrary order was studied
for Boltzmann-BGK equation \cite{NM:2010}
and for high Mach number flow \cite{FAR:2012}.
Based on the observation that the characteristic polynomial of the flux Jacobian in the Grad moment system did not depend on the intermediate moments,  a regularization  was presented in  \cite{1DB:2013,HY:2014,HG:2014} for the one- and multi-dimensional Grad moment systems
  to achieve global hyperbolicity.
 The quadrature based projection methods were used
to derive hyperbolic systems for the solution of the Boltzmann equation
  \cite{QI:2014,HG2:2014}
 by  using the quadrature rule instead of the exact integration.
In the 1D case,  it is similar to the  regularization in  \cite{1DB:2013}. %, and extended to the multi-dimensional case in \cite{QI:2014}.
 % Both methods in \cite{1DB:2013} and \cite{QI:2013} have been generalized to more general cases in \cite{HG:2014} and , respectively
Those contributions  led to  well understanding the hyperbolicity of  the Grad moment systems.
Based on the operator projection,  a general framework of model reduction technique was recently given in \cite{Fan2015Model}.
%to cover all the different methods. This method is a model reduction by the projection operator.
It  projected the time and space derivatives in the kinetic equation   into a finite-dimensional weighted polynomial space synchronously,
and could give most of the existing moment systems in the literature.
Recently, such model reduction method   was also successfully extended to the 1D special relativistic Boltzmann equation
and the globally hyperbolic moment model of arbitrary order  was derived in \cite{Kuang-Tang2016}.
 %, such as the hyperbolic moment equations proposed in \cite{1DB:2013,HY:2014,HY2:2014}, the quadrature-based moment equations \cite{QI:2013}, and Levermore's maximum entropy method \cite{Lev:1996}.
 %
%

The aim of this paper is to extend the model reduction method by the operator projection
to the two-dimensional kinetic description of the Vicsek swarming  model
and
derive corresponding globally hyperbolic moment system of arbitrary order.
The paper is organized as follows.
Section \ref{sect:model} introduces the kinetic and macroscopic equations for the Vicsek model.
Section \ref{sec:moment} gives  a family of orthogonal functions dependent on a parameter and their properties
and  derives the arbitrary order globally hyperbolic moment system
of the  kinetic description for the Vicsek model.
Section \ref{sec:4} investigates the mathematical properties of moment system, including: hyperbolicity,  rotational invariance,   mass-conservation
and relationship between Grad type expansions in different parameter.
Section \ref{sect:numEXP} presents a semi-implicit numerical scheme and conducts a numerical experiment
to check the convergence of the proposed hyperbolic moment system.
Section \ref{sec:conclud} concludes the paper. % with several remarks.

\section{Kinetic and macroscopic equations for the Vicsek model}\label{sect:model}
This section introduces the non-linear kinetic and macroscopic equations for the Vicsek model.

\subsection{Kinetic equation}
The kinetic equation for the Vicsek model can be written as follows \cite{Gamba-Haack2015}
\begin{align} \label{eq:kinetic-1}
	\partial_t f &+\vec \omega\cdot{\nabla_{\bx}} f+{\nabla_{\vec \omega}}\cdot(\vec F[f]f)=\sigma\Delta_{\vec \omega} f,
\end{align}
where $f(t, \bx, \vec \omega)$ is the particle distribution function depending on
 the time $t$,   spatial variable $\bx\in \mathbb R^d$, and the unit velocity vector $\vec \omega$,
  the parameter $\sigma$ is a scaled diffusion constant describing the intensity of the noise with  the Brownian motion,
 the vector   $\vec F[f](t,\vec x,\vec \omega)$ is the
mean-field interaction force between the particles   and
  given by
\begin{align}
 \label{eq:kinetic-2}
	\vec  F[f](t,\bx, \vec \omega)=&(\mathrm{Id}-\vec \omega\otimes\vec \omega)\vec \Omega(\bx),
\quad \vec \Omega(t,\bx)=\dfrac{\vec J(t,\bx)}{\abs{\vec J(t,\bx)}},
\end{align}
the notation $\mathrm{Id}$ is the identity operator,
$\vec \Omega(t,\bx)$ is the mean velocity,
and $\vec J(t,\bx)$ denotes the mean flux at $\vec x$ and is defined by
\begin{align}
	\label{eq:kr}
	\vec J(t,\bx)=&\int_{\bbR^d}\int_{S^{d-1}}K(\by-\bx)\vec \omega f(t,\by,\vec \omega)\dd \by\dd \vec \omega.
\end{align}
Here $K(\by-\bx)$ is the characteristic function of the ball
$B(0, R)=\{\vec x: |\vec x|\leq R\}$, i.e. 	$K(\vec x)={\boldsymbol{1}}_{\abs{\vec x}<R}$, and
$R$ is the radius of  the alignment interactions between the particles.
The vector field $\vec F[f](t,\vec x,\vec\omega)$
tends to align the particles to the direction
$\vec \Omega$ which is the director of the particle flux $\vec J$ and
becomes
spatially local in the large scale limit of space and time so that $\vec J$ can be approximated by
\begin{equation}
	\vec J(t,\bx)=\int_{ S^{d-1}}\vec \omega f(t,\bx, \vec \omega)\dd \vec\omega.
	\label{eq:J-limit}
\end{equation}
In the numerical computations, %when the mean velocity is calculated over each cell,
$\vec J(t,\bx)$ can be   considered within the approximation \eqref{eq:J-limit}
because the spatial mesh stepsize may be smaller than the radius of  ball $B(0, R)$ \cite{Gamba-Haack2015}.

It has been shown \cite{Gamba-Kang2015} that the non-linear kinetic equation \eqref{eq:kinetic-1} with \eqref{eq:kinetic-2} and \eqref{eq:kr} or \eqref{eq:J-limit}
has a  non-negative global weak solution in $C(0, T;L^1(D)\cap L^\infty\big(( 0, T)\times D\big)$ with $D=\bbR^d\times S^{d-1}$,
for any time $T$, given non-negative initial value $f(0,\bx, \vec \omega)$ in $L^1(D)\cap L^\infty(D)$ and
$\vec J(t,\bx)$ which is always not equal to   $\mathrm{0}$ for $t\in[0, T]$.
%{\color{red} Moreover, the operator $\mathrm{Id}-\vec \omega\otimes\vec \omega$ in \eqref{eq:kinetic-2} enforces
% the size of particle velocity $\vec \omega_k$ is always 1}

Throughout the paper, we will only consider the case of \eqref{eq:J-limit} and
the direction of mean velocity of $f$ becomes % from \eqref{eq:kinetic-2} and \eqref{eq:J-limit}
	\begin{equation*}
		\vec\Omega=\dfrac{\int_{ S^{d-1}} \vec \omega f \dd\vec\omega}{\abs{\int_{ S^{d-1}}
 \vec \omega f \dd\vec \omega}}.
	\end{equation*}
Moreover, the kinetic equation \eqref{eq:kinetic-1} is rewritten as following
\begin{align}\label{eq:kinetic-1-0000}
	  \partial_t f + \vec \omega\cdot \nabla_{\bx} f = Q(f),
\end{align}
where the collision term $Q(f)$ is defined by
\begin{align}
	  Q(f)=-{\nabla_{\vec \omega}}\cdot(\vec F[f]f)+\sigma\Delta_{\vec \omega} f.
\end{align}

\begin{lemma}[\cite{Gamba-Haack2015}]
The collision term $Q(f)$  can be rewritten as follows % Fokker-Planck operator
	\begin{equation}
		Q(f)=\sigma\nabla_{\vec\omega}\cdot\left(M_{\vec \Omega}\nabla_{\vec\omega}\big(\dfrac{f}{M_{\vec\Omega}}\big)\right),
		\label{eq:QfM}
	\end{equation}
and satisfies
	\begin{equation}
	\int_{ S^{d-1}}  	Q(f)f \frac{\dd\vec \omega}{M_{\vec \Omega}}\leq 0,
	\end{equation}
where
\begin{equation}
	M_{\vec \Omega}(\vec \omega)=C_0\exp(\dfrac{\vec\omega\cdot\vec \Omega}{\sigma}),
\label{eq:VM}
\end{equation}
is the equilibrium function, also known as the Von Mises distribution, and $C_0$  is a constant of normalization.
\end{lemma}

\begin{remark}
The equilibria of operator $Q$ are given by the set $\{\rho M_{\vec \Omega}| \rho\in \mathbb R, \vec \Omega\in S^{d-1}\}$
and forms a set of dimension $d$, but the collisional invariants of $Q$ are only of dimension 1.
Particularly, $Q$ does preserve  the mass but cannot preserve the flux, that is, for a general $f$, it holds
$$
\int_{ S^{d-1}}  	Q(f) \dd\vec \omega=0, \quad \int_{S^{d-1}}  	Q(f) \vec \omega \dd\vec \omega\neq 0.
$$
\end{remark}

\begin{remark}
In the  case of $d=2$, if letting $\vec\Omega=(\cos\bt, \sin\bt)^T$,
then the equilibrium can be expressed as follows
\begin{equation}
	M(\theta-\bt)=C_0\exp\left(\dfrac{\cos(\theta-\bt)}{\sigma}\right),
\end{equation}
and the kinetic equation \eqref{eq:kinetic-1-0000} becomes
\begin{equation}
	\label{eq:kinetic2D}
\begin{aligned}
	f_t&+\cos\theta f_x+\sin\theta f_y=Q(f),  \\
	Q(f) =& \sigma\partial_{\theta}\left(
M_{\bar \theta}\partial_{\theta}\big(\dfrac{f}{M_{\bar\theta}}\big)\right)=\partial_\theta(\sin(\theta-\bt)f)+\sigma\partial_\theta^2f,
\end{aligned}
\end{equation}
where \ $\theta$  and $\bt$ denote the
angles of microscopic and macroscopic velocities, respectively.
\end{remark}

\subsection{Macroscopic equations}
The kinetic equation \eqref{eq:kinetic-1-0000} is written at the microscopic level, i.e. at time and length
scales which are characteristic of the dynamics of the individual particles.
When investigating the dynamics of the system at large time and length scales compared with the scales
of the individuals, a set of new   variables  $\tilde{t}=\varepsilon t$ and $\tilde{\bx}=\varepsilon \bx$
has to be introduced \cite{Degond-Motsch2008}, where $\varepsilon$  denotes the ratio between micro and macro variables,
 $\varepsilon\ll 1$.
 In this new set of variables, the kinetic equation \eqref{eq:kinetic-1-0000} is written (after dropping the tildes for clarity) as follows
%If introducing a hydrodynamic scaling
%then after dropping the primes for the sake of clarity, the evolution of $f^\epsilon$ Eq. \eqref{eq:kinetic-1} is governed by
\begin{eqnarray}
	&&\partial_t f^\varepsilon+\vec \omega\cdot{\nabla_{\bx}} f^\varepsilon
	=\dfrac{1}{\varepsilon}Q( f^\varepsilon).
	\label{eq:macroscopic}
\end{eqnarray}
In the limit $\varepsilon\rightarrow0$, \ $f^\varepsilon$ converges locally in space to an
equilibrium state in the local space, that is
\begin{equation}
	f^\varepsilon\rightarrow f(t,\bx, \vec\omega)=\rho(t,\bx)M_{\vec \Omega}(\vec \omega).
\end{equation}
At this moment, the evolution of macroscopic density   $\rho=\int_{S^{d-1}} f\dd\vec \omega$ and
 velocity $\vec \Omega$
is described by the following equations
\begin{align}\label{eq:hongguan-a}
	& \partial_t\rho+\nabla_{\bx}\cdot(c_1\rho \vec \Omega)=0,  \\
	& \rho(\partial_t \vec \Omega + c_2(\vec \Omega\cdot \nabla_{\bx})\vec \Omega)+\lambda (\mathrm{Id}-\vec \Omega\otimes \vec \Omega)\nabla_{\bx}\rho=0, \\
	& \abs{\vec \Omega}=1, \label{eq:hongguan-c}
\end{align}
where $c_1, c_2$, and $\lambda$ are three constants  depending on $\sigma$.
Such macroscopic system is hyperbolic but non-conservative, and the operator
$\mathrm{Id}-\vec \Omega\otimes \vec \Omega$ ensures the constraint \eqref{eq:hongguan-c}.

\begin{remark}
In the 2D case, the macroscopic density $\rho$ and angle of velocity $\bt$
are related to the particle distribution function $f$ by
\begin{subequations}\label{eq:btdingyi}
\begin{eqnarray}
&&\rho(t,\bx)=\int_0^{2\pi}f(t,\bx, \theta)\dd \theta,
    \label{eq:rho}
 \\
&&\vec J(t,\bx)=\int_0^{2\pi}\begin{pmatrix}
		\cos\theta \\ \sin\theta   \end{pmatrix}f(t,\bx, \theta)\dd \theta, \label{eq:btguanxi2}
\\
&&\vec\Omega= (\cos\bt, \sin\bt)^T
	=\dfrac{\vec J(t,\bx)}{\abs{\vec J(t,\bx)}}. \label{eq:btguanxi1}
\end{eqnarray}
\end{subequations}
\end{remark}

%%%%%%%%%%%%%%%%%%%%%%%%%%%%%%%%%%%%%%%%%%%%%%\
%%%%%%%%%%%%%%%
\section{Derivation of moment system}\label{sec:moment}
This section derives the moment system for the 2D kinetic description
of  Vicsek swarming model by using the operator projection \cite{Fan2015Model,Kuang-Tang2016}.
For the sake of simplicity, we assume $\bt=0$
in Sections \ref{ssec:1} and \ref{ssec:2}. In fact, the general case of $\bt\not=0$  can
be converted into such simple case by operating
a translation transformation $\theta\mapsto \theta-\bt$.

\subsection{Orthogonal functions}\label{ssec:1}
%In 1949, Grad [8] assumed that the distribution function is close to a local Maxwellian and
%expanded the distribution function f into Hermite series to obtain the Grad 13 and Grad 20
%moment systems.
The Hermite polynomials are used to derive the Grad's moment system \cite{GRAD:1949,GRAD2:1949}, where
the distribution function is assumed to close to a local Maxwellian.
%一组关于$\exp(-x^2)$正交的Hemite 多项式,
Here we   need to find the orthogonal functions with respect to the weight $M(\theta)$ defined in $[0, 2\pi]$,
because $\theta \in [0, 2\pi]$.
Those  functions, denoted by
$$H^c_0(\theta), H^c_1(\theta), \cdots, H^c_{N}(\theta), \cdots,
H^s_1(\theta), \cdots, H^s_N(\theta), \cdots,$$
are built on the trigonometric functions
$$\{1, \cos\theta, \cdots, \cos(N\theta), \cdots;  \sin\theta, \cdots,
\sin(N\theta), \cdots\}$$
by using the Schmit orthogonal process,
where the superscript $c$ (resp. $s$) denotes the  functions
consisting of a linear combination of $\cos(k\theta)$ (resp.
 $\sin(k\theta)$).
 %It is the result of the parity of the function
 %and will be explained in the following.
Because the function $M(\theta)$ is  even,
it holds
\begin{equation*}
	\int_0^{2\pi}\sin(m\theta)\cos(n\theta)M(\theta)\dd\theta=0, \quad \forall m, n\in\bbN,
\end{equation*}
so in the process of Schimidt orthogonalization,
the linear combination of  $\cos(k\theta)$ is orthogonal to that of $\sin(k\theta)$,
so that there are two sets of ``independent'' orthogonal  functions.
The first $2N+1$  functions are
$H^c_0(\theta), H^c_1(\theta), \cdots, H^c_{N}(\theta)$
and $H^s_1(\theta), \cdots, H^s_N(\theta)$,
and can be expressed in terms of the trigonometric functions
$\{1, \cos\theta, \cdots, \cos(N\theta)$, \ $\sin\theta, \cdots, \sin(N\theta)\}$
as follows
\begin{align} \label{eq:H1}
	&\begin{pmatrix}
		H^c_0(\theta), H^c_1(\theta),  \cdots ,  H^c_N(\theta) \\
	\end{pmatrix}^\mathrm{T}=
		A_{N}
	\begin{pmatrix}
1,  \cos\theta,  \cdots,  \cos(N\theta)\\
	\end{pmatrix}^\mathrm{T},  \\
	&\begin{pmatrix} \label{eq:H2}
		H^s_{1}(\theta) ,  \cdots,  H^s_{N}(\theta) \\
	\end{pmatrix}^\mathrm{T}=
		B_{N}
	\begin{pmatrix}
\sin\theta,  \cdots,  \sin(N\theta) \\
	\end{pmatrix}^\mathrm{T},
\end{align}
where $A_{N}\in\bbR^{(N+1)\times(N+1)}, B_{N}\in\bbR^{N\times N}$,
and both $A_{N}$ and $B_{N}$ are invertible lower triangular matrix.
This set of  functions  satisfy
\begin{equation*}
	\int_0^{2\pi}H^l_m(\theta)H^l_n(\theta)M(\theta)\dd\theta=\delta_{m, n}, \quad m, n\in\{0, 1, \cdots, N\},  l\in\{c, s\}.
\end{equation*}
When $l=s$,  $m, n\not=0$.
For the sake of simplicity, define
\begin{align*}
\begin{pmatrix} 1,  \cos\theta,  \cdots,  \cos(N\theta)\\ \end{pmatrix}^\mathrm{T}\triangleq \bE^c_{N}(\theta),  \ \
&\begin{pmatrix}\sin\theta,  \cdots,  \sin(N\theta) \\ \end{pmatrix}^\mathrm{T}\triangleq \bE^s_{N}(\theta),  \\
\begin{pmatrix} H^c_0(\theta), H^c_1(\theta), \cdots, H^c_{N}(\theta)\\ \end{pmatrix}^\mathrm{T}\triangleq \bH^c_{N}(\theta),  \ \
&\begin{pmatrix} H^s_1(\theta), H^s_2(\theta), \cdots, H^s_{N}(\theta)\\ \end{pmatrix}^\mathrm{T}\triangleq \bH^s_{N}(\theta),  \\
\begin{pmatrix} H^c_0(\theta), H^c_1(\theta), \cdots\\ \end{pmatrix}^\mathrm{T}\triangleq \bH^c(\theta),  \ \
&\begin{pmatrix} H^s_1(\theta), H^s_2(\theta), \cdots\\ \end{pmatrix}^\mathrm{T}\triangleq \bH^s(\theta),
\end{align*}
then \eqref{eq:H1} and \eqref{eq:H2} can be rewritten as follows
\begin{equation}\label{eq:3.3}
	\bH^c_{N}(\theta) = A_N\bE^c_N(\theta), \quad
	\bH^s_{N}(\theta) = B_N\bE^s_N(\theta).
\end{equation}
It is worth noting that the first polynomial is 1, and will be applied to the calculation of density $\rho$.

\begin{remark} The coefficients $A_N$ and $B_N$ in \eqref{eq:3.3}
can be calculated by using
	the following regular modified cylindrical Bessel function of order $n$
	\begin{equation}
    I_n(x)=\dfrac{1}{\pi} \int_0^{\pi} \exp(x\cos\theta)\cos(n\theta)\dd\theta.
	\end{equation}
In fact,	because of the identities
	\begin{align*}
		\sin(m\theta)\sin(n\theta)&=-\dfrac{1}{2}[\cos((m+n)\theta)-\cos((m-n)\theta)], \\
		\cos(m\theta)\cos(n\theta)&=\dfrac{1}{2}[\cos((m+n)\theta)+\cos((m-n)\theta)],
	\end{align*}
 the Schmidt process is carried out with the aid of   $I_n(x)$,
and	then the calculation of $A_N$ and $B_N$ is completed by the existing packages.
\end{remark}

\subsection{Hilbert space and orthonormal basis}\label{ssec:2}
On the interval $[0, 2\pi]$,  define Hilbert space $\mathcal H$ and
inner product with respect to the weight $M(\theta)$ by
\begin{equation}
	{\mathcal H}=\left\{f\Big|\int_0^{2\pi}f^2(\theta)\frac{\dd\theta}{M(\theta)}<\infty\right\},
\end{equation}
\begin{equation}\label{eq:neiji}
\pro{f(\theta)}{g(\theta)}_{M(\theta)}\triangleq \int_0^{2\pi}f(\theta)g(\theta)\dfrac{\dd\theta}{M(\theta)}.
\end{equation}
Such
inner product is symmetric, i.e.
\begin{equation*}
	\pro{f}{g}_{M(\theta)}=\pro{g}{f}_{M(\theta)}.
\end{equation*}
Moreover, due to \eqref{eq:QfM},  the inner product still satisfies
\begin{equation}\label{eq:Qneiji}
	\pro{Q(f)}{g}_{M(\theta)}
	=-\sigma\int_0^{2\pi}M(\theta)\partial_\theta\left(\dfrac{f}{M}\right)\partial_\theta\left(\dfrac{g}{M}\right){\dd \theta}
	=\pro{f}{Q(g)}_{M(\theta)},
\end{equation}
and
\begin{equation}\label{eq:Q<0}
	\pro{Q(f)}{f}_{M(\theta)}
	=-\sigma\int_0^{2\pi}M(\theta)\partial_\theta^2\left(\dfrac{f}{M}\right){\dd \theta}\leqslant 0.
\end{equation}

Take a basis of $\mathcal H$ as $P^c_0(\theta), \cdots, P^c_{N}(\theta)\cdots$, \ $P^s_1(\theta), \cdots, P^s_N(\theta)\cdots$, and denote
\begin{subequations}
\begin{eqnarray}
&&\begin{pmatrix} P^c_0(\theta), P^c_1(\theta), \cdots\\ \end{pmatrix}^\mathrm{T}\triangleq \bP^c(\theta),  \\
&&\begin{pmatrix} P^s_1(\theta), P^s_2(\theta), \cdots\\ \end{pmatrix}^\mathrm{T}\triangleq \bP^s(\theta).
\end{eqnarray}
\end{subequations}
Such basis can be generated by the previous orthogonal functions $\bH^c(\theta)$ and $\bH^s(\theta)$
as follows
\begin{equation}
	\bP^c(\theta)=\bH^c(\theta)M(\theta), \quad
	\bP^s(\theta)=\bH^s(\theta)M(\theta).
\end{equation}

\begin{lemma}
The functions $\bP^c(\theta)$ and $\bP^s(\theta)$  form
 an orthonormal basis of $\mathcal H$, and satisfy the following properties
		\begin{equation}\label{eq:xingzhi}
			\begin{aligned}
				&\pro{\cos(k\theta)M(\theta)}{P^c_n(\theta)}_{M(\theta)}=0, \quad k\leqslant n-1,\\
				&\pro{\sin(k\theta)M(\theta)}{P^s_n(\theta)}_{M(\theta)}=0, \quad k\leqslant n-1,\\
				&\pro{\cos(k\theta)M(\theta)}{P^s_n(\theta)}_{M(\theta)}=0, \quad \forall k, n,   \\
				&\pro{\sin(k\theta)M(\theta)}{P^c_n(\theta)}_{M(\theta)}=0, \quad \forall k, n.
			\end{aligned}
		\end{equation}
\end{lemma}

Because $\{\bP^c, \bP^s\}$ is an orthonormal basis, one has
\begin{equation*}
\mathcal H=span\{P^c_0, P^c_1,  P^c_2,\cdots; P^s_1, \cdots\},
\end{equation*}
and  any function $f\in\mathcal{H}$ may be expressed as follows
\begin{equation}\label{eq:fk}
	f=\sum_{k=0}^\infty f^c_kP^c_k+\sum_{k=1}^\infty f^s_kP^s_k,
\end{equation}
where 	$f^c_k=\pro{f}{P^c_k}$ and $f^s_k=\pro{f}{P^s_k}$.

Take a subspace of $\mathcal H$ as
\begin{equation*}
\mathcal H_N=span\{P^c_0, P^c_1, \cdots, P^c_N; P^s_1, \cdots, P^s_N\}.
\end{equation*}
It is obvious that $\{\bP^c_N(\theta), \bP^s_N(\theta)\}$ forms an orthonormal basis
of $\mathcal H_N$, where
\begin{subequations}
\begin{eqnarray}
&&\begin{pmatrix} P^c_0(\theta), P^c_1(\theta), \cdots,  P^c_N(\theta)\\ \end{pmatrix}^\mathrm{T}\triangleq \bP^c_N(\theta),  \\
&&\begin{pmatrix} P^s_1(\theta), P^s_2(\theta), \cdots,  P^s_N(\theta)\\ \end{pmatrix}^\mathrm{T}\triangleq \bP^s_N(\theta).
\end{eqnarray}
\end{subequations}

For any $f\in \mathcal H$, it is expanded as follows
\begin{equation*}
	f(t, \bx, \theta)=\sum_{k=0}^\infty f^c_k(t, \bx, \bt)P^c_k(\theta-\bt)+\sum_{k=1}^\infty f^s_k(t, \bx, \bt)P^s_k(\theta-\bt),
\end{equation*}
so one can define projection operator $\Pi_N[\bt] : \mathcal H \mapsto \mathcal H_N$ by
\begin{equation}\label{def:pi}
	\Pi_N[\bt]f =\sum_{k=0}^N f^c_k(t, \bx, \bt)P^c_k(\theta-\bt)+\sum_{k=1}^N f^s_k(t, \bx, \bt)P^s_k(\theta-\bt).
\end{equation}
Let $\tilde {\vec f}=\begin{pmatrix}
	f^c_0, \cdots, f^c_N, f^s_1, \cdots, f^s_N
\end{pmatrix}^\mathrm{T}$,
the above equation can be rewritten as follows
\begin{equation}
  \Pi_Nf=\tilde {\vec f}^\mathrm{T}\mcP_N(\theta-\bt).
\end{equation}
In the case of no confusion, the symbol $\bt$ will be ignored in the projection operator, that is,
 $\Pi_N=\Pi_N[\bt]$.

%\subsection{乘速度关系}

In  the following let us  projecting the product of $\cos\theta$ or $\sin\theta$ and orthonormal basis.
Let
\begin{equation}
	\mathcal{\bP}_N(\theta)=
\begin{pmatrix}
	\bP^c_N(\theta) \\
	\bP^s_N(\theta) \\
\end{pmatrix},
\end{equation}
and the $i$th component of  $\mathcal{\bP}_N$  is denoted by $({\mathcal{\bP}_N})_i$,
 	equal to $\bP^c_i(\theta)$, if $i\leqslant N$, otherwise $\bP^s_{i-N}(\theta)$.

\begin{lemma}[Projecting product of velocity and basis]\label{prop:recc}
The result of the operator $\Pi_N$ acting on the product of velocity and basis is
	\begin{subequations}
	\begin{eqnarray}\label{eq:dituix}
    && \Pi_N[\bt]\cos\theta\mcP_N(\theta-\bt)=J^c(\bt)\mcP_{N}(\theta-\bt),  \\
	\label{eq:dituiy}
    && \Pi_N[\bt]\sin\theta\mcP_N(\theta-\bt)=J^s(\bt)\mcP_{N}(\theta-\bt),
\end{eqnarray}
	\end{subequations}
where   $J^c(\bt)\in\bbR^{(2N+1)\times(2N+1)}$ and $J^s(\bt)\in\bbR^{(2N+1)\times(2N+1)}$.
Both matrices $J^c(\bt)$ and $J^s(\bt)$
 are symmetric and thus can be real diagonalizable.
For any eigenvalue $\lambda$,  $\abs{\lambda}\leqslant1$,  and $J^c(\bt), J^s(\bt)$
has the following form
  \begin{equation}
    J^c(\bt)=\cos\bt J_1-\sin\bt J_2, \quad J^s(\bt)=\sin\bt J_1+\cos\bt J_2,
    \label{eq:recc2}
  \end{equation}
where $J_1$ and $J_2$ satisfy
  $$
  \Pi_N[0]\cos\theta\mcP_N(\theta)=J_1\mcP_{N}(\theta), \quad
  \Pi_N[0]\sin\theta\mcP_N(\theta)=J_2\mcP_{N}(\theta).
  $$
\end{lemma}

\begin{proof}
	According to \eqref{eq:fk} and \eqref{def:pi},
the $(i, j)$th component of matrix $J^c(\bt)$ is calculated as follows
	\begin{equation}\label{eq:J^c}
		J^c_{ij}(\bt)=\pro{\cos\theta(\mcP_N)_i(\theta-\bt)}{(\mcP_N)_j(\theta-\bt)}_{M(\theta-\bt)}.
	\end{equation}
Using the definition of inner product \eqref{eq:neiji} gives
	\begin{equation*}
	J^c_{ij}(\bt)=J^c_{ji}(\bt),
	\end{equation*}
so $J^c(\bt)$ is symmetric and can be written as follows
	$J^c(\bt)=\int_0^{2\pi}\cos\theta\mcP_N\mcP_N^{\mathrm T}\dfrac{\dd\theta}{M}$.

Because $\mcP_N$ is an orthonormal basis,
 $\int_0^{2\pi}\mcP_N\mcP_N^{\mathrm T}\dfrac{\dd\theta}{M}=I$.
For any $\lambda\in\bbR$, and non-zero vector $\bx\in\bbR^{2N+1}$,
one has
	\begin{align*}
		\bx^\mathrm{T}(\lambda I-J^c(\bt))\bx=&
		\bx^\mathrm{T}\left(\int_0^{2\pi}(\lambda-\cos\theta)\mcP_N\mcP_N^{\mathrm T}\dfrac{\dd\theta}{M}\right)\bx
		=\int_0^{2\pi}(\lambda-\cos\theta)(\bx^\mathrm{T}\mcP_N)^2\dfrac{\dd\theta}{M}.
	\end{align*}
When $\lambda>1$ or $\lambda<-1$,  $\lambda-\cos\theta>0$ or $\lambda-\cos\theta<0$,
so that the above is greater than 0 or less than 0. Thus the matrix $\lambda I-J^c(\bt)$
is positive definite or negative definite,
	$\lambda I-J^c(\bt)$ is non singular,
so $\lambda>1$ or $\lambda<-1$ is not eigenvalue. Therefore, $\abs{\lambda}\leqslant1$.
For $J^s (\bt) $, the conclusion can be similar  proved.

Because $\cos\theta=\cos\bt\cos(\theta-\bt)-\sin\bt\sin(\theta-\bt)$,
  $\sin\theta=\sin\bt\cos(\theta-\bt)+\cos\bt\sin(\theta-\bt)$,
using the definition of $J_1, J_2$ gives
	\begin{align*}
    J^c_{ij}(\bt)=&\pro{\cos\theta(\mcP_N)_i(\theta-\bt)}{(\mcP_N)_j(\theta-\bt)}_{M(\theta-\bt)} \\
    =&\cos\bt(\pro{\cos\theta(\mcP_N)_i(\theta)}{(\mcP_N)_j(\theta)}_{M(\theta)})
    -\sin\bt(\pro{\sin\theta(\mcP_N)_i(\theta)}{(\mcP_N)_j(\theta)}_{M(\theta)}) \\
    =&\cos\bt (J_1)_{ij}-\sin\bt (J_2)_{ij}.
  \end{align*}
Thus one gets
  \begin{equation*}
    J^c(\bt)=\cos\bt J_1-\sin\bt J_2.
  \end{equation*}
Similarly,  one has
  \begin{equation*}
    J^s(\bt)=\sin\bt J_1+\cos\bt J_2.
  \end{equation*}
\end{proof}

%\subsection{微分关系}
%再考虑微分关系, 即对基函数$\mcP_N(\theta-\bt)$求导后再做投影.

\begin{lemma}[Projecting derivatives of basis] \label{prop:weifen}
The projection of   derivatives of basis function $\mcP_N(\theta-\bt)$
is
	\begin{eqnarray}\label{eq:weifen}
	\Pi_N\dd(\mcP_N(\theta-\bt))=\tilde D\mcP_{N}(\theta-\bt)\dd(\theta-\bt),
\end{eqnarray}
where $\tilde D\in\bbR^{(2N+1)\times(2N+1)}$ is a constant matrix,
and its $(i, j)$th element is
$$-\int_0^{2\pi}{ ({\mcP_N})_i(\theta)}\dd\dfrac{({\mcP_N})_j(\theta)}{M(\theta)}.$$
\end{lemma}

\begin{proof}
Similar to the calculation of \eqref{eq:J^c},  the $(i, j)$th component of matrix $\tilde D$
is calculated as follows
	\begin{align*}
		\tilde D_{ij}=&\pro{\dd({\mcP_N}_i(\theta-\bt))}{{\mcP_N}_j(\theta-\bt)}_{M(\theta-\bt)}\\
		=&\int_0^{2\pi}{(\mcP_N)_j(\theta)}\dfrac{{\dd((\mcP_N)_i(\theta))}}{M(\theta)}
		=-\int_0^{2\pi}{(\mcP_N)_i(\theta)}\dd\dfrac{(\mcP_N)_j(\theta)}{M(\theta)},
	\end{align*}
where the periodicity of the basis functions and $M(\theta)$ has been used in the	second equal sign
while the integration by parts is used in the third equal sign.
\end{proof}

%\subsection{Collision term}
%Substituting the basis functions into the expression of collision term gives
\begin{lemma}[Projecting collision term]\label{prop:coll}
The result of the operator $\Pi_N$ acting on the collision term is
\begin{eqnarray}\label{eq:coll}
	\Pi_NQ(\mcP_N(\theta-\bt))=Q_N\mcP_{N}(\theta-\bt),
\end{eqnarray}
where $Q_N\in\bbR^{(2N+1)\times(2N+1)}$ is a symmetric and negative semidefinite matrix,
and its first column and first row elements are zeros.
\end{lemma}

\begin{proof}
Using the property of inner product \eqref{eq:Qneiji}, one has
	\begin{align*}
		(Q_N)_{ij}=&\pro{Q((\mcP_N)_i(\theta-\bt))}{(\mcP_N)_j(\theta-\bt)}_{M(\theta-\bt)}\\
		=&\pro{(\mcP_N)_i(\theta-\bt)}{Q((\mcP_N)_j(\theta-\bt))}_{M(\theta-\bt)}\\
		=&\pro{Q((\mcP_N)_j(\theta-\bt))}{(\mcP_N)_i(\theta-\bt)}_{M(\theta-\bt)}
		=(Q_N)_{ji}.
	\end{align*}
Because ${\bP^c_0}=M(\theta-\bt)$ and $Q({\bP^c_0})=0$,
the first column and first row elements of $Q_N$ are zeros.

On the other hand, for any non-zero vector $\bx\in\bbR^{2N+1}$, using \eqref{eq:Q<0}, one has
\begin{align*}
	\bx^\mathrm{T}Q_N\bx=&
	\bx^\mathrm{T}\left(\int_0^{2\pi}Q(\mcP_N)\mcP_N^{\mathrm T}\dfrac{\dd\theta}{M}\right)\bx
	=\int_0^{2\pi}Q(\bx^\mathrm{T}\mcP_N)(\bx^\mathrm{T}\mcP_N)\dfrac{\dd\theta}{M} \\
	=&\pro{Q(\bx^\mathrm{T}\mcP_N)}{\bx^\mathrm{T}\mcP_N}\leqslant 0,
\end{align*}
which implies that the matrix $Q_N$ is negative semidefinite.

\end{proof}

\subsection{Moment system}
%现在我们有了$[0, 2\pi]$区间上Hilbert空间$\mathcal H$的一组权函数为$M$的标准正交基,
%以上3节分别考虑了乘速度算子, 微分算子和碰撞项算子作用在这组基上并投影的形式,
%按照\cite{Fan2015Model}所介绍的投影方法, 我们可以完成矩方程组的推导.

%首先将密度分布函数$f$按基函数展开,
%\begin{equation*}
%	f(t, \bx, \theta)=\sum_{k=0}^\infty f^c_k(t, \bx, \bt)P^c_k(\theta-\bt)+\sum_{k=1}^\infty f^c_k(t, \bx, \bt)P^s_k(\theta-\bt),
%\end{equation*}
%做投影,
%\begin{equation*}
%	\Pi_Nf=\sum_{k=0}^Nf^c_kP^c_k(\theta-\bt)+\sum_{k=1}^Nf^c_kP^s_k(\theta-\bt),
%\end{equation*}

The moment system of kinetic equation \eqref{eq:kinetic-1-0000} is derived by the following steps.

\begin{enumerate}
\item[(i)] Calculate the partial derivatives of $\Pi_Nf$ with respect to $t$ and
$x,y$, and then project it onto ${\mathcal H}_N$. Using \eqref{eq:weifen} gives
\begin{subequations}
\begin{align}
	\Pi_N\partial_t\Pi_Nf&=(\pd{\tilde {\vec f}^\mathrm{T}}{t}-\pd{\bt}{t}\tilde {\vec f}^\mathrm{T}\tilde D)\mcP_N(\theta-\bt)
	\triangleq G_1\mcP_N(\theta-\bt),  \\
\label{eq:pianx}
\Pi_N\partial_{x}\Pi_Nf&=(\pd{\tilde {\vec f}^\mathrm{T}}{x}-\pd{\bt}{x}\tilde {\vec f}^\mathrm{T}\tilde D)\mcP_N(\theta-\bt)
\triangleq G_2\mcP_N(\theta-\bt),   \\
\label{eq:piany}
\Pi_N\partial_{y}\Pi_Nf&=(\pd{\tilde {\vec f}^\mathrm{T}}{y}-\pd{\bt}{y}\tilde {\vec f}^\mathrm{T}\tilde D)\mcP_N(\theta-\bt)
	\triangleq G_3\mcP_N(\theta-\bt).
\end{align}
\end{subequations}
\item[(ii)] Multiply \eqref{eq:pianx} and \eqref{eq:piany} by $\cos\theta$ and $\sin\theta$ respectively,
and then project them onto ${\mathcal H}_N$.
	Using \eqref{eq:dituix} and \eqref{eq:dituiy} gives
	\begin{subequations}
\begin{eqnarray}
	\Pi_N\cos\theta\Pi_N\partial_x\Pi_Nf=G_2J^c(\bt)\mcP_N(\theta-\bt),  \\
	\Pi_N\sin\theta\Pi_N\partial_y\Pi_Nf=G_3J^s(\bt)\mcP_N(\theta-\bt).
\end{eqnarray}
\end{subequations}
\item[(iii)]
	Substituting $\Pi_Nf$ into the collision term and projecting it, and using \eqref{eq:coll}
gives
\begin{align}
  \Pi_NQ(\Pi_Nf)&=\tilde {\vec f}^\mathrm{T}Q_N\mcP_N(\theta-\bt).
\end{align}
\item[(iv)] Substituting above all into the
kinetic equation,   and comparing the coefficients of each basis function
gives the following moment system
\begin{equation}\label{eq:moment3}
\begin{aligned}
	\Pi_N[\bt]\partial_t(\Pi_N[\bt]f)
	+& \Pi_N[\bt](\cos\theta\Pi_N[\bt](\partial_x(\Pi_N[\bt]f))) \\
	+& \Pi_N[\bt](\sin\theta\Pi_N[\bt](\partial_y(\Pi_N[\bt]f)))
  = \Pi_N[\bt]Q(\Pi_N[\bt]f),
\end{aligned}
\end{equation}
i.e.
\begin{equation}\label{eq:moment1}
  G_1+G_2J^c(\bt)+G_3J^s(\bt)=\tilde {\vec f}^\mathrm{T}Q_N.
\end{equation}
\end{enumerate}

Because the moment system has $2N+1$ equations but $2N+2$ unknowns
$$\{f^c_0, \cdots, f^s_{N}, \bt\},$$ it needs an additional relationship between $\bt$ and $\{f^c_0, \cdots, f^s_{N}\}$.

\begin{lemma}\label{prop:a0a1}
One has  $\rho=f^c_0$, \ $f^s_1=0$, and $a_0f^c_0+a_1f^c_1>0$,
	where $a_0=\int_0^{2\pi}\cos(\theta-\bt)P^c_0(\theta-\bt)\dd\theta$,
	$a_1=\int_0^{2\pi}\cos(\theta-\bt)P^c_1(\theta-\bt)\dd\theta$.
\end{lemma}

\begin{proof}
Using \eqref{eq:rho} gives
\begin{align*}
  \rho=&\int_0^{2\pi}1\cdot \tilde {\vec f}^\mathrm{T}\mcP_N(\theta-\bt)\dd\theta
  = \int_0^{2\pi}(\mcP_N)_0(\theta-\bt)\cdot \tilde {\vec f}^\mathrm{T}\mcP_N(\theta-\bt)\dfrac{\dd\theta}{M(\theta-\bt)} \\
	=&\sum_{k=0}^{2N}{\tilde f_{k}}\pro{(\mcP_N)_0}{(\mcP_N)_k}_{M(\theta-\bt)}
	=f^c_0.
\end{align*}
where $\tilde f_k$ denotes the $k$th component of vector $\tilde {\vec f}$.

Similarly, using \eqref{eq:btguanxi1}, \eqref{eq:btguanxi2}, and property \eqref{eq:xingzhi} yields
\begin{align*}
	\int_0^{2\pi}\cos\theta f\dd\theta
	=&\cos\bt\sum_{k=0}^{2N}{\tilde f_{k}}\pro{\cos(\theta-\bt)M(\theta-\bt)}{(\mcP_N)_k}_{M(\theta-\bt)} \\
	&-\sin\bt\sum_{k=0}^{2N}{\tilde f_{k}}\pro{\sin(\theta-\bt)M(\theta-\bt)}{(\mcP_N)_k}_{M(\theta-\bt)}\\
	=&\cos\bt(a_0f^c_0+a_1f^c_1)-\sin\bt a_2f^s_1
	= \cos\bt\abs{j(\bx)},
\\
	\int_0^{2\pi}\sin\theta f\dd\theta
	=&\sin\bt(a_0f^c_0+a_1f^c_1)+\cos\bt a_2f^s_1
	=\sin\bt\abs{j(\bx)},
\end{align*}
where
$a_0=\int_0^{2\pi}\cos(\theta-\bt)P^c_0(\theta-\bt)\dd\theta$,
  $a_1=\int_0^{2\pi}\cos(\theta-\bt)P^c_1(\theta-\bt)\dd\theta$,
and $a_2=\int_0^{2\pi}\sin(\theta-\bt)P^s_1(\theta-\bt)\dd\theta$.

 From the above two equations, one can  obtain
$\abs{j(\bx)}=\sqrt{(a_0f^c_0+a_1f^c_1)^2+(a_2f^s_N)^2}$,
thus one has
\begin{align}
	\label{eq:fs1}
	\sqrt{(a_0f^c_0+a_1f^c_1)^2+(a_2f^s_{1})^2}\cos\bt =\cos\bt(a_0f^c_0+a_1f^c_1)-\sin\bt a_2f^s_{1},  \\
	\label{eq:fs2}
	\sqrt{(a_0f^c_0+a_1f^c_1)^2+(a_2f^s_{1})^2}\sin\bt =\sin\bt(a_0f^c_0+a_1f^c_1)+\cos\bt a_2f^s_{1}.
\end{align}
Multiplying
\eqref{eq:fs1} and \eqref{eq:fs2}by $\sin\bt$ and $\cos\bt$ respectively
and subtracting them gives $a_2f^s_1=0$.
Because $a_2$ is larger than zero, $f^s_{1}=0$.
%根据定义计算大于0.
And by hypothesis $\abs{j(\bx)}\not=0$,  one has $a_0f^c_0+a_1f^c_1>0$.
\end{proof}

Based on the above lemma,
the variable $\bt$ can be used to replace $f^s_{1}$, and unknowns in
the system \eqref{eq:moment1} become $\{f^c_0, \cdots, f^c_N, \bt, f^s_{2}, \cdots, f^s_{N}\}$,
thus the variable number  of the moment equations is equal to the equation number.
%我们可以在变量$f^c_0, \cdots, f^s_{N}$ 中消去$f^s_{1}$, 加入变量$\bt$,
Let $\bF=[f^c_0, \cdots, f^c_N, \bt, f^s_{2}, \cdots, f^s_{N}]^{\mathrm T}$,
and rewrite the moment system  \eqref{eq:moment1} in the matrix-vector form
\begin{equation}\label{eq:moment2}
	D\bF_t+ {J^c}D\bF_x+ {J^s}D\bF_y= \tilde Q_N^\mathrm{T}\bF.
\end{equation}
where $D=(I-\vec e_{N+2}\vec e^{\mathrm T}_{N+2}-\tilde{\tilde D}^{\mathrm T}\bF \vec e^{\mathrm T}_{N+2})$,
$\vec e_{N+2}$ denotes the  $(N+2)$th column vector of $2N+1$ order unit matrix,
$\tilde{\tilde D}$ and $\tilde Q_N$ are the matrices defined respectively
by replacing  the $(N+2)$th row of matrices $\tilde D$ and $Q_N$ with zero.
%至此, 矩方程组的形式已经完全给出. 表示将矩阵$\tilde D$第$N+2$行全置为0,

\section{Properties of moment system}\label{sec:4}
This section investigates the mathematical properties of moment system \eqref{eq:moment2}.

\subsection{Hyperbolicity}
In order to prove the hyperbolicity of the moment system \eqref{eq:moment2},
one should show that the matrix $D$ is invertible, and the matrix $\alpha J^c+\beta J^s$
can be really diagonalizable for any $\alpha, \beta\in\bbR$.
\begin{lemma}\label{prop:keni}
	The matrix $D$ is invertible.
\end{lemma}

\begin{proof}
	From the definition of $D$, one has
	\begin{equation*}
		D=\begin{pmatrix}
			1 & 0 & 0 & \cdots & 0 & -\tilde{\tilde D}^\mathrm{T}_1{\vec F} & 0 & \cdots & 0 & 0\\
			0 & 1 & 0 & \cdots & 0 & -\tilde{\tilde D}^\mathrm{T}_2{\vec F} & 0 & \cdots & 0 & 0\\
      \vdots & \vdots &\vdots &\ddots &\vdots          &\vdots &\vdots &\ddots & \vdots & \vdots \\
			0 & 0 & 0 & \cdots & 0 & -\tilde{\tilde D}^\mathrm{T}_{N+2}{\vec F} & 0 & \cdots & 0 & 0\\
      \vdots & \vdots &\vdots &\ddots &\vdots &\vdots &\vdots &\ddots & \vdots & \vdots\\
			0 & 0 & 0 & \cdots & 0 & -\tilde{\tilde D}^\mathrm{T}_{2N+1}{\vec F} & 0 & \cdots & 0 & 1 \\
		\end{pmatrix},
	\end{equation*}
  and	$\det{D}=-\tilde{\tilde D}^\mathrm{T}_{N+2}{\vec F}$,
where $\tilde{\tilde D}^\mathrm{T}_{j}$
denotes the $j$th row of matrix $\tilde{\tilde D}^\mathrm{T}$.
According to Lemma \ref{prop:weifen}, one has
	\begin{align*}
		\tilde D_{i, N+2}=&-\int_0^{2\pi}{(\mcP_N)_i(\theta)}\dd\dfrac{(\mcP_N)_{N+2}(\theta)}{M(\theta)}
		=-\int_0^{2\pi}{(\mcP_N)_i(\theta)}\dd (B_{1, 1}\sin\theta)\\
		=&-b\int_0^{2\pi}{(\mcP_N)_i(\theta)}\cos\theta\dd\theta
		=-b\pro{\cos\theta M(\theta)}{(\mcP_N)_i(\theta)}_{M(\theta)}\\
		=&\begin{cases}
			-ba_0,  &\quad i=0, \\
			-ba_1,  &\quad i=1, \\
			0,  &\quad \text{otherwise}, \\
		\end{cases}
	\end{align*}
where $b\not=0$ is the element of matrix $B_N$ at $(1, 1)$,
and the definition of  $a_0, a_1$ are given in Lemma \ref{prop:a0a1}.
	Thus, \ $\det{D}=b(a_0f^c_0+a_1f^c_1)$. Using Lemma \ref{prop:a0a1} gives
 $\det{D}\not=0$.
\end{proof}

\begin{lemma}\label{prop:duicheng}
	$\alpha J^c+\beta J^s$ can be really diagonalizable  for all   $\alpha, \beta\in\bbR$.
\end{lemma}

\begin{proof}
	Thanks to Lemma \ref{prop:recc},  both matrices $J^c$ and $J^s$
are really symmetric so that their linear combination is really symmetric too,
	and thus really diagonalizable.
\end{proof}

Combing Lemma \ref{prop:keni} with Lemma \ref{prop:duicheng} gives the following
 conclusion.

\begin{lemma}
	The moment system \eqref{eq:moment2} is hyperbolic in time.
\end{lemma}

\subsection{Rotational invariance}
\begin{lemma}
	Under the rotation coordinate transformation
	\begin{align*}
		\begin{pmatrix} x' \\ y' \\ \end{pmatrix}=&
		\begin{pmatrix} \cos\alpha & \sin\alpha \\ -\sin\alpha & \cos\alpha \\ \end{pmatrix}
		\begin{pmatrix} x \\ y \\ \end{pmatrix},  \\
		\theta'=&\theta-\alpha,  \ \
		\bt'= \bt-\alpha,
	\end{align*}
the moment system \eqref{eq:moment2} keeps invariant.
\end{lemma}
\begin{proof}
	Because the moment system \eqref{eq:moment2} is derived by
changing unknowns of \eqref{eq:moment1} and the coordinate transformation is not involved,
our proof will be completed for \eqref{eq:moment1}.

	If using the identity $f(t, x, y, \theta)=f'(t, x', y', \theta')$, and
	$\tilde {\vec f}=(
		f^c_0, \cdots, f^c_N, f^s_1, \cdots, f^s_N)^\mathrm{T}$
and
	$\tilde {\vec f}'=(
		(f^c_0)', \cdots, (f^c_N)', (f^s_1)', \cdots, (f^s_N)')^\mathrm{T}$
to denote the expansion coefficients    of $f(t, x, y, \theta)$ and
$f'(t, x', y', \theta')$, respectively,
then
	\begin{align*}
		&\pd{\tilde {\vec f}^\mathrm{T}}{t}=\pd{\tilde {\vec f'}^\mathrm{T}}{t}, \quad \pd{\bt}{t}=\pd{\bt'}{t},  \
		\pd{\tilde {\vec f}^\mathrm{T}}{x}=\cos\alpha\pd{\tilde {\vec f'}^\mathrm{T}}{x'}
		-\sin\alpha\pd{\tilde {{\vec f}'}^\mathrm{T}}{y'},  \\
		&\pd{\bt}{x}=\cos\alpha\pd{\bt'}{x'}
		-\sin\alpha\pd{\bt'}{y'},  \
		\pd{\tilde {\vec f}^\mathrm{T}}{y}=\sin\alpha\pd{\tilde {\vec f'}^\mathrm{T}}{x'}
		+\cos\alpha\pd{\tilde {{\vec f}'}^\mathrm{T}}{y'},  \\
		&\pd{\bt}{y}=\sin\alpha\pd{\bt'}{x'}
		+\cos\alpha\pd{\bt'}{y'}.
	\end{align*}
	On the other hand, the matrices $J^c(\bt), J^s(\bt), J^c(\bt'), J^s(\bt')$ satisfy
	\begin{align*}
		&J^c(\bt)=\cos\bt J_1-\sin\bt J_2, \quad J^s(\bt)=\sin\bt J_1+\cos\bt J_2,  \\
		&J^c(\bt')=\cos\bt' J_1-\sin\bt' J_2, \quad J^s(\bt')=\sin\bt' J_1+\cos\bt' J_2,
	\end{align*}
	thus one has
	\begin{align*}
		\pd{\tilde {\vec f}^\mathrm{T}}{x}J^c(\bt)+
		\pd{\tilde {\vec f}^\mathrm{T}}{y}J^s(\bt)=&
		\cos\bt'\pd{\tilde {{\vec f}'}^\mathrm{T}}{x}J_1(\bt)-
		\sin\bt'\pd{\tilde {{\vec f}'}^\mathrm{T}}{x}J_2(\bt) \\
		&+\sin\bt'\pd{\tilde {{\vec f}'}^\mathrm{T}}{y}J_1(\bt)+
		\cos\bt'\pd{\tilde {{\vec f}'}^\mathrm{T}}{y}J_2(\bt) \\
		=&\pd{\tilde {{\vec f}'}^\mathrm{T}}{x'}J^c(\bt')+
		\pd{\tilde {{\vec f}'}^\mathrm{T}}{y'}J^s(\bt'),  \\
		\pd{\bt}{x}J^c(\bt)+ \pd{\bt}{y}J^s(\bt)=&
		\pd{\bt'}{x'}J^c(\bt')+
		\pd{\bt'}{y'}J^s(\bt'),  \\
		\tilde {\vec f}^\mathrm{T}Q_N=&\tilde {\vec f'}^\mathrm{T}Q_N,
	\end{align*}
Substituting them into	\eqref{eq:moment1} completes the proof.
\end{proof}

\subsection{Relationship between Grad type expansions  in different $\bt$}
Transformation of a density distribution function under different $\bt$ basis functions.
For the purpose of numerical computations, let us establish the relationship between Grad type expansions of density distribution in different $\bt$.
\begin{lemma}
If assuming that	$\Pi_N[\bt_1]f=\sum_{k=0}^{2N}\tilde f^1_k(\mcP_N)_k(\theta-\bt_1)$,
	$\Pi_N[\bt_2]f=\sum_{k=0}^{2N}\tilde f^2_k(\mcP_N)_k(\theta-\bt_2)$, then it holds
\begin{equation}
	[\tilde f^2_0, \cdots, \tilde f^2_{2N}]^{\mathrm T}=
	T(\bt_1-\bt_2)[\tilde f^1_0, \cdots, \tilde f^1_{2N}]^{\mathrm T},
	\label{eq:trans}
\end{equation}
where
\begin{equation}
	T(\bt)=\begin{pmatrix} A_N & O \\ O & B_N  \end{pmatrix}
	 X(\bt)
	 \begin{pmatrix} A_N^{-1} & O \\ O & B_N^{-1}  \end{pmatrix},
\end{equation}
and
	$$X(\theta)=
	\begin{pmatrix}
		1 & 0 & \cdots & 0 &0 & \cdots & 0 \\
		0 & \cos(\theta) & \cdots &0& -\sin(\theta) & \cdots & 0 \\
		\vdots &\vdots &\ddots&\vdots&\vdots&\ddots&\vdots\\
		0 &0& \cdots & \cos (N\theta) &0& \cdots & -\sin (N\theta) \\
		0 & \sin(\theta) & \cdots     &  0                      & \cos(\theta) & \cdots & 0\\
	 \vdots &    \vdots                             &\ddots   &   \vdots                       &                      \vdots           &\ddots & \vdots\\
		0 &  0                        & \cdots & \sin (N\theta) & 0               &\cdots & \cos (N\theta)
	\end{pmatrix}.   $$
\end{lemma}

\begin{proof}
First estibalish the relationship between basis functions.
 Because
\begin{align*}
	&\bE^c_{N}(\theta-\bt_2)= \\
	&\begin{pmatrix}
		1 & 0 & \cdots & 0&0 & \cdots & 0 \\
		0 & \cos(\bt_1-\bt_2) & \cdots &0& -\sin(\bt_1-\bt_2) & \cdots & 0 \\
		\vdots&\vdots&\ddots&\vdots&\vdots&\ddots&\vdots\\
		0 &0& \cdots & \cos (N(\bt_1-\bt_2)) &0& \cdots & -\sin (N(\bt_1-\bt_2))
	\end{pmatrix} \\
	&\cdot \begin{pmatrix}
		\bE^c_{N}(\theta-\bt_1) \\
		\bE^s_{N}(\theta-\bt_1)
	\end{pmatrix},
\end{align*}
and
\begin{align*}
	&\bE^s_{N}(\theta-\bt_2)= \\
	&\begin{pmatrix}
		0 & \sin(\bt_1-\bt_2) & \cdots     &         0               & \cos(\bt_1-\bt_2) & \cdots &0
\\
	 \vdots &         \vdots                        &\ddots   &     \vdots                     &             \vdots                    &\ddots &\vdots
\\
		0  &  0                       & \cdots & \sin (N(\bt_1-\bt_2)) & 0                & \cdots & \cos (N(\bt_1-\bt_2))
	\end{pmatrix} \\
	&\cdot \begin{pmatrix}
		\bE^c_{N}(\theta-\bt_1) \\
		\bE^s_{N}(\theta-\bt_1)
	\end{pmatrix},
\end{align*}
then one has
\begin{align}
	\nonumber
	\begin{pmatrix}
	\bH_N^c(\theta-\bt_2) \\
	\bH_N^s(\theta-\bt_2) \\
\end{pmatrix}
	=& \begin{pmatrix} A_N & O \\ O & B_N \\ \end{pmatrix}
	 X(\bt_1-\bt_2)
	 \begin{pmatrix} A_N^{-1} & O \\ O & B_N^{-1} \\ \end{pmatrix}
	\begin{pmatrix}
	\bH_N^c(\theta-\bt_1) \\
	\bH_N^s(\theta-\bt_1) \\
\end{pmatrix} \\
	=&T(\bt_1-\bt_2)
	\begin{pmatrix}
	\bH_N^c(\theta-\bt_1) \\
	\bH_N^s(\theta-\bt_1) \\
\end{pmatrix}.
\end{align}
If denoting $\mcP_N(\theta-\bt_1)=\tilde T(\bt_1-\bt_2)\mcP_N(\theta-\bt_2)$,
then one has
\begin{align*}
	\tilde T_{ij}(\bt_1-\bt_2)=&\int_0^{2\pi}(\mcP_N)_i(\theta-\bt_1)(\mcP_N)_j(\theta-\bt_2)/M(\theta-\bt_2)\dd\theta \\
	=&\int_0^{2\pi}(\bH_N)_i(\theta-\bt_1)(\bH_N)_j(\theta-\bt_2)M(\theta-\bt_1)\dd\theta
={T}_{ji}(\bt_1-\bt_2),
\end{align*}
where $({\bH_N})_i$ denotes the $i$th component of $\begin{pmatrix} \bH_N^c \\ \bH_N^s \\ \end{pmatrix}$.
So one has $\tilde T(\bt_1-\bt_2)=T'(\bt_1-\bt_2)$ and
\begin{align*}
\Pi_Nf=&[\tilde f^1_0, \cdots, \tilde f^1_{2N}]{\mcP_N}(\theta-\bt_1) \\
=&[\tilde f^1_0, \cdots, \tilde f^1_{2N}]\tilde T(\bt_1-\bt_2){\mcP_N}(\theta-\bt_2) \\
=&[\tilde f^2_0, \cdots, \tilde f^2_{2N}]{\mcP_N}(\theta-\bt_2).
\end{align*}
Hence one has
\begin{equation*}
	[\tilde f^2_0, \cdots, \tilde f^2_{2N}]^{\mathrm T}=
	T(\bt_1-\bt_2)[\tilde f^1_0, \cdots, \tilde f^1_{2N}]^{\mathrm T}.
\end{equation*}
\end{proof}

\subsection{Mass conservation}
\begin{lemma}
	The moment system \eqref{eq:moment2} preserves the mass-conservation.
\end{lemma}
\begin{proof}
	Because $\rho=f^c_0$,
one just needs to see whether $f^c_0$ is changed.

It is obvious that
the convective term does not change $f^c_0$.
On the other hand, according to Lemma \ref{prop:coll},
the first row and first column of matrix $Q_N$ are zeros,
so the first element of $\tilde Q_N^\mathrm{T}\vec F$ is equal to 0,
and does not change the value of $\bF_0$ or $f^c_0$.
In summary, The moment system \eqref{eq:moment2} is mass-conservative.
\end{proof}

%\subsection{线性稳定性}
%矩方程组\eqref{eq:moment2}显然存在一个局部平衡态解
%$\bF^{(0)}=[\rho_0, 0, \cdots, 0, \bt_0, 0, \cdots, 0]^\mathrm{T}$,
%这里$\rho_0, \theta_0$为常数,
%对于这个矩方程组, 它具有以下线性稳定性.
%\begin{proposition}
	%矩方程组\eqref{eq:moment2}是线性稳定的, 即若
	%\begin{equation}
		%\bF(t, \bx)=\bF^{(0)}\exp(i\omega t-ik_1x-ik_2y),
		%\label{eq:lstable}
	%\end{equation}
	%是\eqref{eq:moment2}的解, 其中$i=\sqrt{-1}, \omega\in\bbC, k\in\bbR, \norm{\bF^{(0)}}_2<\infty$,
	%则$\omega$的虚部不小于0, \ $Im(w)\geqslant 0$.
%\end{proposition}
%\begin{proof}
	%将\eqref{eq:lstable}带入\eqref{eq:moment2}, 得到,
	%\begin{equation*}
		%\bF^{(0)}(i\omega D-ik_1J^c(\bt)D-ik_2J^s(\bt)D-\tilde Q_N^\mathrm{T})|_{\bF^{(0)}}=0,
	%\end{equation*}
	%因为$\bF^{(0)}\not=\b0$, 所以$\omega$是以下方程的解,
	%\begin{equation*}
		%\det(i\omega D-ik_1J^c(\bt)D-ik_2J^s(\bt)D-\tilde Q_N^\mathrm{T})|_{\bF^{(0)}}=0,
	%\end{equation*}
	%由于$D$可逆,
	%\begin{equation*}
		%\det(D)\det(i\omega I-ik_1J^c(\bt)-ik_2J^s(\bt)-\tilde Q_N^\mathrm{T}D^{-1})|_{\bF^{(0)}}=0,
	%\end{equation*}
	%令$U=ik_1J^c(\bt)+ik_2J^s(\bt)+\tilde Q_N^\mathrm{T}D^{-1}$, 则$i\omega$ 是$U$的特征值.
	%要证$Im(\omega)\geqslant 0$, 只需证$Re(\lambda(U))\leqslant 0$,
%\end{proof}

\section{Numerical experiments}\label{sect:numEXP}
This section conducts  numerical experiments to check the behavior of our hyperbolic
moment system \eqref{eq:moment3} or \eqref{eq:moment2}.
%\begin{equation}\label{eq:moment3}
%\begin{aligned}
	%\Pi_N[\bt]\partial_t(\Pi_N[\bt]f)
	%+& \Pi_N[\bt](\cos\theta\Pi_N[\bt](\partial_x(\Pi_N[\bt]f))) \\
	%+& \Pi_N[\bt](\sin\theta\Pi_N[\bt](\partial_y(\Pi_N[\bt]f)))
  %= \Pi_N[\bt]Q(\Pi_N[\bt]f).
%\end{aligned}
%\end{equation}

\subsection{Numerical scheme}
 The spatial grid $\{(x_i,y_j), i,j\in \mathbb N\}$ considered here is uniform
 %$[x_{i-\frac{1}{2}},x_{i+\frac{1}{2}}]$  denotes the $i$th cell,
 so that the stepsizes $\Delta x=x_{i+1}-x_{i}$ and $\Delta y=y_{j+1}-y_j$ are constant.
The grid in $t$-direction $\{t_{n+1}=t_n+\Delta t, n\in \mathbb N\}$
 is also  given with the stepsize  $\Delta t=C_{\mbox{\tiny CFL}} \Delta x$,
where $C_{\mbox{\tiny CFL}} $ denotes the
CFL (Courant-Friedrichs-Lewy) number.  Use $f_{i,j}^n$ and $\bt^n_{i, j}$  to denote the approximations of $f, \bt$ at $t=t_n$ and $(x_i, y_j)$ respectively.
Denote $(\Pi f)^n_{i, j}=\Pi_N[\bt_i^n]f^n_{i, j}$.
 For the purpose of  checking the behavior of our  hyperbolic moment system,
 similar to  \cite{NM:2010},
 we only consider a first-order accurate semi-implicit operator-splitting type  numerical scheme
 for the system \eqref{eq:moment3} or \eqref{eq:moment2},
 which is formed  into the   convection  and collision steps:
\begin{subequations}
\begin{eqnarray}
	\label{scheme:flux_x}
	&&\Pi_N[\bt^n_{i, j}](\Pi f)^{n*}_{i, j}
	= (\Pi f)^n_{i, j}-\dfrac{\Delta t}{\Delta x}[(\Pi F^-)^n_{i+\frac12, j}-(\Pi F^+)^n_{i-\frac12, j}],  \\
	&&\Pi_N[\bt^{n*}_{i, j}](\Pi f)^{n**}_{i, j}
	= (\Pi f)^{n*}_{i, j}-\dfrac{\Delta t}{\Delta y}[(\Pi F^-)^{n*}_{i, j+\frac12}-(\Pi F^+)^{n*}_{i, j-\frac12}],  \\
	%&-&\dfrac{\Delta t}{\Delta y}[(\Pi F^-)^{n}_{i, j+\frac12}-(\Pi F^+)^{n}_{i, j-\frac12}],  \\
	&&\Pi_N[\bt^{n**}_{i, j}]\left(\dfrac{(\Pi f)^{n+1}_{i, j}-(\Pi f)^{n**}_{i, j}}{\Delta t}\right)
  = \Pi_N[\bt^{n**}_{i, j}]Q((\Pi f)^{n+1}_{i, j}), \label{eq:collQ}
\end{eqnarray}
\end{subequations}
where the numerical fluxes are chosen as the nonconservative HLL flux \cite{Rhebergen2008Discontinuous}.
As an example, the flux in $x$ direction can be expressed as follows
{\small \begin{align*}
	&(\Pi F^-)^n_{i+\frac12, j}= \\
	&\begin{cases}
		\Pi_N[\bt^n_{i, j}](\cos\theta(\Pi f)^n_{i, j}),  & 0\leqslant\lambda^L_{i+\frac12, j}, \\
		\dfrac{\lambda^R_{i+\frac12, j}\Pi_N[\bt^n_{i, j}](\cos\theta(\Pi f)^n_{i, j})
		-\lambda^L_{i+\frac12, j}\Pi_N[\bt^n_{i, j}](\cos\theta\Pi_N[\bt^n_{i, j}](\Pi f)^{n}_{i+1, j})}
		{\lambda^R_{i+\frac12, j}-\lambda^L_{i+\frac12, j}} \\
		+\dfrac{\lambda^L_{i+\frac12, j}\lambda^R_{i+\frac12, j}(\Pi_N[\bt^n_{i, j}](\Pi f)^{n}_{i+1, j}-(\Pi f)^n_{i, j})}
    {\lambda^R_{i+\frac12, j}-\lambda^L_{i+\frac12, j}},  & \lambda^L_{i+\frac12, j}<0<\lambda^R_{i+\frac12, j}, \\
		\Pi_N[\bt^n_{i, j}](\cos\theta\Pi_N[\bt^n_{i, j}](\Pi f)^{n}_{i+1, j}),  & 0\geqslant\lambda^R_{i+\frac12, j}, \\
	\end{cases}
\end{align*}
}
and
{\small  \begin{align*}
	&(\Pi F^+)^n_{i-\frac12, j}= \\
	&\begin{cases}
		\Pi_N[\bt^n_{i, j}](\cos\theta\Pi_N[\bt^n_{i, j}](\Pi f)^{n}_{i-1, j}),  & 0\leqslant\lambda^L_{i-\frac12, j}, \\
		\dfrac{\lambda^R_{i-\frac12, j}\Pi_N[\bt^n_{i, j}](\cos\theta\Pi_N[\bt^n_{i, j}](\Pi f)^{n}_{i-1, j})
		-\lambda^L_{i-\frac12, j}\Pi_N[\bt^n_{i, j}](\cos\theta(\Pi f)^n_{i, j})}
		{\lambda^R_{i-\frac12, j}-\lambda^L_{i-\frac12, j}} \\
		+\dfrac{\lambda^L_{i-\frac12, j}\lambda^R_{i-\frac12, j}((\Pi f)^{n}_{i, j}-\Pi_N[\bt^n_{i, j}](\Pi f)^n_{i-1, j})}
		{\lambda^R_{i-\frac12, j}-\lambda^L_{i-\frac12, j}},  & \lambda^L_{i-\frac12, j}<0<\lambda^R_{i-\frac12, j}, \\
		\Pi_N[\bt^n_{i, j}](\cos\theta(\Pi f)^n_{i, j}),  & 0\geqslant\lambda^R_{i-\frac12, j}, \\
	\end{cases}
\end{align*}
}
where $\lambda^L_{i\pm\frac12, j}=\min\{\lambda^{\mathrm {min}}_{i, j}, \lambda^{\mathrm {min}}_{i\pm 1, j}\}$,
$\lambda^R_{i\pm\frac12, j}=\max\{\lambda^{\mathrm {max}}_{i, j}, \lambda^{\mathrm {max}}_{i\pm 1, j}\}$,
$\lambda^{\mathrm {min}}_{i, j}$ and $\lambda^{\mathrm {max}}_{i, j}$
denotes the minimum and maximum eigenvalues of $J^c(\bt)$ at $(i, j)$.

\begin{lemma}
For any $\bt_1, \bt_2$, it holds
	\begin{equation}
		\Pi_N[\bt_1]f=\Pi_N[\bt_1]\Pi_N[\bt_2]f.
	\end{equation}
\end{lemma}

\begin{proof}
If assuming
 $$\Pi_N[\bt_1]f=\sum_{k=0}^{2N}\tilde f^1_k(\mcP_N)_k(\theta-\bt_1), \
\Pi_N[\bt_2]f=\sum_{k=0}^{2N}\tilde f^2_k(\mcP_N)_k(\theta-\bt_2),$$
then using \eqref{eq:trans} gives
\begin{equation*}
	[\tilde f^2_0, \cdots, \tilde f^2_{2N}]^{\mathrm T}=
	T(\bt_1-\bt_2)[\tilde f^1_0, \cdots, \tilde f^1_{2N}]^{\mathrm T}.
\end{equation*}
The term $\Pi_N[\bt_1]\Pi_N[\bt_2]f$ is equivalent to
transferring  the expansion $\Pi_N[\bt_2]f$ in the basis with the parameter $\bt_2$
to the basis with the parameter $\bt_1$,
so the coefficients of expansion  $\Pi_N[\bt_1]\Pi_N[\bt_2]f$ are
\begin{equation*}
	T(\bt_2-\bt_1)[\tilde f^2_0, \cdots, \tilde f^2_{2N}]^{\mathrm T}
	= T(\bt_2-\bt_1)T(\bt_1-\bt_2)[\tilde f^1_0, \cdots, \tilde f^1_{2N}]^{\mathrm T}
	= [\tilde f^1_0, \cdots, \tilde f^1_{2N}]^{\mathrm T},
\end{equation*}
which are the  coefficients of expansion $\Pi_N[\bt_1]f$.
\end{proof}
The above lemma tell us that  in \eqref{scheme:flux_x},
$\Pi_N[\bt^{n*}_{i, j}]=\Pi_N[\bt^{n*}_{i, j}]\Pi_N[\bt^n_{i, j}]$.
Hence, it can first calculate the coefficients of  $f_{i, j}$ in $\Pi_N[\bt^n_{i, j}]$,
then uses the definition to get the value of $\bt^{n*}$. Finally,
 the transition matrix \eqref{eq:trans} for the projection  in different $\bt$
is used to get  the coefficients of  $f_{i, j}$ in $\Pi_N[\bt^{n*}_{i, j}]$.

From what has been discussed above, the following steps are required to complete the numerical scheme.
\begin{enumerate}
	\item[(a)]	solve the $x$-convective step to get $\Pi_N[\bt^n_{i, j}](\Pi f)^{n*}_{i, j}$.
	\item[(b)]  use the definition to calculate $\bt^{n*}$ and give $(\Pi f)^{n*}_{i, j}$.
	\item[(c)]	solve the $y$-convective step to  obtain $\Pi_N[\bt^n_{i, j}](\Pi f)^{n**}_{i, j}$.
	\item[(d)] use the definition to calculate $\bt^{n**}$ and get $(\Pi f)^{n**}_{i, j}$.
	\item[(e)]	solve the collision step  to yield $\Pi_N[\bt^n_{i, j}](\Pi f)^{n+1}_{i, j}$.
	\item[(f)] calculate $\bt^{n+1}$ and $(\Pi f)^{n+1}_{i, j}$.  Set $n=n+1$ and goto (a).
\end{enumerate}

For the collision step, an implicit scheme is used.
Substituting the matrix $Q_N$ into \eqref{eq:collQ}
%\begin{align*}
	%\Pi_N[\bt^{n*}_{i, j}](\Pi f)^{n+1}_{i, j}-\Pi_N[\bt^{n*}_{i, j}](\Pi f)^{n*}_{i, j}
	%= \Delta t\Pi_N[\bt^{n*}_{i, j}]Q((\Pi f)^{n+1}_{i, j})
	%= \Delta t\Pi_N[\bt^{n*}_{i, j}]Q_N(\Pi f)^{n+1}_{i, j},
%\end{align*}
gives
\begin{eqnarray}
	\Pi_N[\bt^{n**}_{i, j}](\Pi f)^{n+1}_{i, j}
	&=& \Pi_N[\bt^{n**}_{i, j}](I-\Delta tQ_N)^{-1}(\Pi f)^{n**}_{i, j}.
\end{eqnarray}

\begin{lemma}
The implicit discretization for the collision step is unconditionally stable.
\end{lemma}

\begin{proof}
	Because the matrix $Q_N$ is semi negative definite,
 any eigenvalue $\lambda$ of $Q_N$ is not larger than $0$.
 On the other hand,
 because $1-\Delta t\lambda$ is the eigenvalue of	$I-\Delta tQ_N$,
 the eigenvalues of $(I-\Delta tQ_N)^{-1}$ satisfy
	$0\leqslant(1-\Delta t\lambda)^{-1}\leqslant 1$.
Therefore, the implicit scheme for the collision step is unconditionally stable.
\end{proof}

%%%%%%%%%%%%%%%%%%%%%%%%%%%%%%%%%%%%%%%%%%%%%
\subsection{Treatment of
reflection boundary}
Here the treatment of reflection boundaries is similar to the upwind scheme.
Take the left boundary in $x$ direction  as an example to illustrate our treatment of
reflection boundary.

For all $j$, when $\cos\theta\leqslant 0$, the  macroscopic numerical flux $\tilde F$
at the left boundary can be defined by $(\Pi f)_{0, j}$,
while when $\cos\theta>0$, it is defined by reflection of $(\Pi f)_{0, j}$, that is
\begin{equation}
	\tilde F=\begin{cases}
		\cos\theta (\Pi f)_{0, j}(\theta), & \quad \cos\theta\leqslant 0,  \\
		\cos\theta (\Pi f)_{0, j}(\pi-\theta), & \quad \cos\theta>0.
	\end{cases}
\end{equation}
If expanding $\tilde F$ in $\mcP_N(\theta-\bt_{0, j})$,
then its $k$th component is
\begin{align*}
	\tilde F_k=&\int_{\cos\theta\leqslant 0}\cos\theta (\Pi f)_{0, j}(\theta)(\mcP_N)_k(\theta-\bt_{0, j})\dfrac{\dd \theta}{M(\theta-\bt_{0, j})} \\
	&+\int_{\cos\theta>0}\cos\theta (\Pi f)_{0, j}(\pi-\theta)(\mcP_N)_k(\theta-\bt_{0, j})\dfrac{\dd \theta}{M(\theta-\bt_{0, j})}.
\end{align*}
Therefore the numerical flux at the left reflection boundary
in $x$ direction is given by
\begin{equation}
  (\Pi \hat F^{+})_{-\frac12, j}=\sum_{k=0}^{2N}\tilde F_k(\mcP_N)_k(\theta-\bt_{0, j}).
\end{equation}

\subsection{Numerical results}
The 2D scheme of moment system is  used to solve three Riemann problems
and a vortex formation problem and will be compared to the spectral method proposed in \cite{Gamba-Haack2015}.
For our computations, the computational domain of Riemann problems in $x$-direction is taken as $[-5, 5]$,
 the CFL number is chosen as $0.5$,  and the value of parameter $\sigma$ is set as 0.2.

\begin{figure}
  \centering
\subfigure[$N=1$]{
\includegraphics[width=0.3\textwidth]{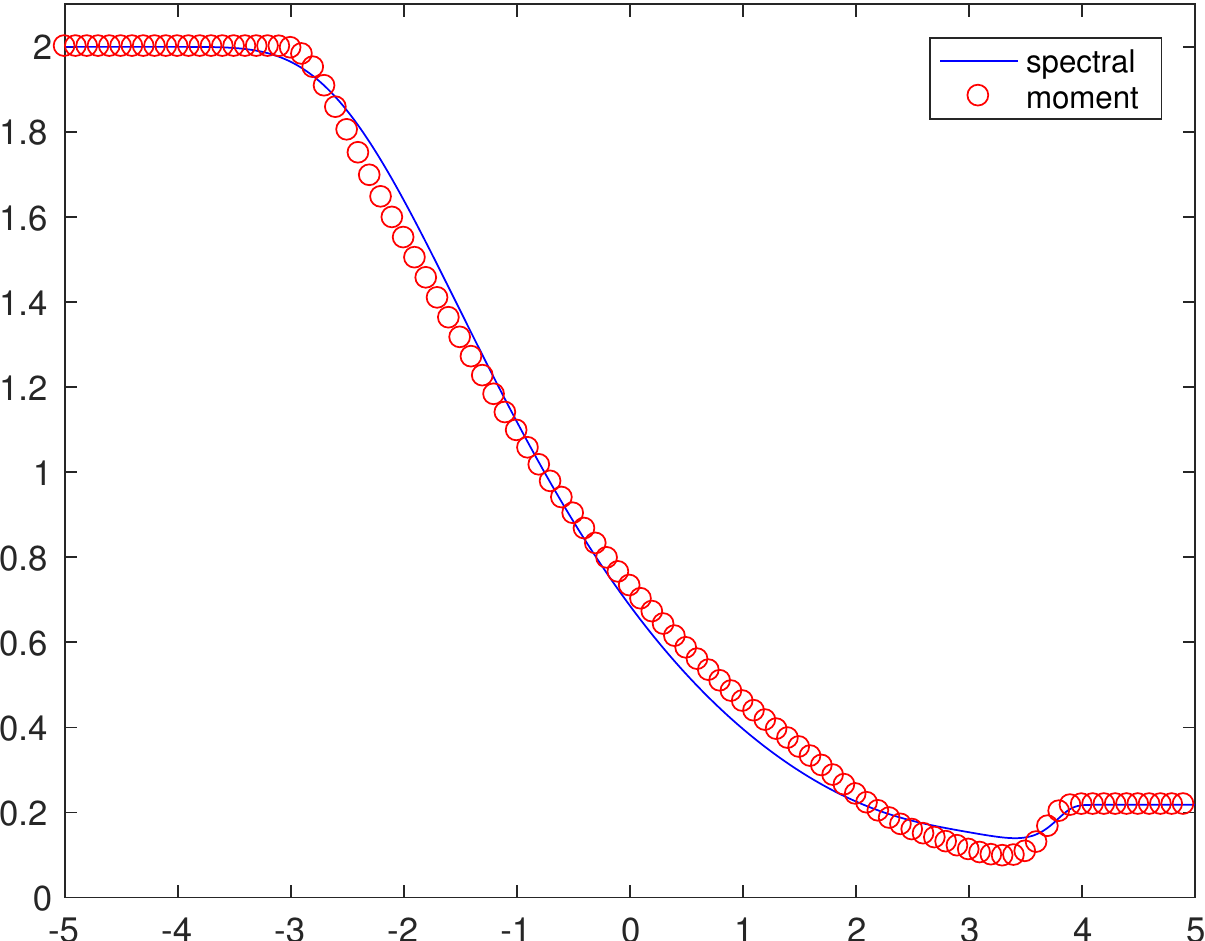}
}
\subfigure[$N=2$]{
\includegraphics[width=0.3\textwidth]{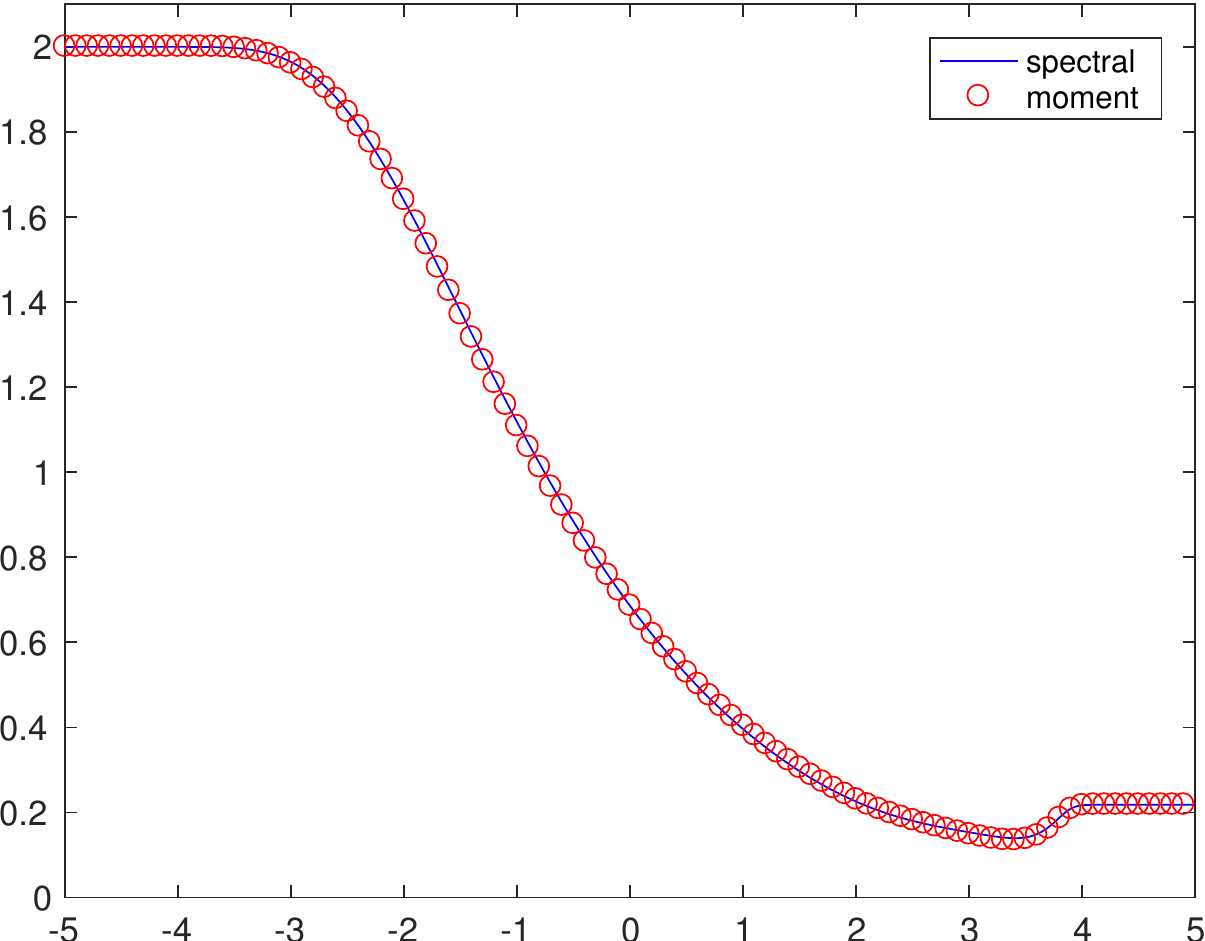}
}
\subfigure[$N=3$]{
\includegraphics[width=0.3\textwidth]{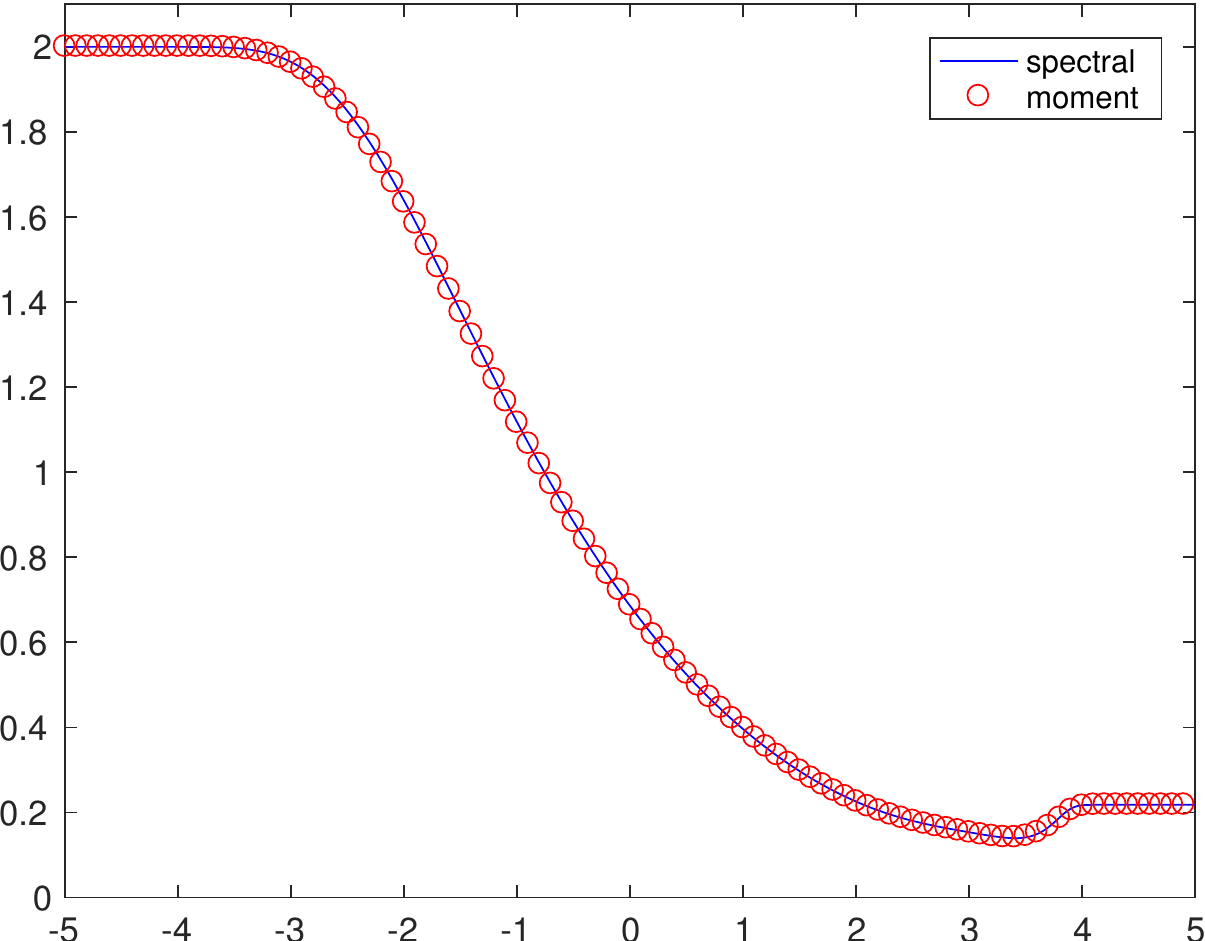}
}

  \centering
\subfigure[$N=4$]{
\includegraphics[width=0.3\textwidth]{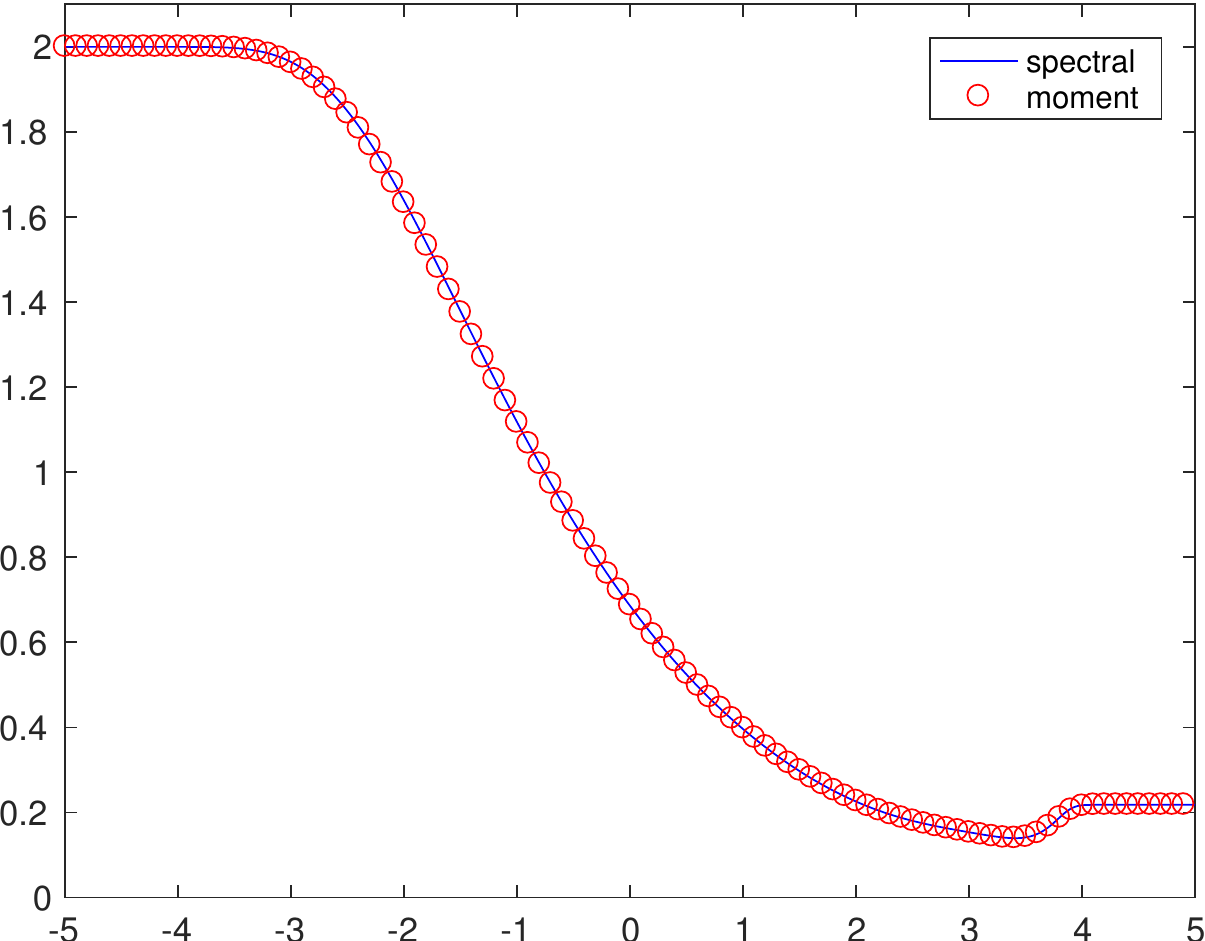}
}
\subfigure[$N=5$]{
\includegraphics[width=0.3\textwidth]{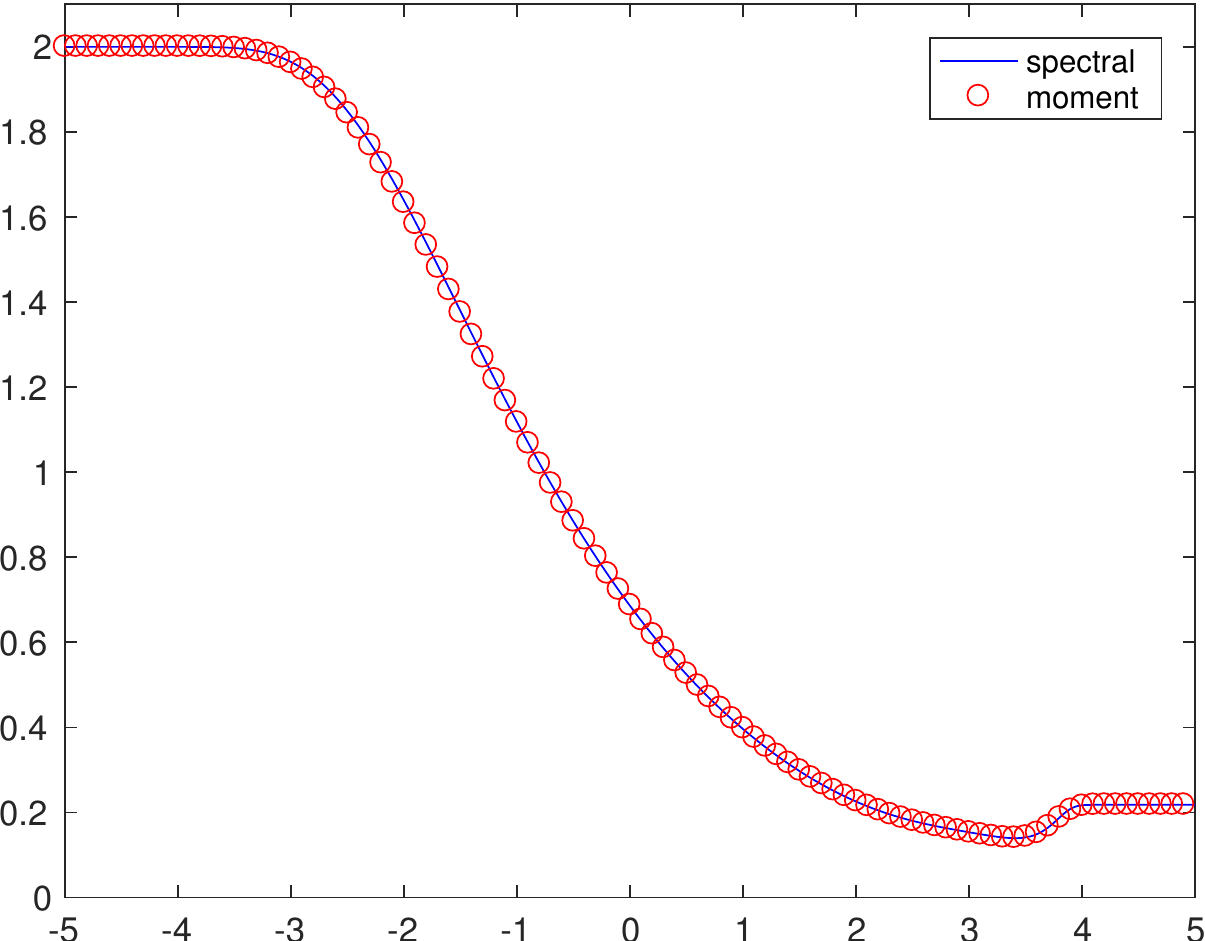}
}
\subfigure[$N=6$]{
\includegraphics[width=0.3\textwidth]{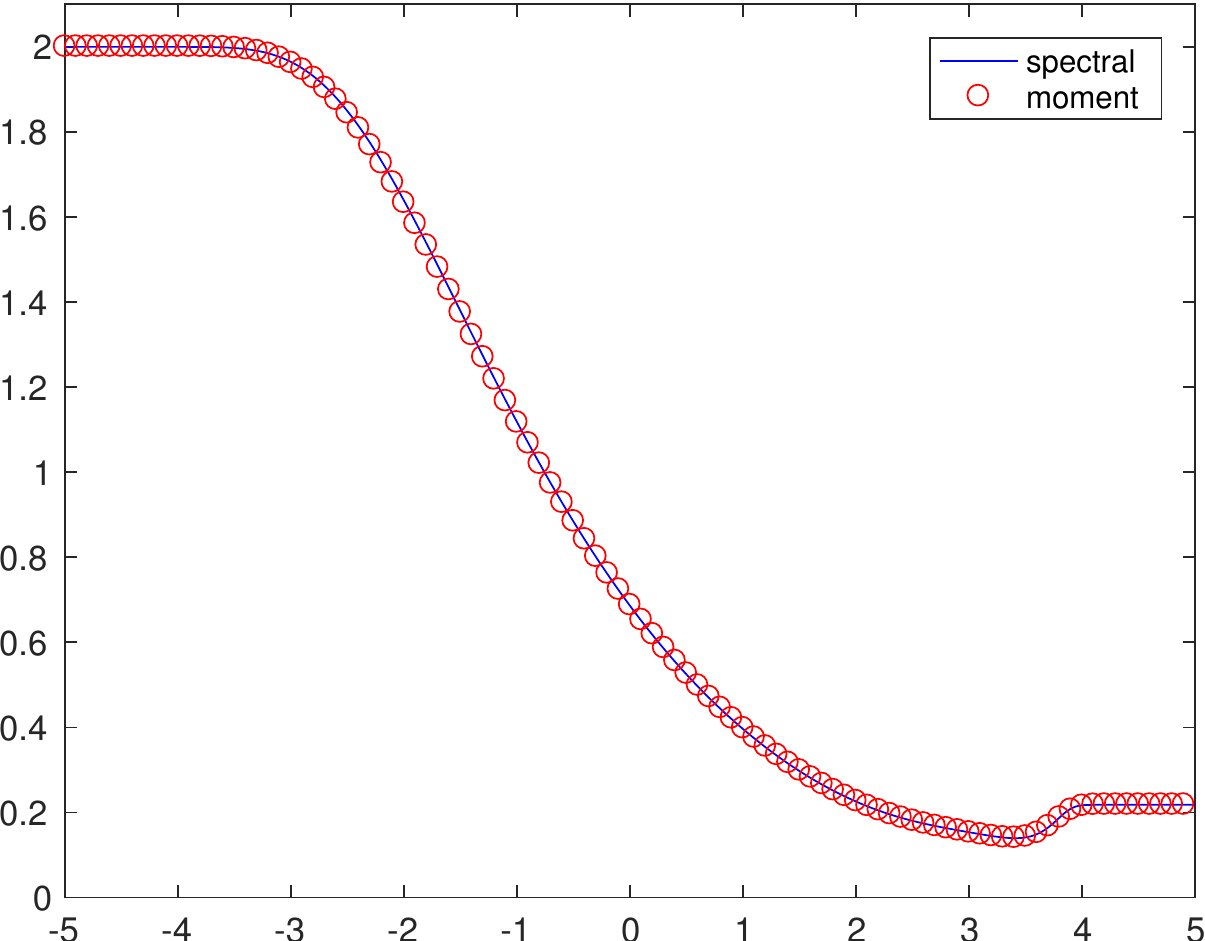}
}

%  \centering
%\subfigure[$N=7$]{
%\includegraphics[width=0.3\textwidth]{./images/riemann1/cmpN_e1000/r1_cmpN_rho_n2000_e1000_N7_circle.pdf}
%}
%\subfigure[$N=8$]{
%\includegraphics[width=0.3\textwidth]{./images/riemann1/cmpN_e1000/r1_cmpN_rho_n2000_e1000_N8_circle.pdf}
%}
%\subfigure[$N=9$]{
%\includegraphics[width=0.3\textwidth]{./images/riemann1/cmpN_e1000/r1_cmpN_rho_n2000_e1000_N9_circle.pdf}
%}
\caption{Example \ref{example:5.1}: The densities at $t=4$ obtained by the moment
method with $N=1,2,\cdots,6$ and 2000 cells. The solid line is the reference solution obtained by using the
spectral method with 4000 cells.  $\varepsilon=1$.}\label{figure:5.1}
\end{figure}

\begin{figure}
  \centering
\subfigure[$N=1$]{
\includegraphics[width=0.3\textwidth]{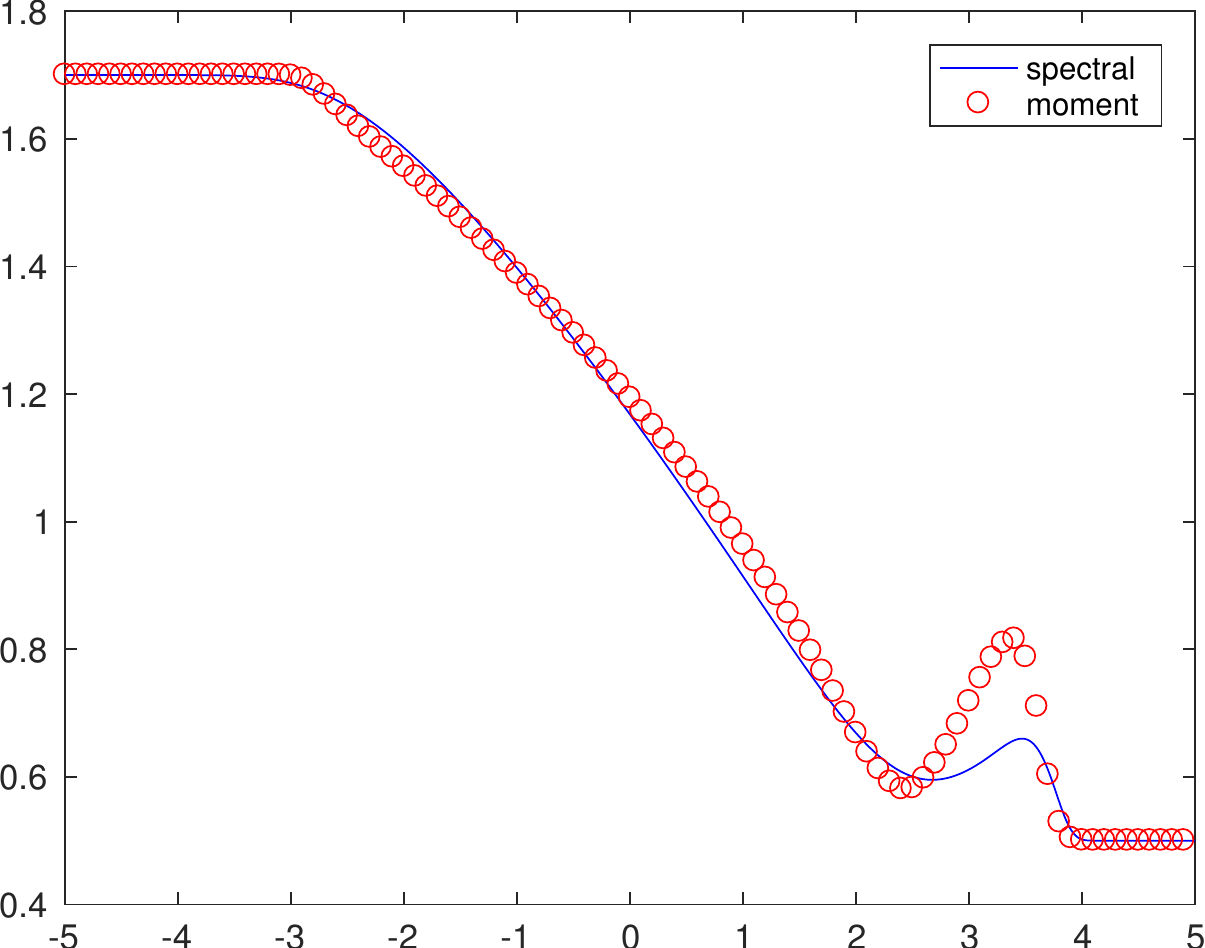}
}
\subfigure[$N=2$]{
\includegraphics[width=0.3\textwidth]{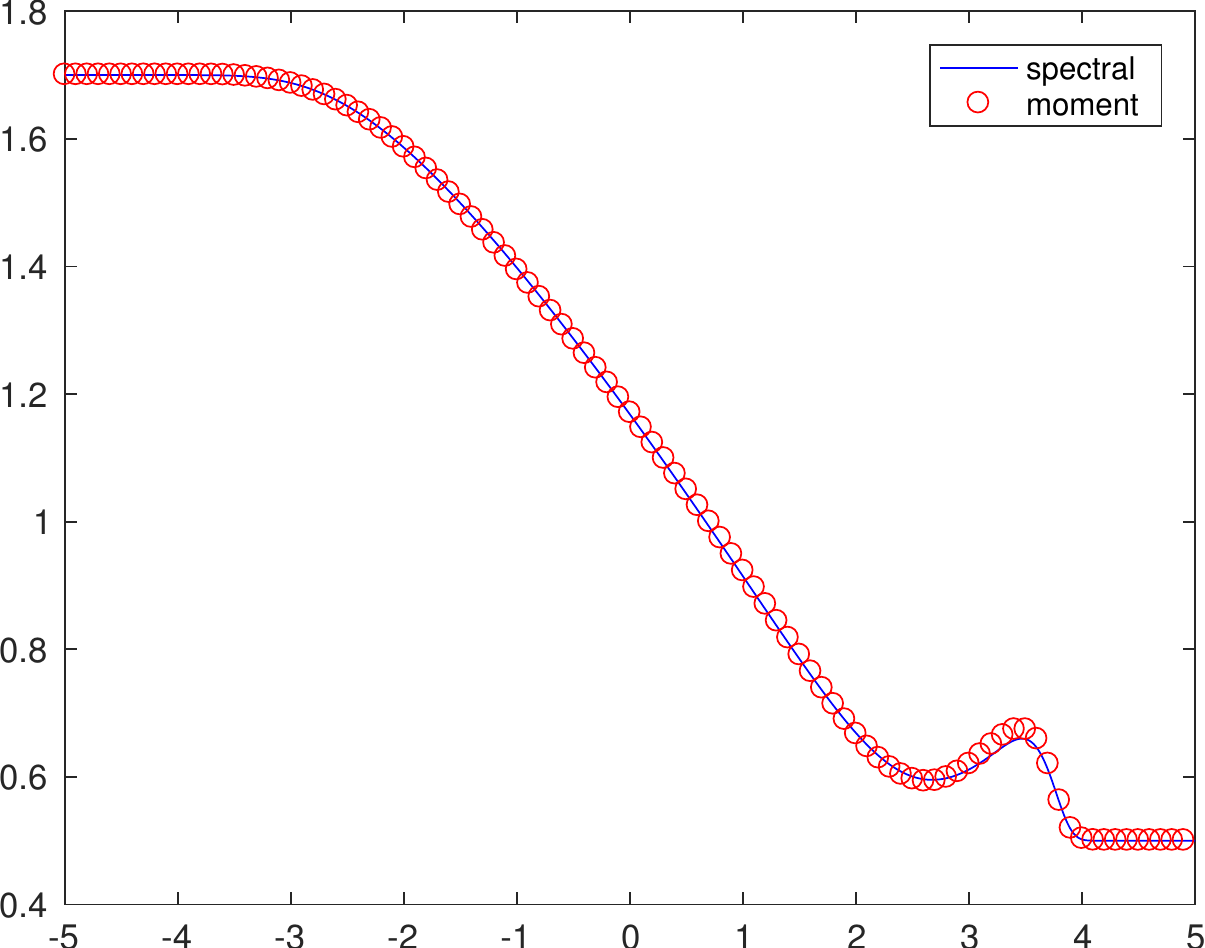}
}
\subfigure[$N=3$]{
\includegraphics[width=0.3\textwidth]{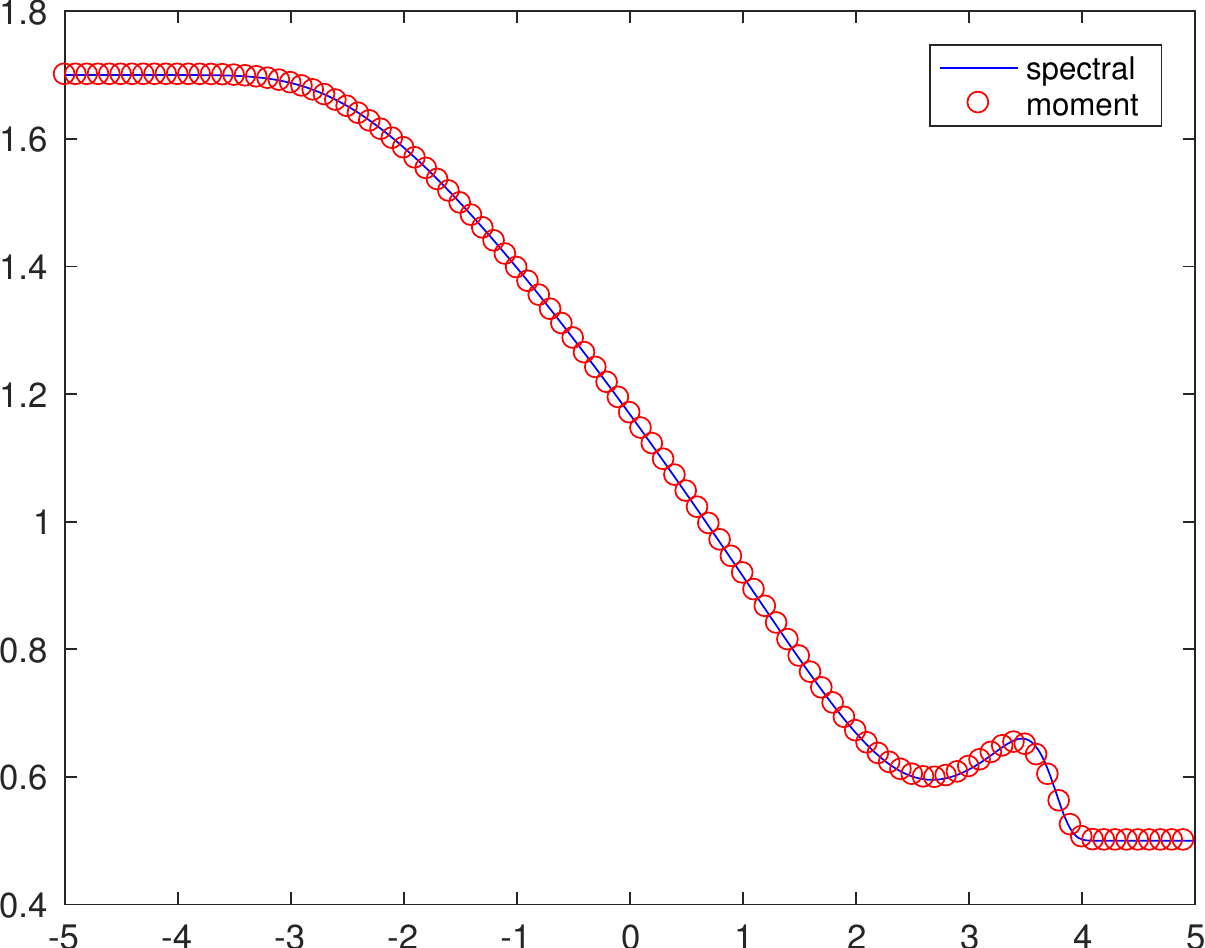}
}

  \centering
\subfigure[$N=4$]{
\includegraphics[width=0.3\textwidth]{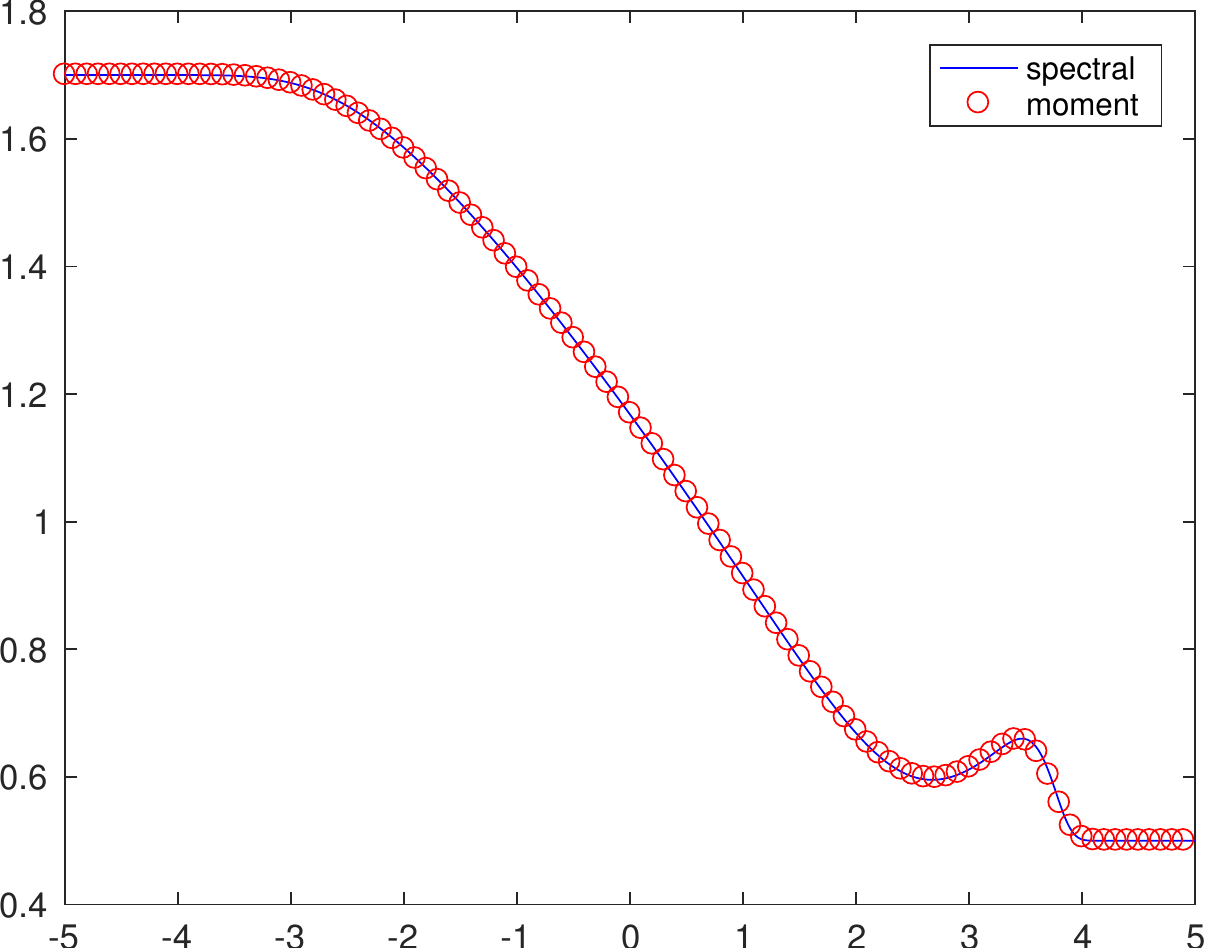}
}
\subfigure[$N=5$]{
\includegraphics[width=0.3\textwidth]{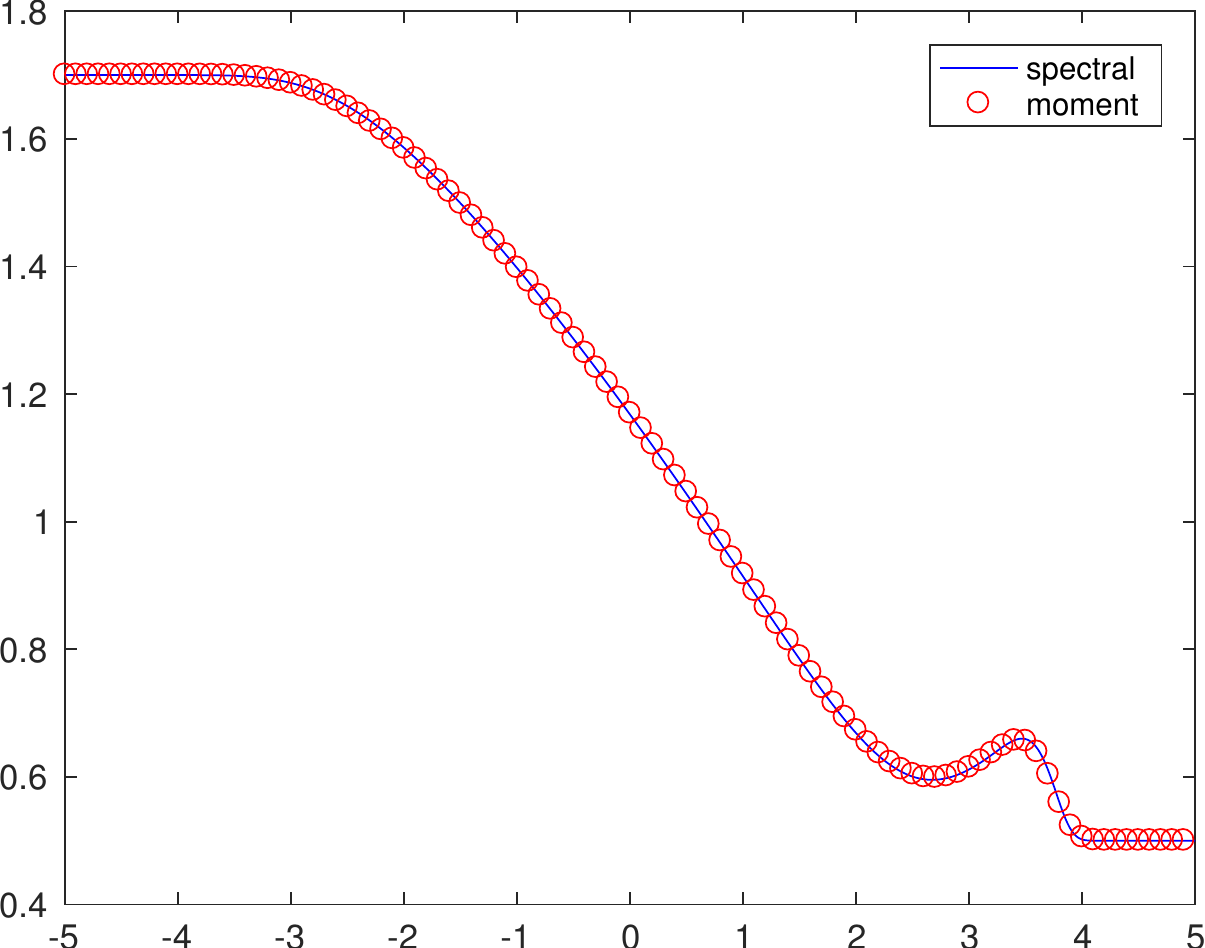}
}
\subfigure[$N=6$]{
\includegraphics[width=0.3\textwidth]{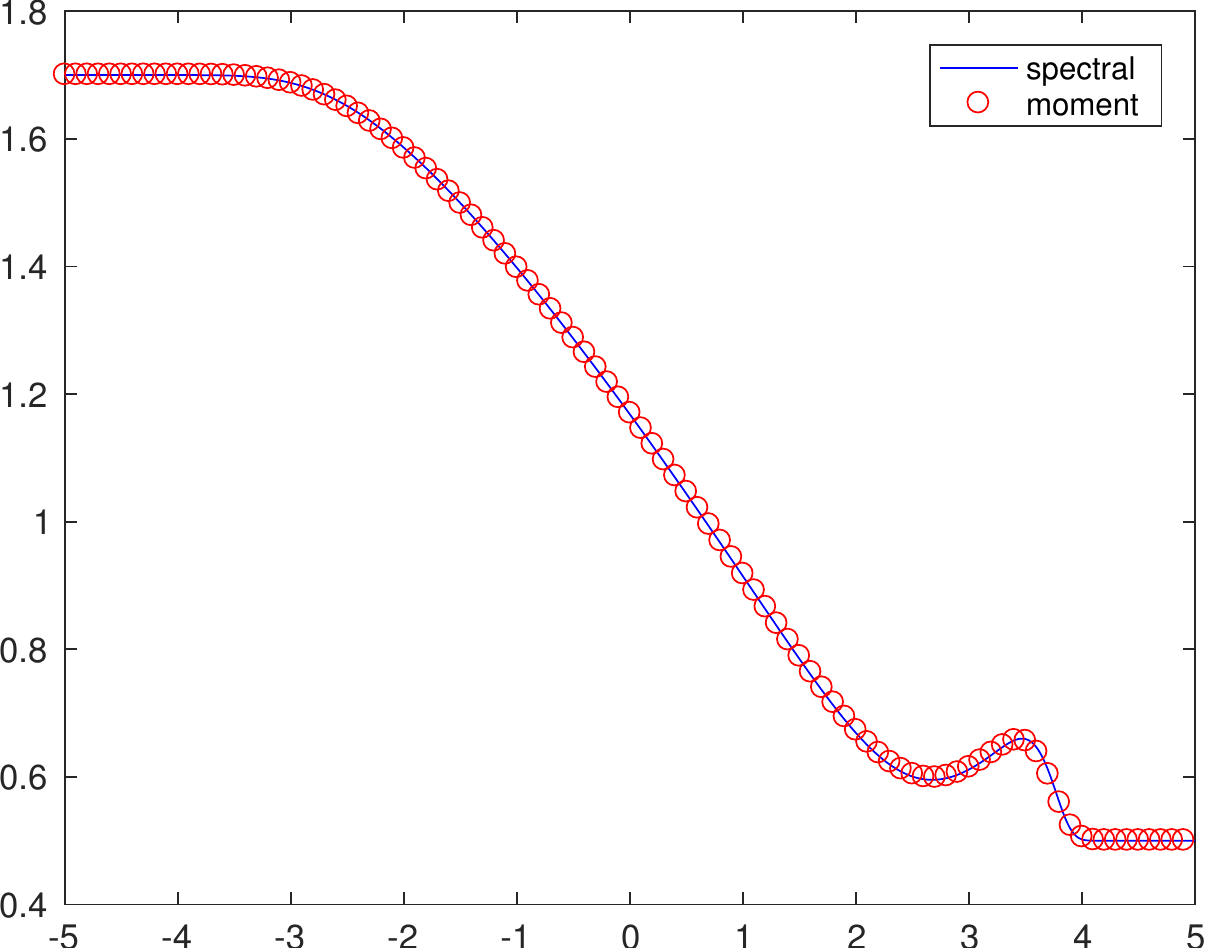}
}

%  \centering
%\subfigure[$N=7$]{
%\includegraphics[width=0.3\textwidth]{./images/riemann1/cmpN_e1000/r1_cmpN_u_n2000_e1000_N7_circle.pdf}
%}
%\subfigure[$N=8$]{
%\includegraphics[width=0.3\textwidth]{./images/riemann1/cmpN_e1000/r1_cmpN_u_n2000_e1000_N8_circle.pdf}
%}
%\subfigure[$N=9$]{
%\includegraphics[width=0.3\textwidth]{./images/riemann1/cmpN_e1000/r1_cmpN_u_n2000_e1000_N9_circle.pdf}
%}
\caption{Same as Fig. \ref{figure:5.1} except for the macroscopic velocity angles.}\label{figure:5.2}
\end{figure}

\begin{figure}
  \centering
\subfigure[$N=1$]{
\includegraphics[width=0.3\textwidth]{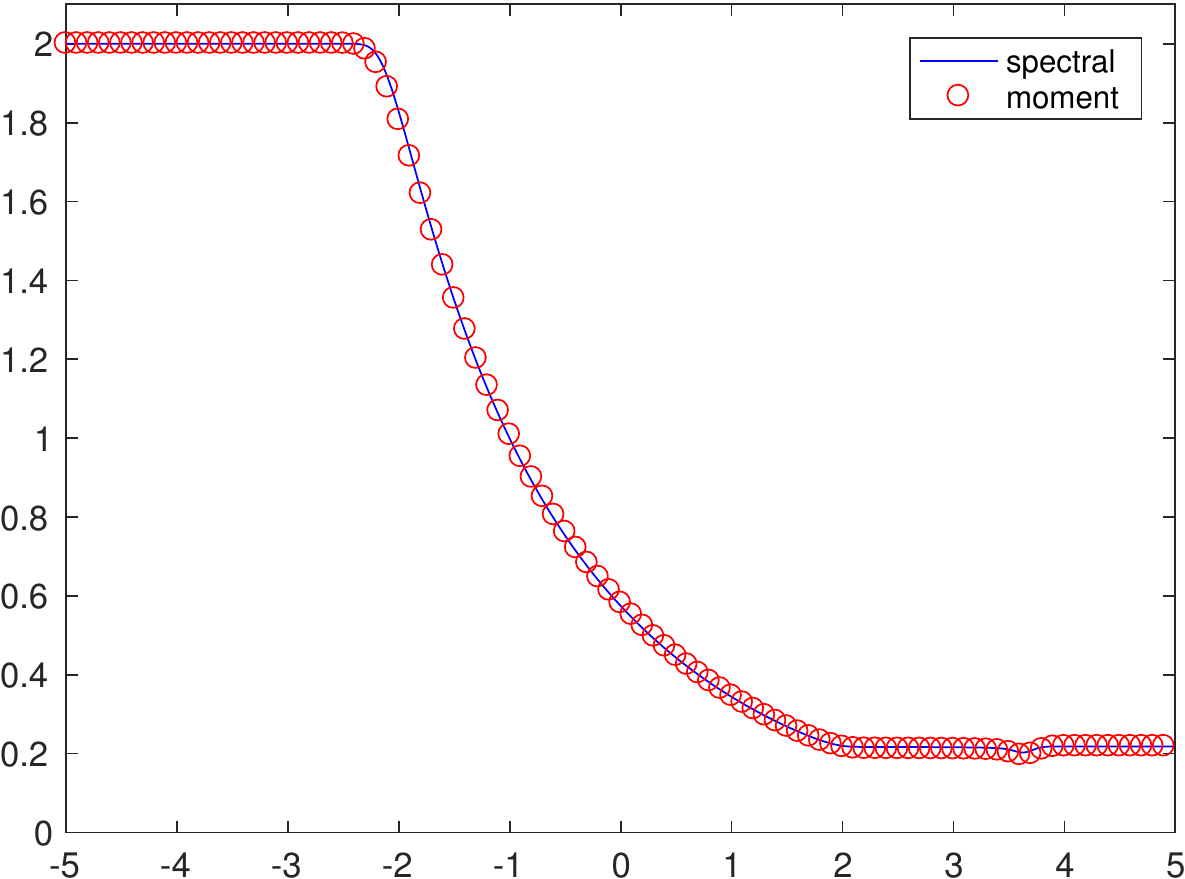}
}
\subfigure[$N=2$]{
\includegraphics[width=0.3\textwidth]{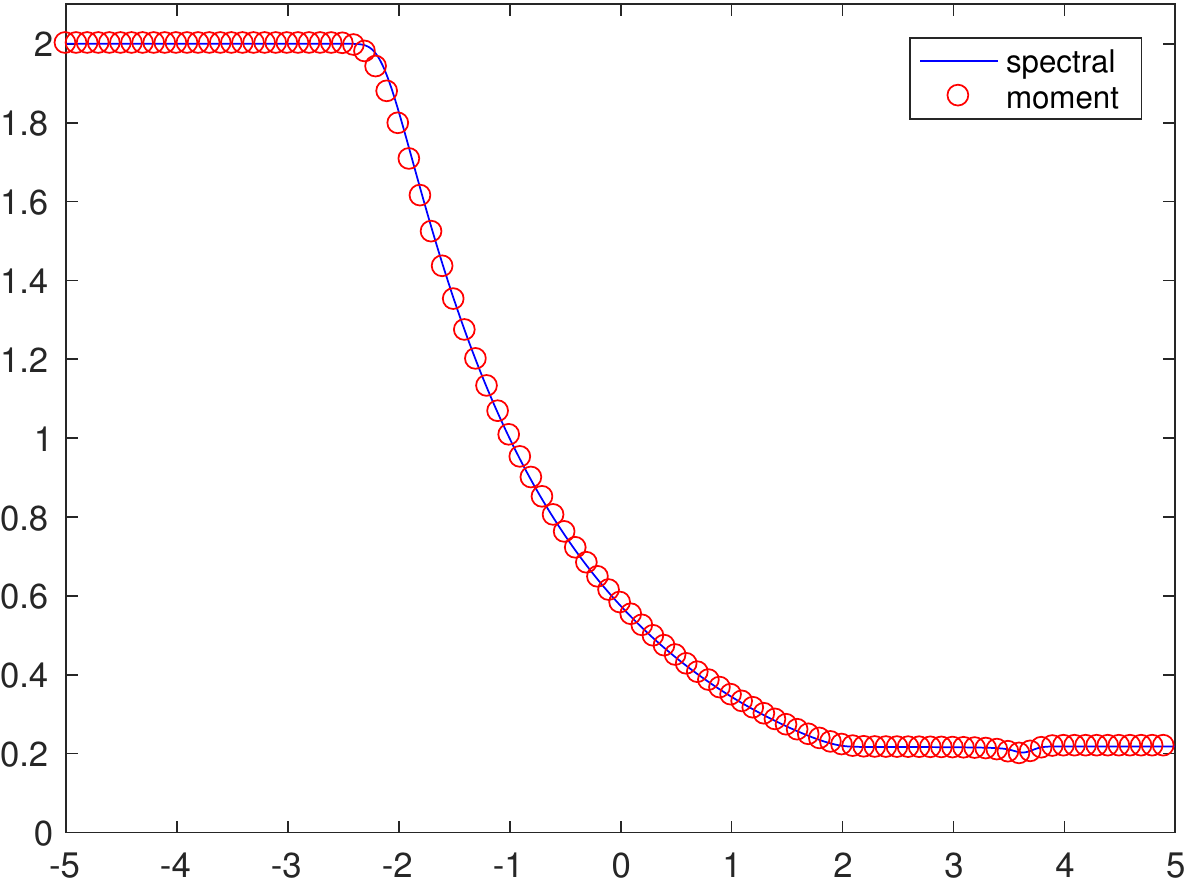}
}
\subfigure[$N=3$]{
\includegraphics[width=0.3\textwidth]{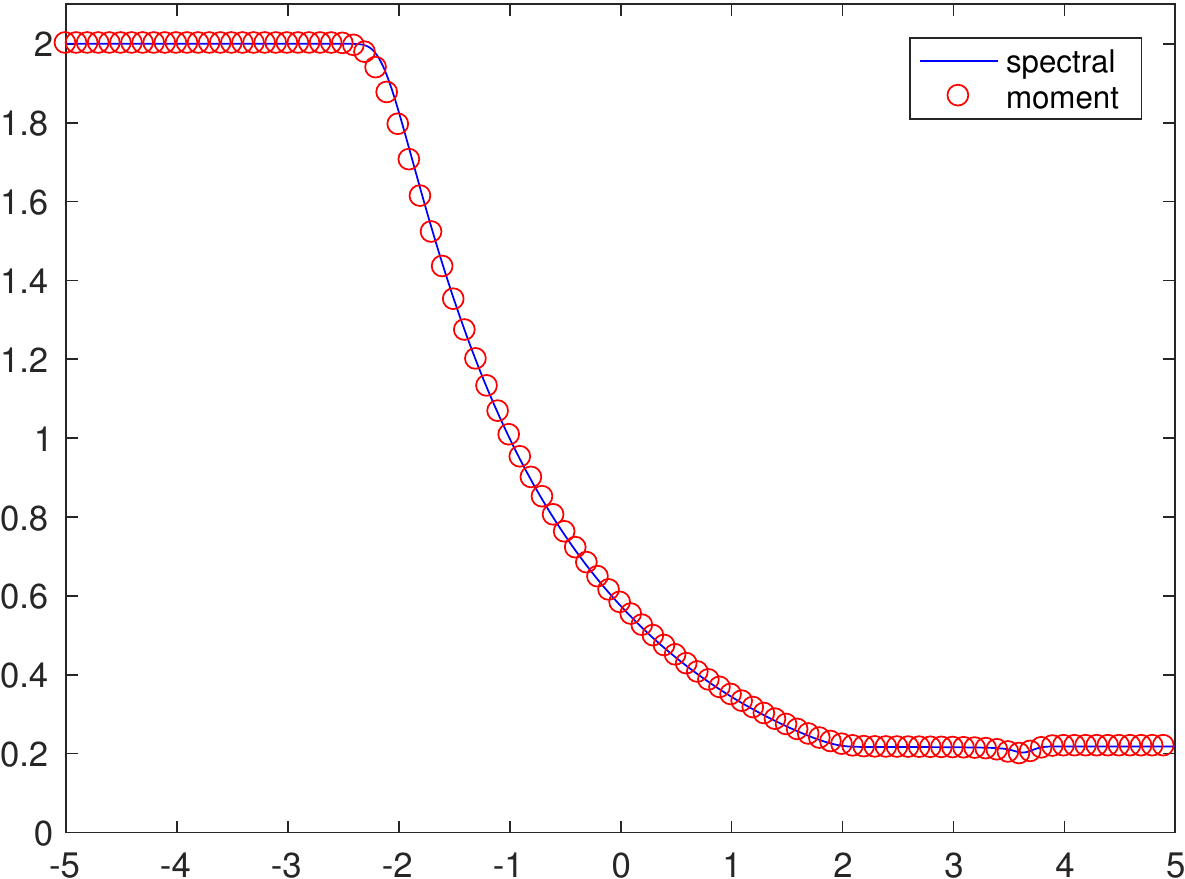}
}

  \centering
\subfigure[$N=4$]{
\includegraphics[width=0.3\textwidth]{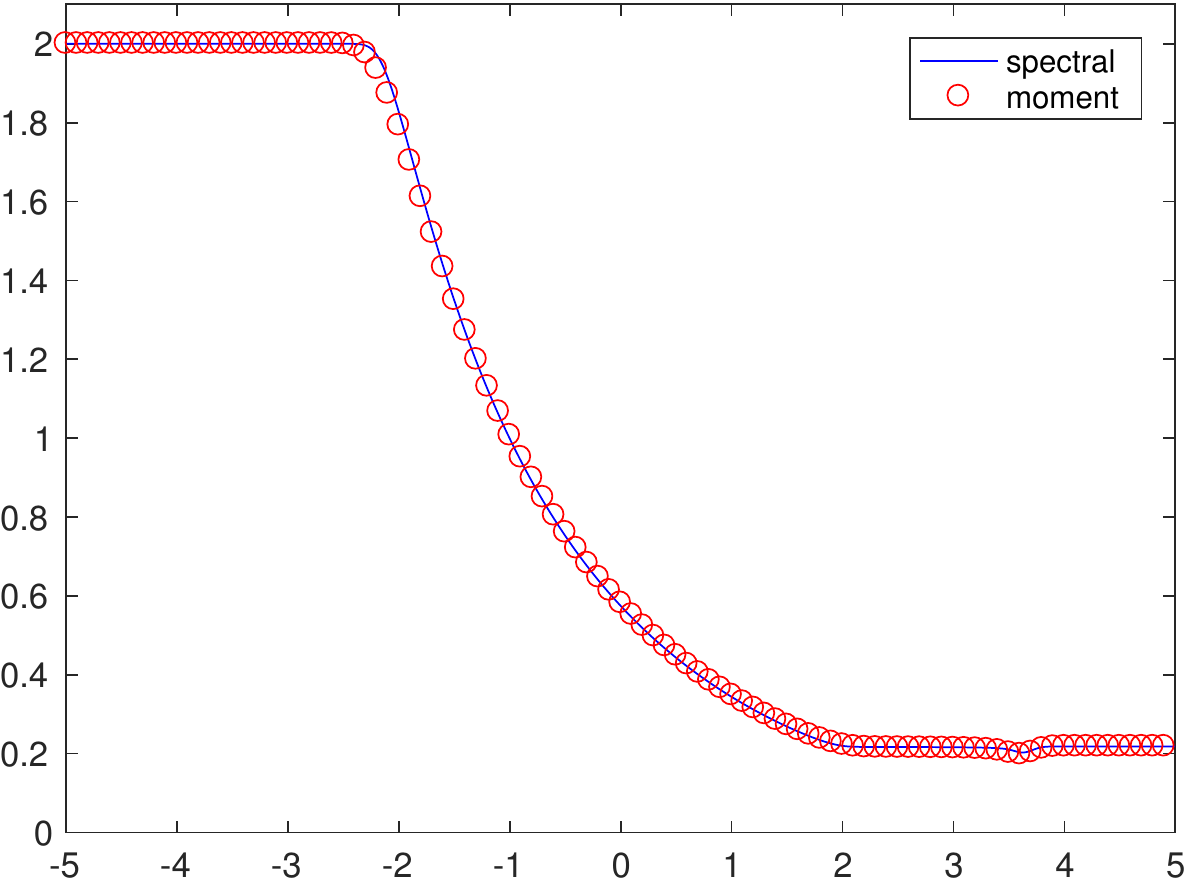}
}
\subfigure[$N=5$]{
\includegraphics[width=0.3\textwidth]{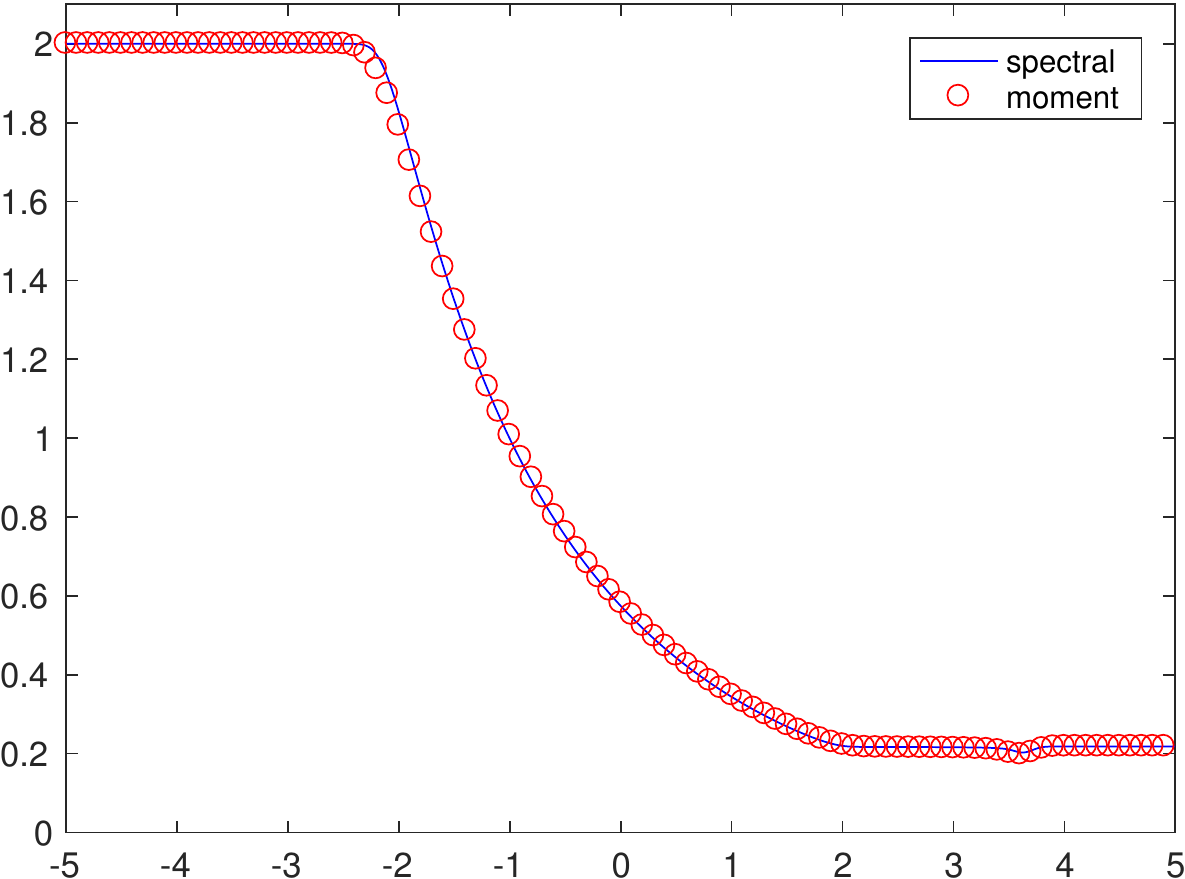}
}
\subfigure[$N=6$]{
\includegraphics[width=0.3\textwidth]{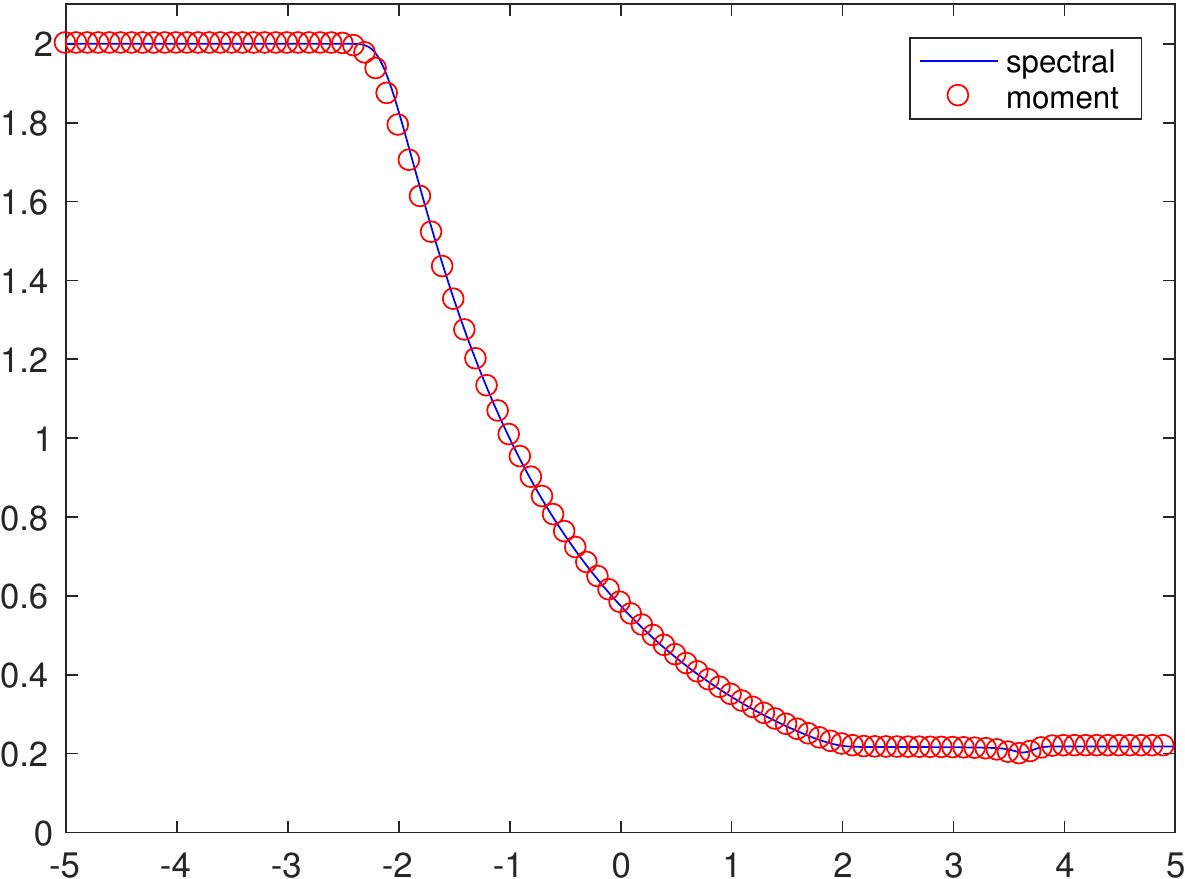}
}

%  \centering
%\subfigure[$N=7$]{
%\includegraphics[width=0.3\textwidth]{./images/riemann1/cmpN_e0010/r1_cmpN_rho_n2000_e0010_N7_circle.pdf}
%}
%\subfigure[$N=8$]{
%\includegraphics[width=0.3\textwidth]{./images/riemann1/cmpN_e0010/r1_cmpN_rho_n2000_e0010_N8_circle.pdf}
%}
%\subfigure[$N=9$]{
%\includegraphics[width=0.3\textwidth]{./images/riemann1/cmpN_e0010/r1_cmpN_rho_n2000_e0010_N9_circle.pdf}
%}
\caption{Same as Fig. \ref{figure:5.1} except for
 $\varepsilon=0.01$.}\label{figure:5.3}
\end{figure}

\begin{figure}
  \centering
\subfigure[$N=1$]{
\includegraphics[width=0.3\textwidth]{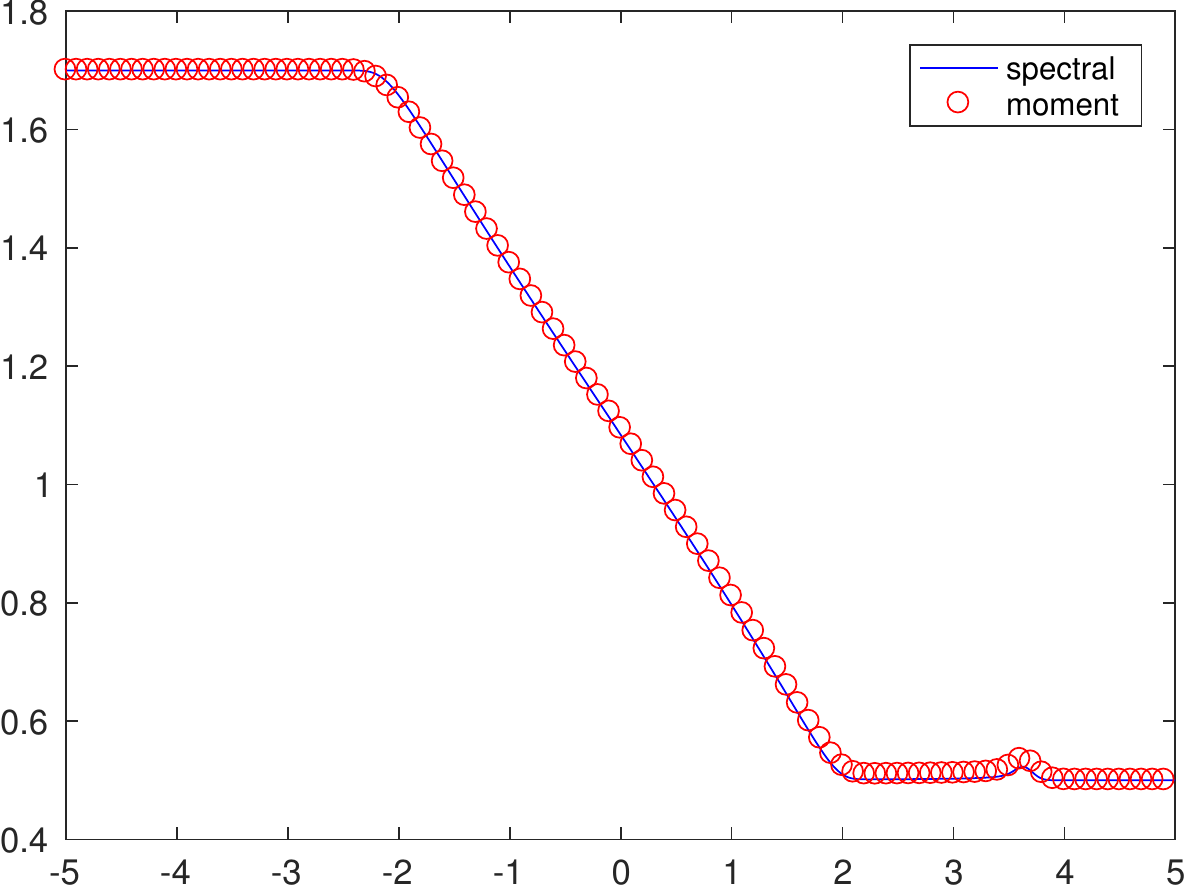}
}
\subfigure[$N=2$]{
\includegraphics[width=0.3\textwidth]{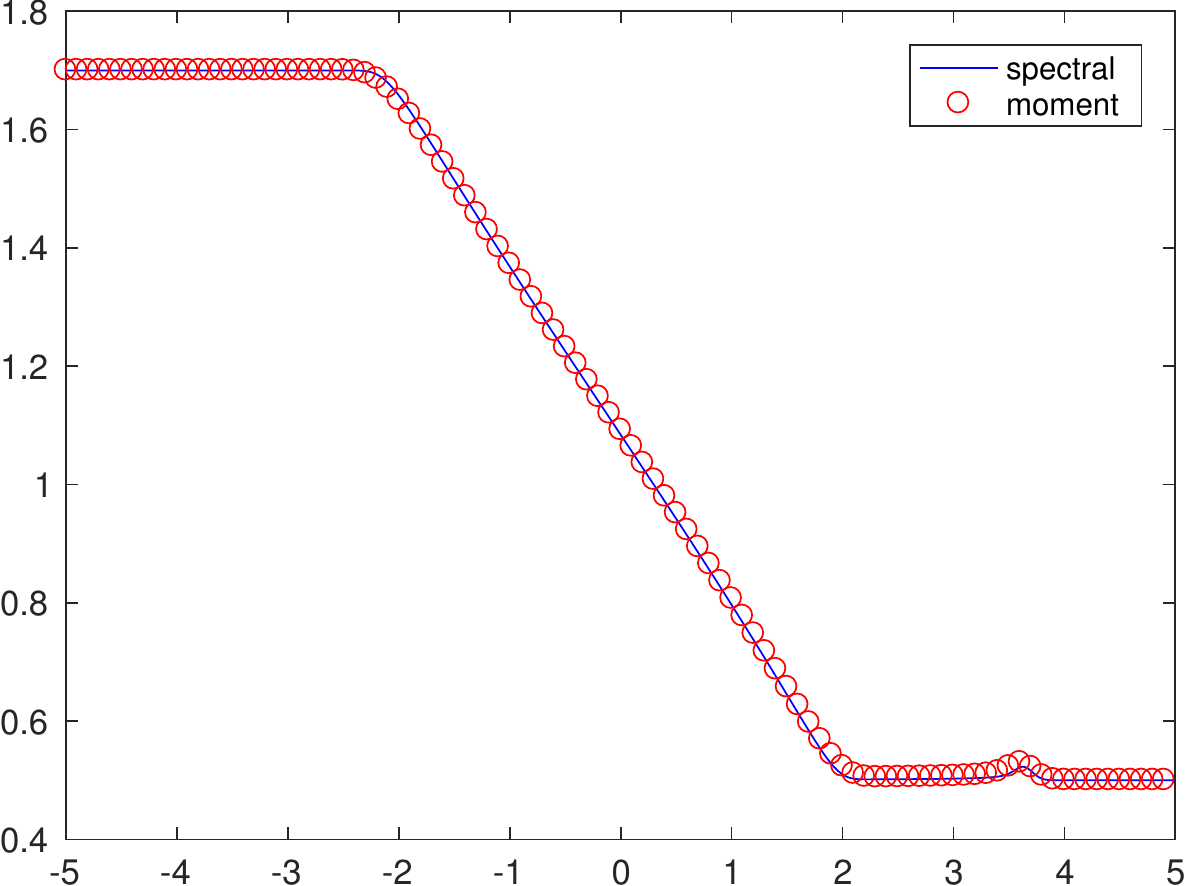}
}
\subfigure[$N=3$]{
\includegraphics[width=0.3\textwidth]{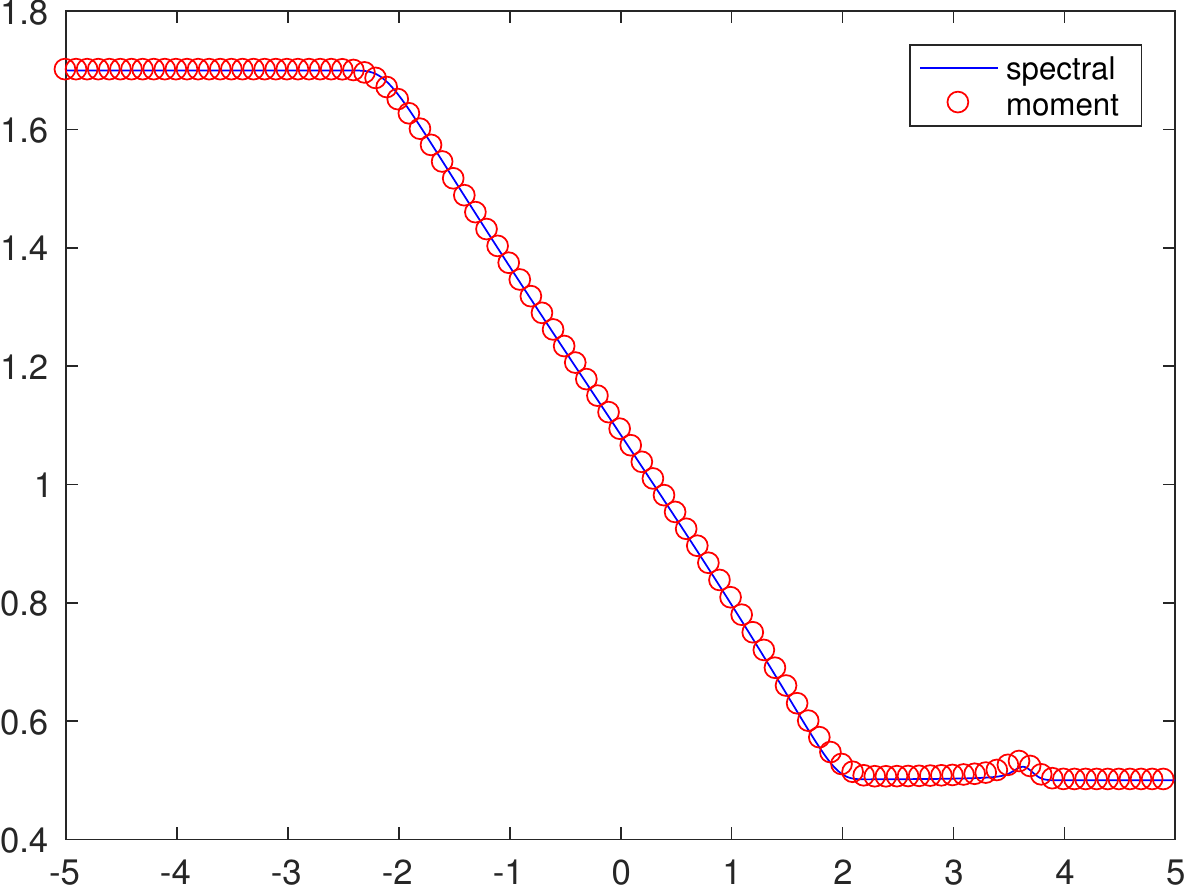}
}

  \centering
\subfigure[$N=4$]{
\includegraphics[width=0.3\textwidth]{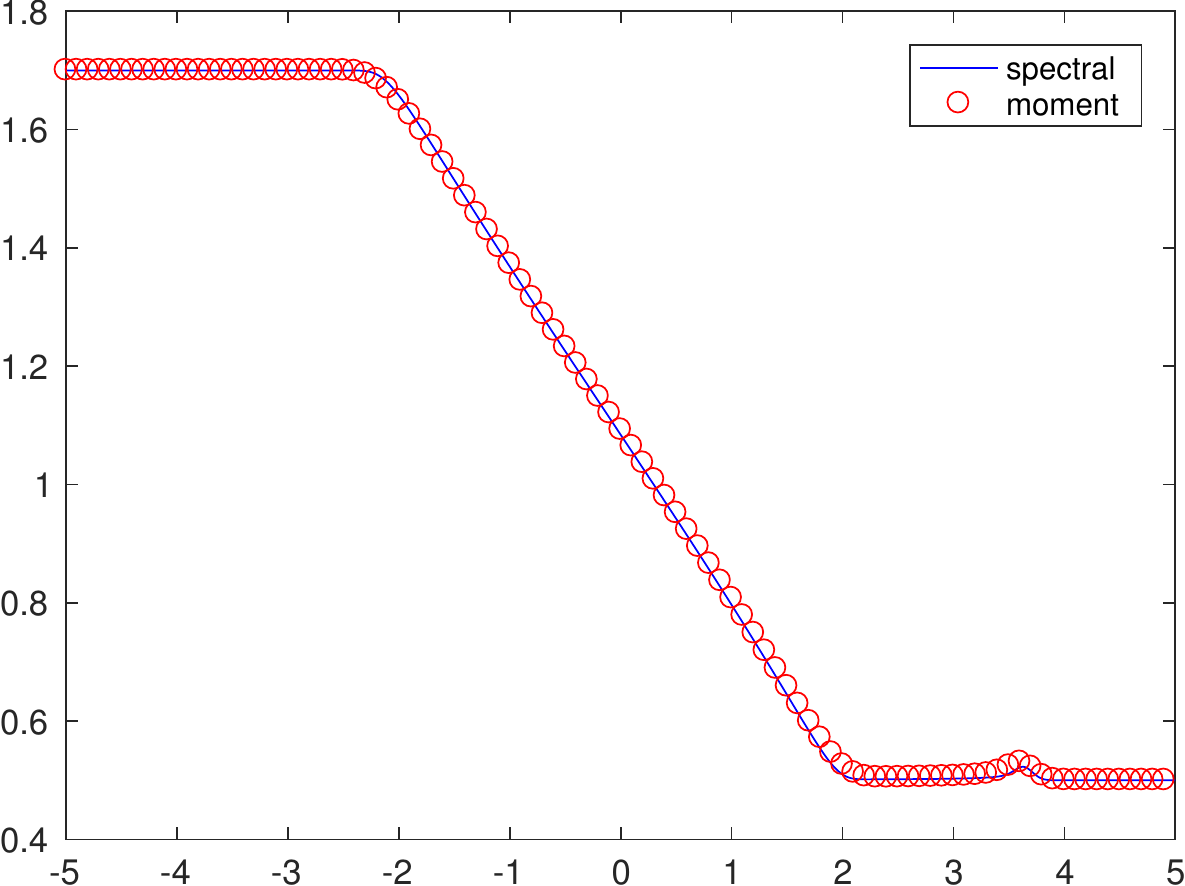}
}
\subfigure[$N=5$]{
\includegraphics[width=0.3\textwidth]{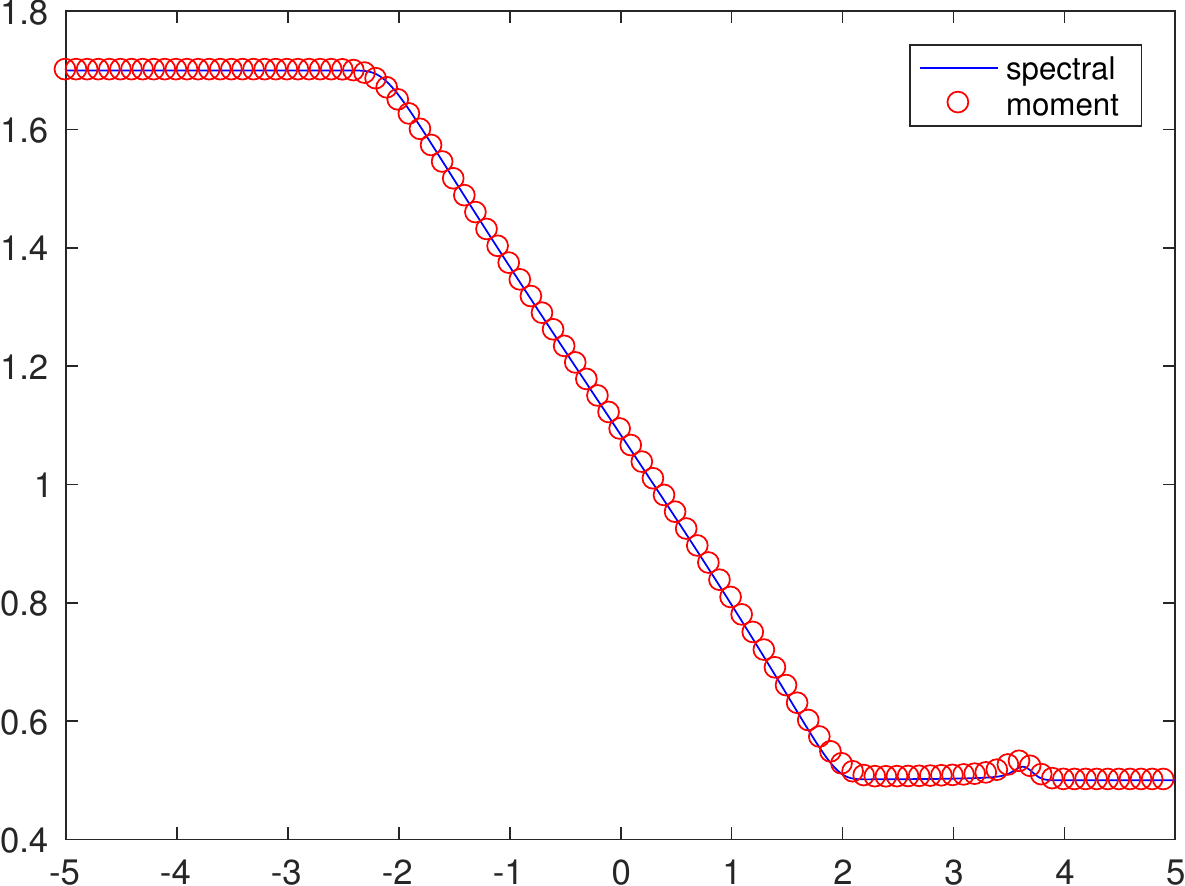}
}
\subfigure[$N=6$]{
\includegraphics[width=0.3\textwidth]{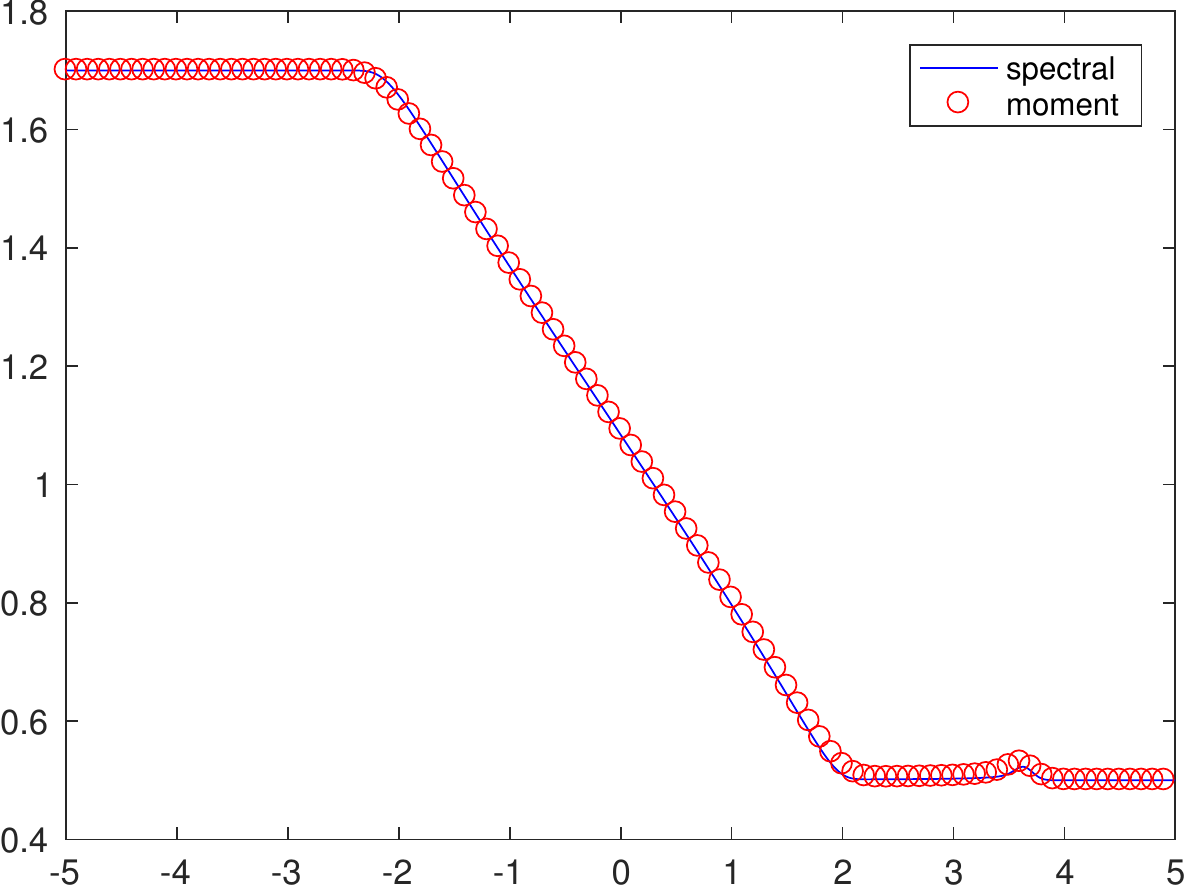}
}

%  \centering
%\subfigure[$N=7$]{
%\includegraphics[width=0.3\textwidth]{./images/riemann1/cmpN_e0010/r1_cmpN_u_n2000_e0010_N7_circle.pdf}
%}
%\subfigure[$N=8$]{
%\includegraphics[width=0.3\textwidth]{./images/riemann1/cmpN_e0010/r1_cmpN_u_n2000_e0010_N8_circle.pdf}
%}
%\subfigure[$N=9$]{
%\includegraphics[width=0.3\textwidth]{./images/riemann1/cmpN_e0010/r1_cmpN_u_n2000_e0010_N9_circle.pdf}
%}
\caption{Same as Fig. \ref{figure:5.3} except for
 the macroscopic velocity angles.}\label{figure:5.4}
\end{figure}

%where the macroscopic density and velocity are determined by \eqref{eq:btdingyi}.
\begin{example}[Rarefaction wave]\label{example:5.1}
The initial data of the first Riemann problem for the density $\rho^\varepsilon$
and velocity angle $\bar\theta^\varepsilon$ are
\begin{equation*}
	%\label{sl:r1}
	(\rho^\varepsilon,\bar\theta^\varepsilon)=\begin{cases} (2, 1.7), & x<0,\\
  (0.218, 0.5), & x>0,
\end{cases}
\end{equation*}
and the initial particle distribution function is set as  the Von Mises distribution
associated with the initial density and velocity angle.
The solutions of this problem are given by a rarefaction wave.
Figs. \ref{figure:5.1} and \ref{figure:5.2} show
the densities $\rho^\varepsilon$  and macroscopic velocity angles
$\bt^\varepsilon$ at $t=4$ obtained by the moment
method with $N=1,2,\cdots,6$, 2000 cells, and $\varepsilon=1$,
where the solid line denotes the reference solution obtained by using the
spectral method  with 4000 cells.  Figs. \ref{figure:5.3} and  \ref{figure:5.4}
display  corresponding solutions for the case of
$\varepsilon=0.01$.
It is seen that
the solutions  of the moment
system well agree with the reference when $N$ is larger than 1,
and for a fixed  $N$, the solutions of moment method
also get closer to the reference as $\varepsilon$ decreases.

%初始化时, 假定空间中每一点处于平衡态, 即每个单元网格上都是Von Mises分布函数,
%在给定初始的宏观密度和宏观速度后每一点的密度分布函数$f(x, \theta)$也就确定了.

\end{example}

%%%%%%%%%%%%%%%%%%%%%%%% cmpN_e1000 %%%%%%%%%%%%%%%%%%%%%%%%
%%%%%%%%%%%%%%%%%%%%%%%% r1 rho circle %%%%%%%%%%%%%%%%%%%%%%%%

%%%%%%%%%%%%%%%%%%%%%%%% cmpN_e1000 %%%%%%%%%%%%%%%%%%%%%%%%
%%%%%%%%%%%%%%%%%%%%%%%% r1 u circle %%%%%%%%%%%%%%%%%%%%%%%%

%%%%%%%%%%%%%%%%%%%%%%%% cmpN_e0010 %%%%%%%%%%%%%%%%%%%%%%%%
%%%%%%%%%%%%%%%%%%%%%%%% r1 rho circle %%%%%%%%%%%%%%%%%%%%%%%%

%%%%%%%%%%%%%%%%%%%%%%%% cmpN_e0010 %%%%%%%%%%%%%%%%%%%%%%%%
%%%%%%%%%%%%%%%%%%%%%%%% r1 u circle %%%%%%%%%%%%%%%%%%%%%%%%

\begin{figure}
  \centering
\subfigure[$N=1$]{
\includegraphics[width=0.3\textwidth]{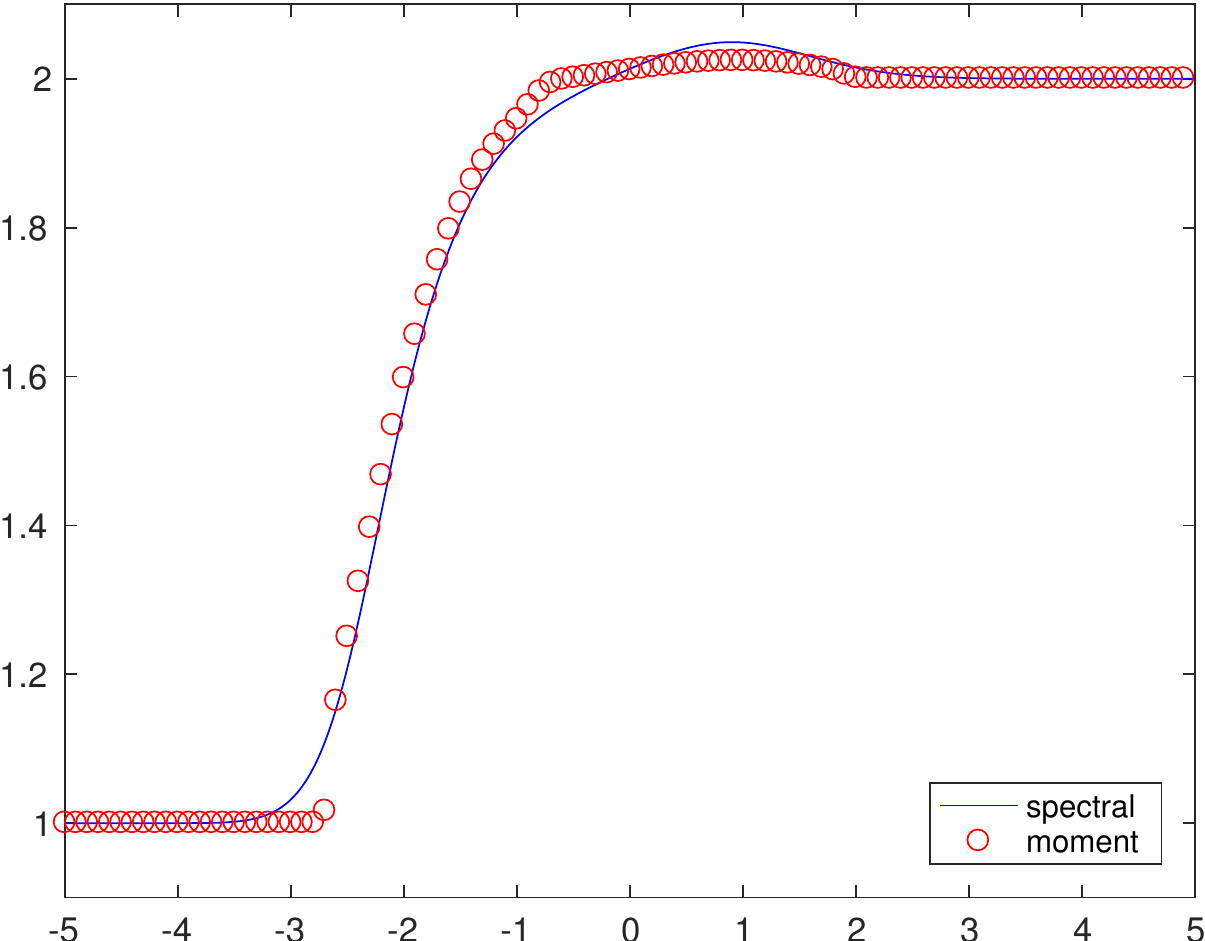}
}
\subfigure[$N=2$]{
\includegraphics[width=0.3\textwidth]{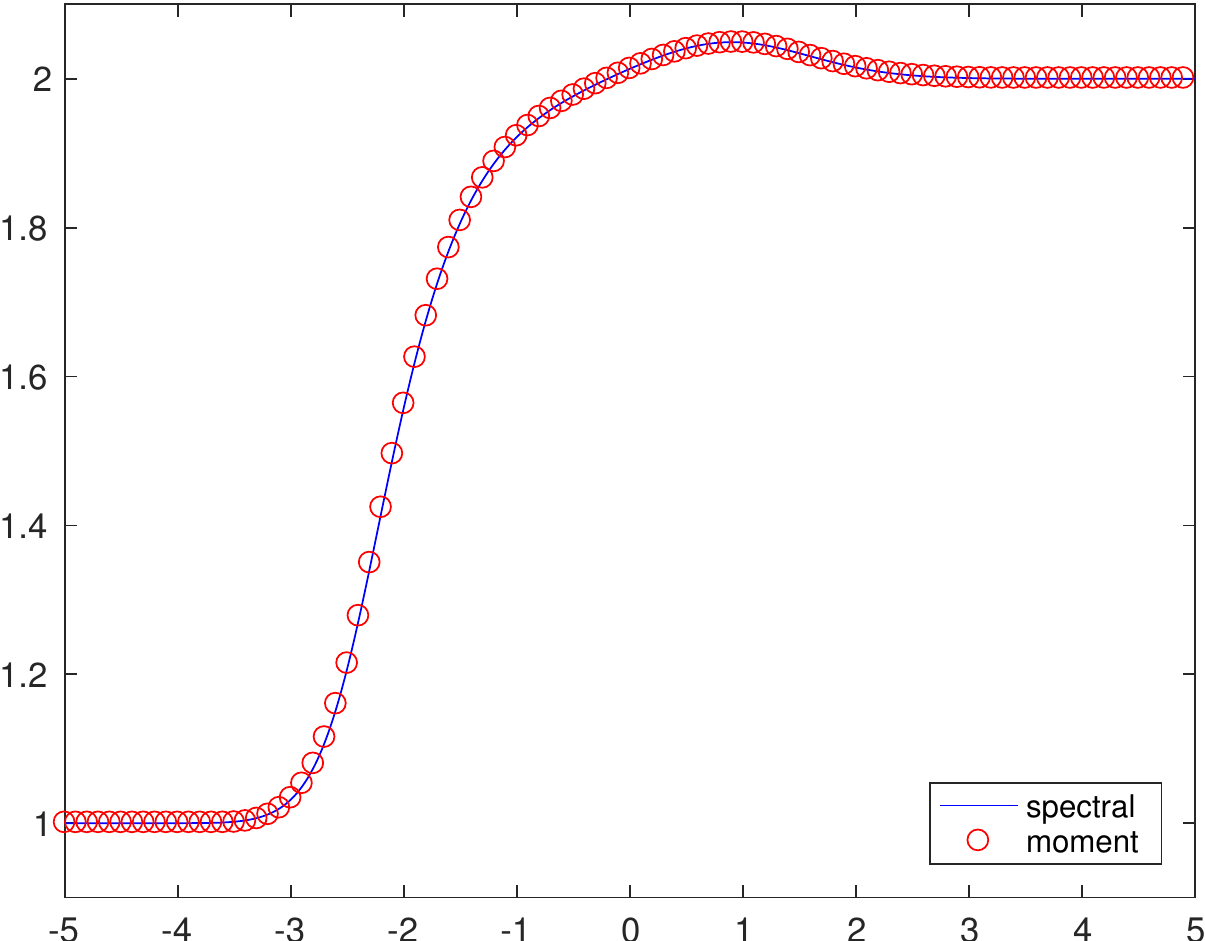}
}
\subfigure[$N=3$]{
\includegraphics[width=0.3\textwidth]{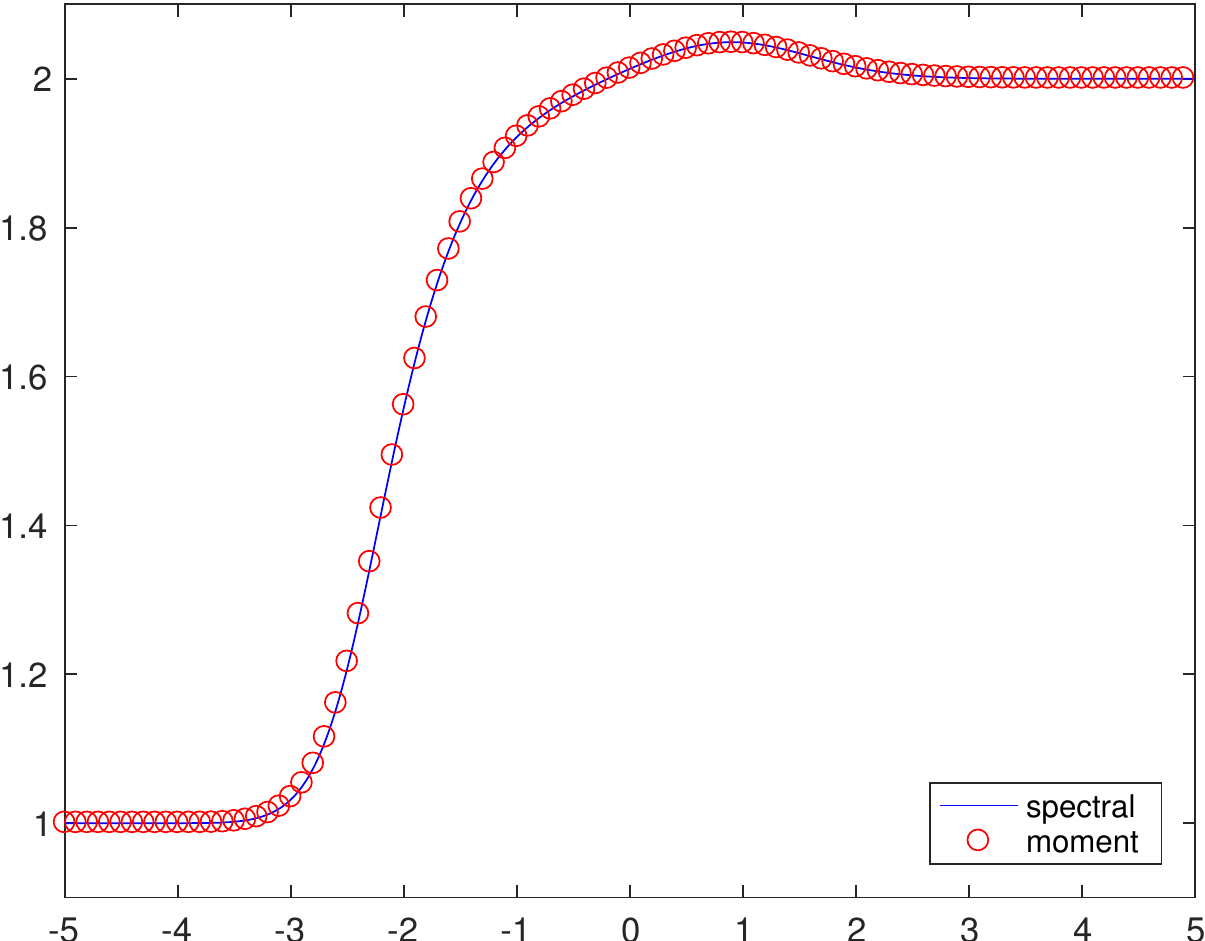}
}

  \centering
\subfigure[$N=4$]{
\includegraphics[width=0.3\textwidth]{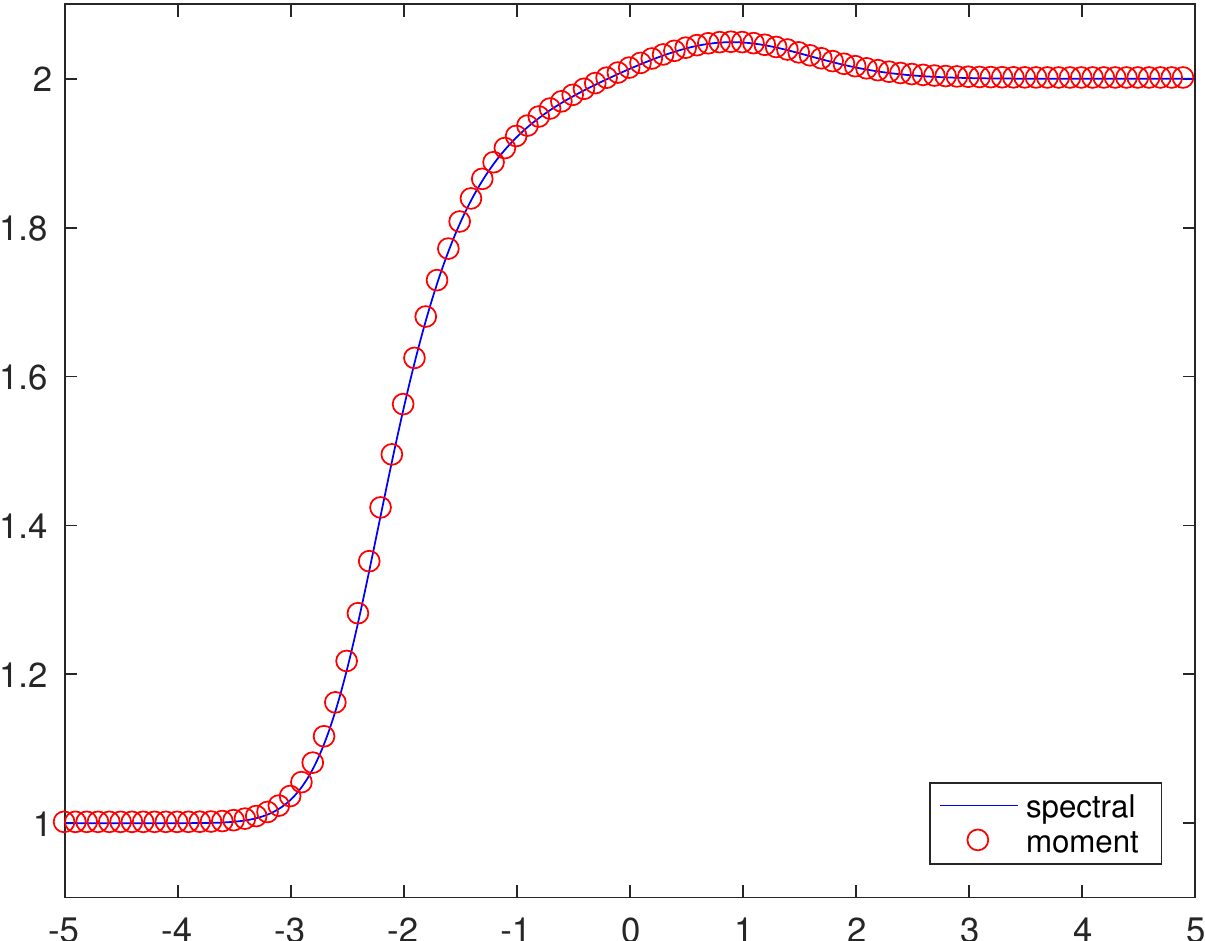}
}
\subfigure[$N=5$]{
\includegraphics[width=0.3\textwidth]{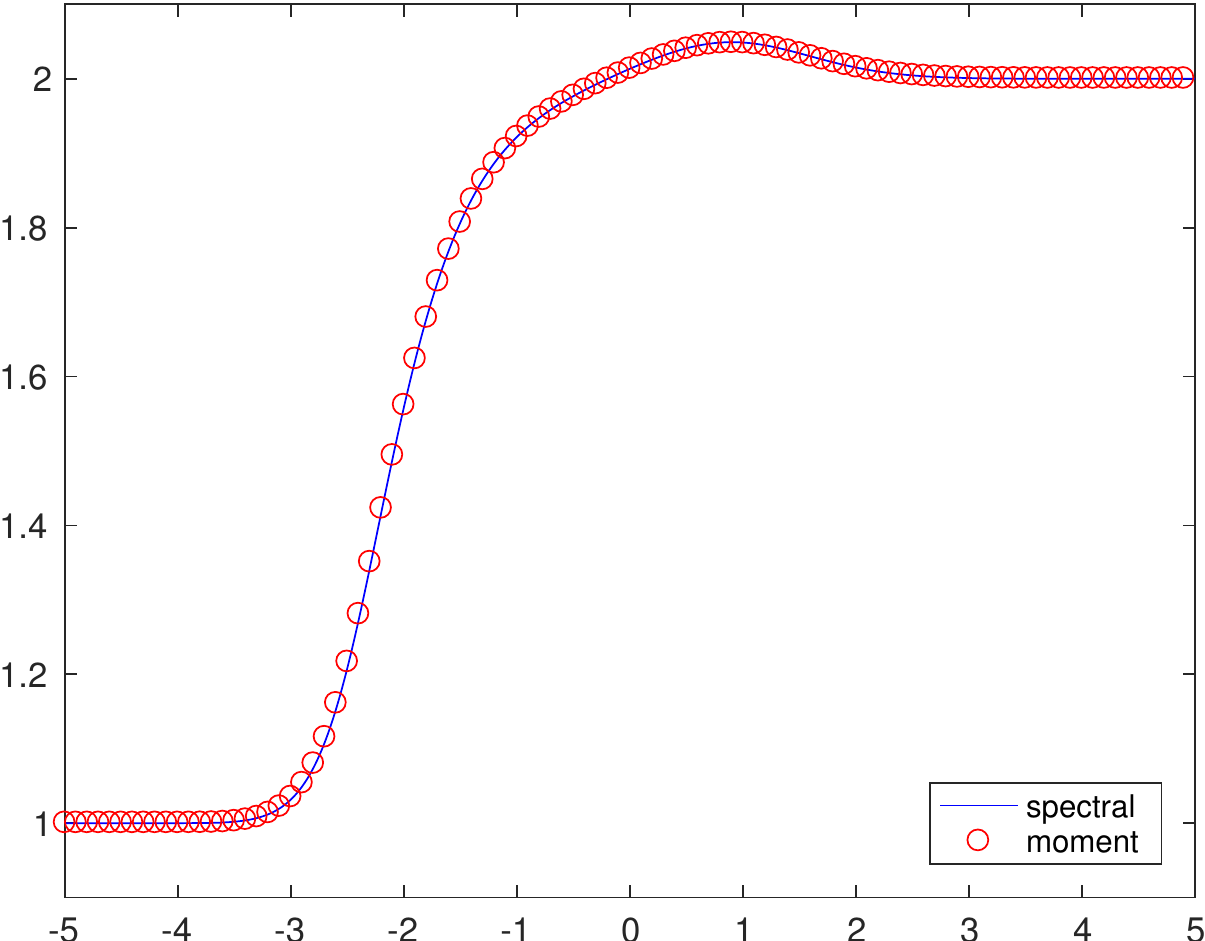}
}
\subfigure[$N=6$]{
\includegraphics[width=0.3\textwidth]{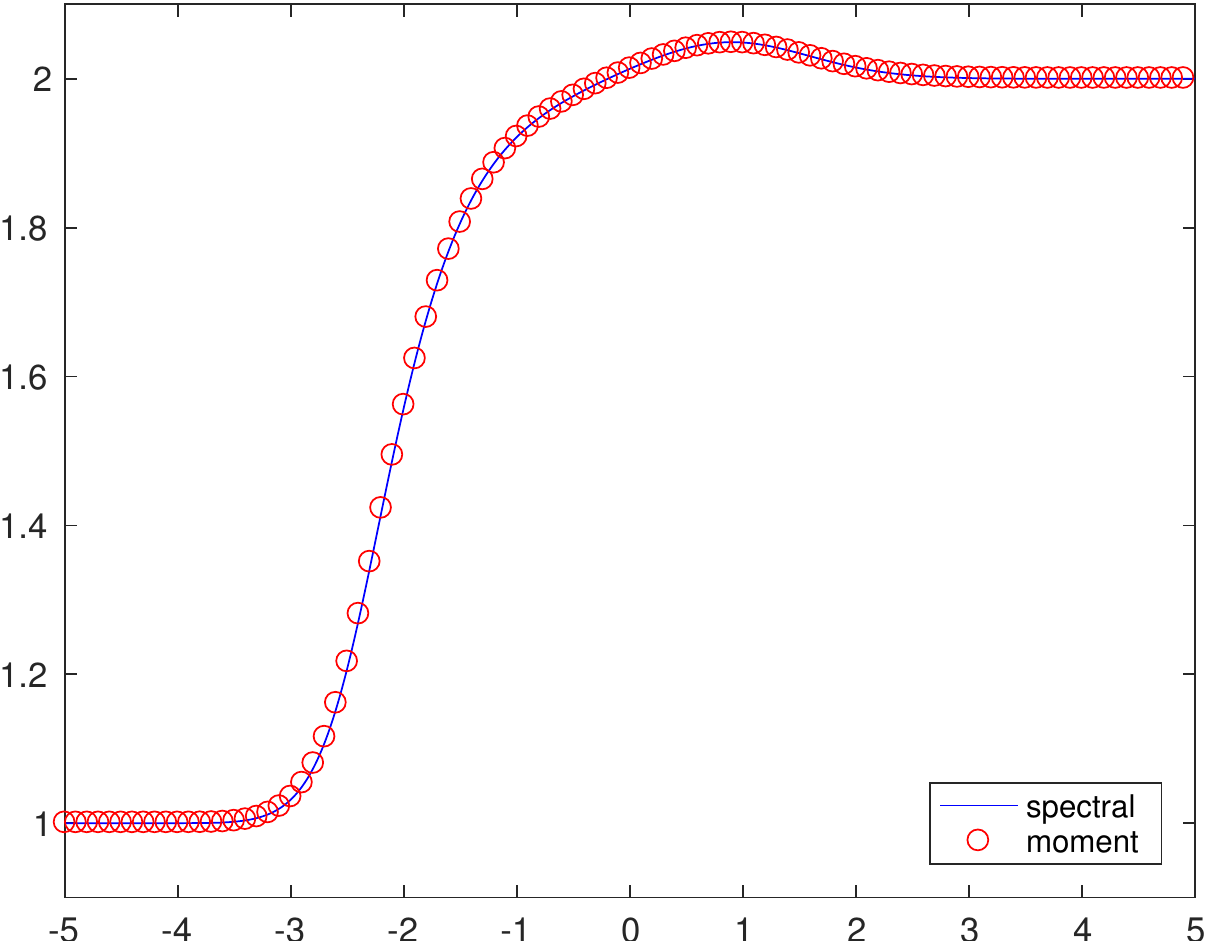}
}

%  \centering
%\subfigure[$N=7$]{
%\includegraphics[width=0.3\textwidth]{./images/riemann2/cmpN_e1000/r2_cmpN_rho_n2000_e1000_N7_circle.pdf}
%}
%\subfigure[$N=8$]{
%\includegraphics[width=0.3\textwidth]{./images/riemann2/cmpN_e1000/r2_cmpN_rho_n2000_e1000_N8_circle.pdf}
%}
%\subfigure[$N=9$]{
%\includegraphics[width=0.3\textwidth]{./images/riemann2/cmpN_e1000/r2_cmpN_rho_n2000_e1000_N9_circle.pdf}
%}
\caption{Example \ref{example:5.2}: The densities at $t=4$ obtained by the moment
method with $N=1,2,\cdots,6$ and 2000 cells. The solid line is the reference solution obtained by using the
spectral method with 4000 cells.  $\varepsilon=1$.} \label{figure-example2-1}
%使用4000个网格Spectral method的结果作为参考,
%和使用2000个网格的矩方法的结果进行比较}
\end{figure}

\begin{figure}
  \centering
\subfigure[$N=1$]{
\includegraphics[width=0.3\textwidth]{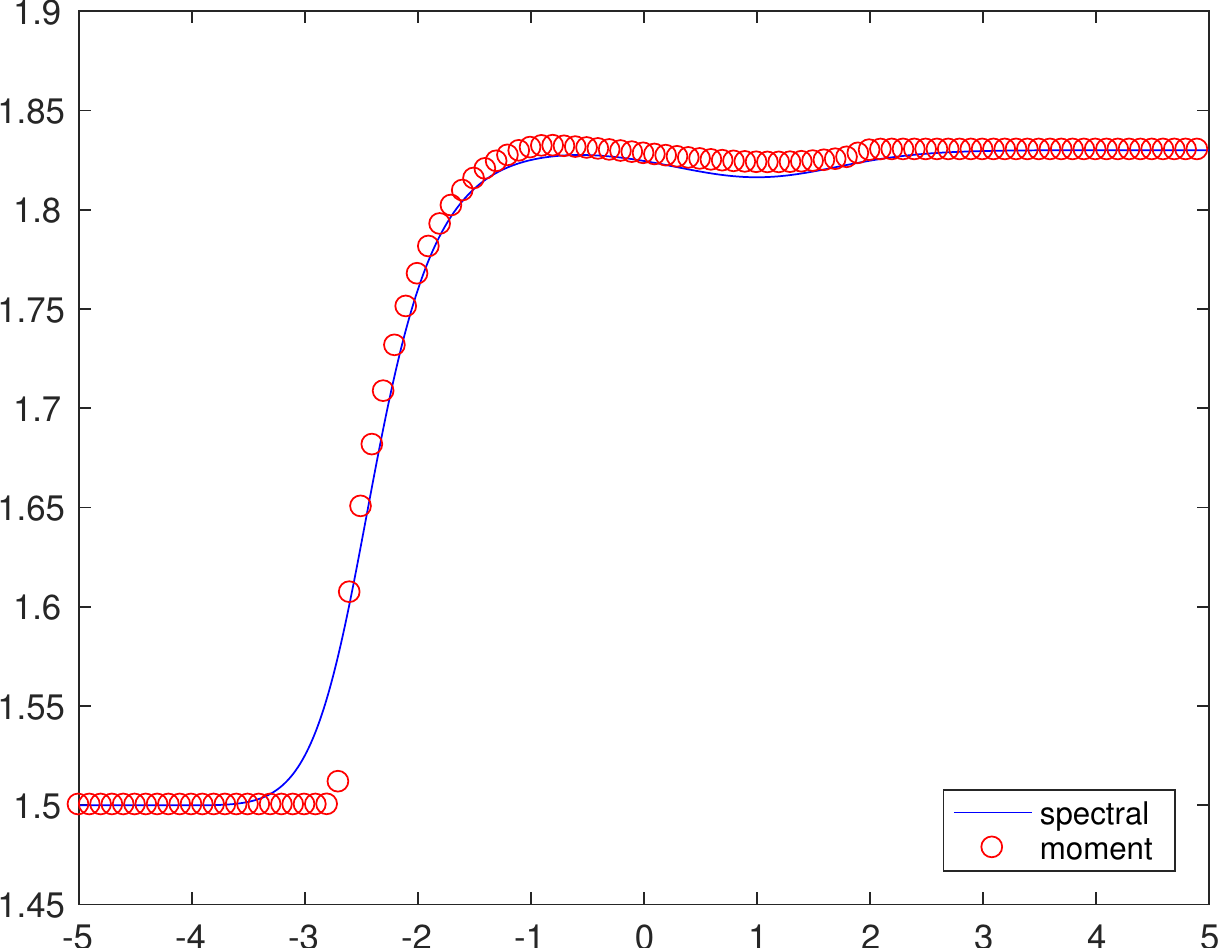}
}
\subfigure[$N=2$]{
\includegraphics[width=0.3\textwidth]{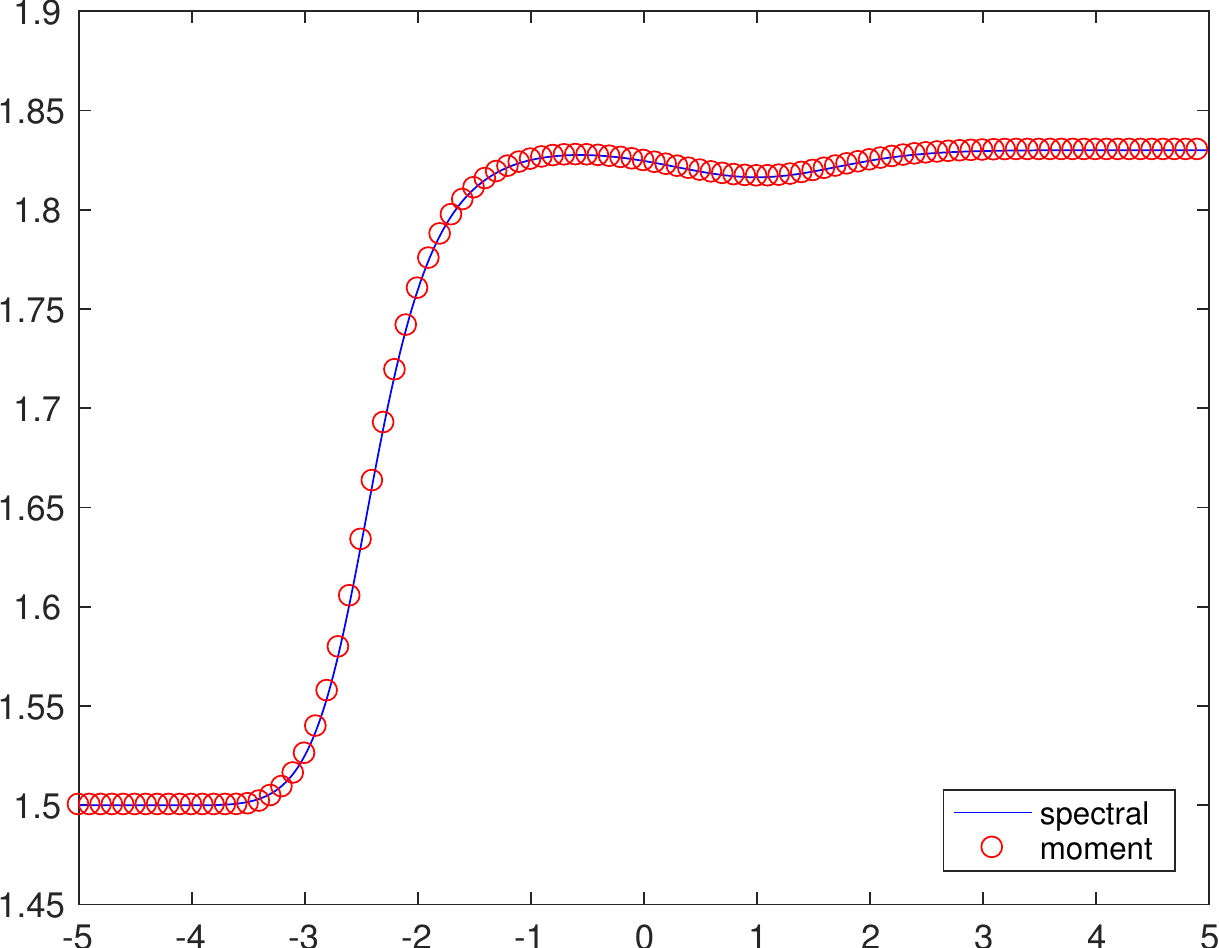}
}
\subfigure[$N=3$]{
\includegraphics[width=0.3\textwidth]{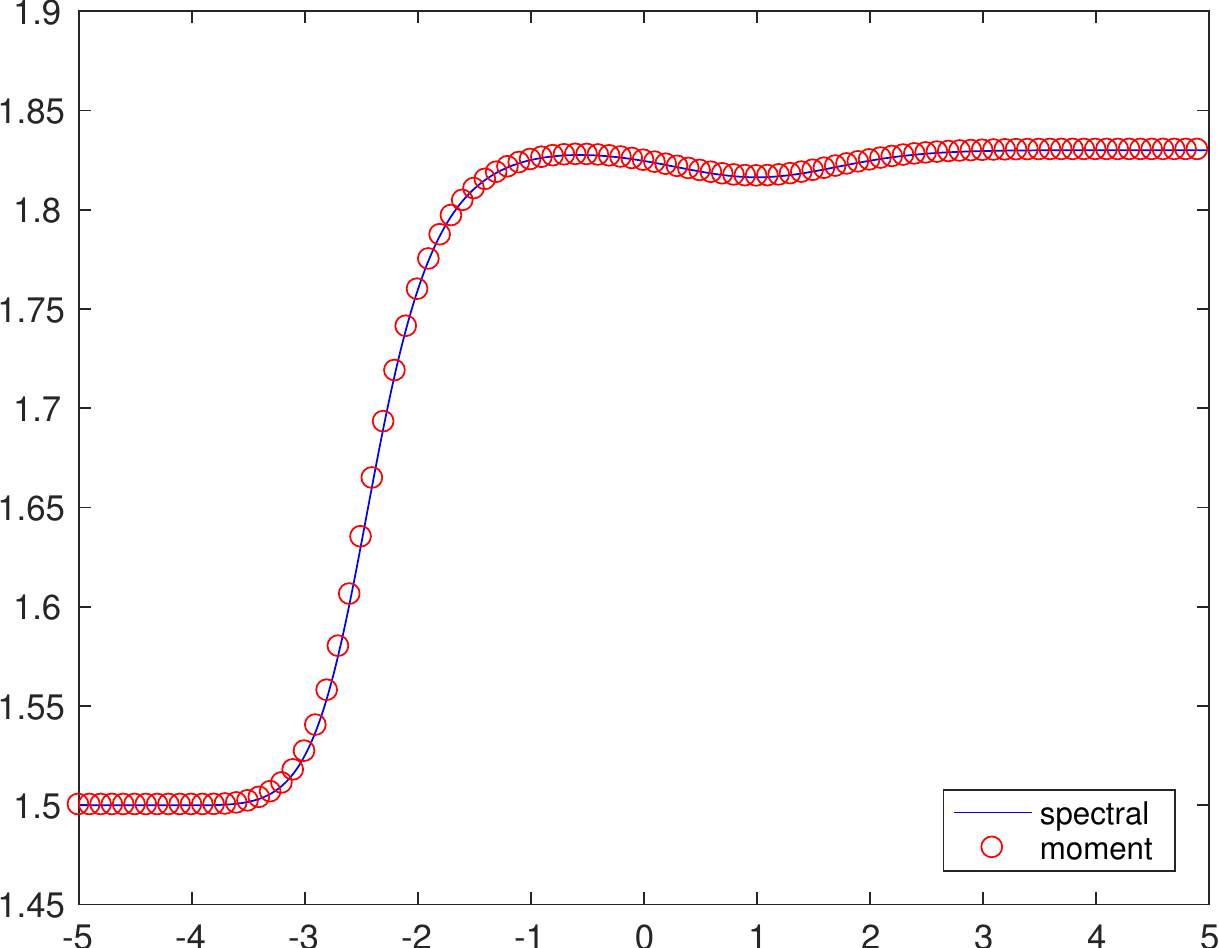}
}

  \centering
\subfigure[$N=4$]{
\includegraphics[width=0.3\textwidth]{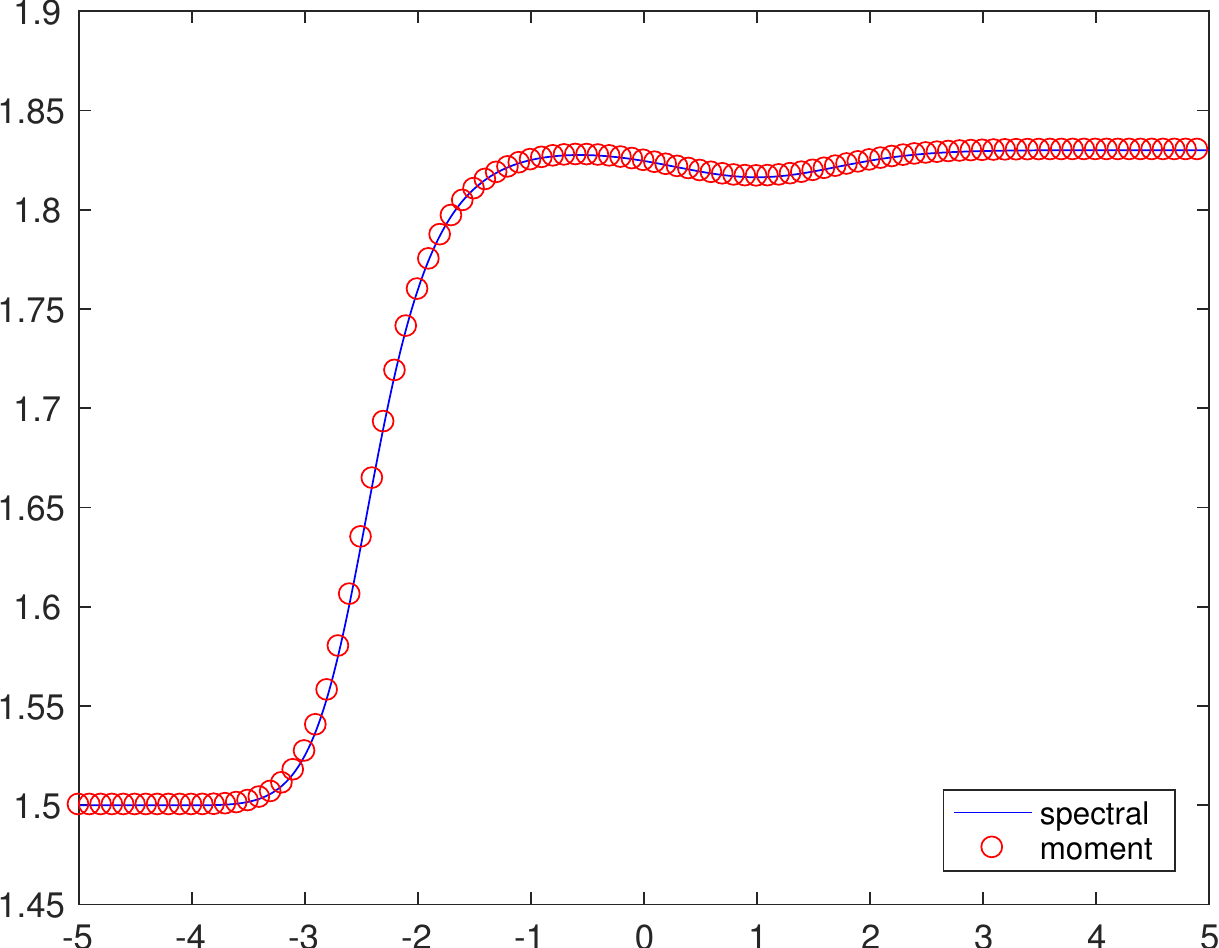}
}
\subfigure[$N=5$]{
\includegraphics[width=0.3\textwidth]{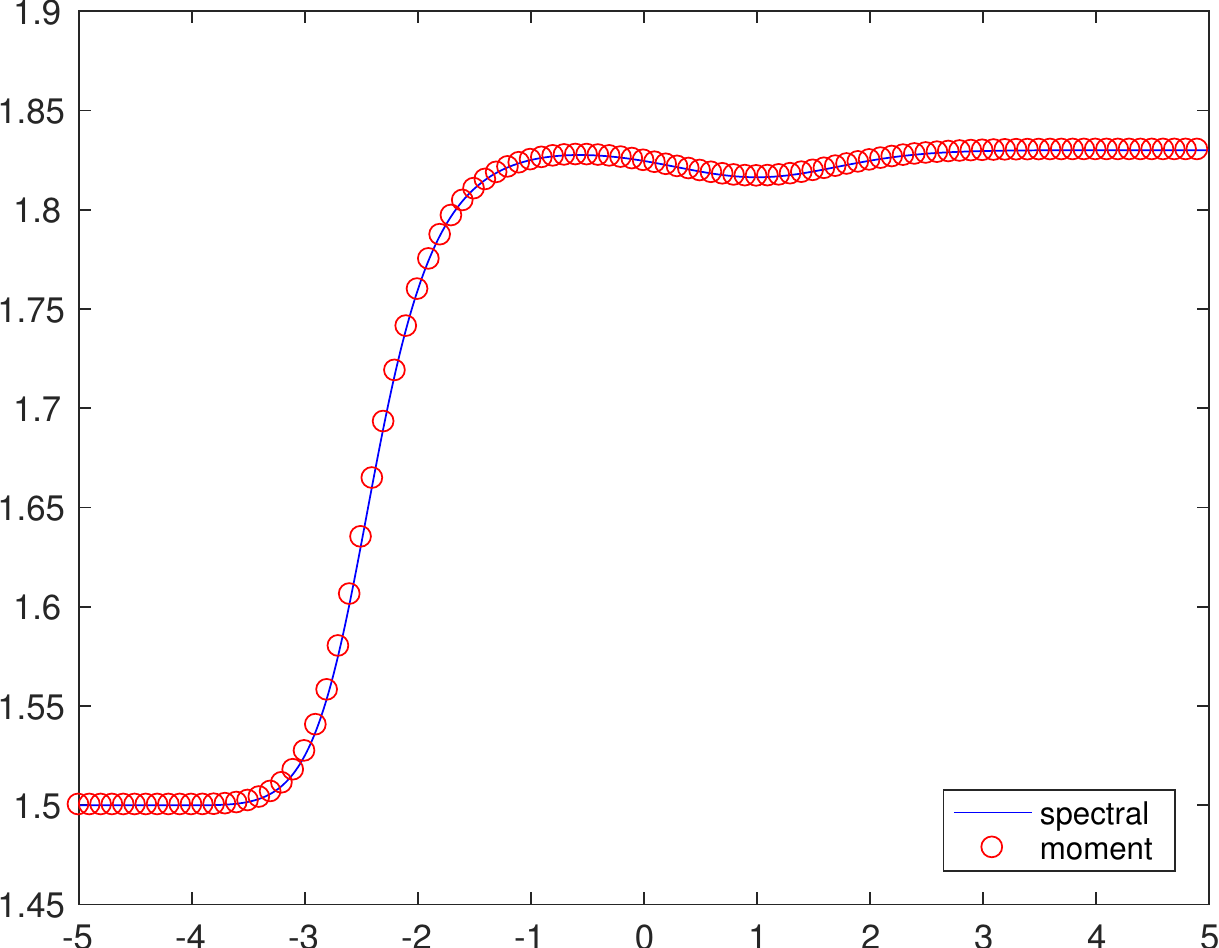}
}
\subfigure[$N=6$]{
\includegraphics[width=0.3\textwidth]{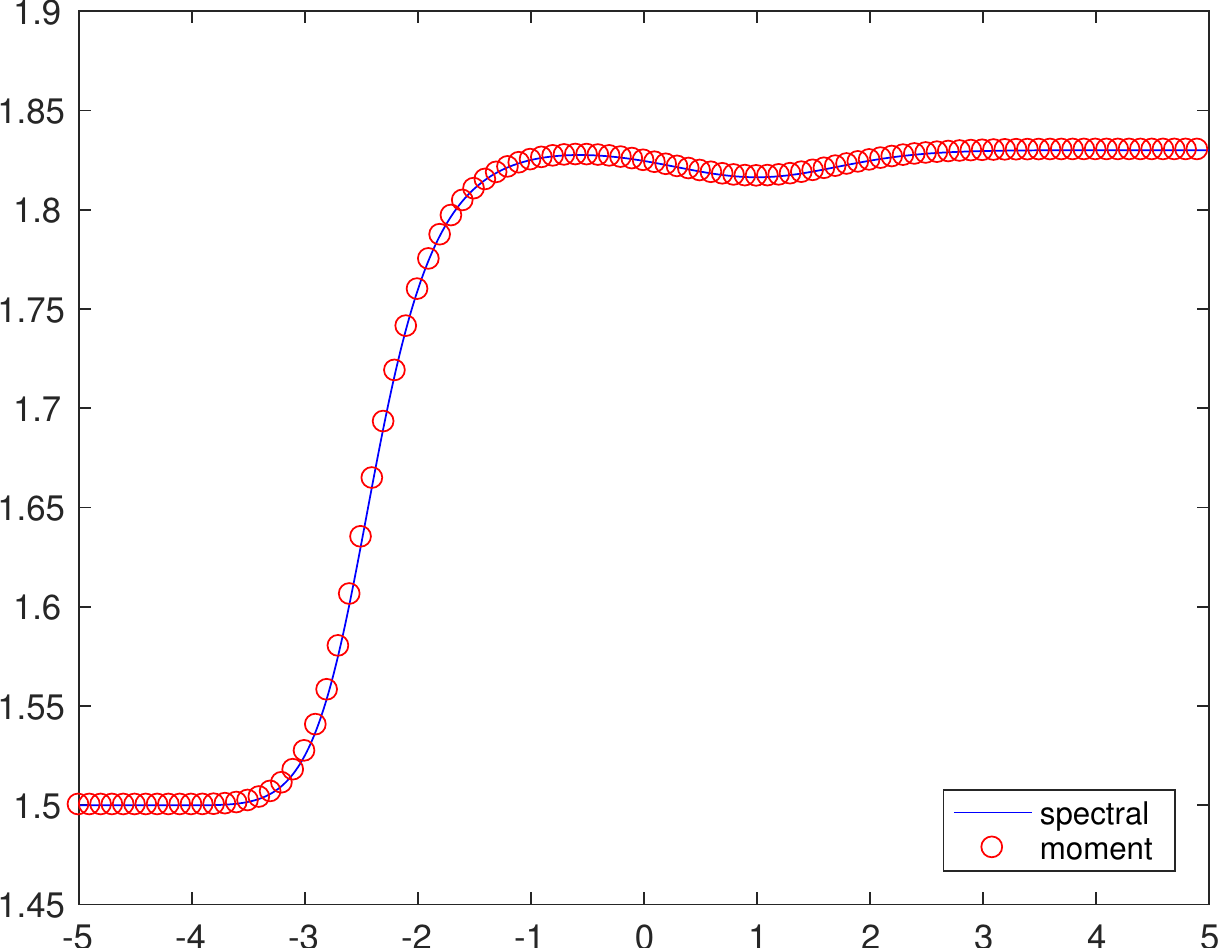}
}

%  \centering
%\subfigure[$N=7$]{
%\includegraphics[width=0.3\textwidth]{./images/riemann2/cmpN_e1000/r2_cmpN_u_n2000_e1000_N7_circle.pdf}
%}
%\subfigure[$N=8$]{
%\includegraphics[width=0.3\textwidth]{./images/riemann2/cmpN_e1000/r2_cmpN_u_n2000_e1000_N8_circle.pdf}
%}
%\subfigure[$N=9$]{
%\includegraphics[width=0.3\textwidth]{./images/riemann2/cmpN_e1000/r2_cmpN_u_n2000_e1000_N9_circle.pdf}
%}
\caption{Same as Fig. \ref{figure-example2-1} except for the macroscopic
velocity angles.}\label{figure-example2-2}
\end{figure}

\begin{figure}
  \centering
\subfigure[$N=1$]{
\includegraphics[width=0.3\textwidth]{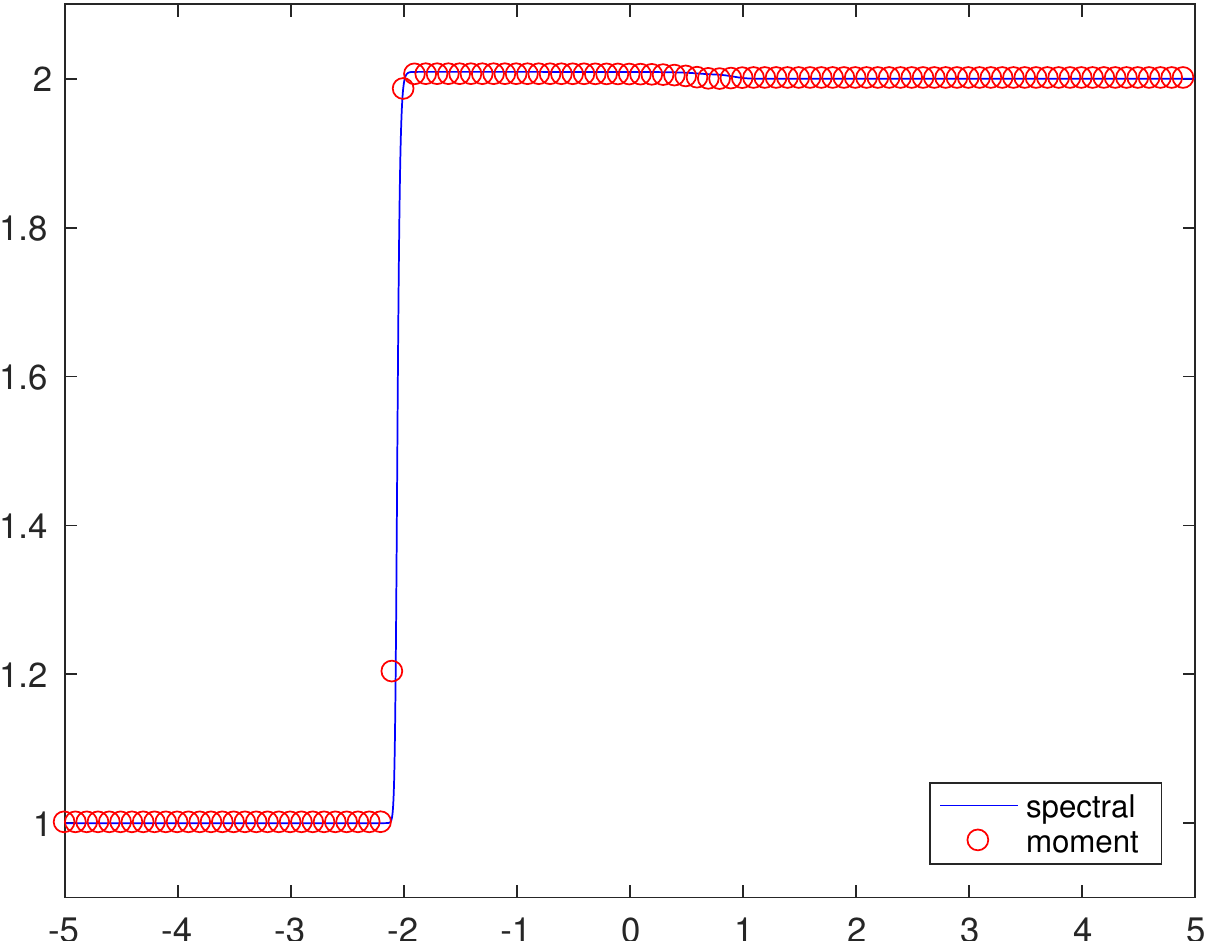}
}
\subfigure[$N=2$]{
\includegraphics[width=0.3\textwidth]{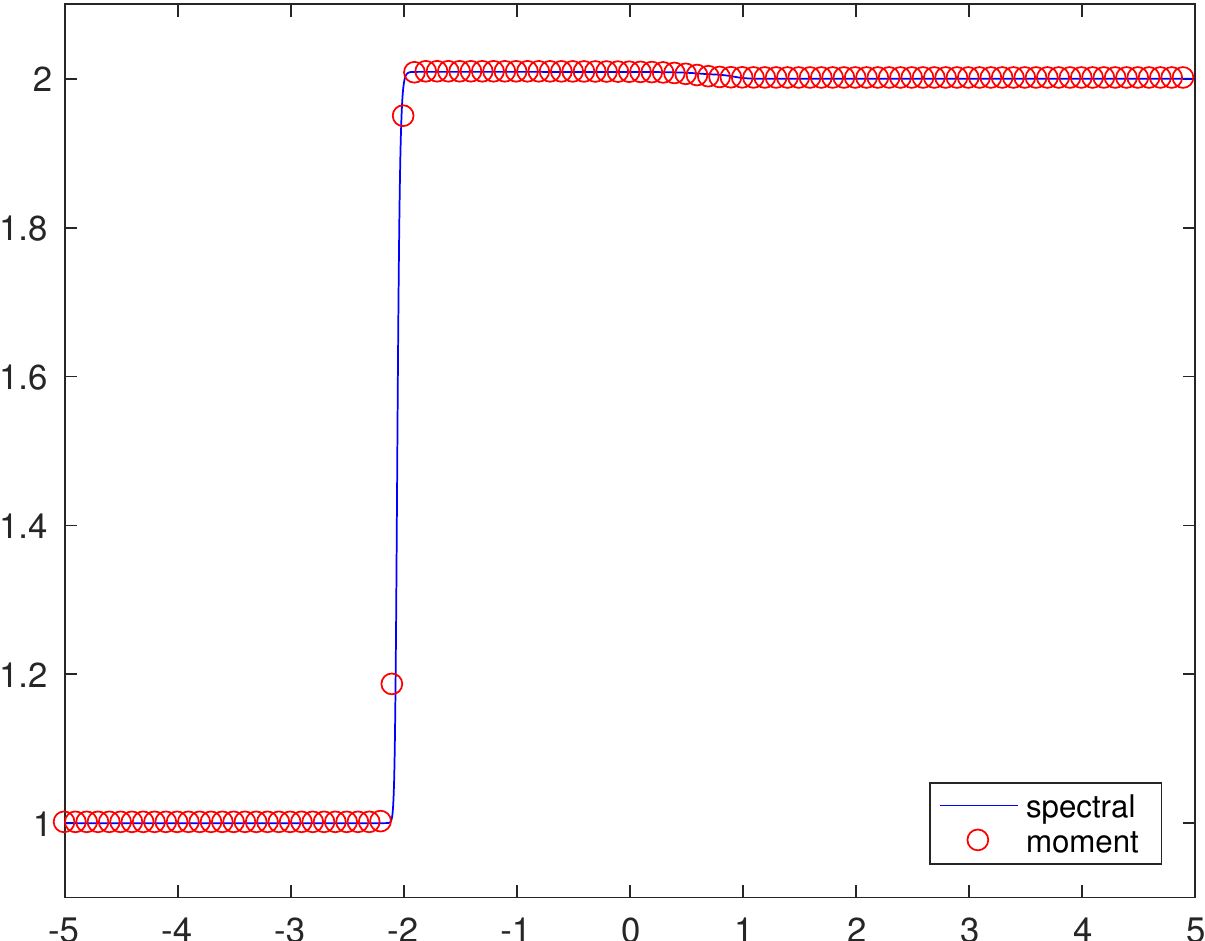}
}
\subfigure[$N=3$]{
\includegraphics[width=0.3\textwidth]{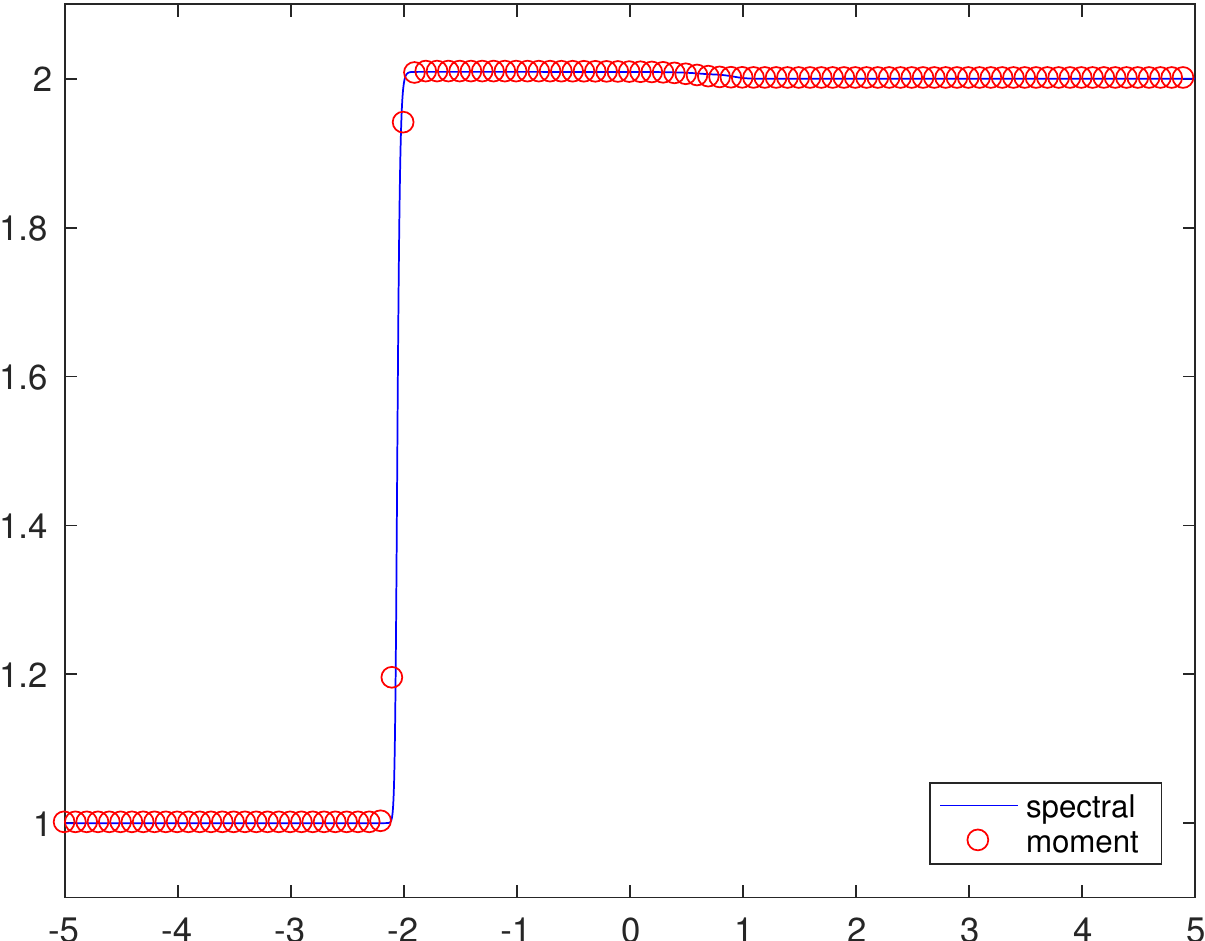}
}

  \centering
\subfigure[$N=4$]{
\includegraphics[width=0.3\textwidth]{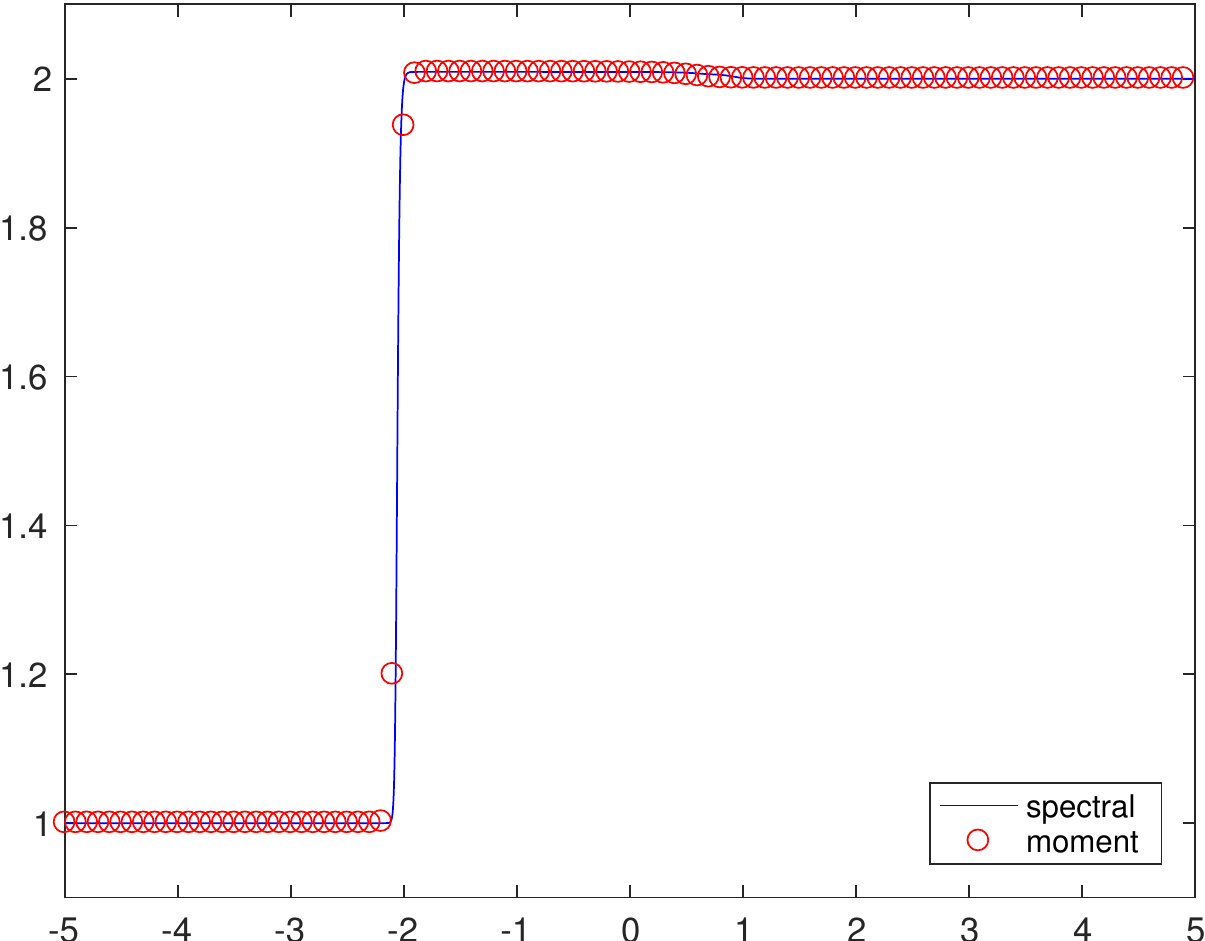}
}
\subfigure[$N=5$]{
\includegraphics[width=0.3\textwidth]{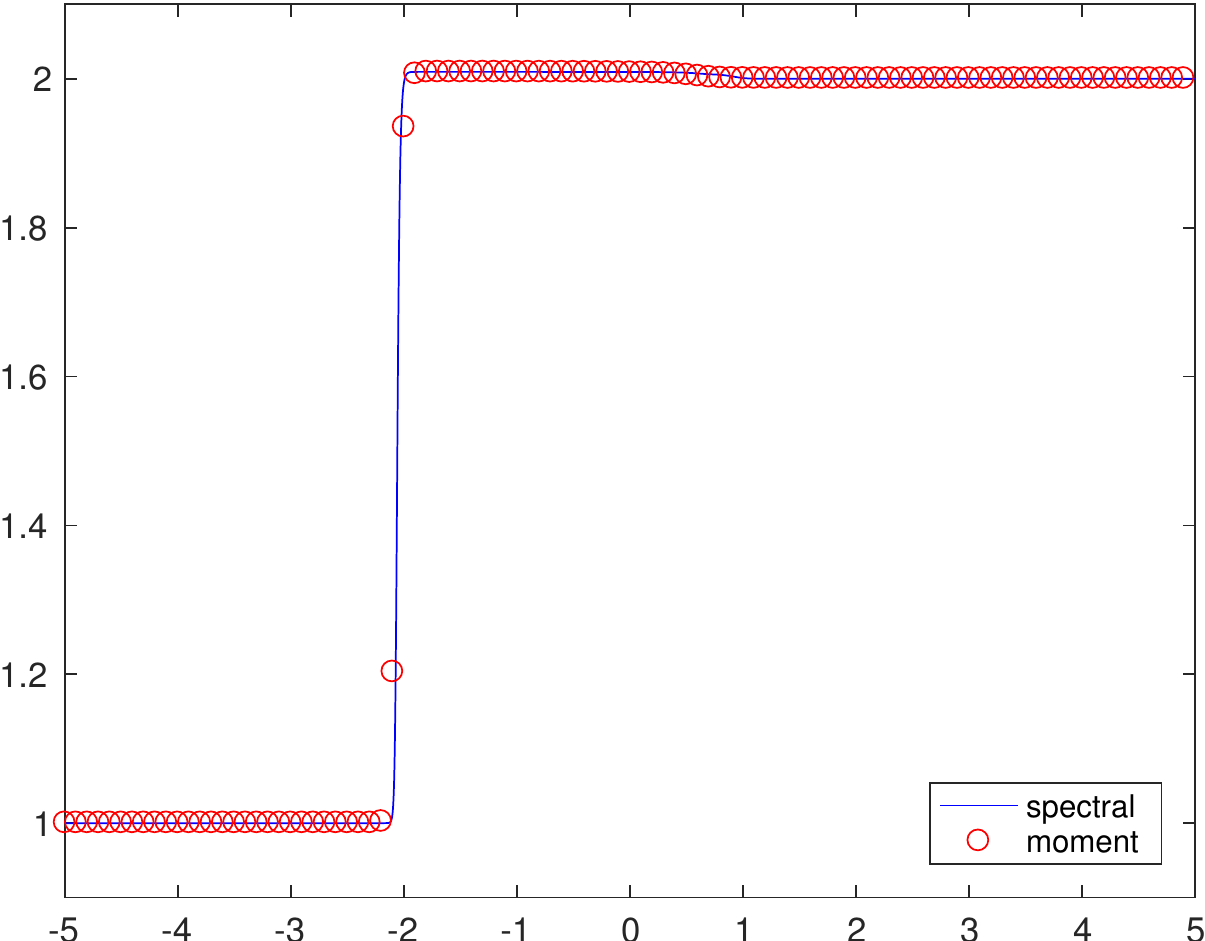}
}
\subfigure[$N=6$]{
\includegraphics[width=0.3\textwidth]{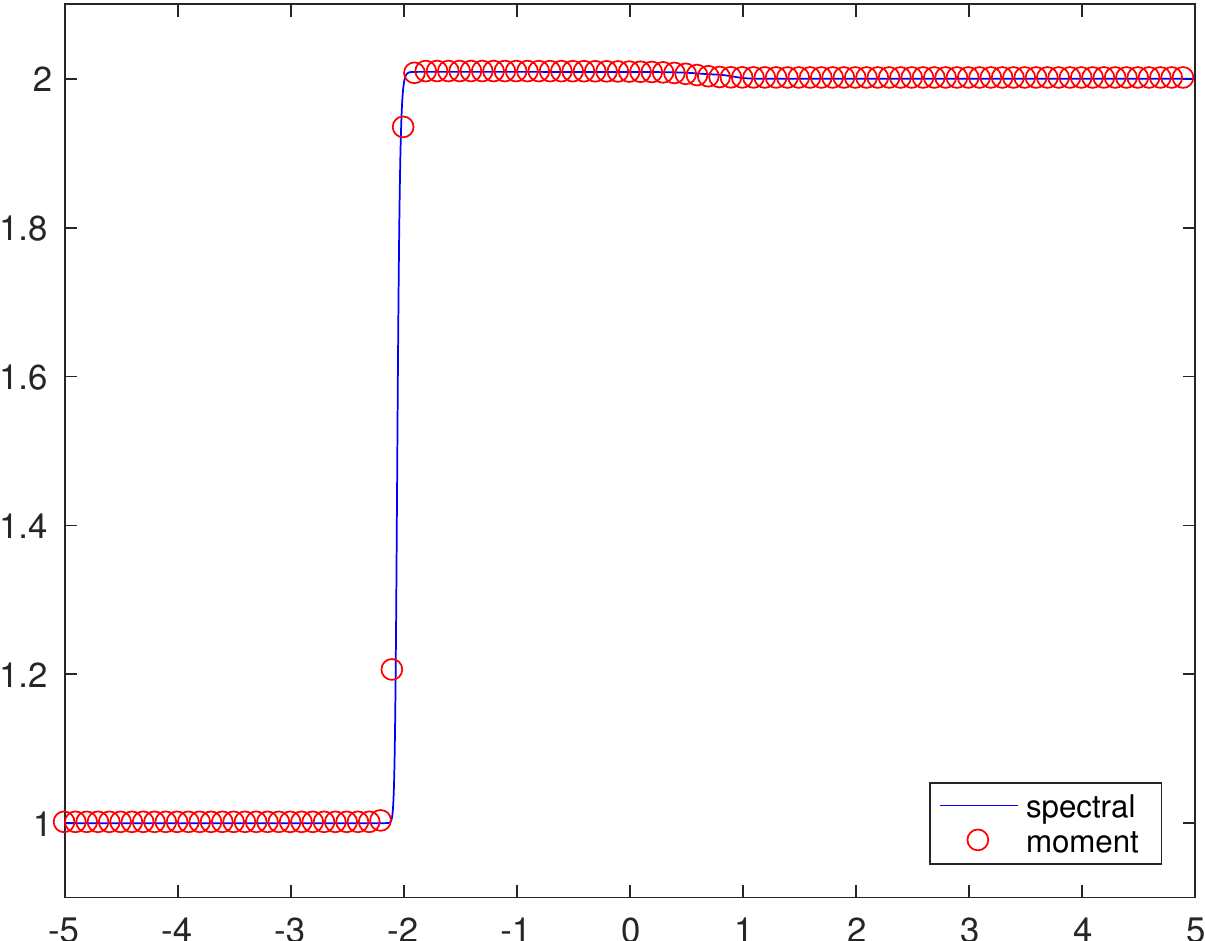}
}
%
%  \centering
%\subfigure[$N=7$]{
%\includegraphics[width=0.3\textwidth]{./images/riemann2/cmpN_e0010/r2_cmpN_rho_n2000_e0010_N7_circle.pdf}
%}
%\subfigure[$N=8$]{
%\includegraphics[width=0.3\textwidth]{./images/riemann2/cmpN_e0010/r2_cmpN_rho_n2000_e0010_N8_circle.pdf}
%}
%\subfigure[$N=9$]{
%\includegraphics[width=0.3\textwidth]{./images/riemann2/cmpN_e0010/r2_cmpN_rho_n2000_e0010_N9_circle.pdf}
%}
\caption{Same as Fig. \ref{figure-example2-1} except for  $\varepsilon=0.01$.}\label{figure-example2-3}
%使用4000个网格Spectral method的结果作为参考,
%和使用2000个网格的矩方法的结果进行比较}
\end{figure}

\begin{figure}
  \centering
\subfigure[$N=1$]{
\includegraphics[width=0.3\textwidth]{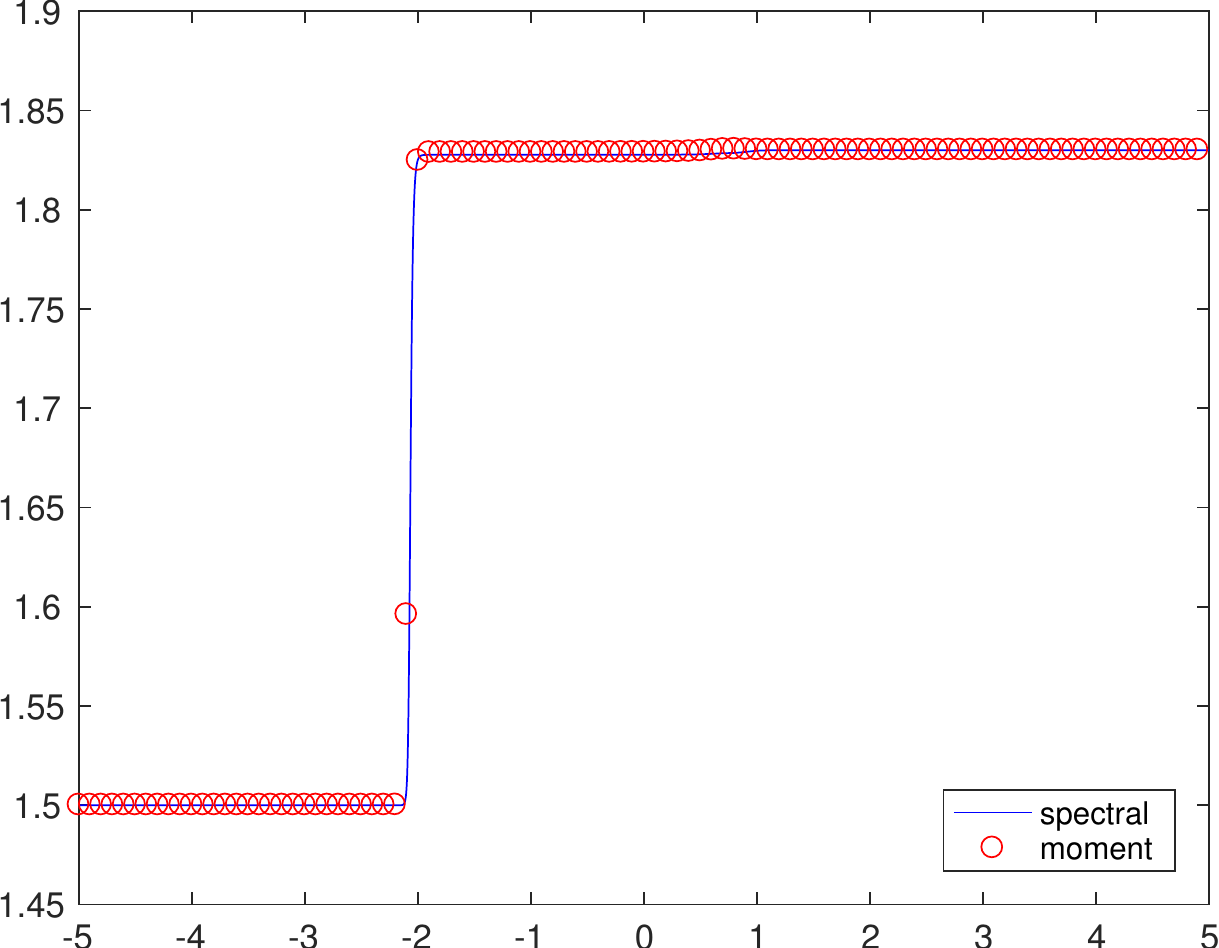}
}
\subfigure[$N=2$]{
\includegraphics[width=0.3\textwidth]{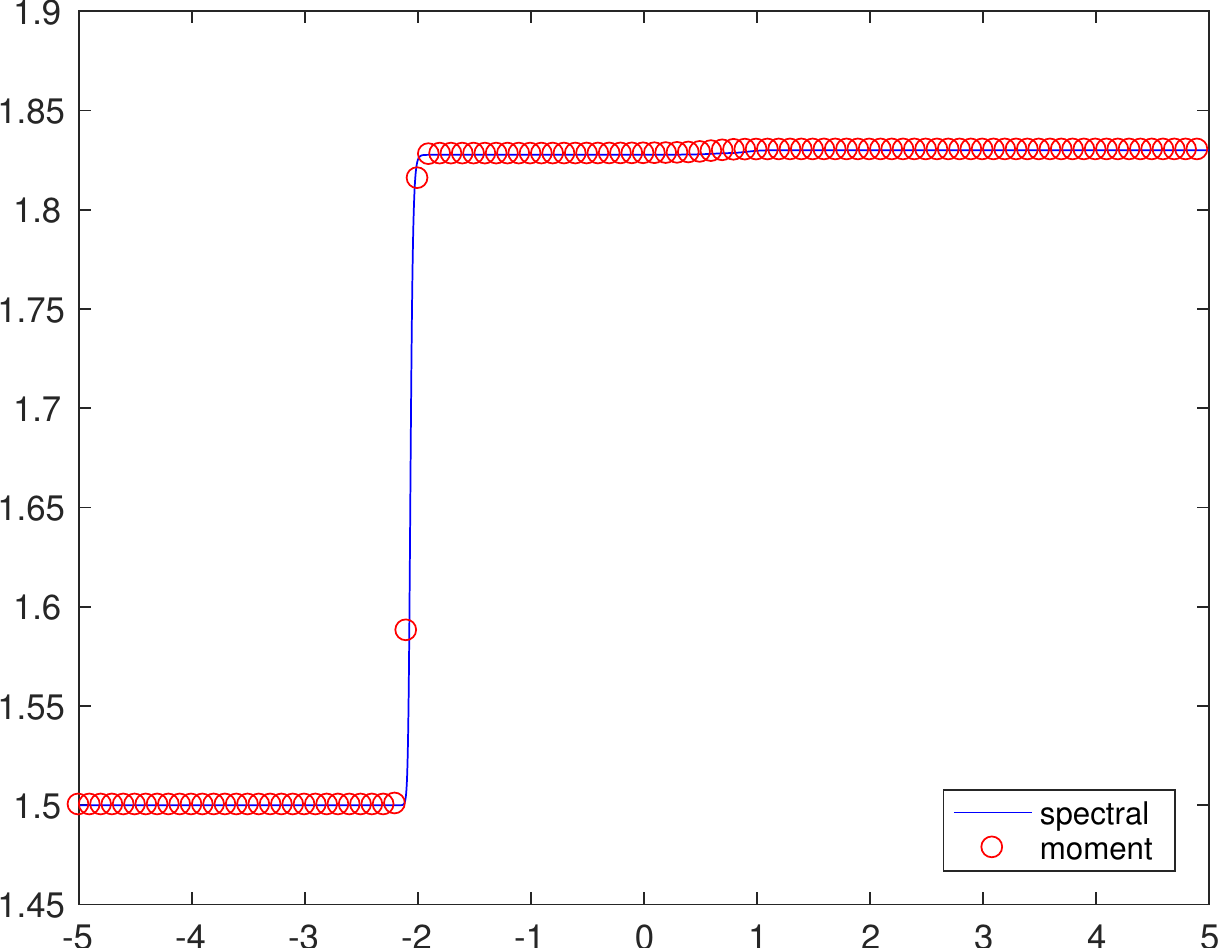}
}
\subfigure[$N=3$]{
\includegraphics[width=0.3\textwidth]{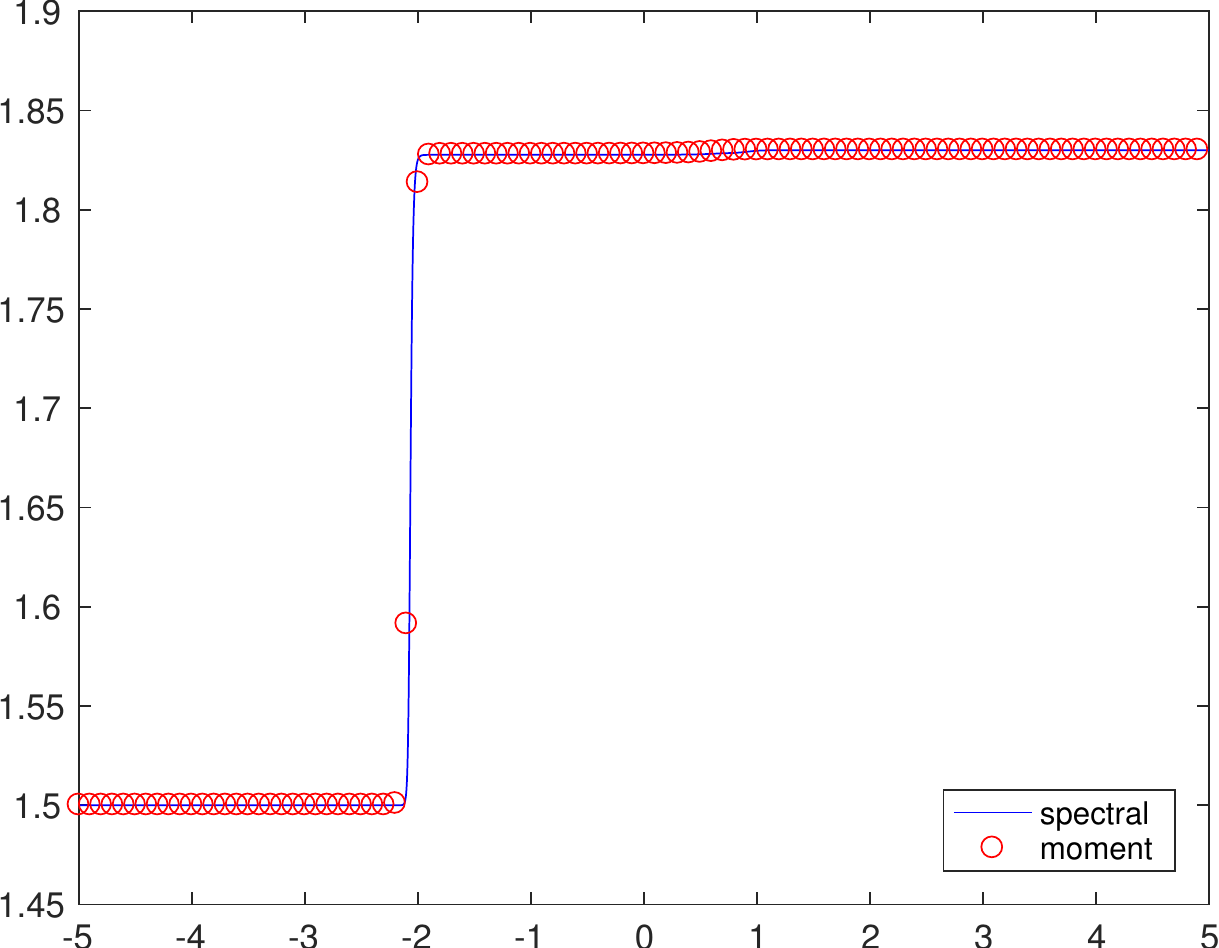}
}

  \centering
\subfigure[$N=4$]{
\includegraphics[width=0.3\textwidth]{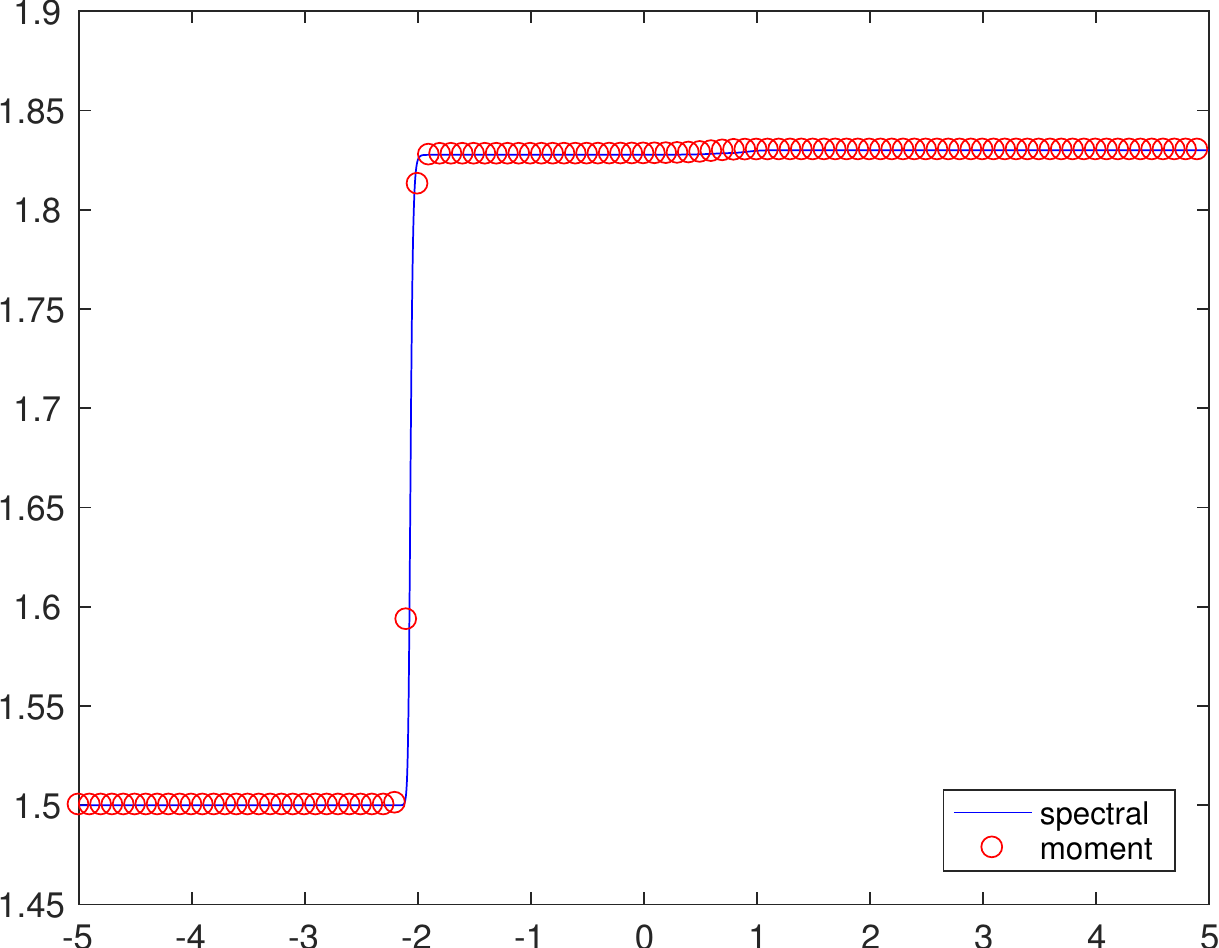}
}
\subfigure[$N=5$]{
\includegraphics[width=0.3\textwidth]{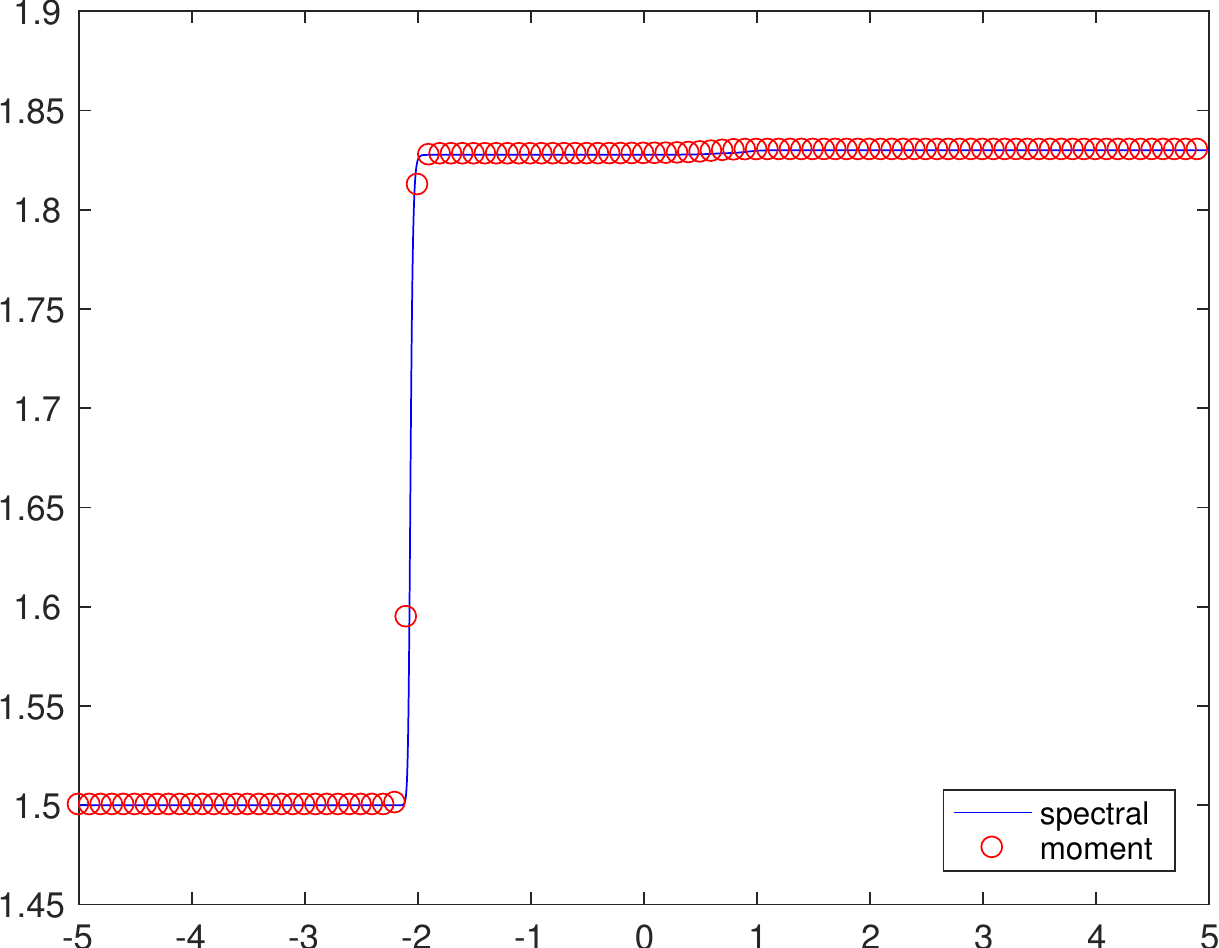}
}
\subfigure[$N=6$]{
\includegraphics[width=0.3\textwidth]{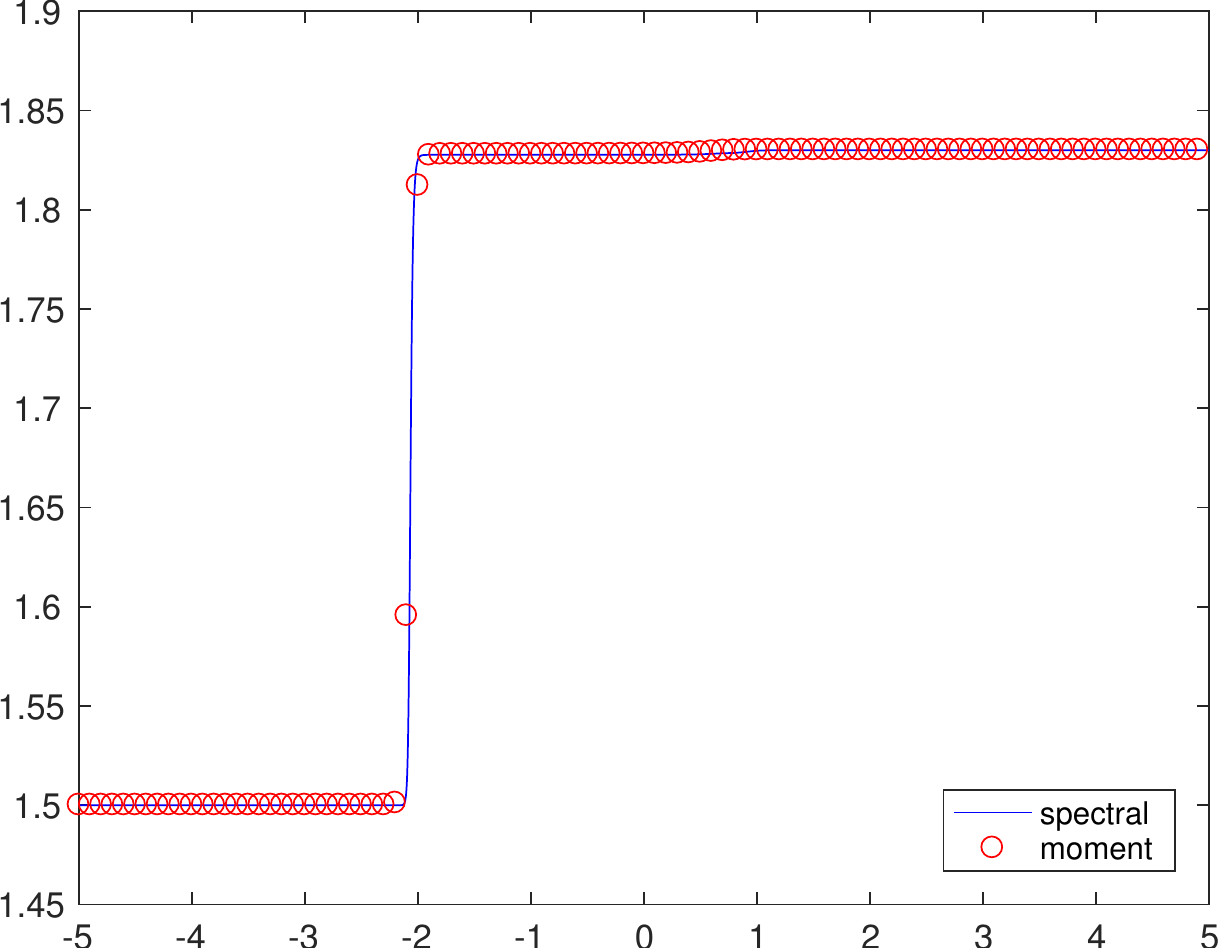}
}

%  \centering
%\subfigure[$N=7$]{
%\includegraphics[width=0.3\textwidth]{./images/riemann2/cmpN_e0010/r2_cmpN_u_n2000_e0010_N7_circle.pdf}
%}
%\subfigure[$N=8$]{
%\includegraphics[width=0.3\textwidth]{./images/riemann2/cmpN_e0010/r2_cmpN_u_n2000_e0010_N8_circle.pdf}
%}
%\subfigure[$N=9$]{
%\includegraphics[width=0.3\textwidth]{./images/riemann2/cmpN_e0010/r2_cmpN_u_n2000_e0010_N9_circle.pdf}
%}
\caption{Same as Fig. \ref{figure-example2-3} except for the macroscopic
velocity angles.}\label{figure-example2-4}
%使用4000个网格Spectral method的结果作为参考,
%和使用2000个网格的矩方法的结果进行比较}
\end{figure}

\begin{example}[Shock wave]\label{example:5.2}
The initial data of the second Riemann problem are
\begin{equation*}
(\rho^\varepsilon,\bar\theta^\varepsilon)=\begin{cases}
 (1, 1.5), & x<0,\\ (2, 1.83),  & x>0. \end{cases}
 \end{equation*}
Figs. \ref{figure-example2-1} and \ref{figure-example2-2} give
the densities $\rho^\varepsilon$  and macroscopic velocity angles
$\bt^\varepsilon$ at $t=4$ obtained by the moment
method with $N=1,2,\cdots,6$, 2000 cells, and $\varepsilon=1$,
where the solid line denotes the reference solution obtained by using the
spectral method  with 4000 cells.  Figs. \ref{figure-example2-3} and \ref{figure-example2-4}
display  corresponding solutions for the case of
$\varepsilon=0.01$. It is observed that a shock wave solution is generated
 and the solutions of moment method $(\rho^\varepsilon, \bt^\varepsilon)$
do converge  the shock profile as $\varepsilon$ becomes small, and
the solutions  of the moment
system well agree with the reference when $N$ is larger than 1.

\end{example}
%%%%%%%%%%%%%%%%%%%%%%%% cmpN_e1000 %%%%%%%%%%%%%%%%%%%%%%%%
%%%%%%%%%%%%%%%%%%%%%%%% r2 rho circle %%%%%%%%%%%%%%%%%%%%%%%%

%%%%%%%%%%%%%%%%%%%%%%%% cmpN_e1000 %%%%%%%%%%%%%%%%%%%%%%%%
%%%%%%%%%%%%%%%%%%%%%%%% r2 u circle %%%%%%%%%%%%%%%%%%%%%%%%

%%%%%%%%%%%%%%%%%%%%%%%% cmpN_e0010 %%%%%%%%%%%%%%%%%%%%%%%%
%%%%%%%%%%%%%%%%%%%%%%%% r2 rho circle %%%%%%%%%%%%%%%%%%%%%%%%

%%%%%%%%%%%%%%%%%%%%%%%% cmpN_e0010 %%%%%%%%%%%%%%%%%%%%%%%%
%%%%%%%%%%%%%%%%%%%%%%%% r2 u circle %%%%%%%%%%%%%%%%%%%%%%%%

\begin{figure}
  \centering
\subfigure[$N=1$]{
\includegraphics[width=0.3\textwidth]{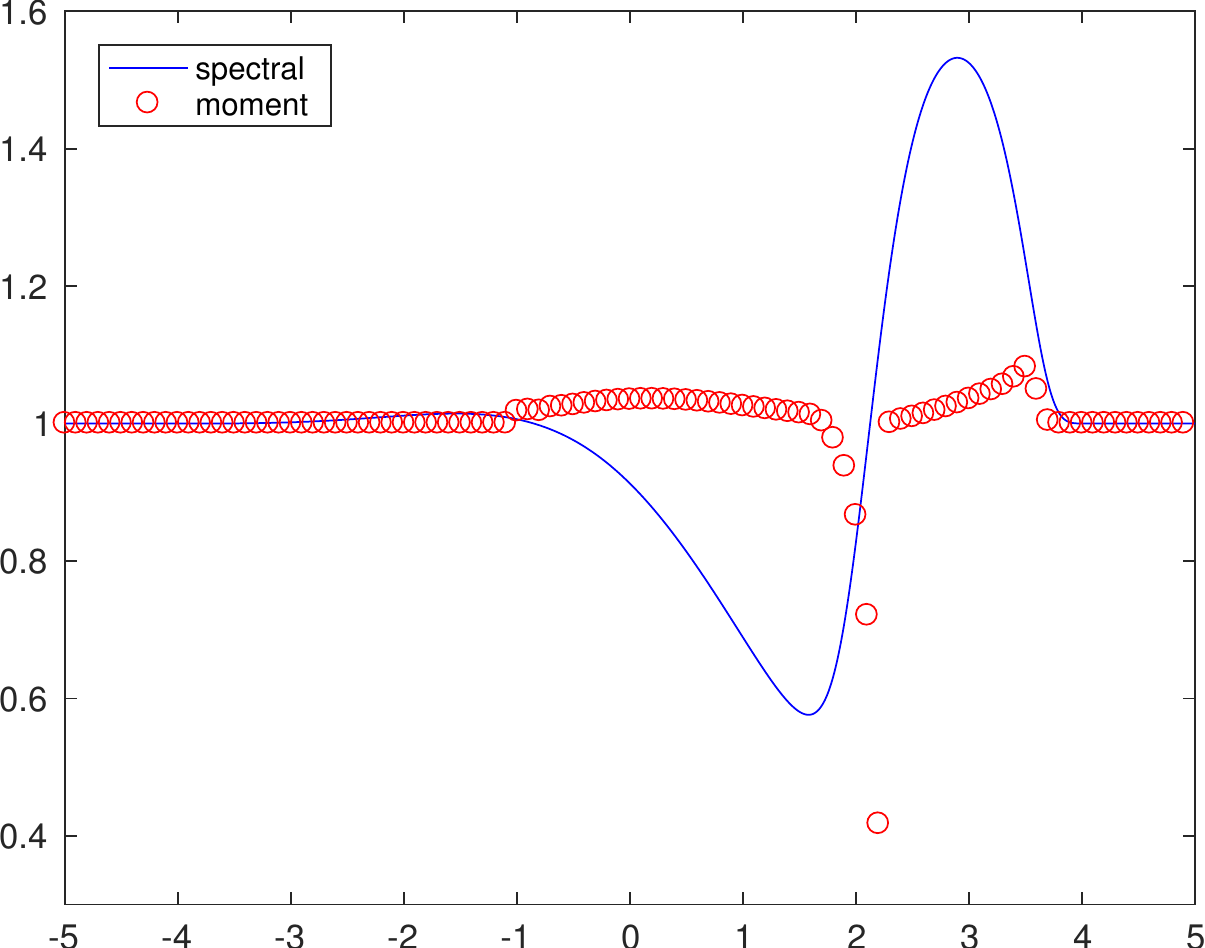}
}
\subfigure[$N=2$]{
\includegraphics[width=0.3\textwidth]{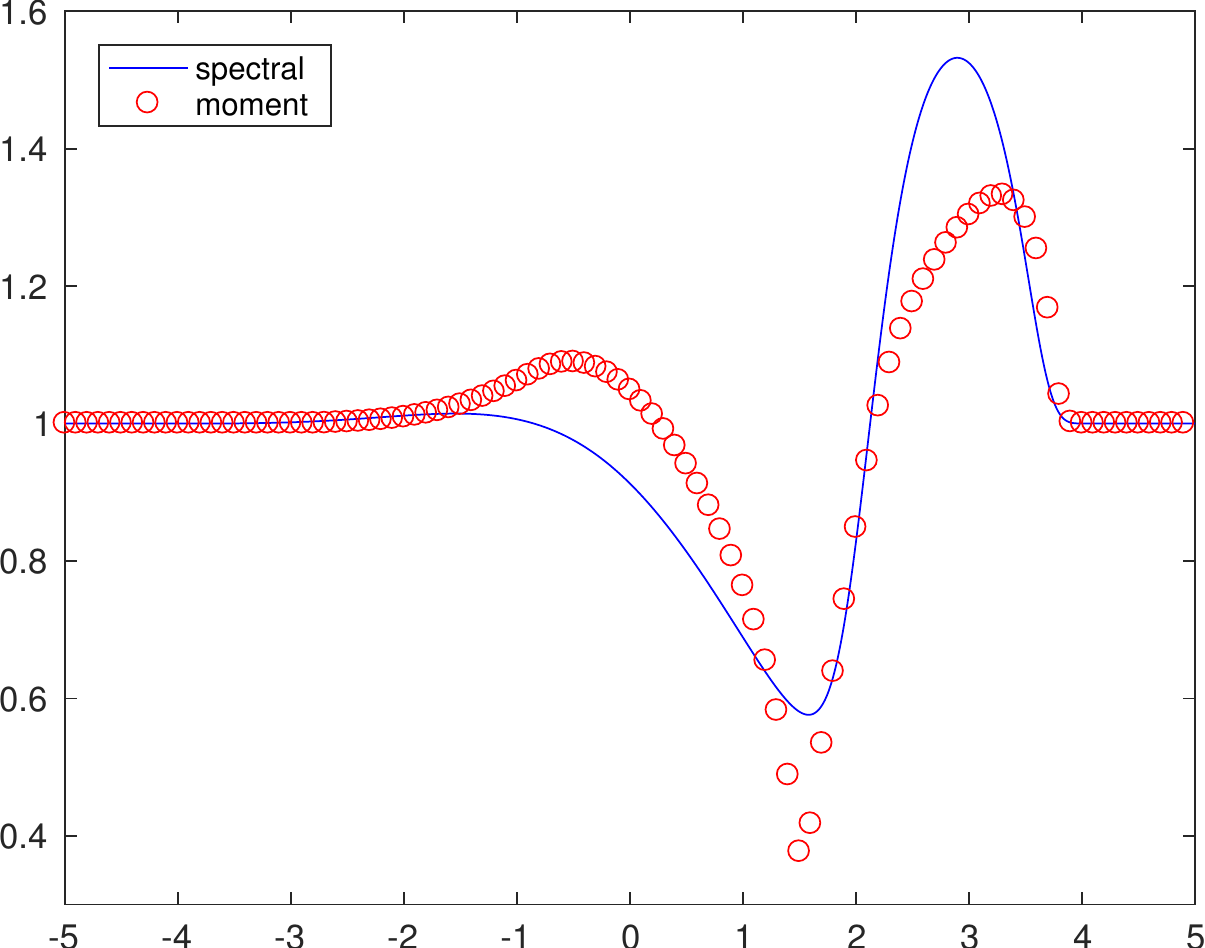}
}
\subfigure[$N=3$]{
\includegraphics[width=0.3\textwidth]{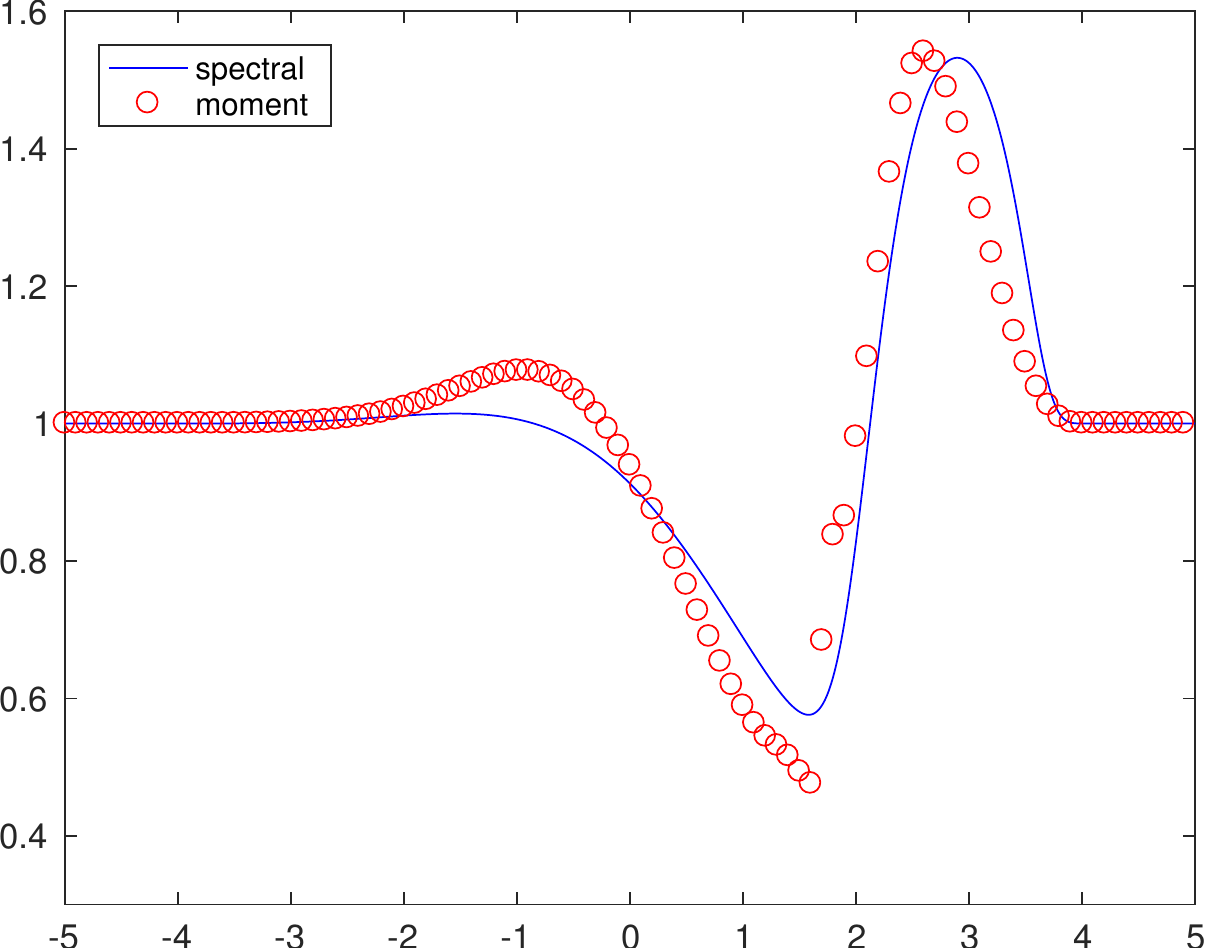}
}

  \centering
\subfigure[$N=4$]{
\includegraphics[width=0.3\textwidth]{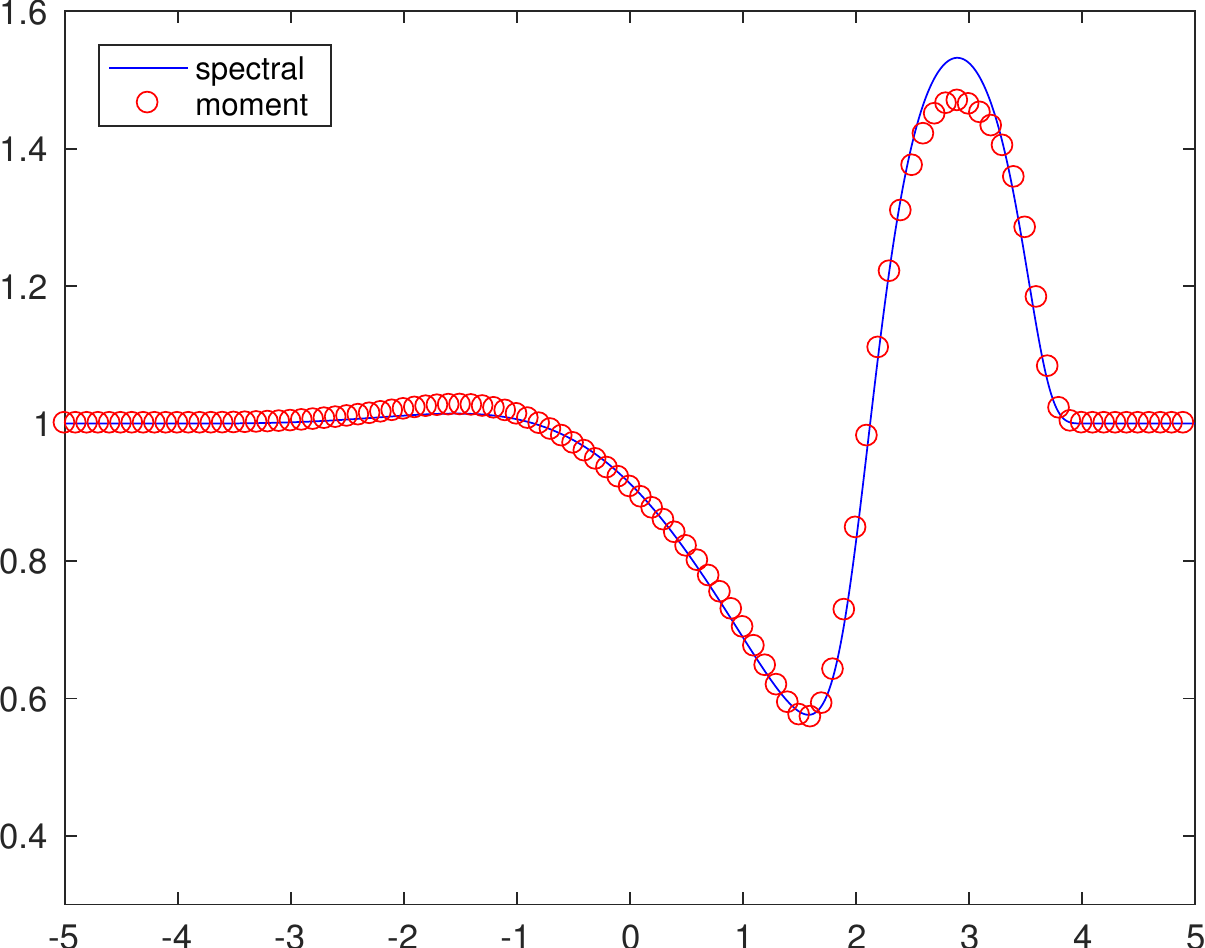}
}
\subfigure[$N=5$]{
\includegraphics[width=0.3\textwidth]{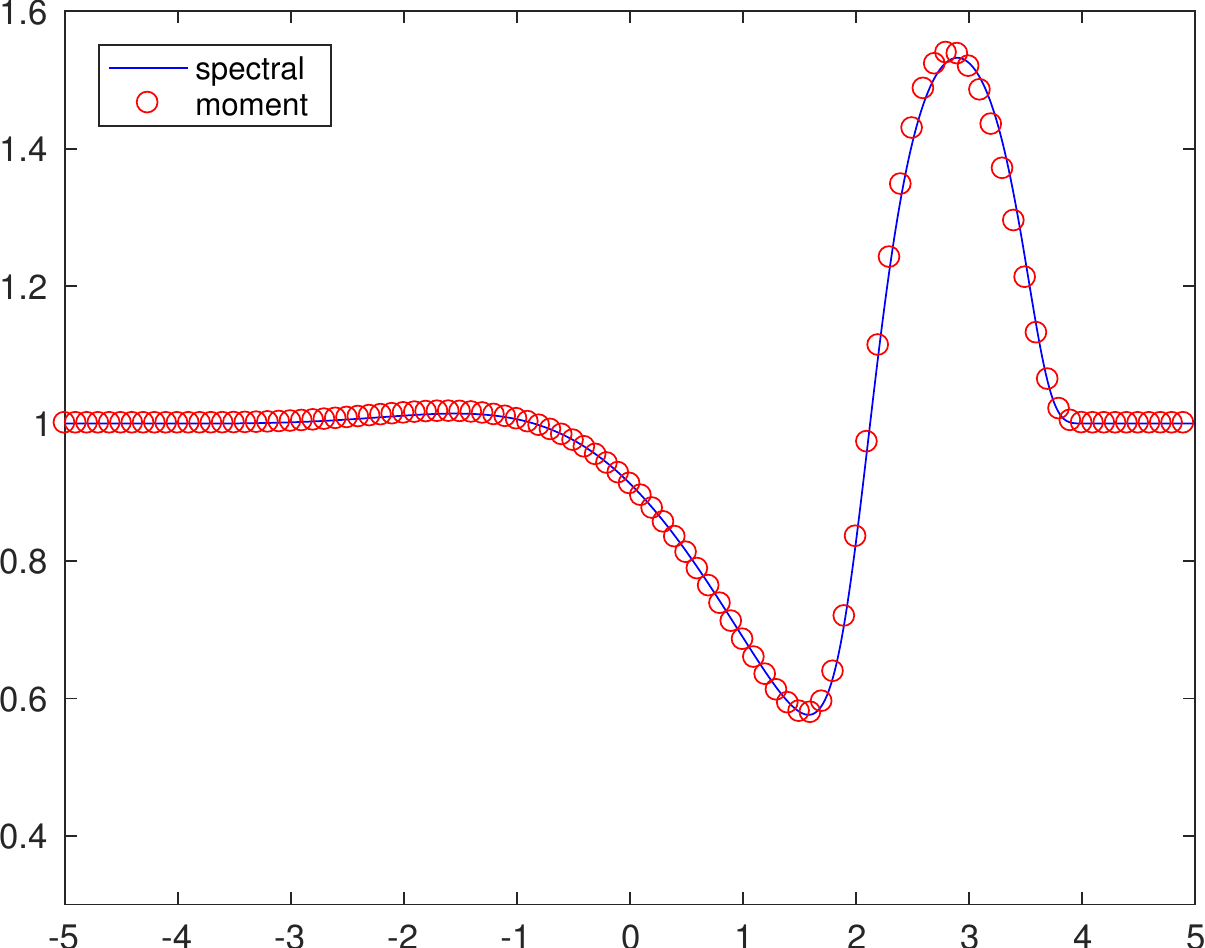}
}
\subfigure[$N=6$]{
\includegraphics[width=0.3\textwidth]{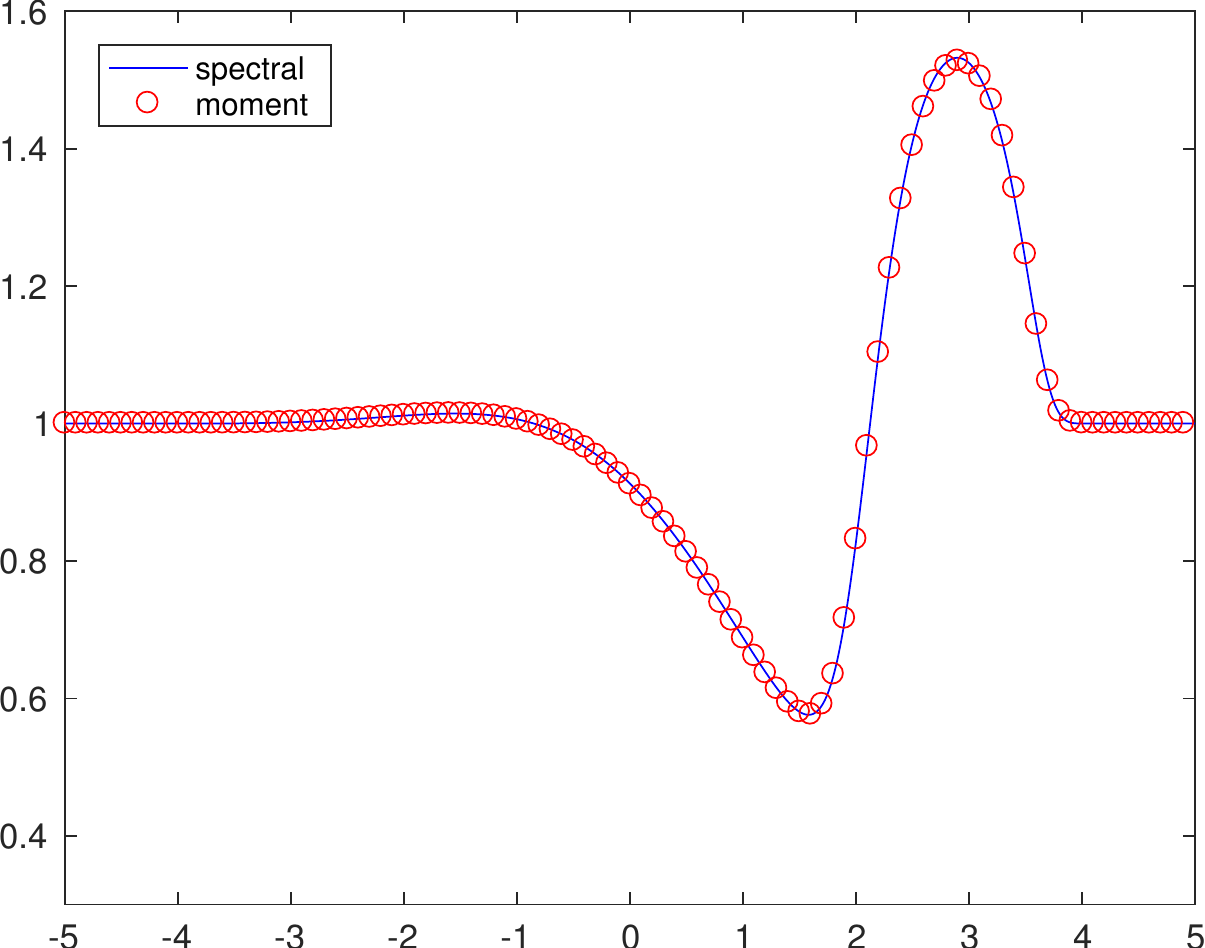}
}

%  \centering
%\subfigure[$N=7$]{
%\includegraphics[width=0.3\textwidth]{./images/riemann3/cmpN_e1000/r3_cmpN_rho_n4000_e1000_N7_circle.pdf}
%}
%\subfigure[$N=8$]{
%\includegraphics[width=0.3\textwidth]{./images/riemann3/cmpN_e1000/r3_cmpN_rho_n4000_e1000_N8_circle.pdf}
%}
%\subfigure[$N=9$]{
%\includegraphics[width=0.3\textwidth]{./images/riemann3/cmpN_e1000/r3_cmpN_rho_n4000_e1000_N9_circle.pdf}
%}
\caption{Example \ref{example:5.3}: The densities at $t=4$ obtained by the moment
method with $N=1,2,\cdots,6$ and 4000 cells. The solid line is the reference solution obtained by using the
spectral method with 8000 cells.  $\varepsilon=1$.} \label{figure-example3-1}
%使用8000个网格Spectral method的结果作为参考,
%和使用4000个网格的矩方法的结果进行比较}
\end{figure}

\begin{figure}
  \centering
\subfigure[$N=1$]{
\includegraphics[width=0.3\textwidth]{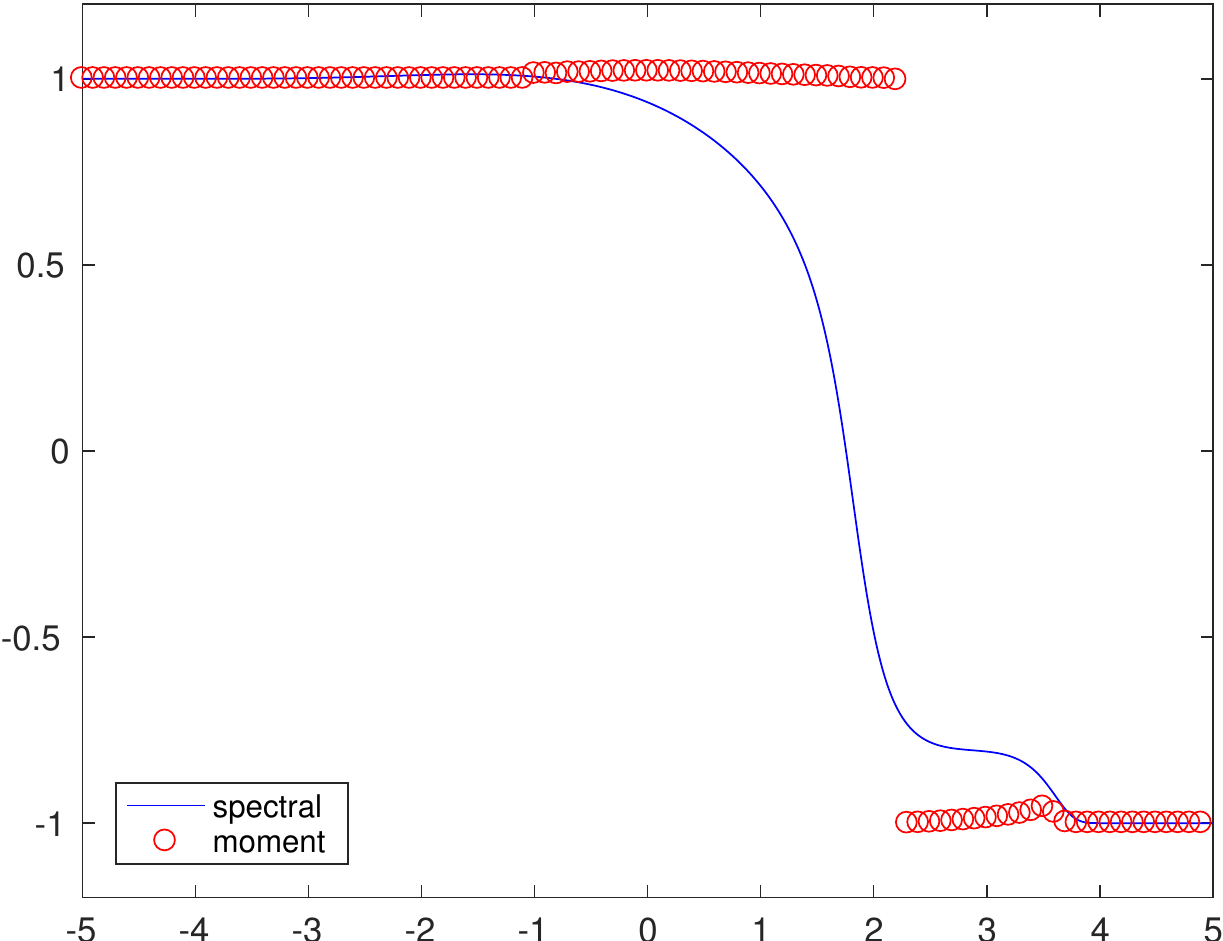}
}
\subfigure[$N=2$]{
\includegraphics[width=0.3\textwidth]{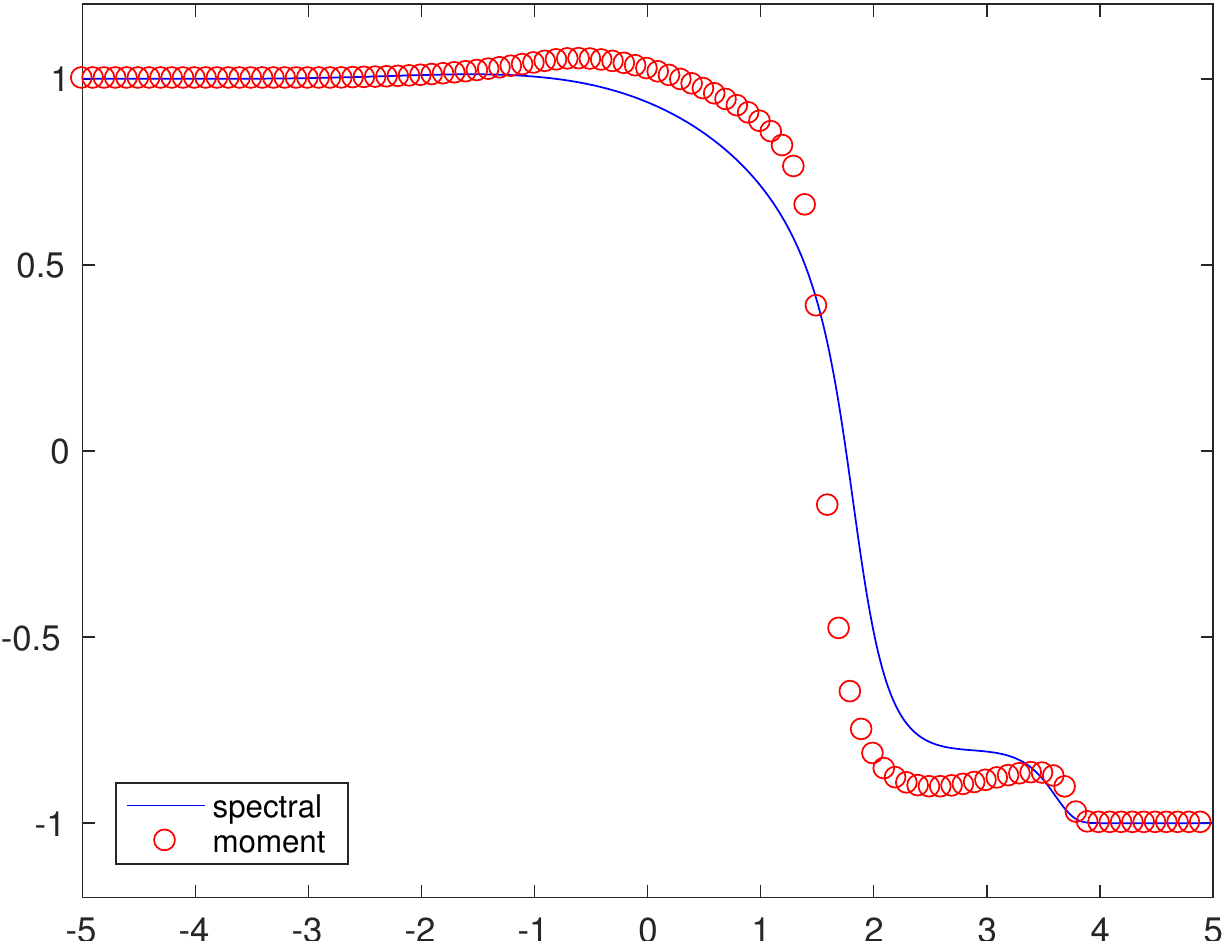}
}
\subfigure[$N=3$]{
\includegraphics[width=0.3\textwidth]{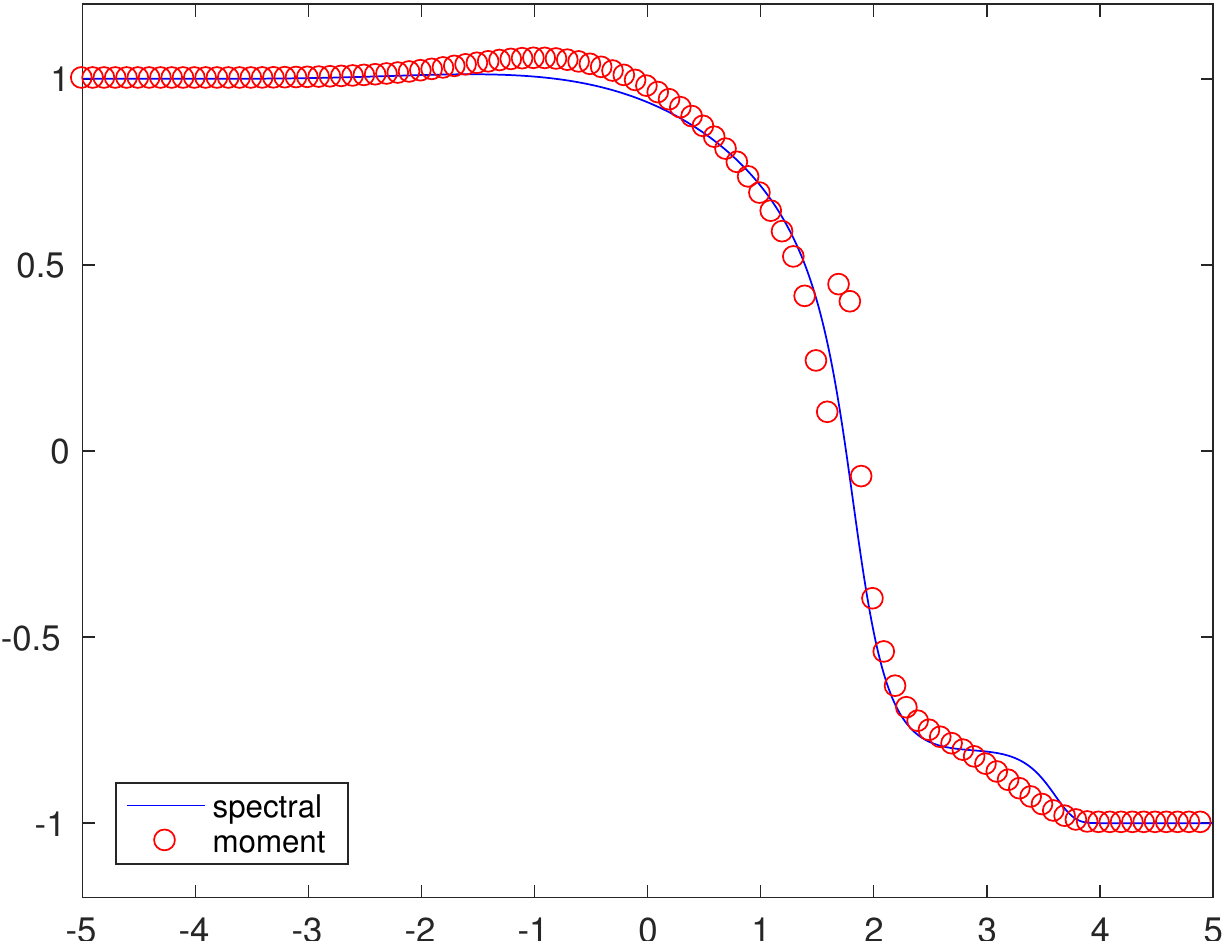}
}

  \centering
\subfigure[$N=4$]{
\includegraphics[width=0.3\textwidth]{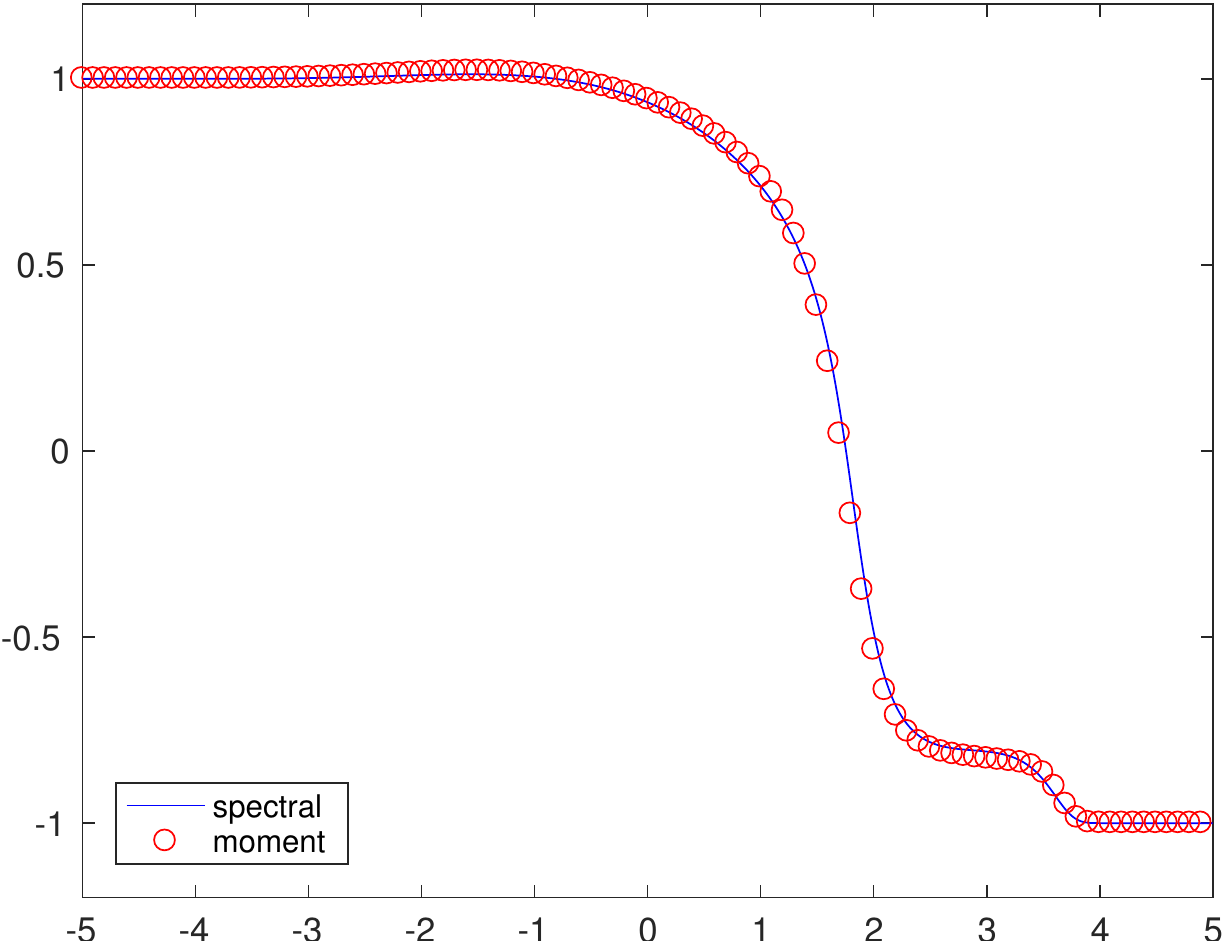}
}
\subfigure[$N=5$]{
\includegraphics[width=0.3\textwidth]{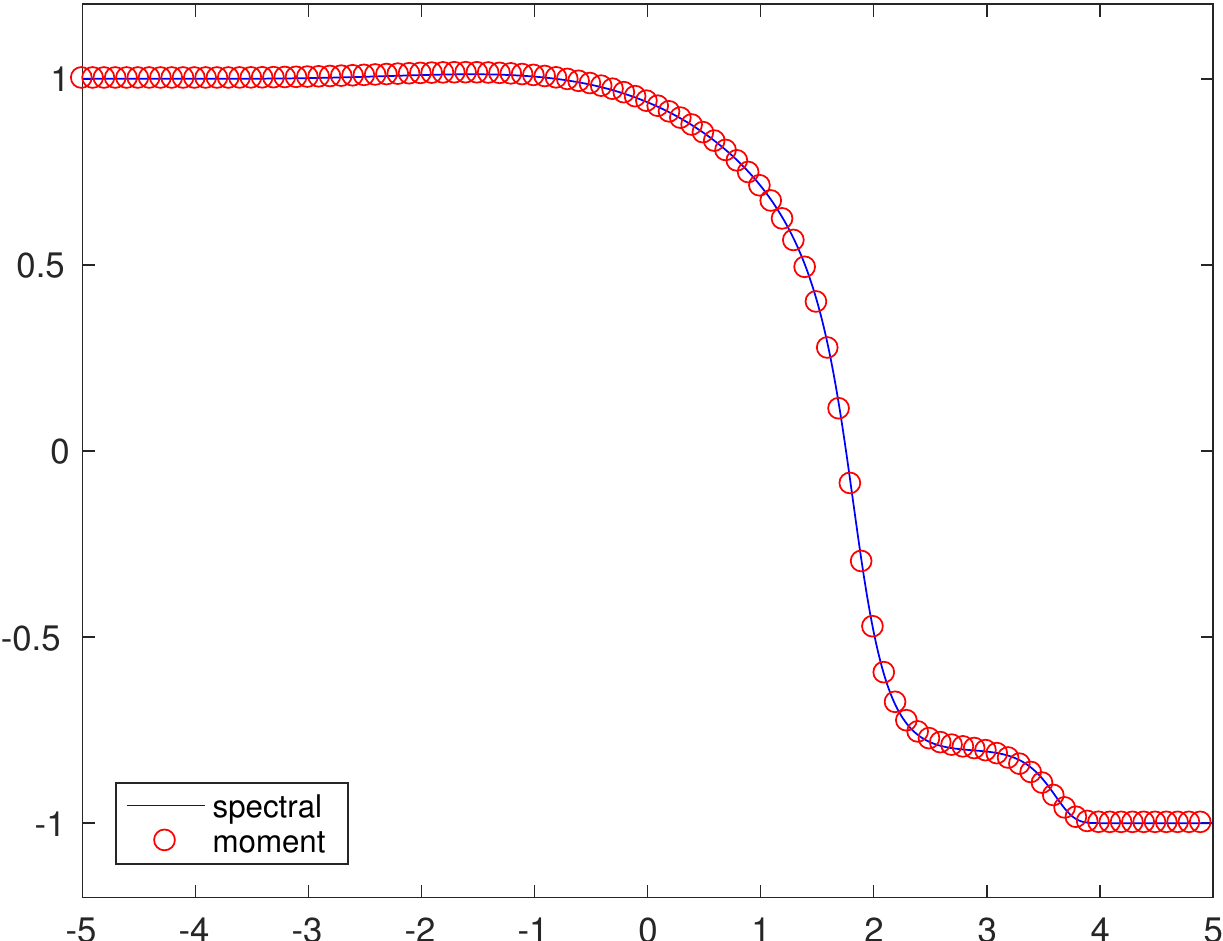}
}
\subfigure[$N=6$]{
\includegraphics[width=0.3\textwidth]{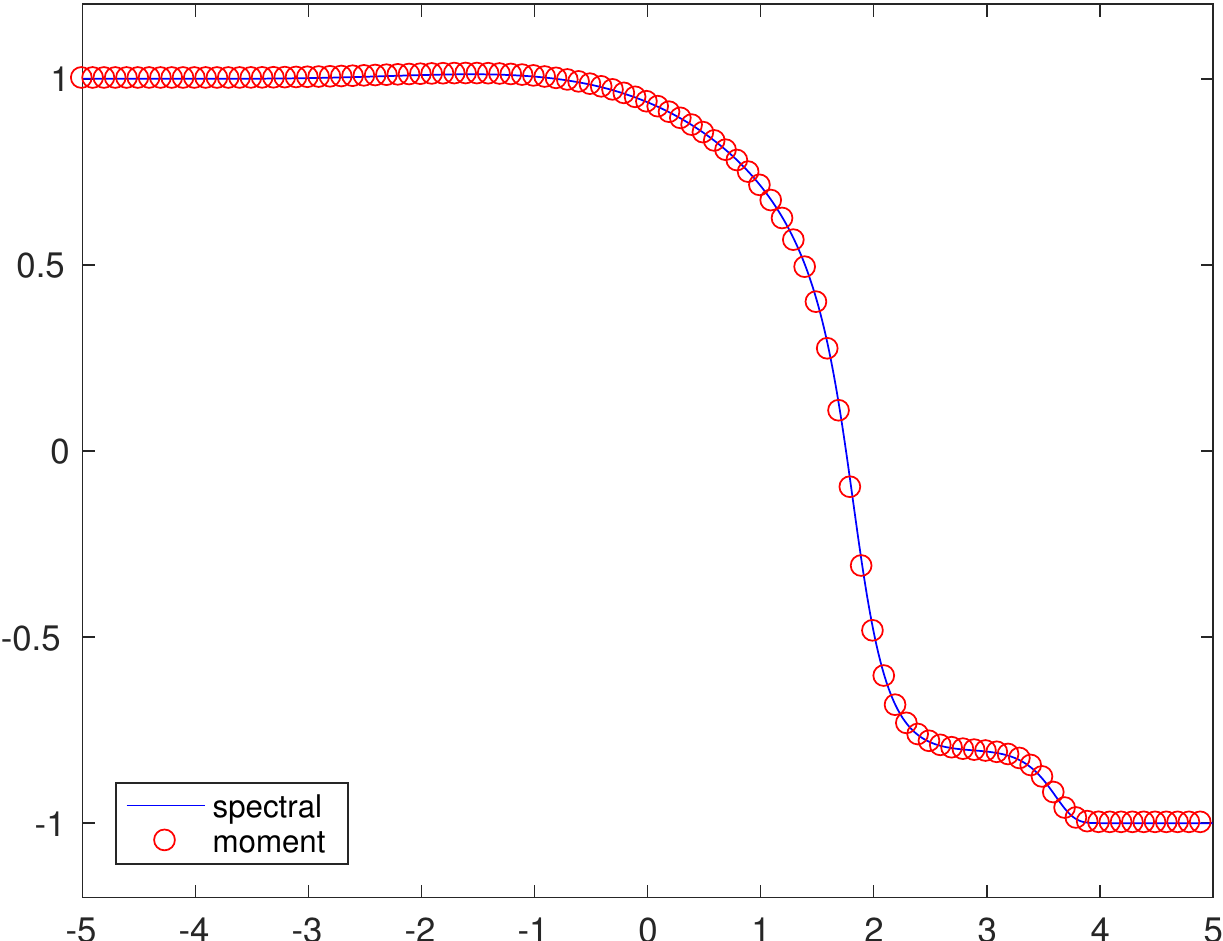}
}
%
%  \centering
%\subfigure[$N=7$]{
%\includegraphics[width=0.3\textwidth]{./images/riemann3/cmpN_e1000/r3_cmpN_u_n4000_e1000_N7_circle.pdf}
%}
%\subfigure[$N=8$]{
%\includegraphics[width=0.3\textwidth]{./images/riemann3/cmpN_e1000/r3_cmpN_u_n4000_e1000_N8_circle.pdf}
%}
%\subfigure[$N=9$]{
%\includegraphics[width=0.3\textwidth]{./images/riemann3/cmpN_e1000/r3_cmpN_u_n4000_e1000_N9_circle.pdf}
%}
\caption{Same as Fig. \ref{figure-example3-1} except for the macroscopic
velocity angles.}\label{figure-example3-2}
%使用8000个网格Spectral method的结果作为参考,
%和使用4000个网格的矩方法的结果进行比较}
\end{figure}

\begin{figure}
  \centering
\subfigure[$N=1$]{
\includegraphics[width=0.3\textwidth]{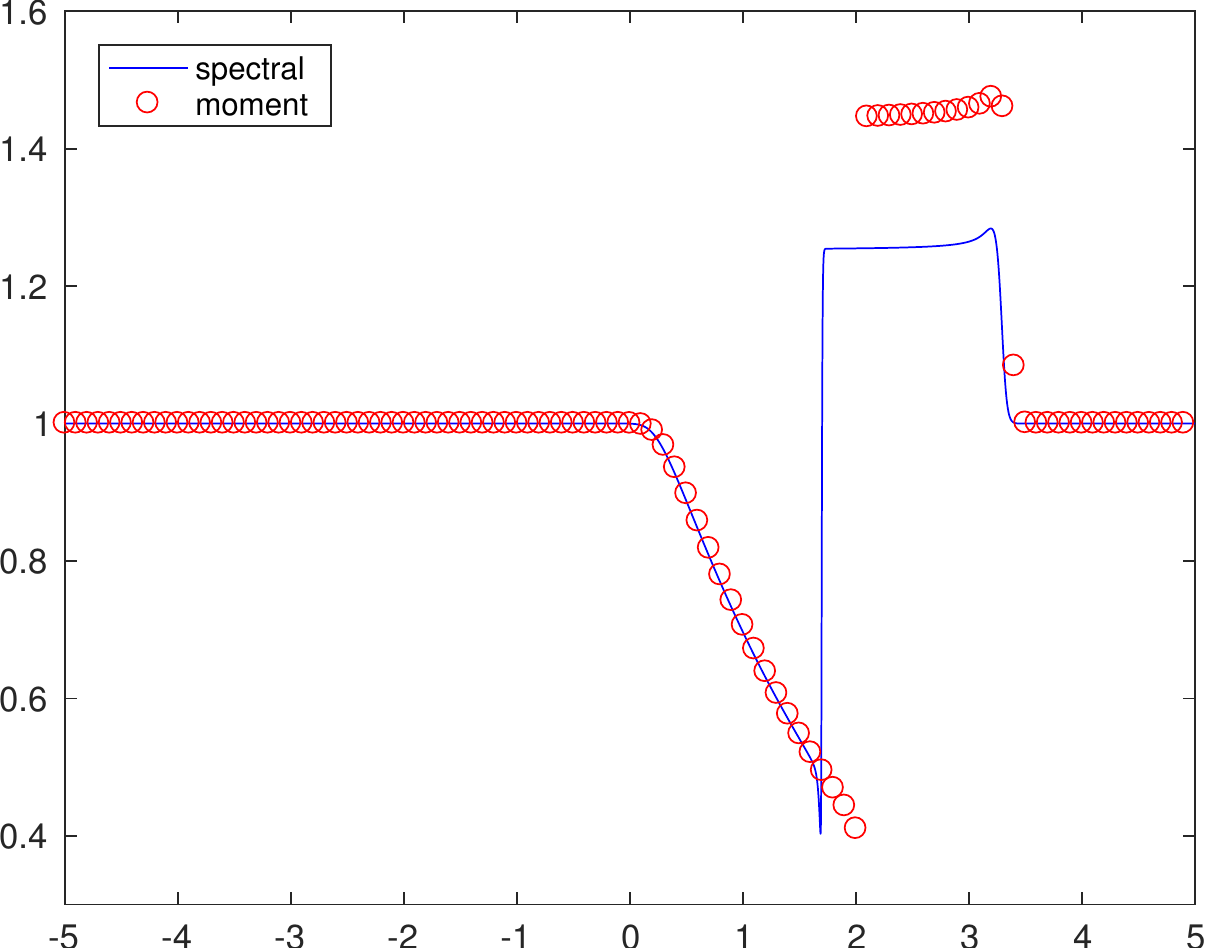}
}
\subfigure[$N=2$]{
\includegraphics[width=0.3\textwidth]{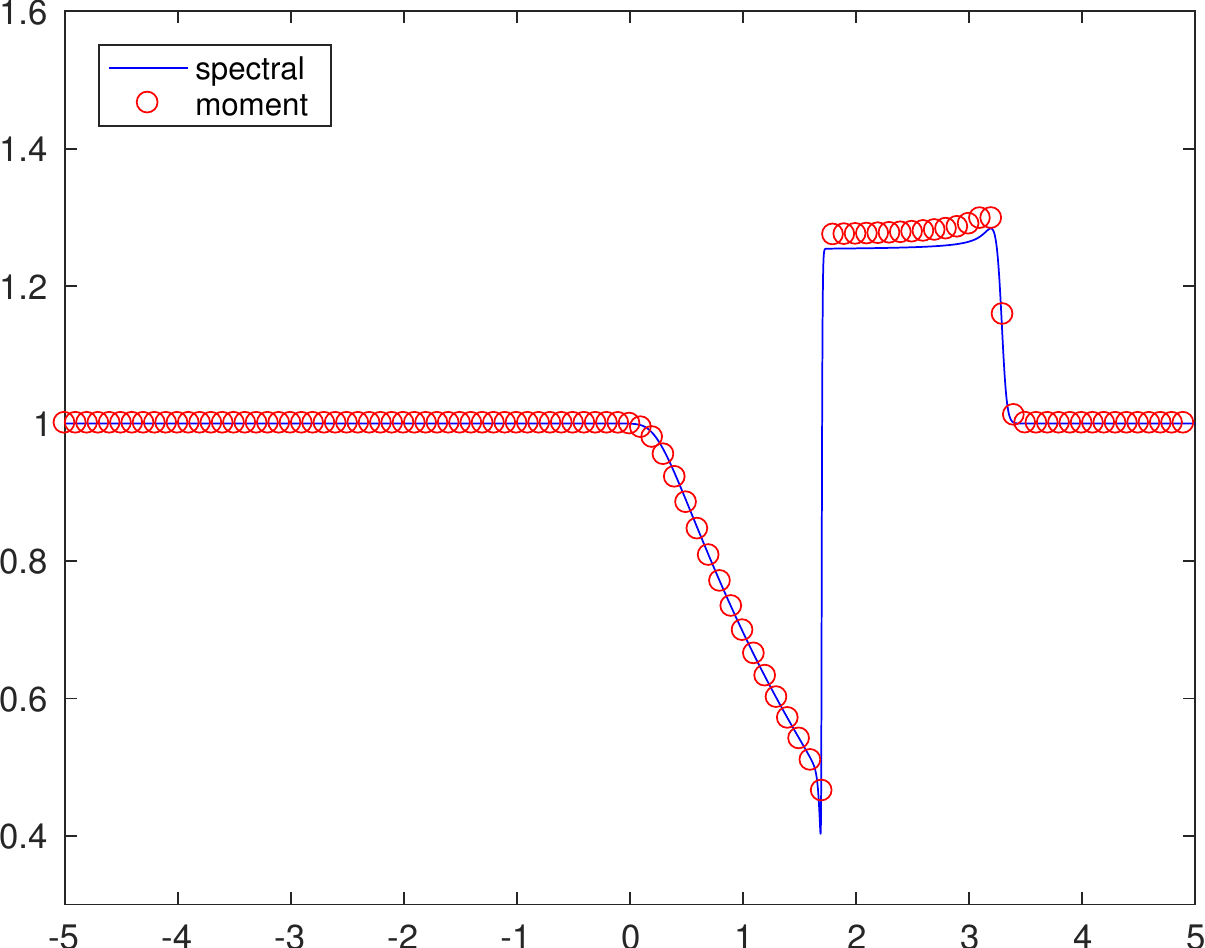}
}
\subfigure[$N=3$]{
\includegraphics[width=0.3\textwidth]{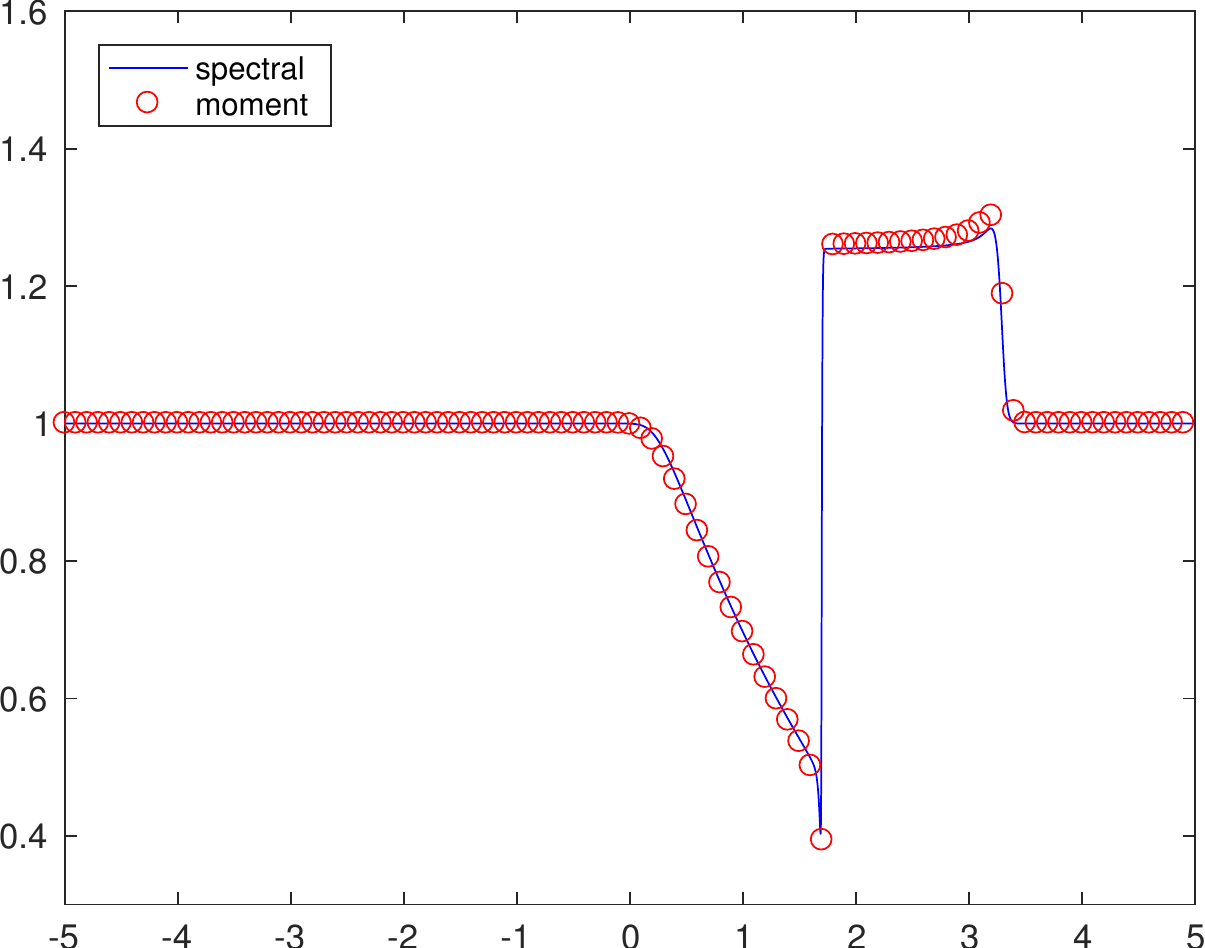}
}

  \centering
\subfigure[$N=4$]{
\includegraphics[width=0.3\textwidth]{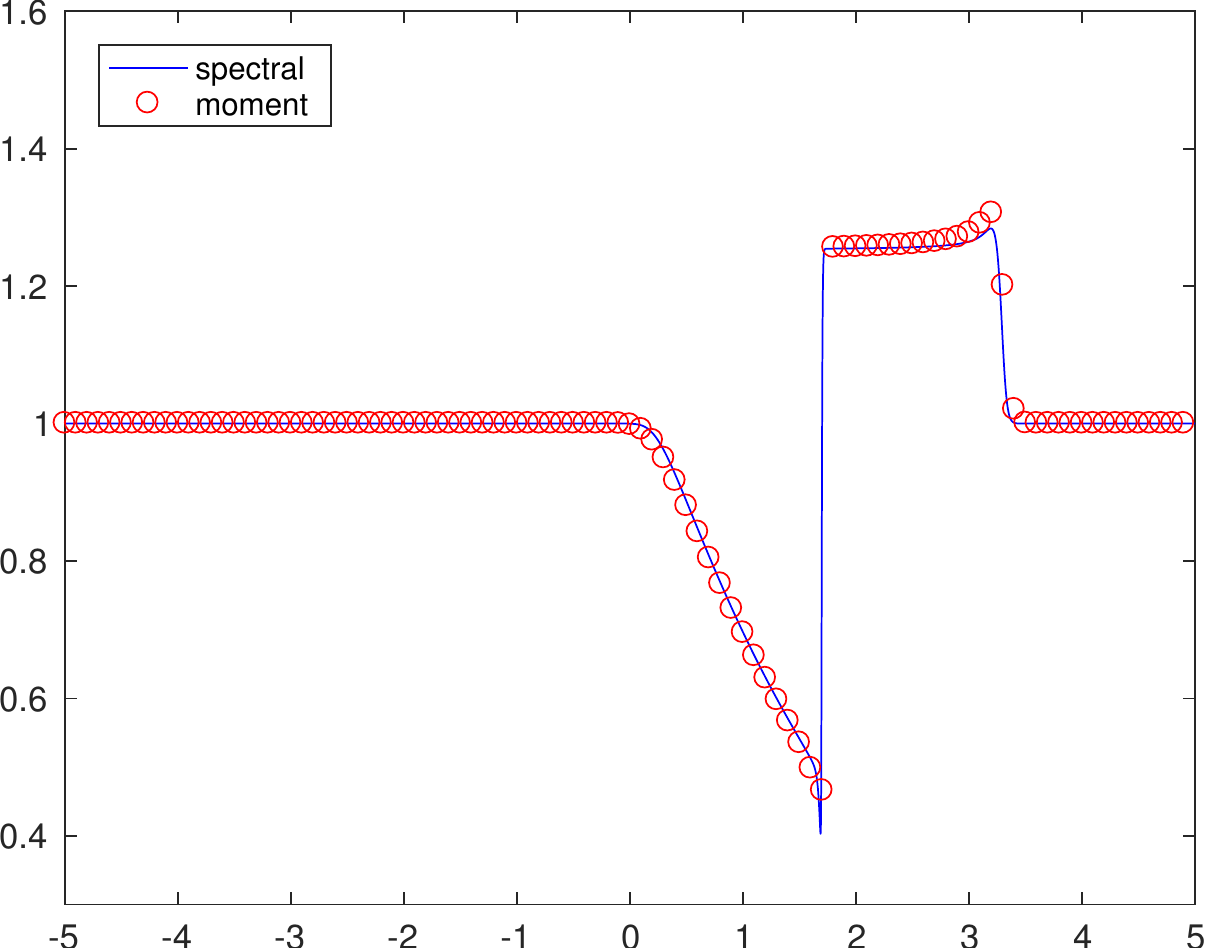}
}
\subfigure[$N=5$]{
\includegraphics[width=0.3\textwidth]{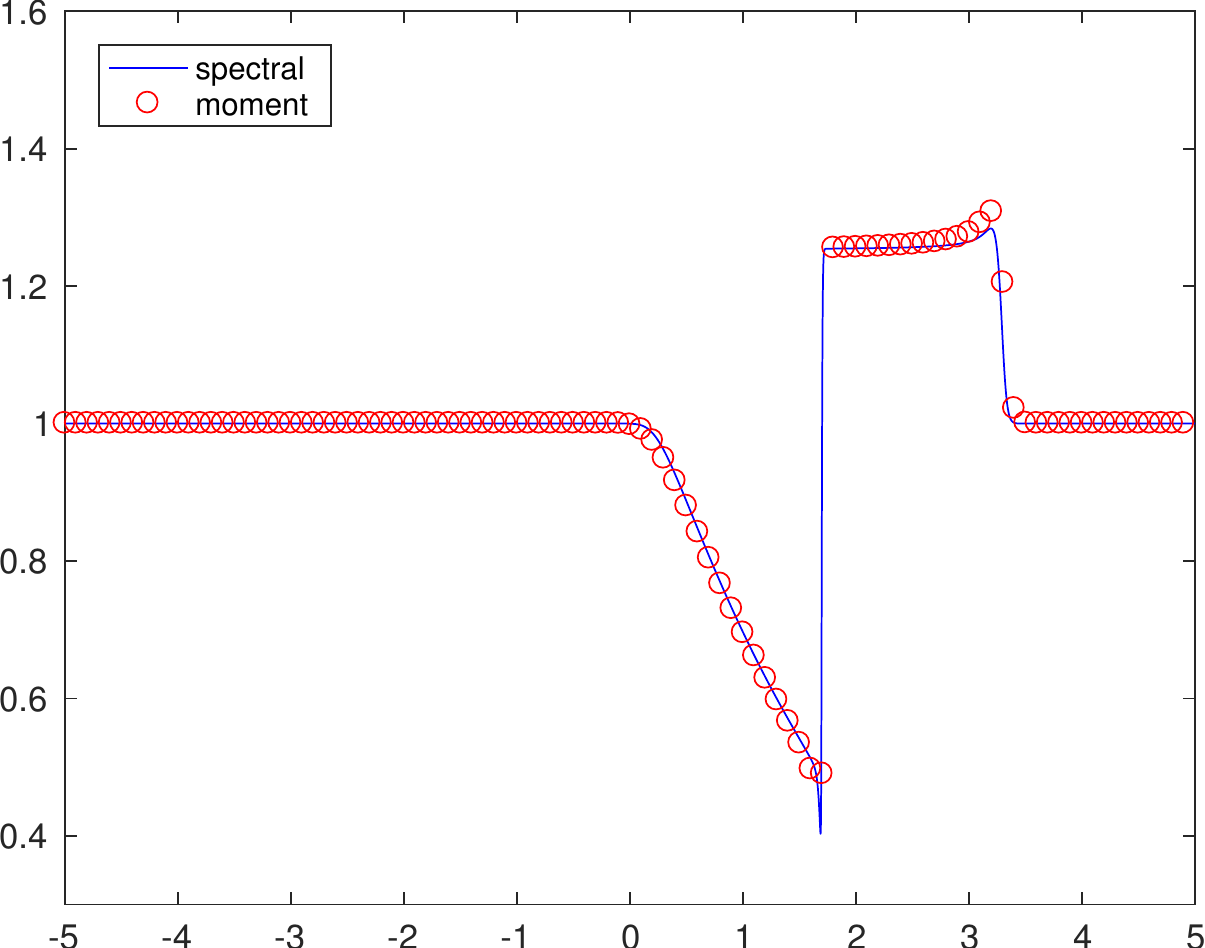}
}
\subfigure[$N=6$]{
\includegraphics[width=0.3\textwidth]{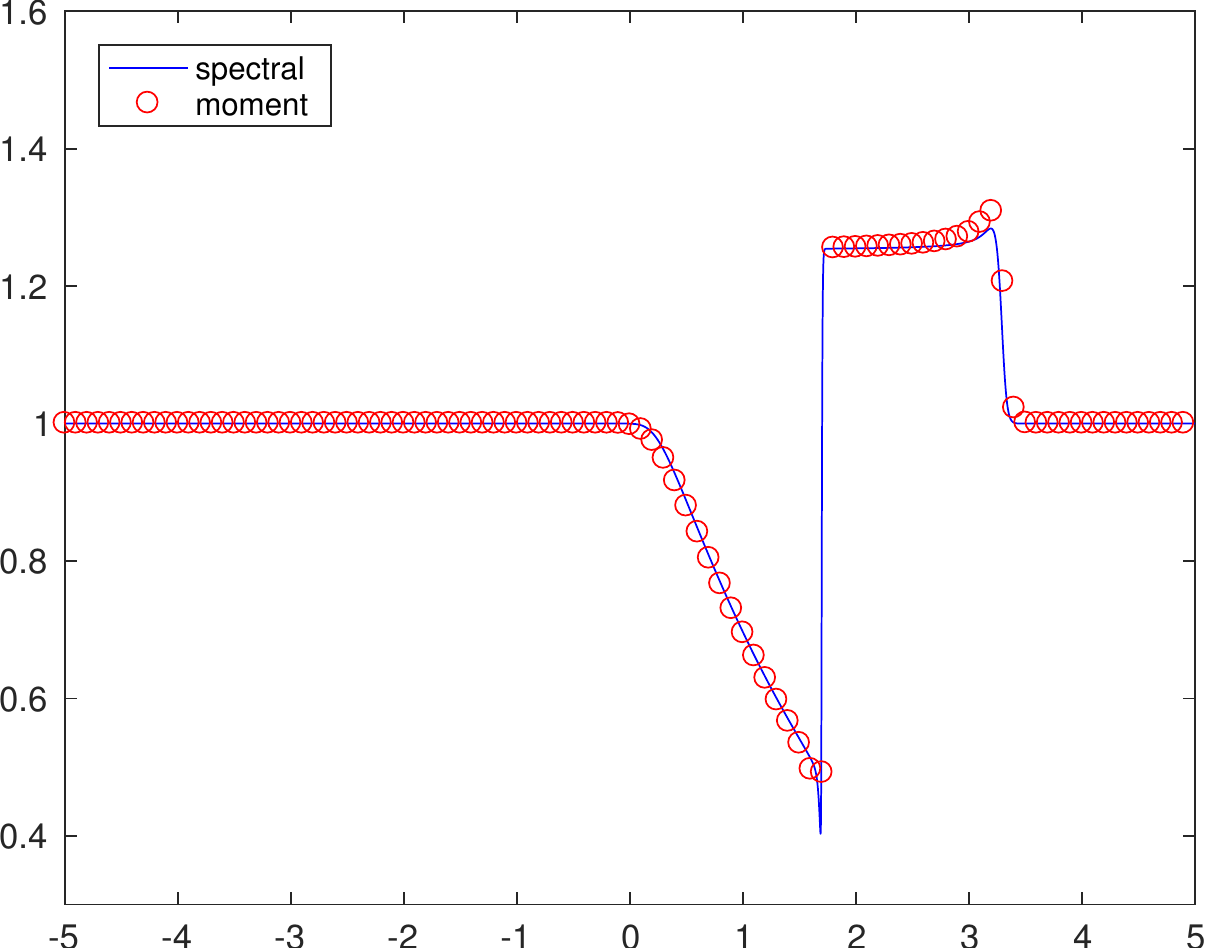}
}

%  \centering
%\subfigure[$N=7$]{
%\includegraphics[width=0.3\textwidth]{./images/riemann3/cmpN_e0010/r3_cmpN_rho_n4000_e0010_N7_circle.pdf}
%}
%\subfigure[$N=8$]{
%\includegraphics[width=0.3\textwidth]{./images/riemann3/cmpN_e0010/r3_cmpN_rho_n4000_e0010_N8_circle.pdf}
%}
%\subfigure[$N=9$]{
%\includegraphics[width=0.3\textwidth]{./images/riemann3/cmpN_e0010/r3_cmpN_rho_n4000_e0010_N9_circle.pdf}
%}
\caption{Same as Fig. \ref{figure-example3-1} except for $\varepsilon=0.01$.}\label{figure-example3-3}
%使用8000个网格Spectral method的结果作为参考,
%和使用4000个网格的矩方法的结果进行比较}
\end{figure}

\begin{figure}
  \centering
\subfigure[$N=1$]{
\includegraphics[width=0.3\textwidth]{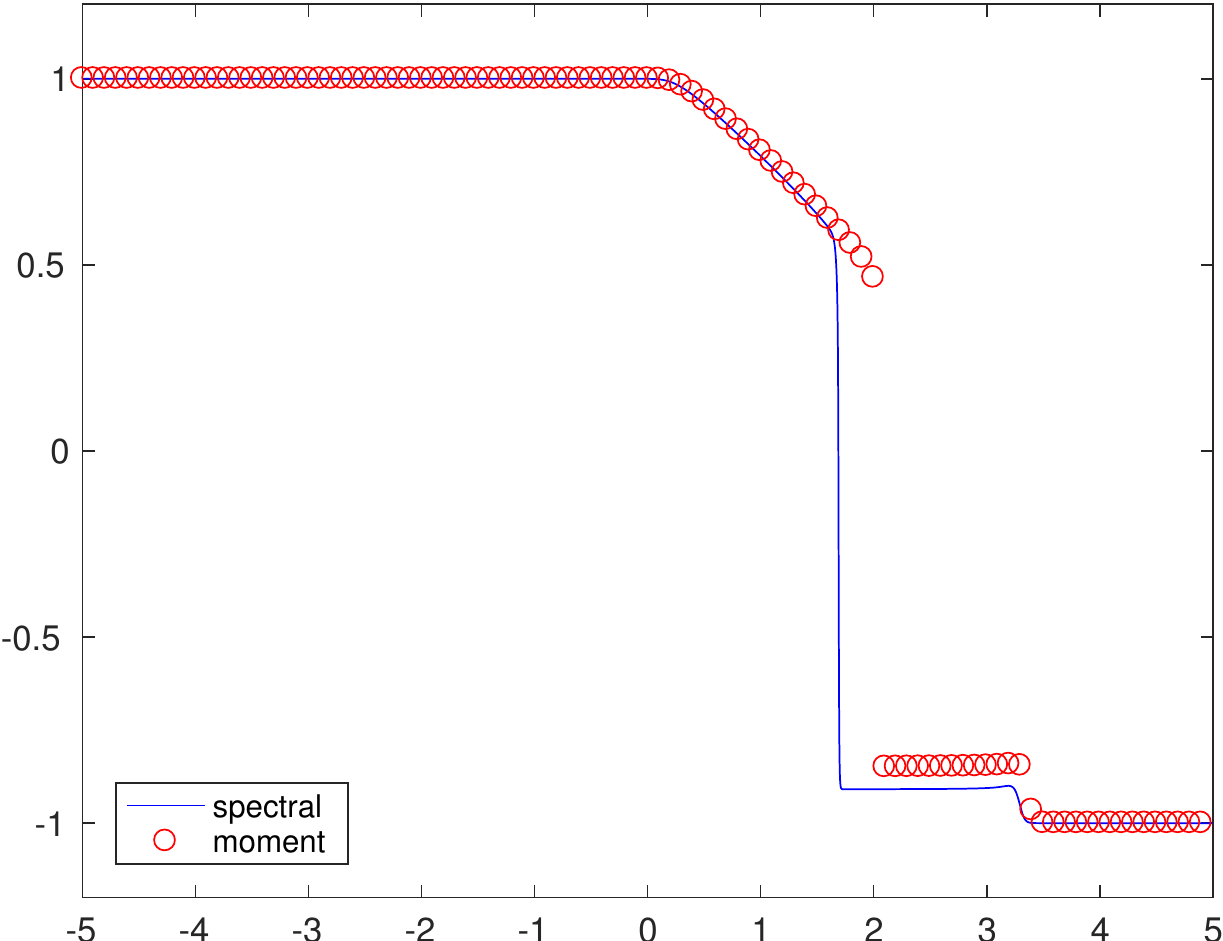}
}
\subfigure[$N=2$]{
\includegraphics[width=0.3\textwidth]{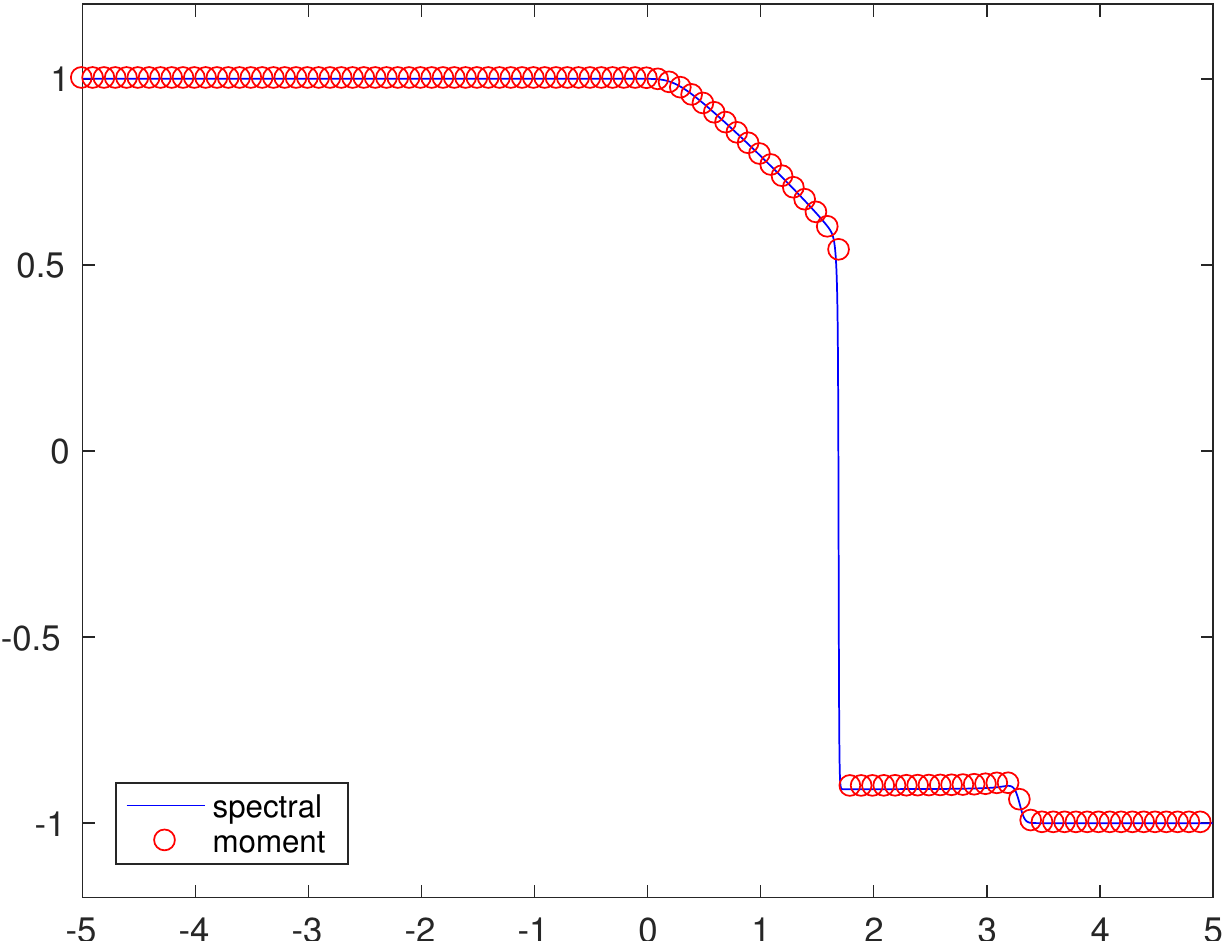}
}
\subfigure[$N=3$]{
\includegraphics[width=0.3\textwidth]{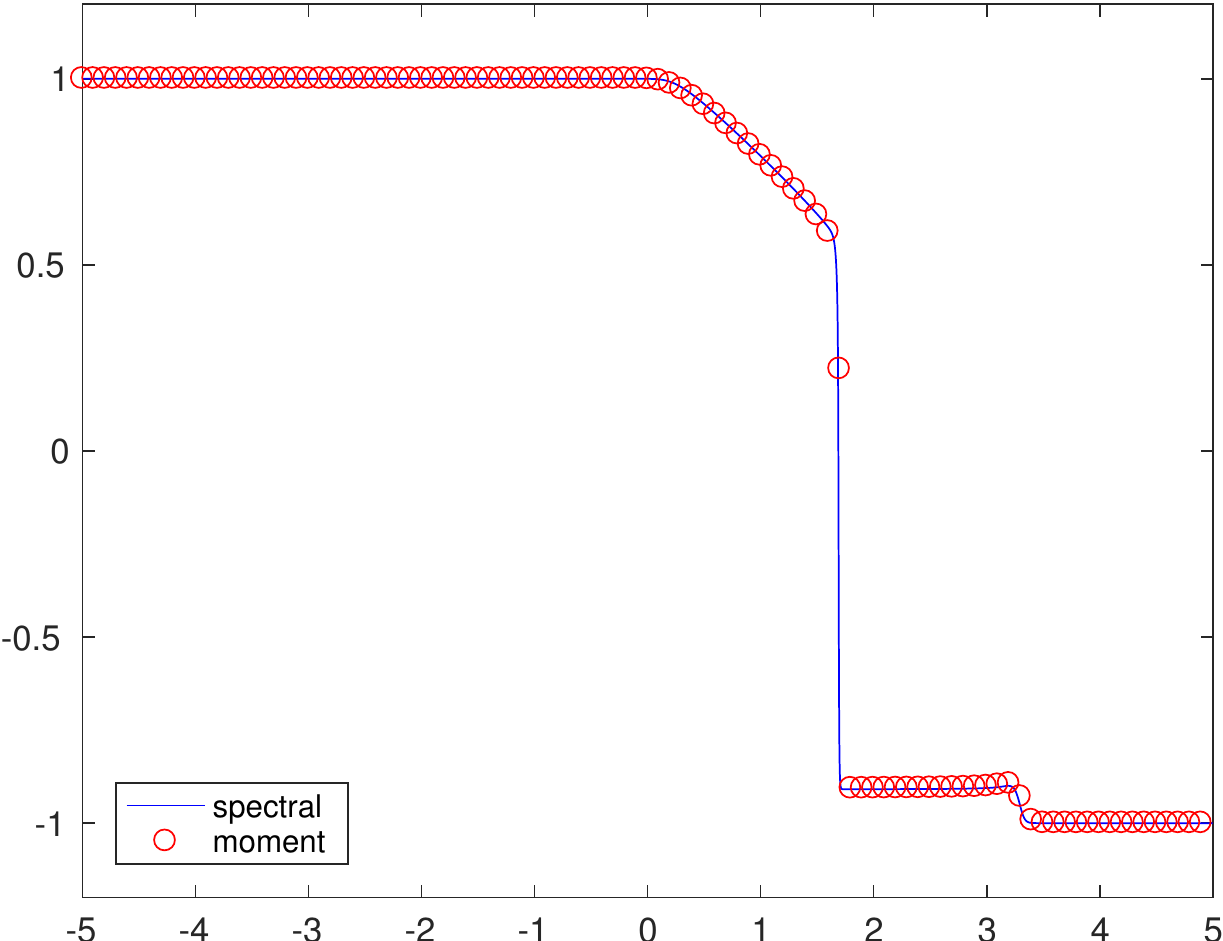}
}

  \centering
\subfigure[$N=4$]{
\includegraphics[width=0.3\textwidth]{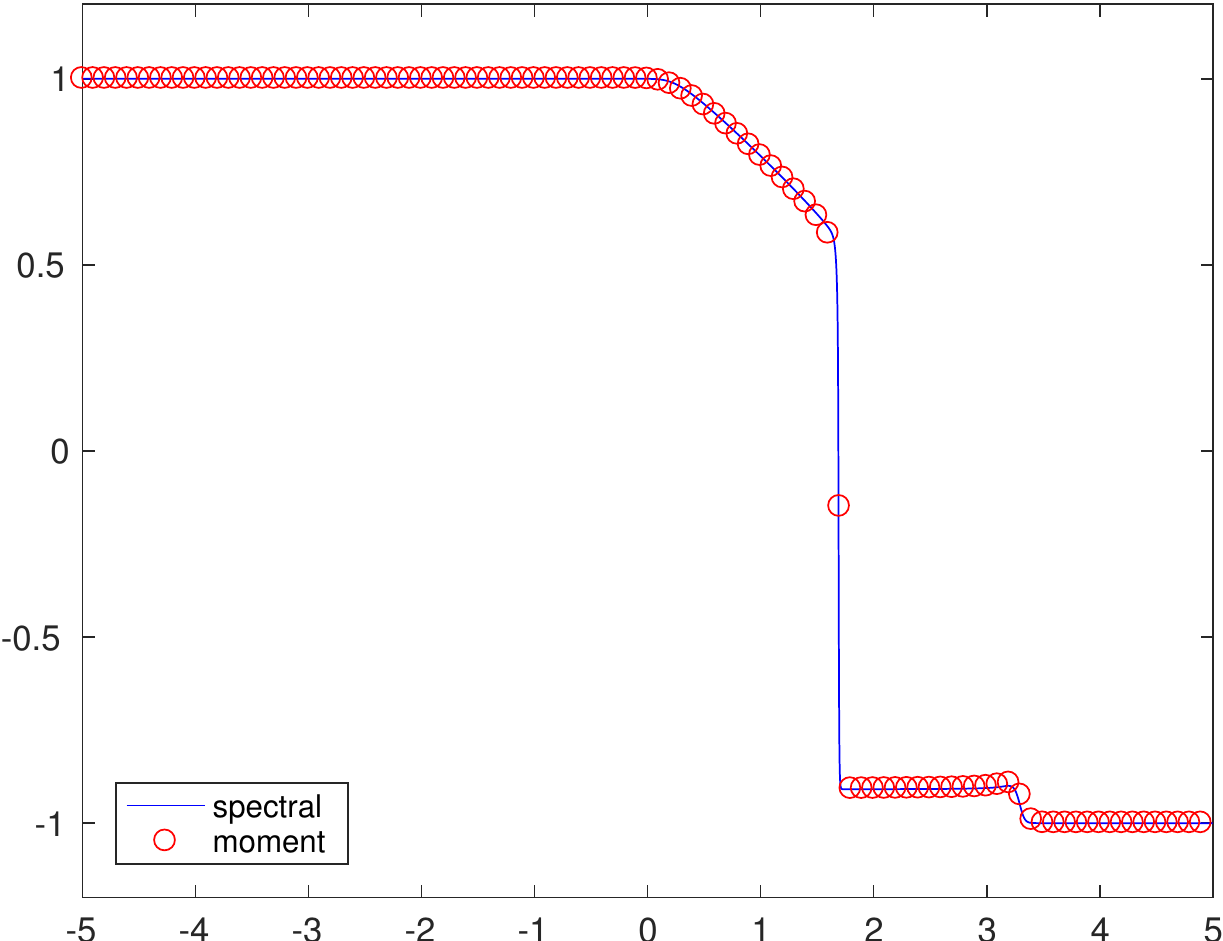}
}
\subfigure[$N=5$]{
\includegraphics[width=0.3\textwidth]{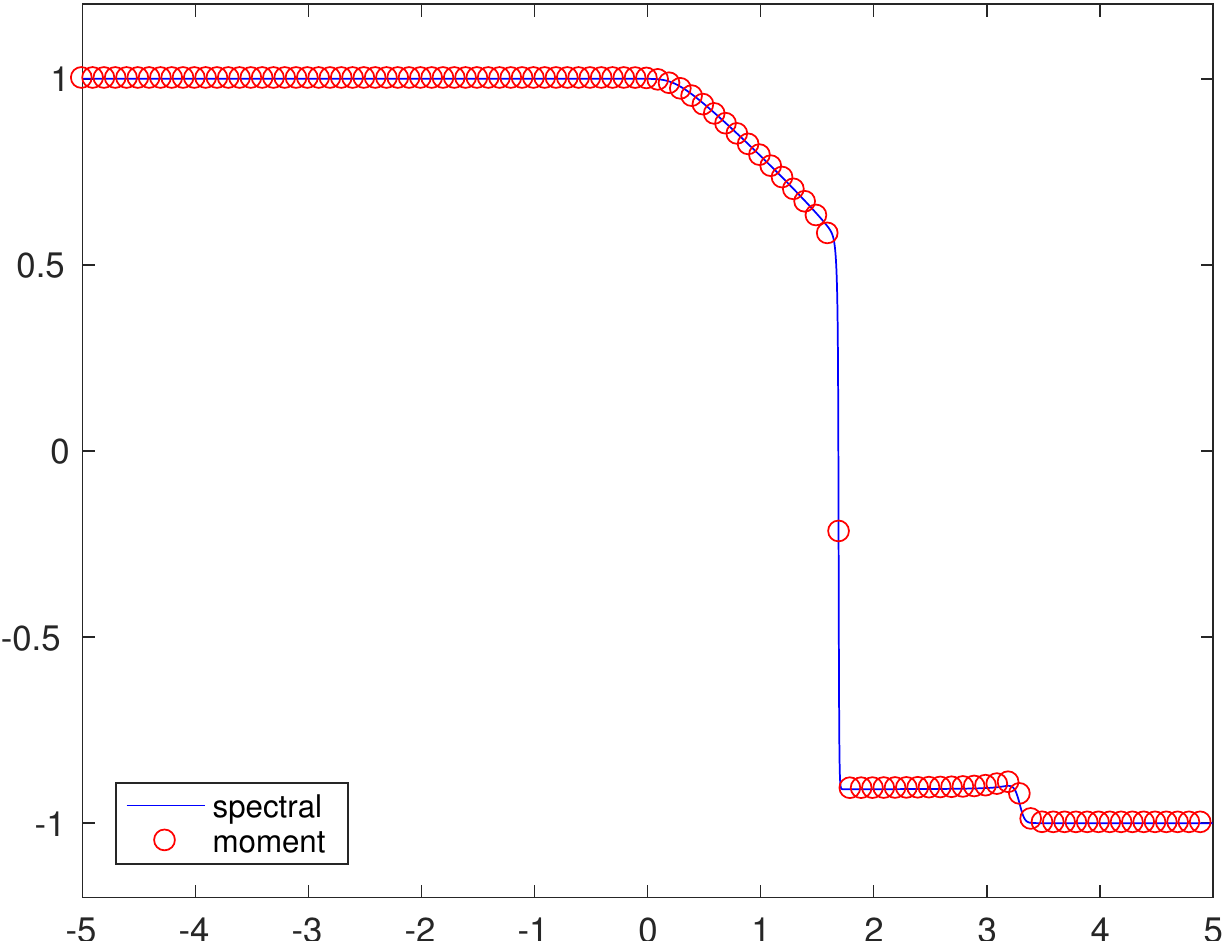}
}
\subfigure[$N=6$]{
\includegraphics[width=0.3\textwidth]{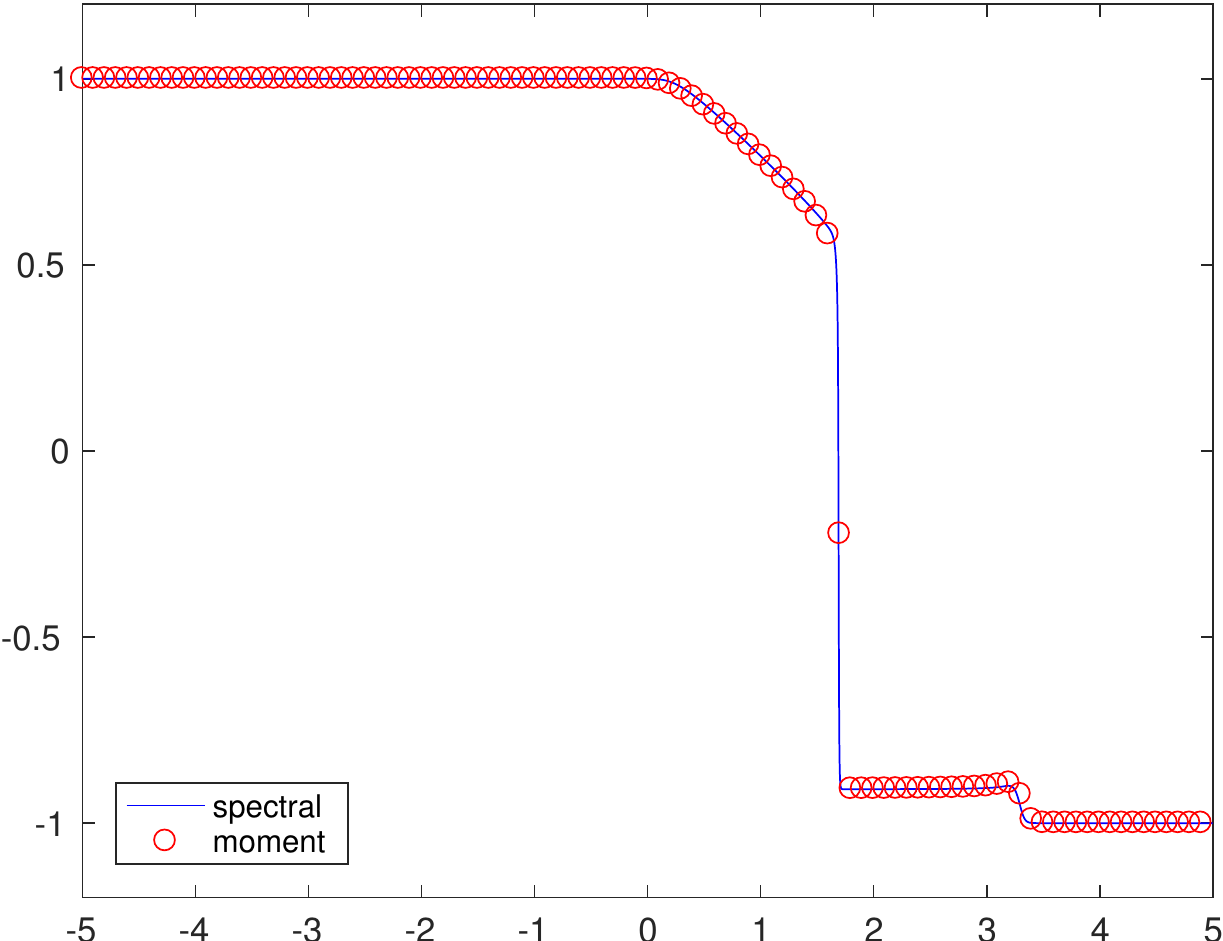}
}
%  \centering
%\subfigure[$N=7$]{
%\includegraphics[width=0.3\textwidth]{./images/riemann3/cmpN_e0010/r3_cmpN_u_n4000_e0010_N7_circle.pdf}
%}
%\subfigure[$N=8$]{
%\includegraphics[width=0.3\textwidth]{./images/riemann3/cmpN_e0010/r3_cmpN_u_n4000_e0010_N8_circle.pdf}
%}
%\subfigure[$N=9$]{
%\includegraphics[width=0.3\textwidth]{./images/riemann3/cmpN_e0010/r3_cmpN_u_n4000_e0010_N9_circle.pdf}
%}
\caption{Same as Fig. \ref{figure-example3-3} except for the macroscopic
velocity angles.}\label{figure-example3-4}
%使用8000个网格Spectral method的结果作为参考,
%和使用4000个网格的矩方法的结果进行比较}
\end{figure}

\begin{example}[Contact discontinuity]\label{example:5.3}
The initial data of the third Riemann problem are
\begin{equation*}
(\rho^\varepsilon,\bar\theta^\varepsilon)=\begin{cases}
	 (1, 1), & x<0,\\  (1, -1), & x>0. \end{cases}
 \end{equation*}
It is a contact discontinuity problem.

Figs. \ref{figure-example3-1} and \ref{figure-example3-2} show the densities $\rho^\varepsilon$  and
macroscopic velocity angles $\bt^\varepsilon$ at $t=4$ obtained by the moment
method with $N=1,2,\cdots,6$, 4000 cells, and $\varepsilon=1$,
where the solid line denotes the reference solution obtained by using the
spectral method  with 8000 cells.  Figs. \ref{figure-example3-3} and  \ref{figure-example3-4}
display  corresponding solutions for the case of
$\varepsilon=0.01$.
It is observed that
  the solutions of moment method $(\rho^\varepsilon, \bt^\varepsilon)$
do converge  the contact profile as $\varepsilon$ becomes small, and
the solutions  of the moment
system   agree with the reference when $N$ is larger than 3.
When the $\varepsilon$ is smaller, the convergence rate of the moment method is faster,
and for a fixed $\varepsilon$, the faster the $N$, the faster the convergence rate.
\end{example}

%%%%%%%%%%%%%%%%%%%%%%%% cmpN_e1000 %%%%%%%%%%%%%%%%%%%%%%%%
%%%%%%%%%%%%%%%%%%%%%%%% r3 u circle %%%%%%%%%%%%%%%%%%%%%%%%

%%%%%%%%%%%%%%%%%%%%%%%% cmpN_e0010 %%%%%%%%%%%%%%%%%%%%%%%%
%%%%%%%%%%%%%%%%%%%%%%%% r3 rho circle %%%%%%%%%%%%%%%%%%%%%%%%

%%%%%%%%%%%%%%%%%%%%%%%% cmpN_e0010 %%%%%%%%%%%%%%%%%%%%%%%%
%%%%%%%%%%%%%%%%%%%%%%%% r3 u circle %%%%%%%%%%%%%%%%%%%%%%%%

\begin{figure}
  \centering
\includegraphics[width=0.30\textwidth,  trim = 40 40 20 40,  clip]{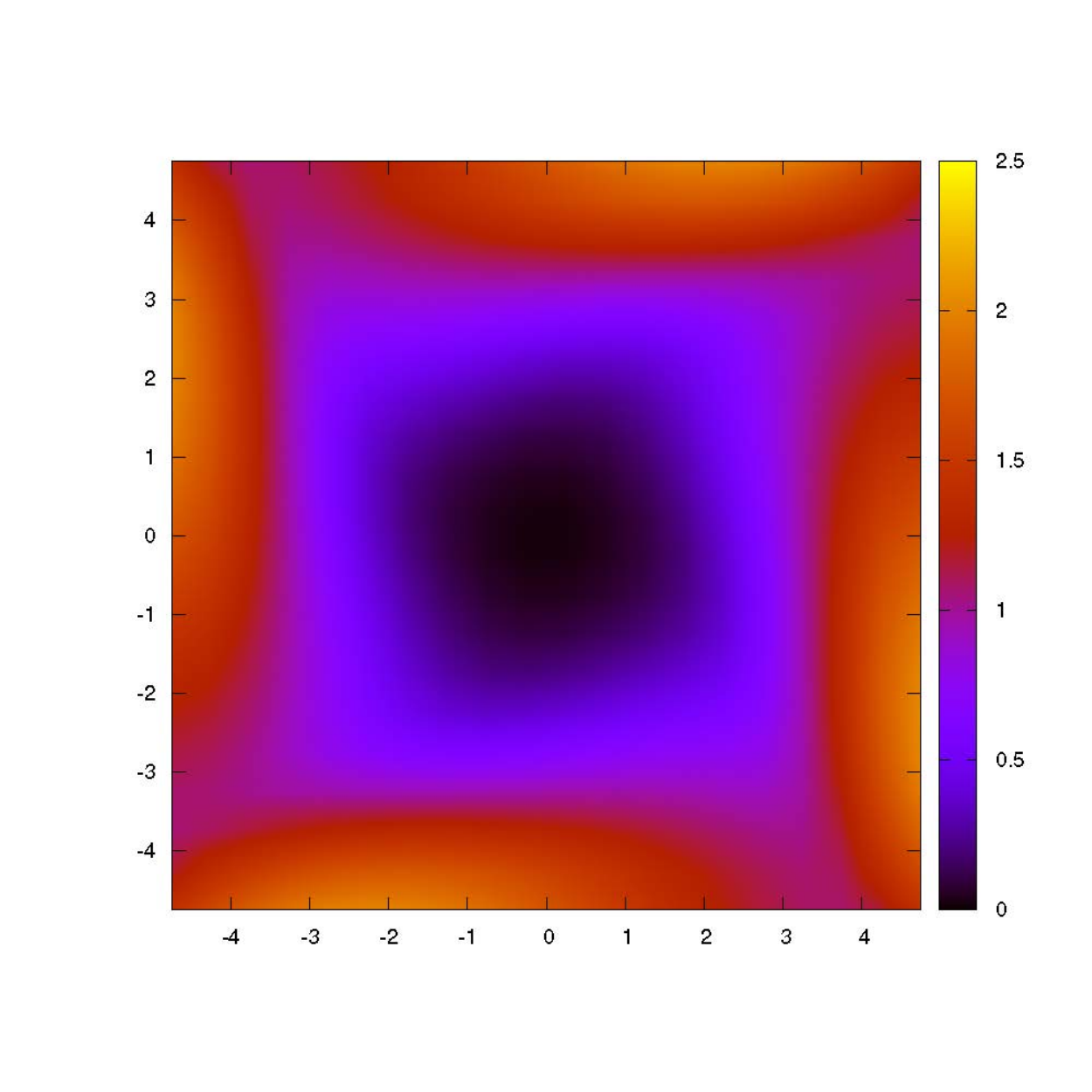}
\includegraphics[width=0.30\textwidth,  trim = 40 40 20 40,  clip]{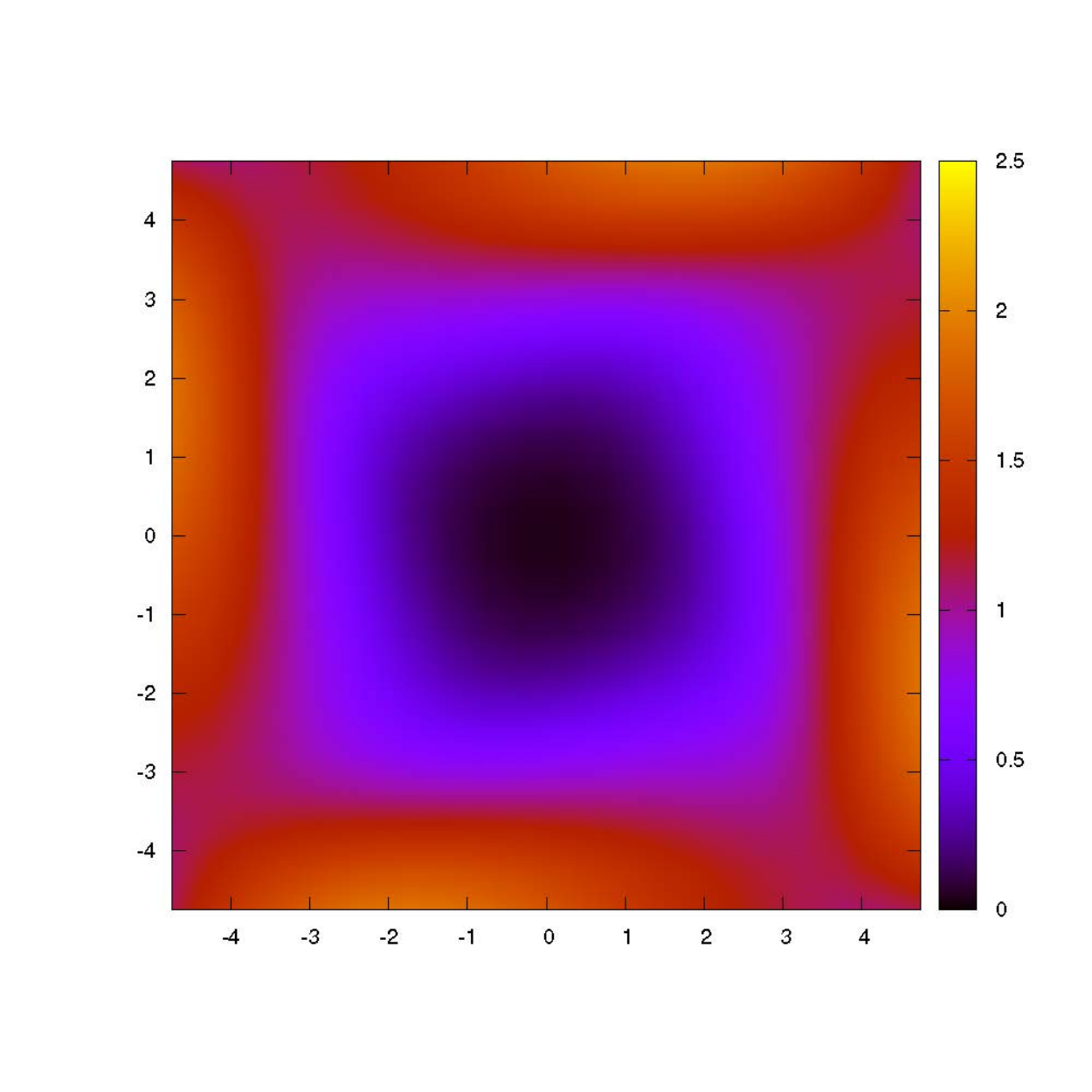}
\includegraphics[width=0.30\textwidth,  trim = 40 40 20 40,  clip]{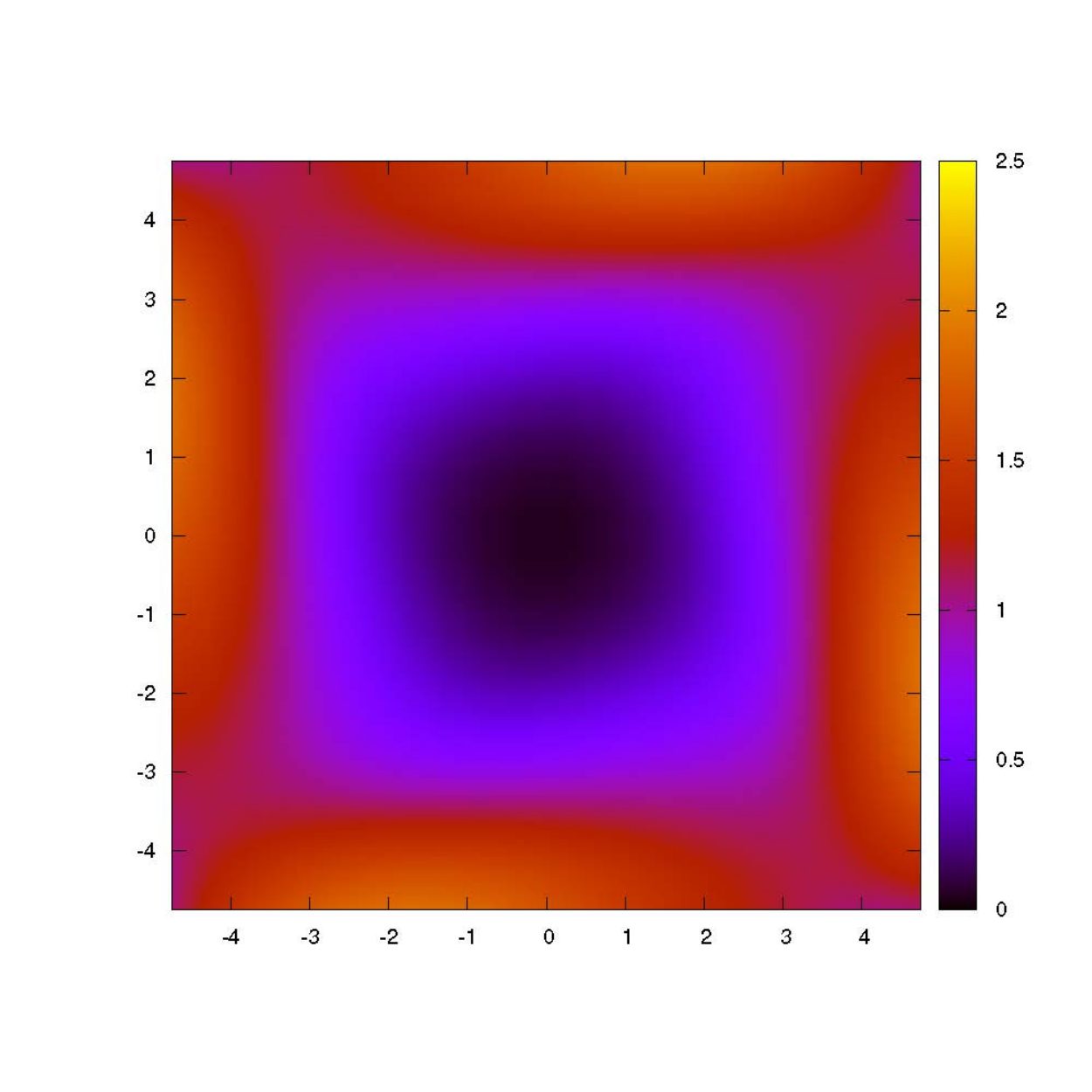}
\caption{Example \ref{2DVortex-formation}: The schlieren images of density obtained by using the moment method with $n=20$. From left to right: $N=2,3,4$.}
\label{fig:5.27}
\end{figure}

\begin{figure}
  \centering
\includegraphics[width=0.30\textwidth]{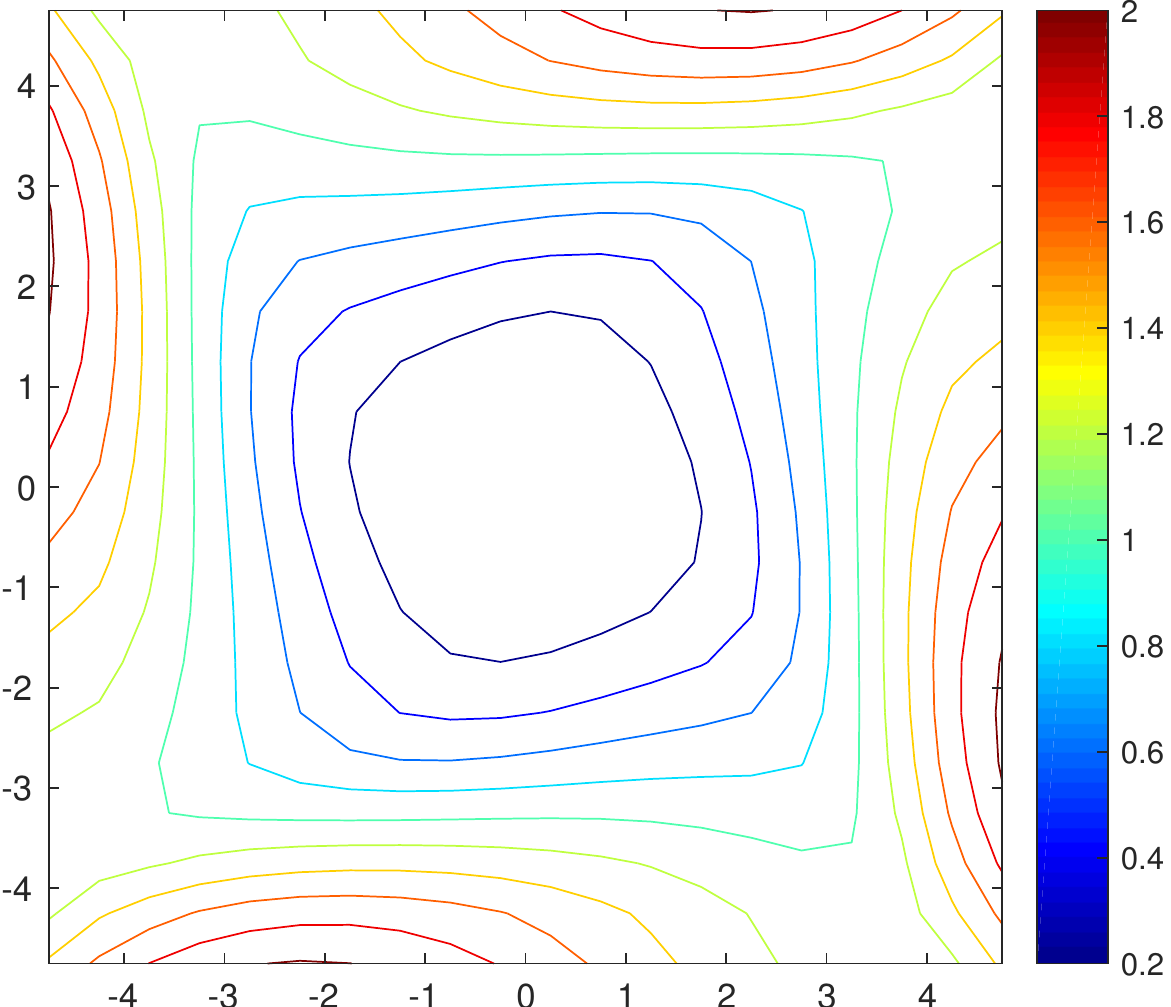}
  \includegraphics[width=0.30\textwidth]{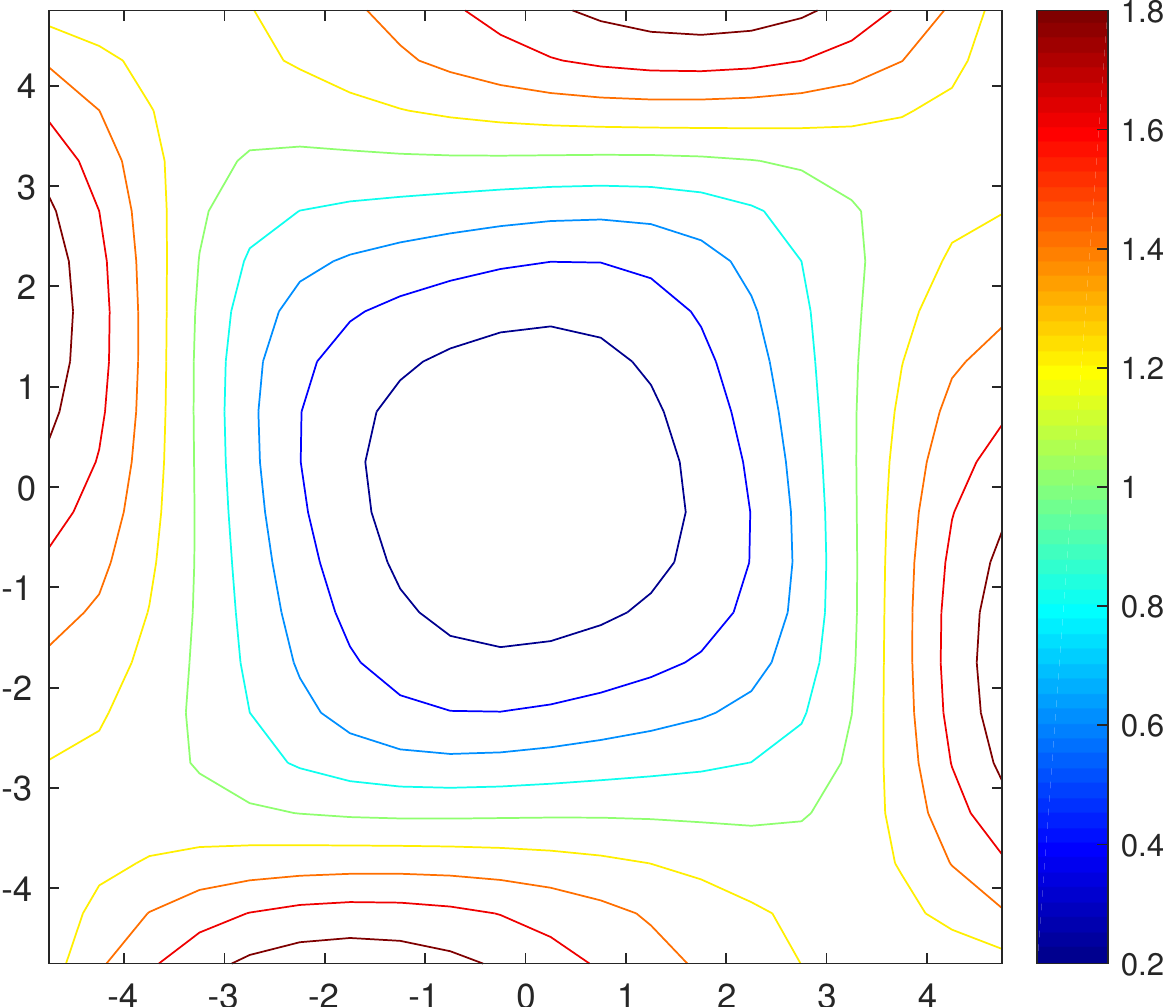}
  \includegraphics[width=0.30\textwidth]{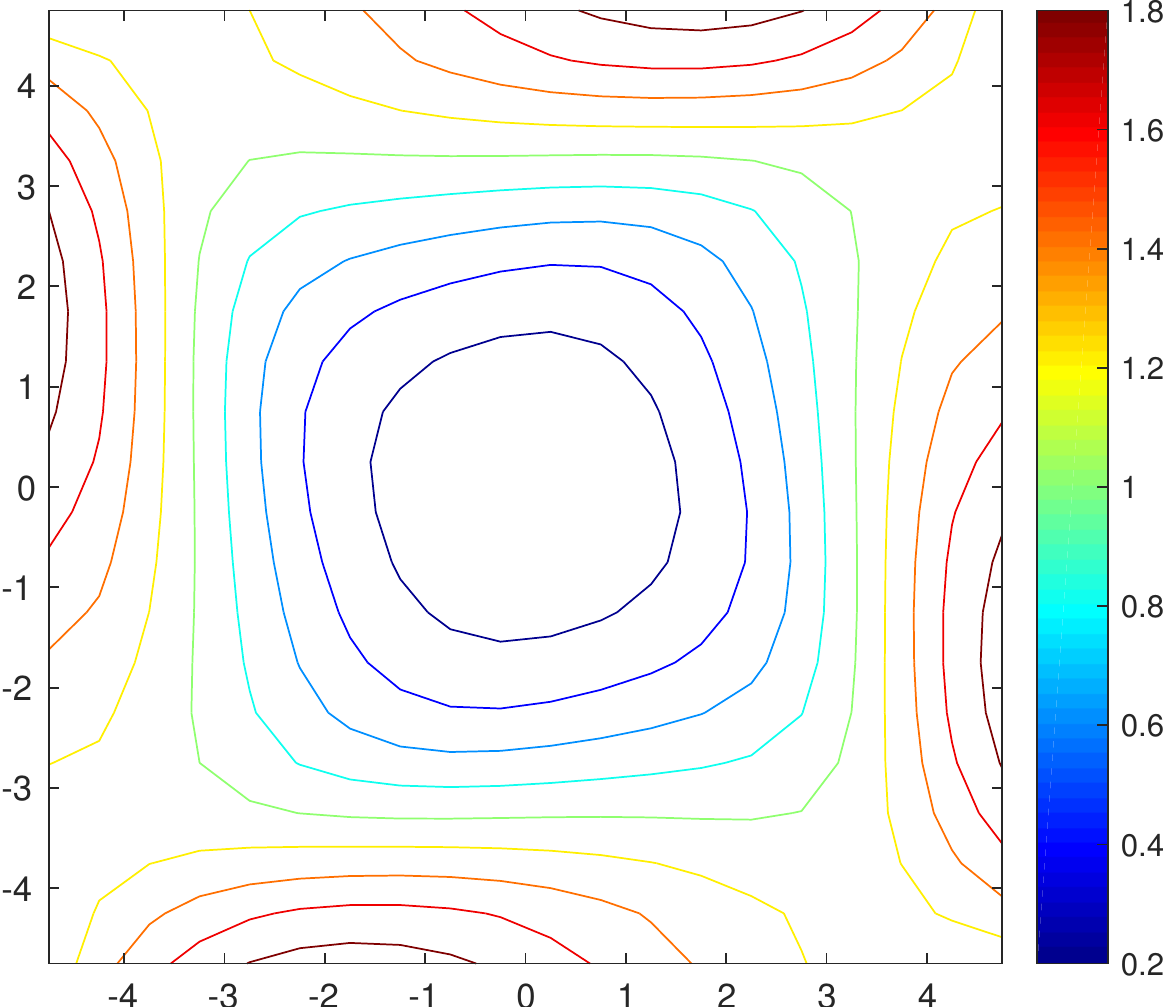}
\caption{Same as Fig.~\ref{fig:5.27} except for the density contours.
}\label{fig:5.28}
% obtained by using the moment method with $n=20$,
%and spectral method with $n=25$, respectively.}
\end{figure}

\begin{figure}
  \centering
\includegraphics[width=0.27\textwidth]{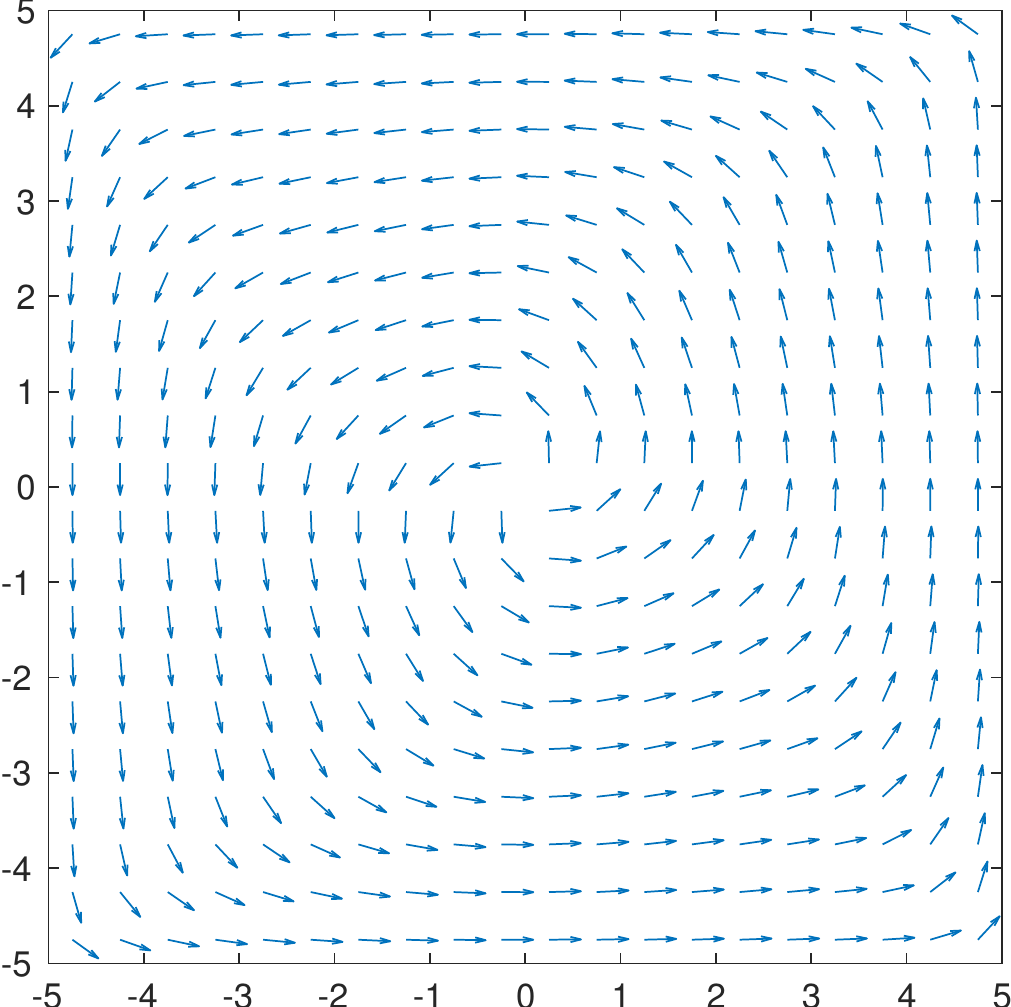}
\includegraphics[width=0.27\textwidth]{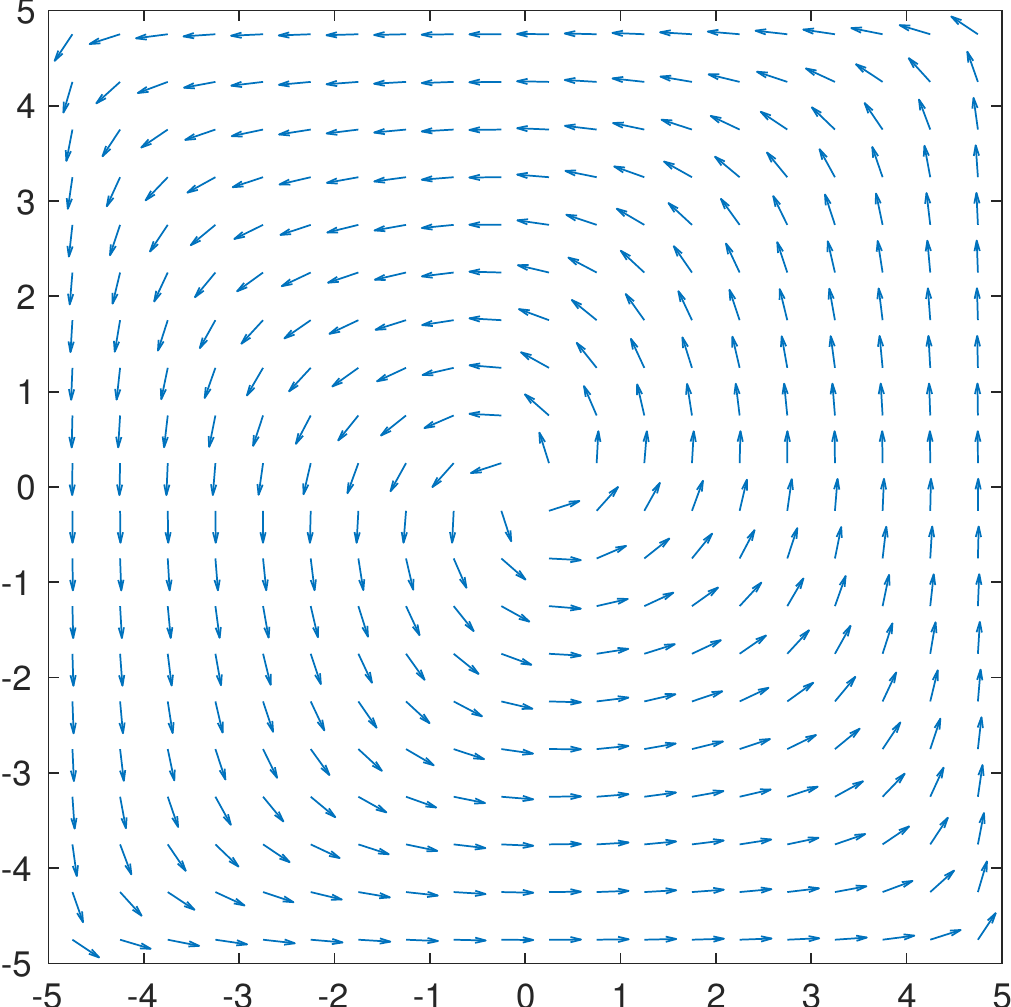}
\includegraphics[width=0.27\textwidth]{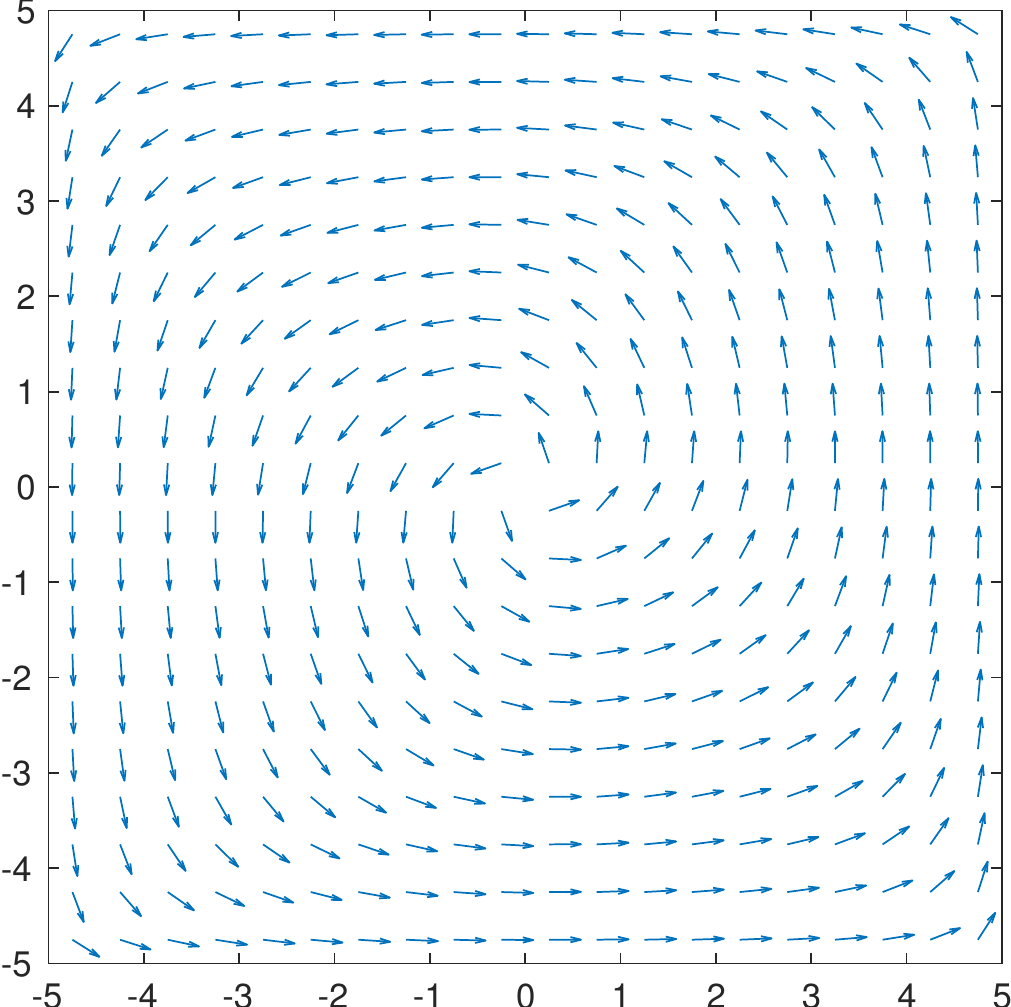}
\caption{Same as Fig.~\ref{fig:5.27} except for the arrow diagrams of velocity.}
\label{fig:5.29}
\end{figure}

\begin{figure}
	\centering
%\subfigure[Spectral method]{
%\includegraphics[width=0.30\textwidth,  trim = 40 40 20 40,  clip]{./images/vortex/gplvor_spc_n25_rho.pdf}
%}
%\subfigure[Moment method, $N=2$]{
\includegraphics[width=0.30\textwidth,  trim = 40 40 20 40,  clip]{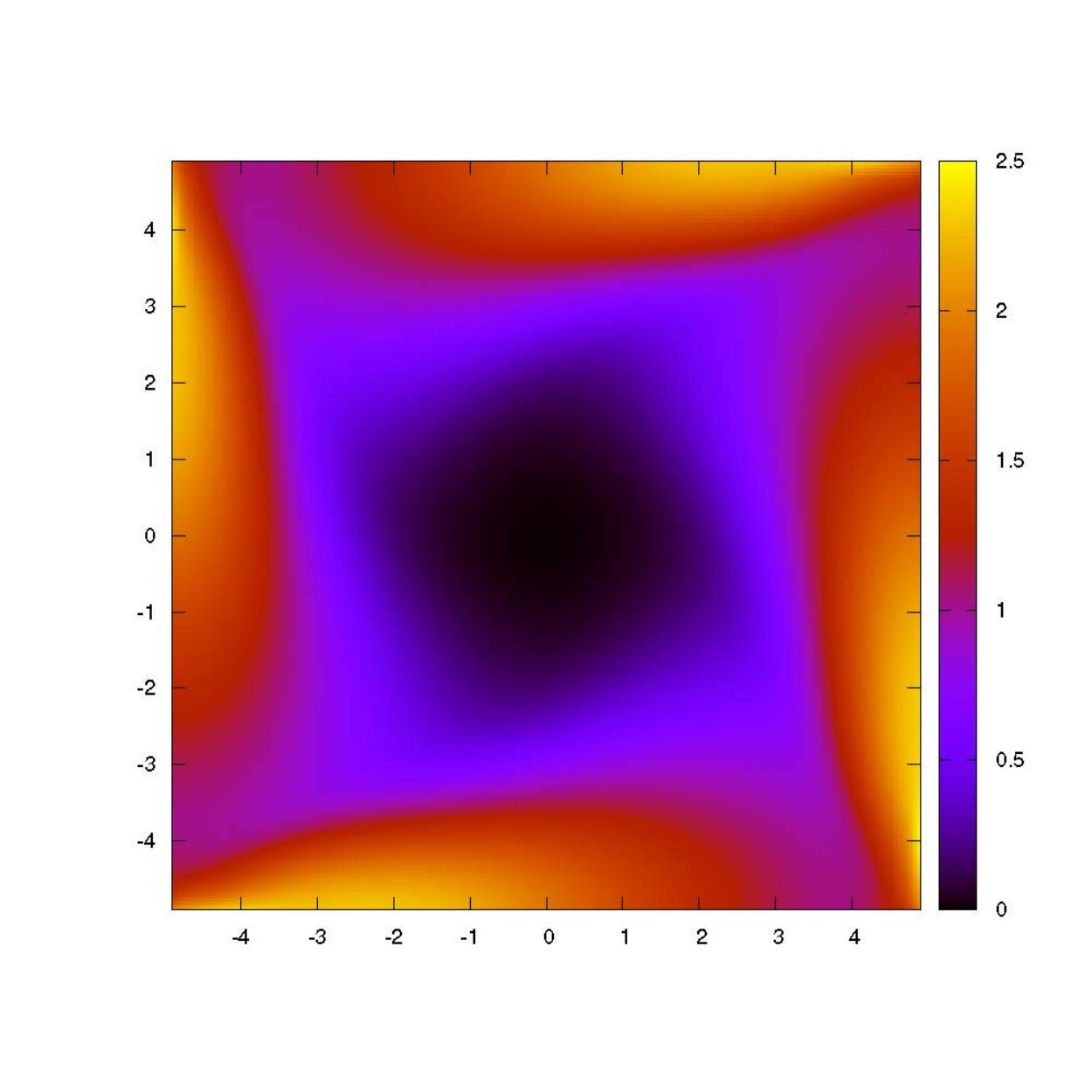}
%}
%\centering
%\subfigure[Moment method, $N=3$]{
	\includegraphics[width=0.30\textwidth,  trim = 40 40 20 40,  clip]{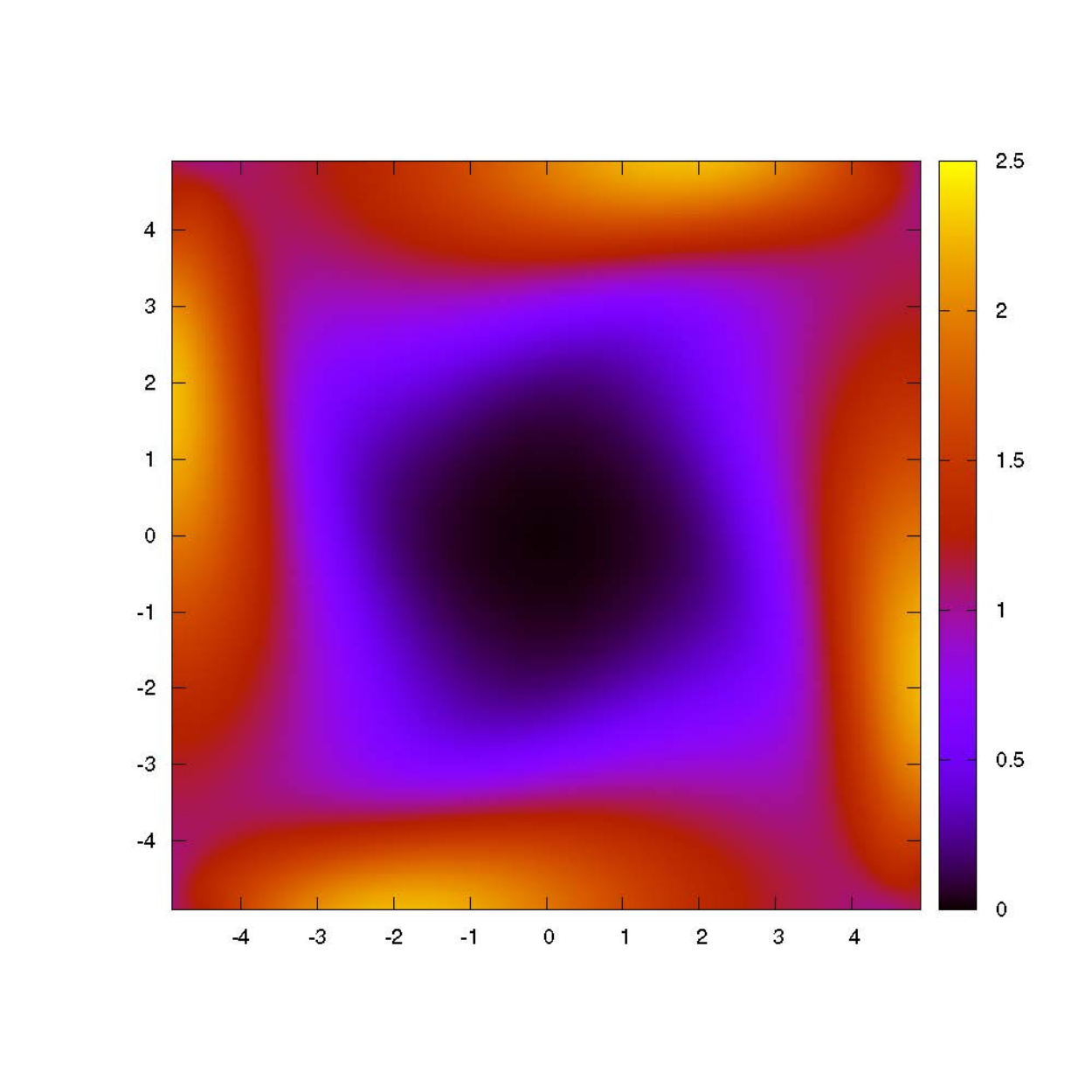}
%}
%\subfigure[Moment method, $N=4$]{
	\includegraphics[width=0.30\textwidth,  trim = 40 40 20 40,  clip]{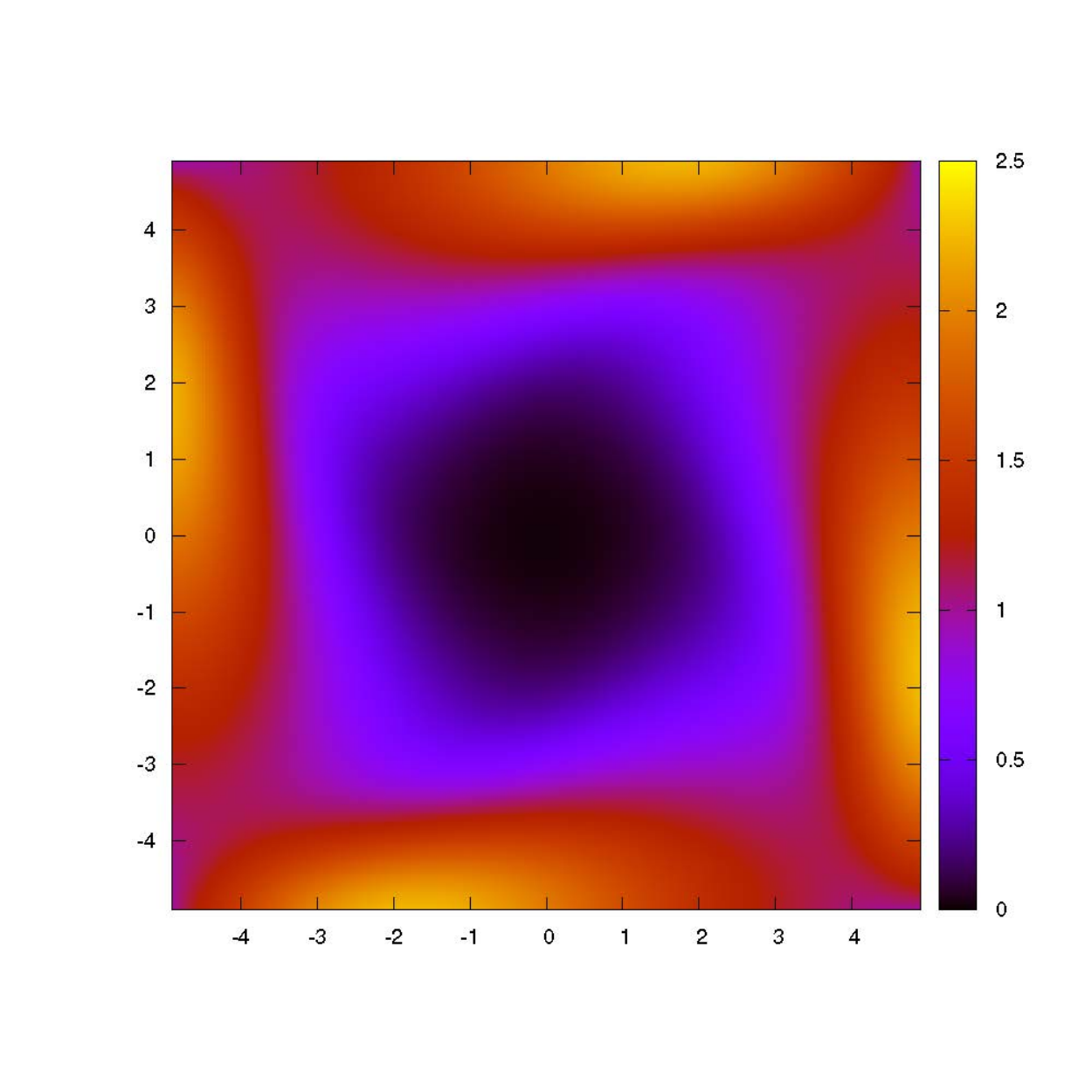}
%}
\caption{Same as Fig.~\ref{fig:5.27} except for
  $n=50$.}
  \label{fig:5.30}
\end{figure}

\begin{figure}
  \centering
%\subfigure[Spectral method]{
%\includegraphics[width=0.30\textwidth]{./images/vortex/spc_n25contour.pdf}
%}
%\subfigure[Moment mentod, $N=2$]{
\includegraphics[width=0.30\textwidth]{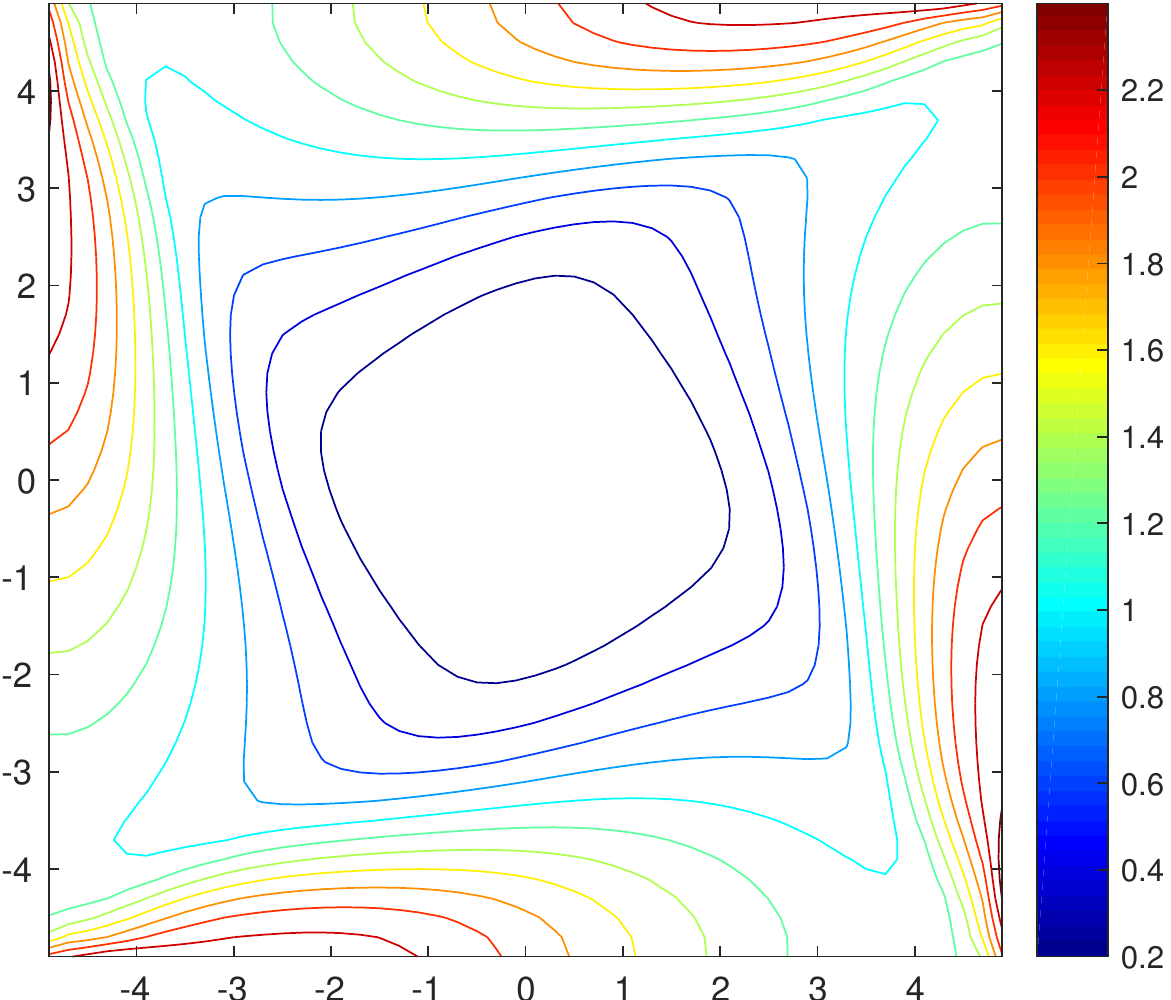}
%}
%\centering
%\subfigure[Moment method, $N=3$]{
  \includegraphics[width=0.30\textwidth]{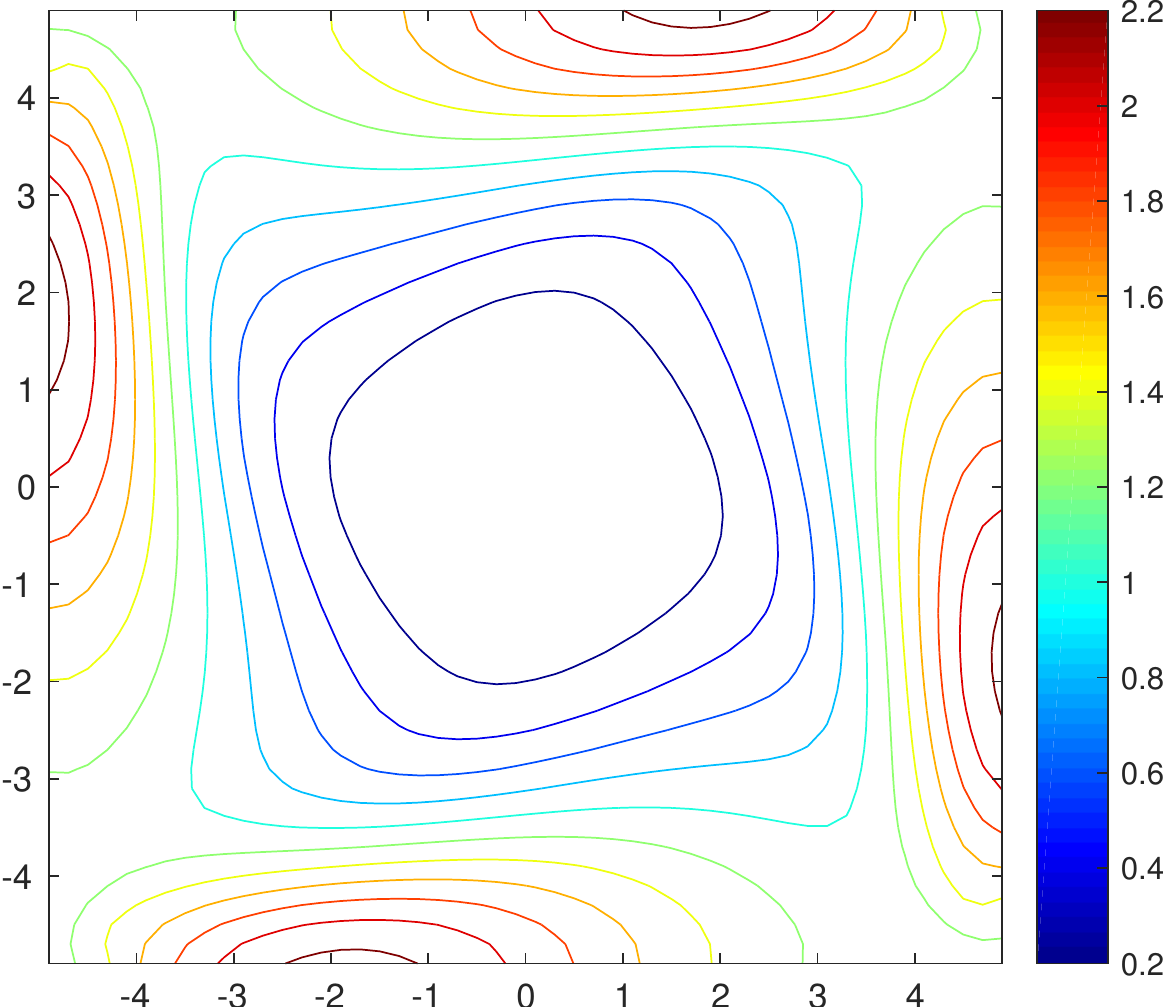}
%}
%\subfigure[Moment method, $N=4$]{
  \includegraphics[width=0.30\textwidth]{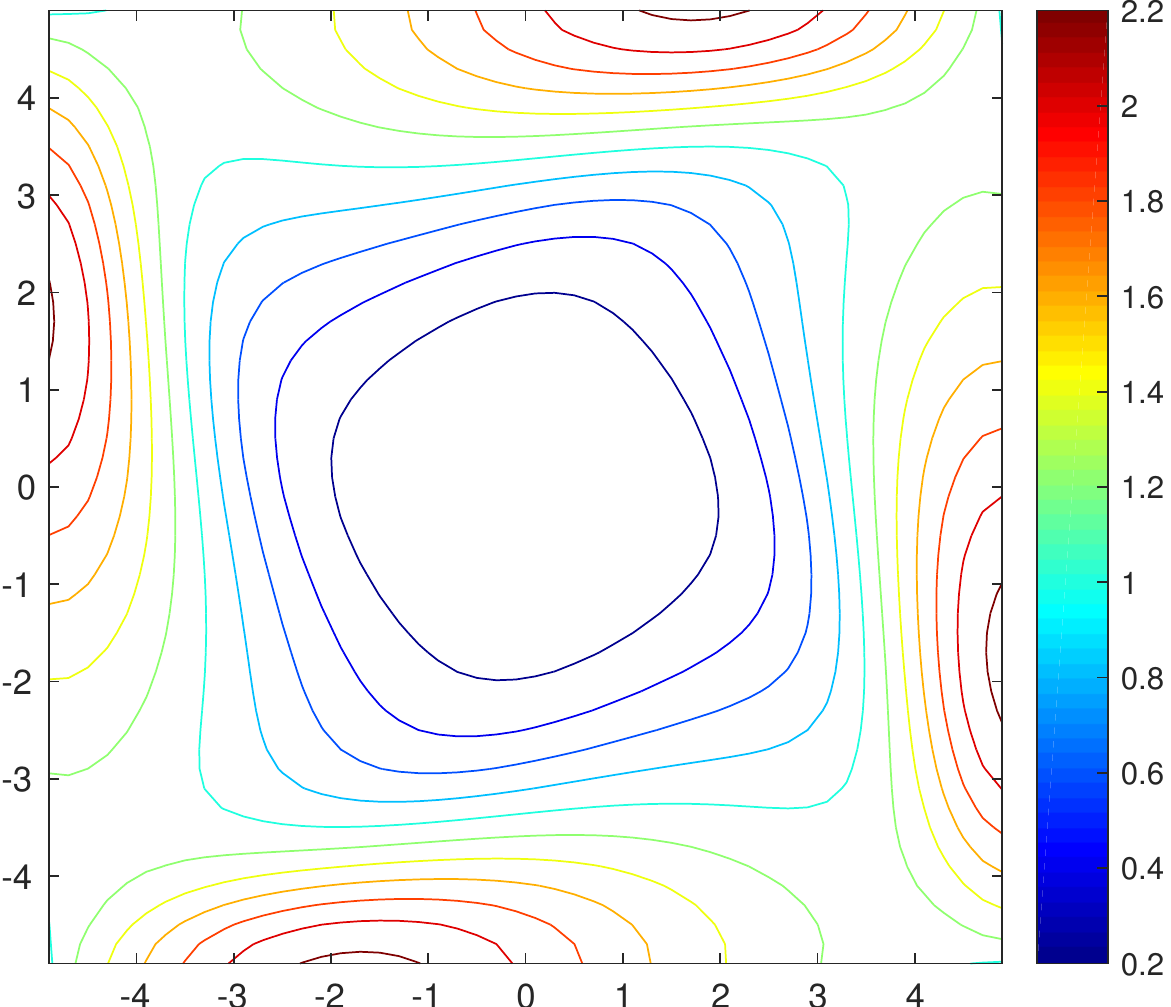}
%}
\caption{Same as Fig.~\ref{fig:5.30} except for  the density contours.}\label{fig:5.32}
\end{figure}

\begin{figure}
  \centering
%\subfigure[Spectral method]{
%\includegraphics[width=0.30\textwidth]{./images/vortex/spc_n25u.pdf}
%}
%\subfigure[Moment method, $N=2$]{
\includegraphics[width=0.27\textwidth]{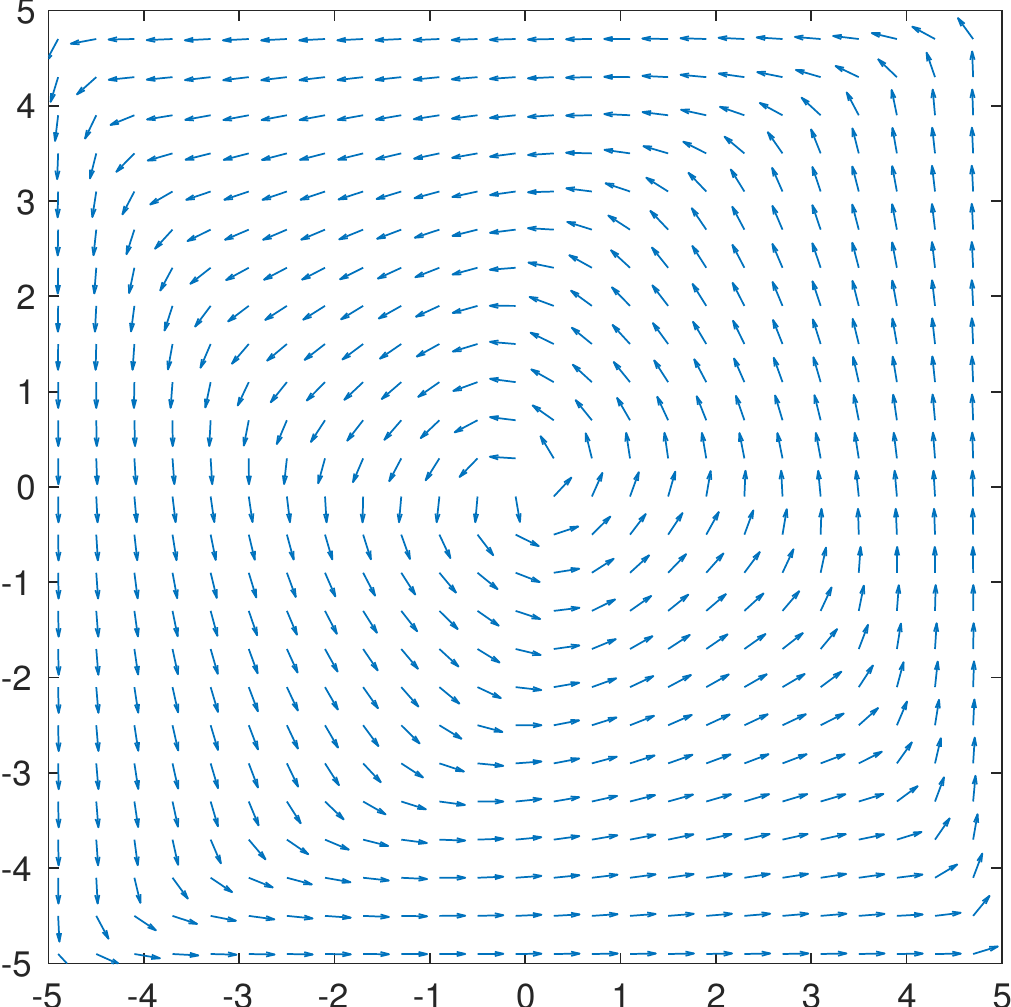}
%}
%\centering
%\subfigure[Moment method, $N=3$]{
  \includegraphics[width=0.27\textwidth]{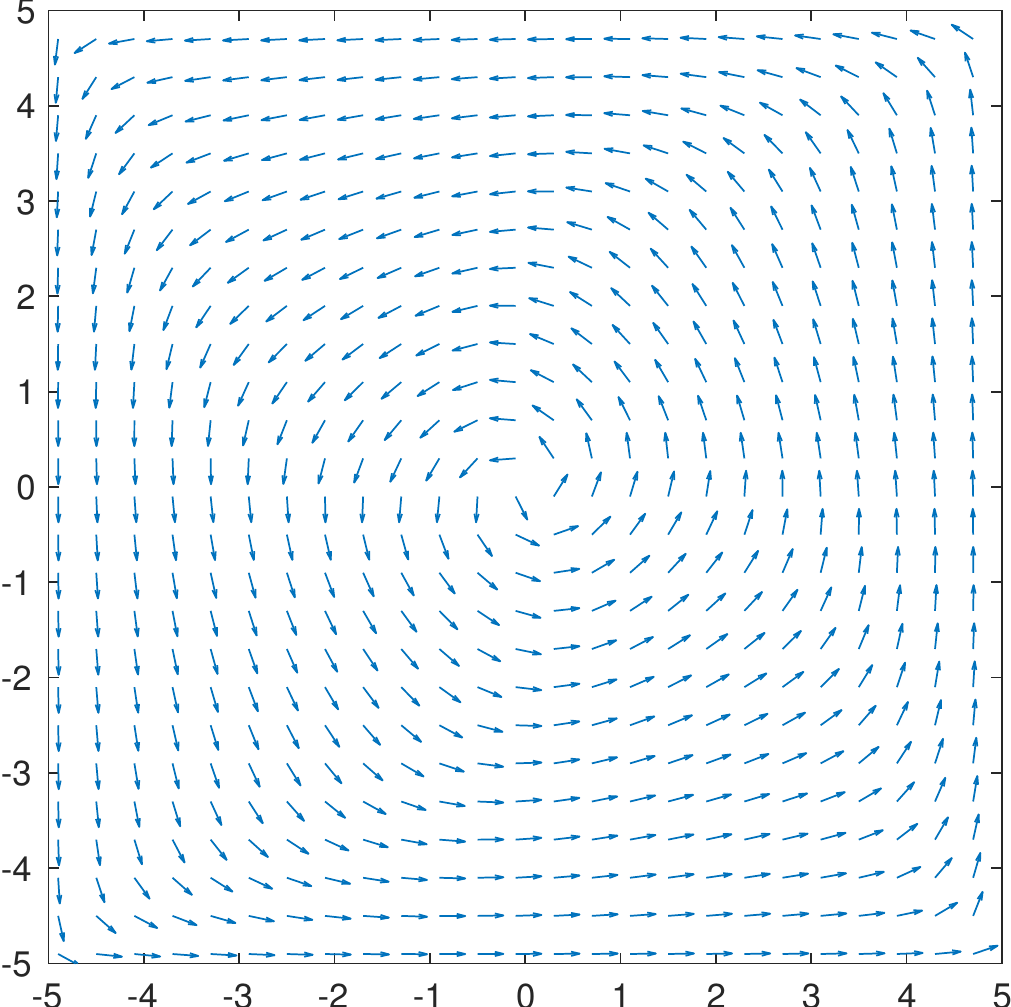}
%}
%\subfigure[Moment method, $N=4$]{
  \includegraphics[width=0.27\textwidth]{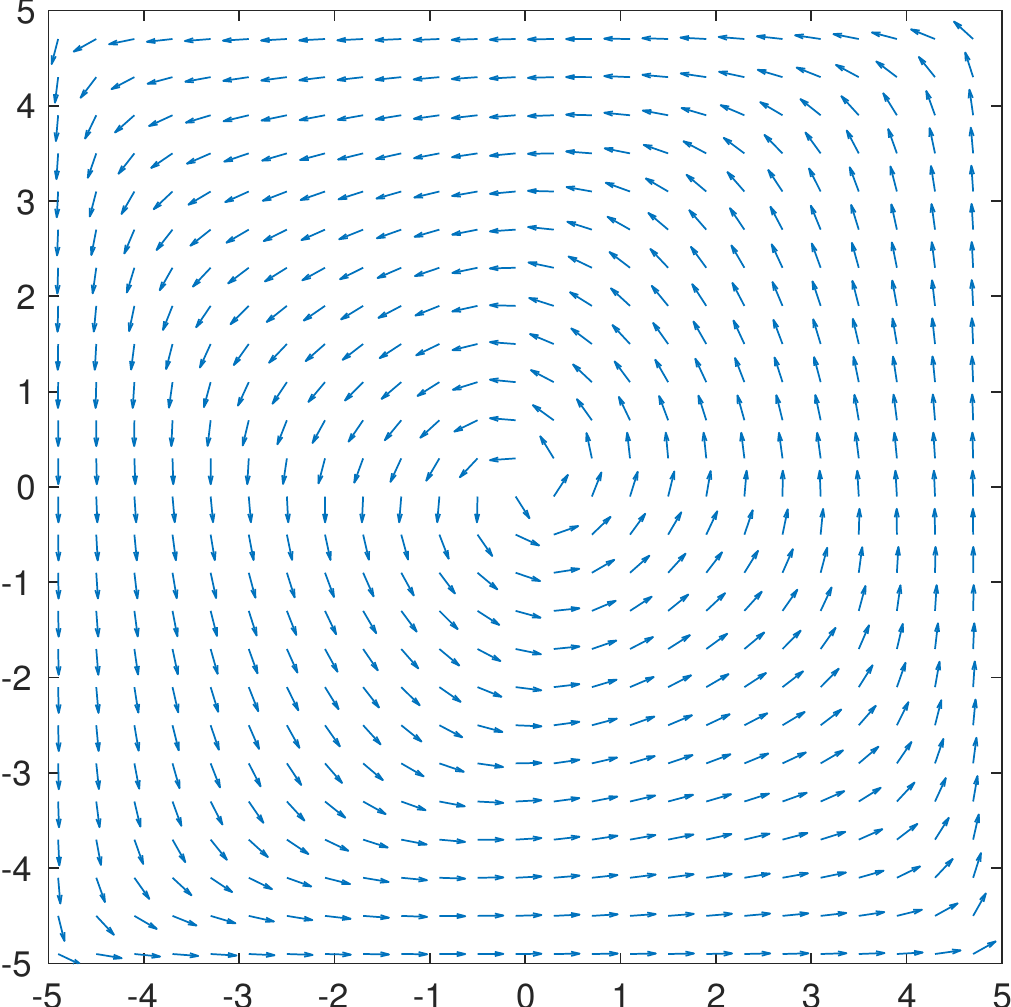}
%}
\caption{Same as Fig.~\ref{fig:5.30} except for the arrow diagram of velocity.}\label{fig:example5-3}
\end{figure}

\begin{figure}
  \centering
\subfigure[Density]{
\includegraphics[width=0.30\textwidth,  trim = 40 40 20 40,  clip]{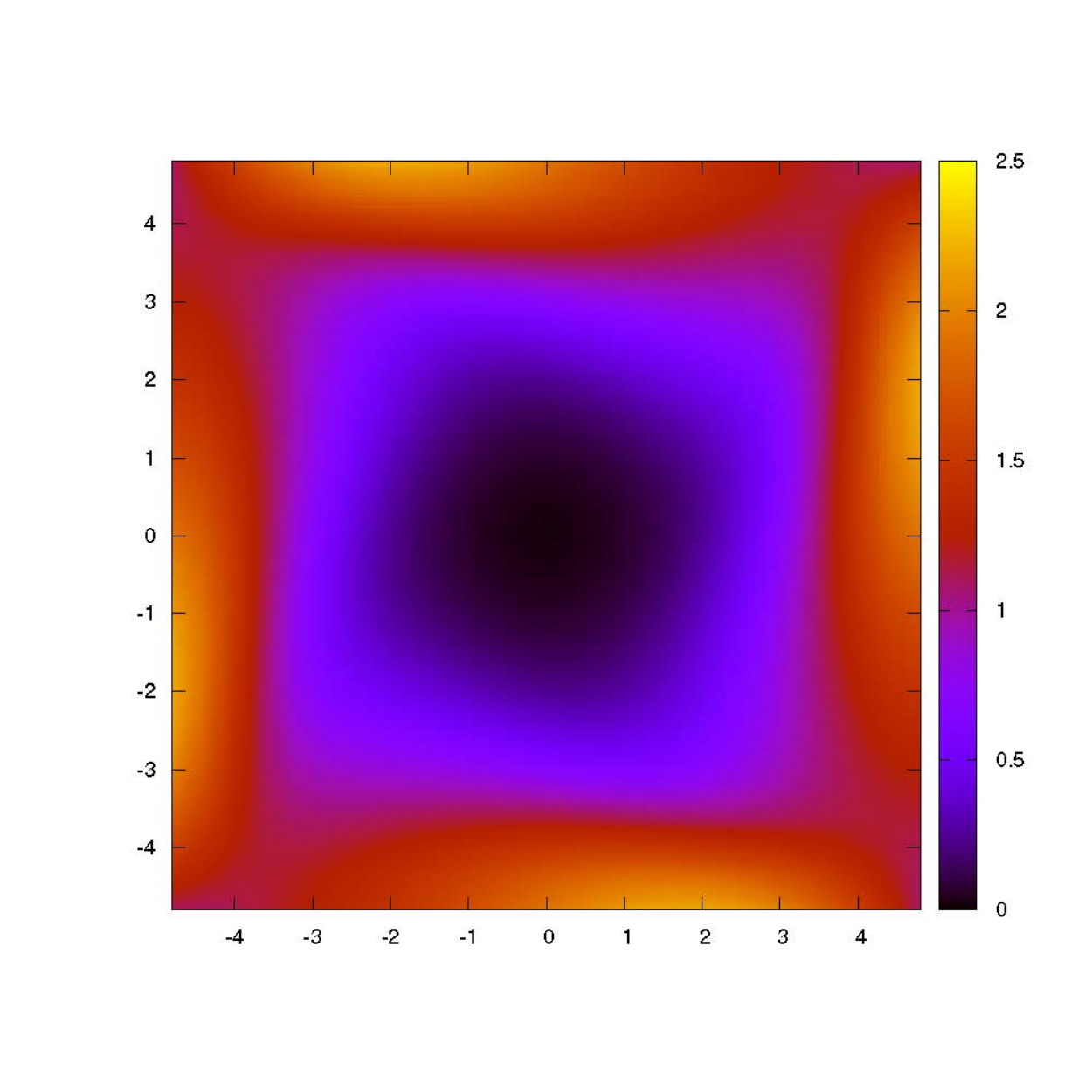}
}
\subfigure[Density]{
\includegraphics[width=0.31\textwidth]{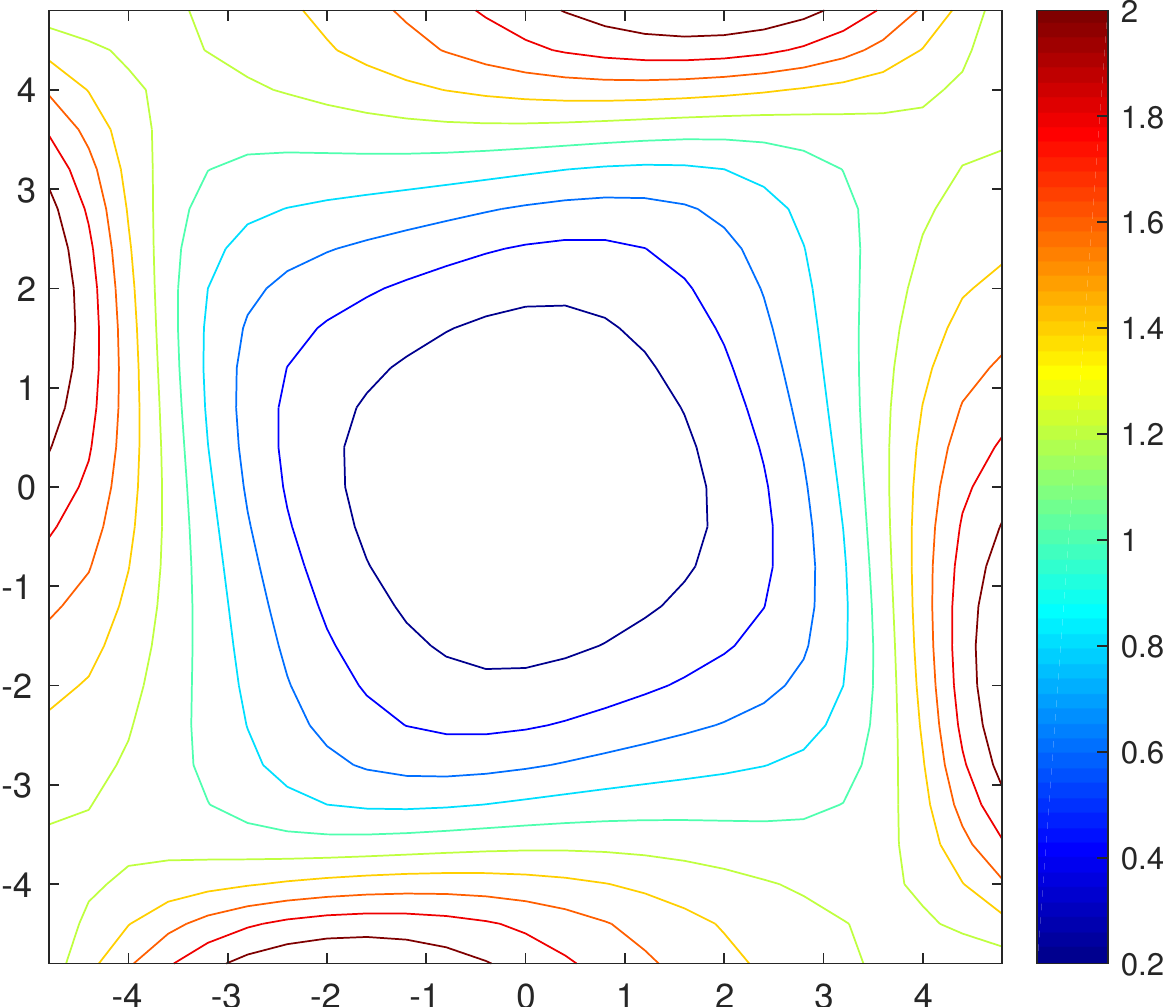}
}
\subfigure[Velocity]{
\includegraphics[width=0.27\textwidth]{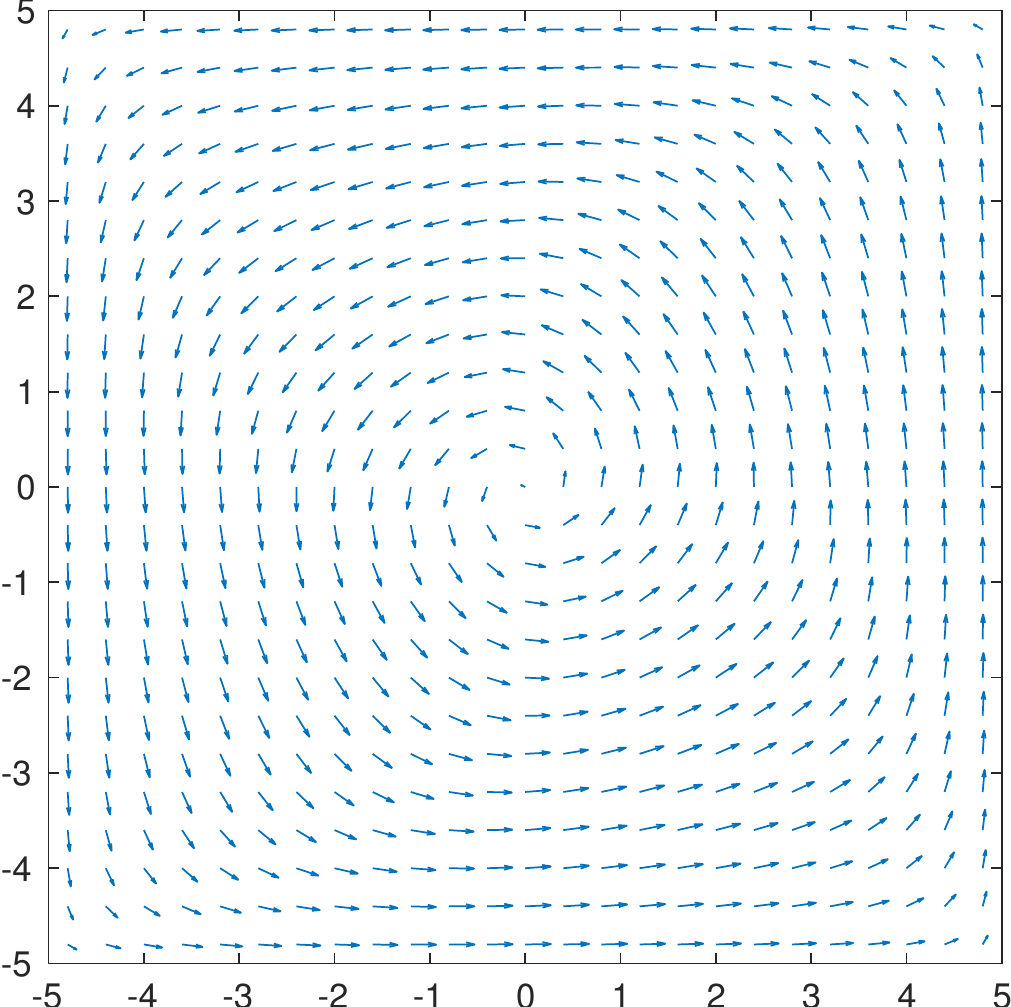}
}
\caption{Same as Figs.~\ref{fig:5.27}-\ref{fig:5.29} except for  the spectral method with $n=25$.}  \label{figure:example4-spectral}
% obtained by using the moment method with $n=20$,
%and spectral method with $n=25$, respectively.}
\end{figure}

\begin{figure}
  \centering
  \subfigure[Total mass $M(l)$]{
  \includegraphics[width=0.45\textwidth]{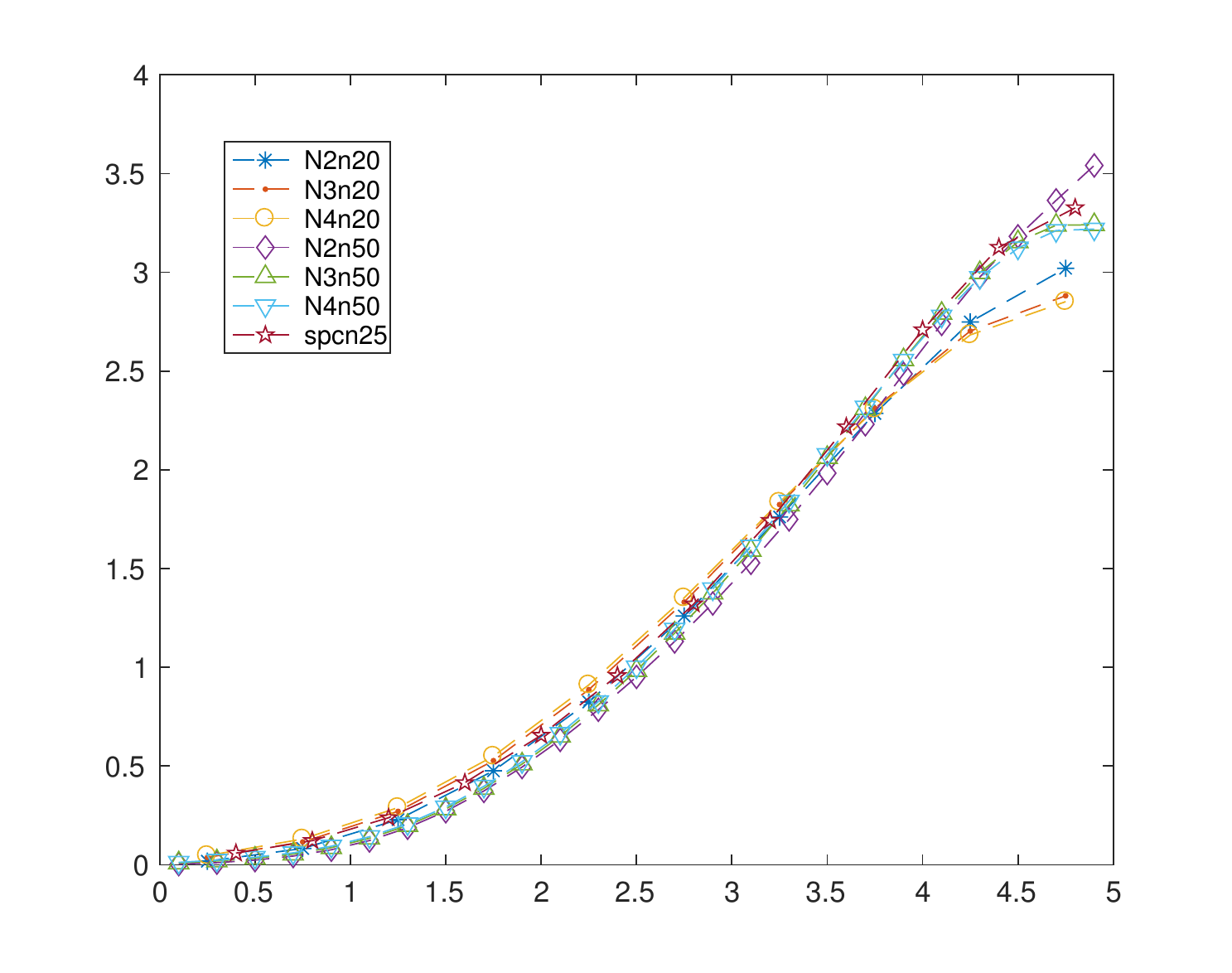}
  }
    \subfigure[Relative $\ell^2$ errors]{
    \includegraphics[width=0.45\textwidth]{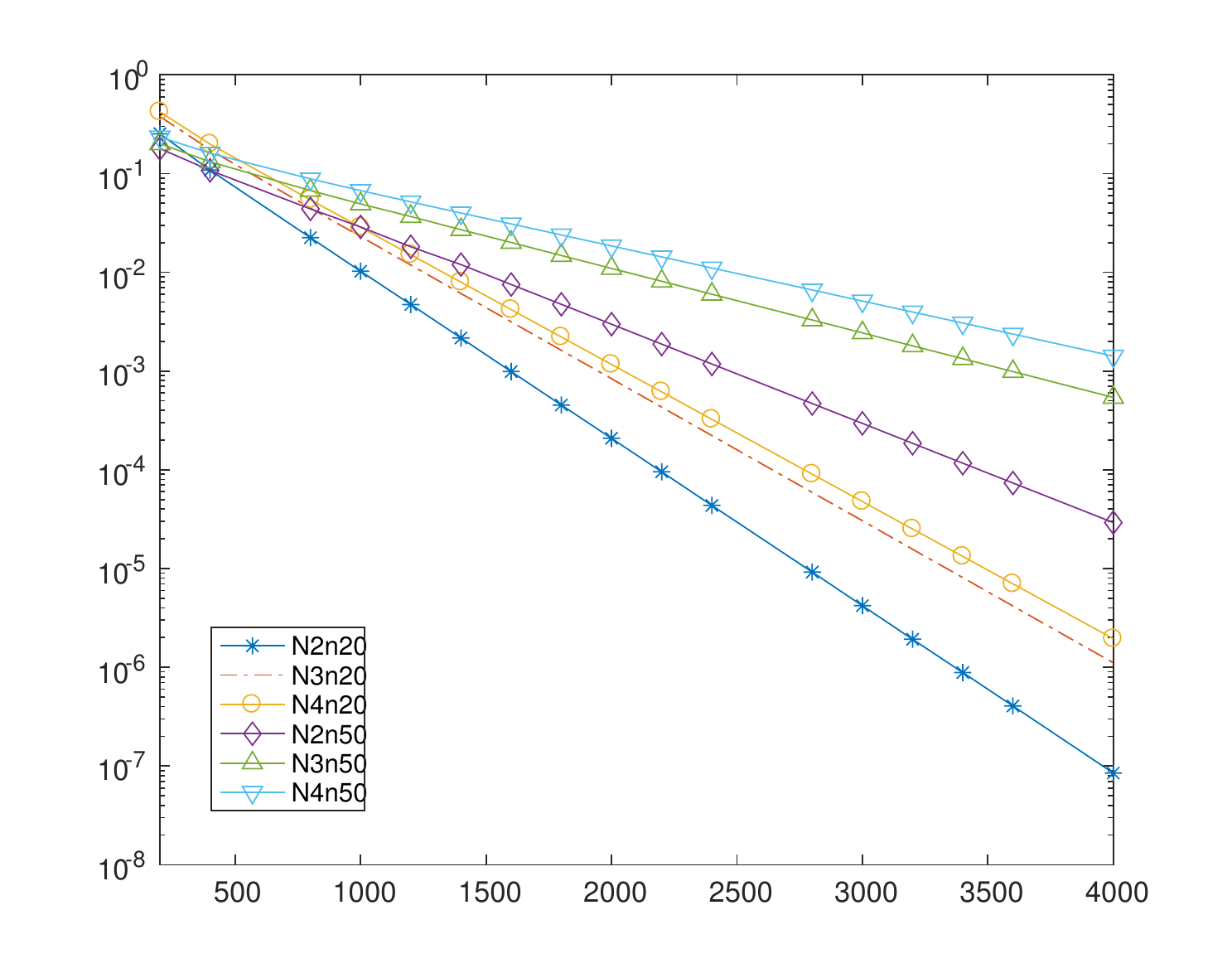}
    }
  \caption{Example \ref{2DVortex-formation}: The total mass $M(l)$ on the square $\Omega (l)$
  and the relative $\ell^2$ errors of density obtained by using moment methods.}
  \label{figure:total-mass2D}
\end{figure}

\begin{example}[Vortex formation]\label{2DVortex-formation}
%下面利用反射边界条件来模拟涡流的形成.
The computational domain is chosen as  the square area $[-5, 5]\times[-5, 5]$
with reflection boundary conditions,
and is divided into a uniform square mesh $\{(x_i,y_j)| x_i=-5+ih,  y_j=-5+jh,
i, j=0, 1, \cdots, n-1\}$.
The initial data are taken as follows
\begin{align*}
  f(0, x_i, y_j,\theta)=\dfrac{1}{2\pi},\ \ \bt(0, x_i, y_j)=\begin{cases}
    0, & x < 4.5,  \\
		\dfrac{\pi}2, & \text{otherwise}.    \end{cases}
\end{align*}
After a transient period, the solution will converge
to a steady state consisting of a vortex-type formation.

In numerical simulation,
a perturbation is added to the initial   velocity direction on the right boundary
in order to ensure that the final steady state is counterclockwise rotation,
and the solutions are output  when the relative $\ell^2$   error  of the  density  between two adjacent iterations
is less than $1.5\times10^{-3}$. %Moreover, epsilon=1,CFL=0.5

Figs. \ref{fig:5.27}-\ref{fig:5.29} show the densities and the velocities
obtained by using the moment methods with
$N=2,3,4$, $\varepsilon=1$,  and the mesh of $n=20$, where 13 equally spaced contour lines
are chosen from  0 to 2.4 with stepsize 0.2.
Figs. \ref{fig:5.30}-\ref{fig:example5-3} shows corresponding results for $n=50$.
For the sake of comparison, Fig. \ref{figure:example4-spectral} also
gives the results obtained by using the spectral method
with $n=25$.
%上面的数值结果展示了计算到平衡态时, 正方形区域上形成了一个漩涡,
%由于我们在初值的边界处做了一个“扰动”, 旋转的方向均为逆时针, 这与我们的预期相符.
%观察平均密度图, 矩方法在右端点处的平均密度是随着$N$增大收敛的,
%$N=3, 4$的平均密度较为接近且与Spectral method的结果较为一致.
%
%最后算例epsilon=1,CFL=0.5，而且epsilon=0.01的时候，这个稳态是算不出来的;
Fig. \ref{figure:total-mass2D} shows the total mass $M(l)$ on the square $\Omega (l)$
and the relative $\ell^2$ errors of density,
where
\begin{equation*}
	M(l)=\int_{\Omega (l)} \rho\dd s,
\end{equation*}
and
\begin{equation*}
	\Omega (l)=\{ (x, y)| \max\{\abs{x}, \abs{y}\}=l \},
\end{equation*}
and ``N2n20'' etc. in the legend represent ``$N=2, n=20$'' etc., while ``specn20'' denotes
the spectral method with $n=20$.
The distributions of $M(l)$  agree well with each other with some discrepancy   near the boundary of the domain ($l\approx 5$).
The discrepancies reduce as the number of moment and mesh cell increases.
%(For interpretation of the references to color in this figure legend, the reader is referred to the web version of this article.)
By observing the numerical error plots, it can be seen that with the increase of time step  number,
the errors are decreased, and  the speed of convergence to the steady state  solution becomes slow
as the mesh number  $n$ or the moment number $N$ is
increasing.

\end{example}

%%%%%%%%%%%5 n=20 rho u contour %%%%%%%%%%%%%%

%%%%%%%%%%%5 n=50 rho u contour %%%%%%%%%%%%%%

%%%%%%%%%%%%%%%%%% avrho %%%%%%%%%%%%%%%%%%%

%%%%%%%%%%%%%%%%%%% avrho %%%%%%%%%%%%%%%%%%%
%\begin{figure}[H]
%  \centering
%  \includegraphics[width=0.5\textwidth]{./images/vortex/error.pdf}
%  \caption{The relative $\ell^2$ errors of density obtained by using moment methods.}
%    \label{figure:2-errors2D}
%\end{figure}

\section{Conclusions}\label{sec:conclud}
 The paper extended the model reduction method by the operator projection to
a non-linear kinetic description of the Vicsek swarming model.
First, a family of the complicate
Grad type orthogonal functions  depending on a parameter (angle of macroscopic velocity)
were carefully studied in the regard of calculating their derivatives
and projection of those derivatives and the product of velocity and basis and collision term.
Next, building on those discussions and the operator projection,
arbitrary order globally hyperbolic moment system of the kinetic description of the Vicsek swarming model
was derived and their mathematical properties such as
  hyperbolicity,      rotational invariance,  mass-conservation
  and relationship between Grad type expansions in different parameter
were also investigated.
Finally, a semi-implicit numerical scheme was presented to solve a Cauchy problem of our
hyperbolic moment system in order to verify the convergence behavior of the moment method.
It was also compared to the spectral method for the kinetic equation.
 It was seen that the solutions of our hyperbolic moment system could {converge}
 to the  solutions of the kinetic equation for the Vicsek swarming model as the order of the  moment system increases, and
 the moment method could successfully capture key features such as
 shock wave, contact discontinuity, rarefaction wave, and vortex formation.

%For acknowledgements section, please don't number the section, please begin it with \section*{Acknowledgements}

%\section*{Acknowledgements}
%This work was partially supported by
% the Special Project on High-performance Computing under the National Key R\&D Program (No. 2016YFB0200603),
%Science Challenge Project (No. JCKY2016212A502),  and
%the National Natural Science Foundation
%of China (Nos.  91330205,  91630310,  \& 11421101).

%\bibliographystyle{siam}
%\bibliography{Paper1.bib}

%\end{CJK*}
\end{document}